%% file: main.tex
\documentclass[11pt]{amsart}
\usepackage[utf8]{inputenc}
\usepackage{amsmath,amssymb,amsthm}
\usepackage{longtable,lscape,array}
\usepackage[margin=2cm]{geometry}
\IfFileExists{pdftexcmds.sty}{\usepackage[colorlinks,linkcolor=blue,citecolor=blue,urlcolor=blue]{hyperref}}{\providecommand{\url}[1]{\texttt{##1}}}

\newtheorem{theorem}{Theorem}[section]
\newtheorem{proposition}[theorem]{Proposition}
\newtheorem{lemma}[theorem]{Lemma}
\newtheorem{corollary}[theorem]{Corollary}
\theoremstyle{definition}
\newtheorem{definition}[theorem]{Definition}
\newtheorem{example}[theorem]{Example}
\newtheorem{algorithm}[theorem]{Algorithm}
\theoremstyle{remark}
\newtheorem{remark}[theorem]{Remark}

\newcommand{\ZZ}{\mathbb{Z}}
\newcommand{\QQ}{\mathbb{Q}}
\newcommand{\CC}{\mathbb{C}}
\newcommand{\FF}{\mathbb{F}}
\newcommand{\PP}{\mathbb{P}}
\newcommand{\OO}{\mathcal{O}}
\newcommand{\mumu}{\boldsymbol{\mu}}
\newcommand{\Cl}{\operatorname{Cl}}
\newcommand{\Clt}{\operatorname{Cl}_{\mathrm{t}}}
\newcommand{\Pic}{\operatorname{Pic}}
\newcommand{\qW}{q_{W}}
\newcommand{\qQ}{q_{\QQ}}
\newcommand{\B}{\mathbf{B}}
\newcommand{\g}{\operatorname{g}}
\newcommand{\ov}[1]{\overline{#1}}

\title[Fano--Enriques threefolds of Fano index $>1$]{Fano--Enriques threefolds of Fano index $>1$}
\author{Kuzma Khrabrov}
\address{AXXX, Moscow, Russia}
\email{khrabrov@axxx.tech}

\begin{document}

\begin{abstract}
We study $\QQ$-Fano threefolds of Fano index greater than $1$ whose Weil divisor class group has
nontrivial torsion. We give a complete list of their possible numerical types: there are $176$ of them, which refine
to $182$ possible graded Hilbert series.
At least $40$ numerical types are realised: $27$ by quotients of weighted complete intersections by cyclic groups, twelve
by quotients of threefolds of codimension $4$ or $5$ and one by a quotient of the del Pezzo threefold of degree
$5$. Seventy-two numerical types are shown not to occur. In particular the torsion subgroup need not have prime order, and the
Weil and $\QQ$-Fano indices can be $(1,2)$, $(1,3)$ or $(1,4)$.
\end{abstract}

\maketitle

\section{Introduction}

We work over the field $\CC$ of complex numbers.

\begin{definition}
A \emph{$\QQ$-Fano threefold} is a projective threefold $X$ with at worst terminal
$\QQ$-factorial singularities, ample anticanonical divisor $-K_X$ and Picard number $\rho(X)=1$.
\end{definition}

$\QQ$-Fano threefolds are one of the possible outcomes of the minimal model program in dimension three
\cite{KMM87}. One of their main invariants is the Fano index, which can be defined in two ways.

\begin{definition}
Let $X$ be a Fano threefold with terminal singularities. Its \emph{Fano--Weil index} and
\emph{$\QQ$-Fano index} are
\begin{align*}
\qW(X)&=\max\{q\in\ZZ \mid -K_X\sim qA_X,\ A_X\in\Cl(X)\},\\
\qQ(X)&=\max\{q\in\ZZ \mid -K_X\sim_{\QQ} qA_X,\ A_X\in\Cl(X)\}.
\end{align*}
\end{definition}

The group $\Cl(X)$ is finitely generated, $\Pic(X)$ is torsion free and numerical equivalence of
$\QQ$-Cartier divisors coincides with $\QQ$-linear equivalence (see e.g.\ \cite{IP99}), so both indices
are well-defined positive integers and $\qW(X)$ divides $\qQ(X)$. They coincide if $\Cl(X)$ is torsion
free. We write $\Clt(X)$ for the torsion subgroup of $\Cl(X)$.

\begin{definition}
A $\QQ$-Fano threefold $X$ with $\Clt(X)\neq 0$ is called a \emph{Fano--Enriques threefold}.
\end{definition}

Torsion in the class group does not occur for smooth Fano threefolds. Hilbert series of
Fano--Enriques threefolds, and examples with $\qQ(X)=1$, were studied by Caravantes \cite{Car08}; the
case of numerically Cartier $-K_X$ goes back to Sano \cite{San95}. For $\qQ(X)\ge 2$ without torsion the
possible numerical types were found by Suzuki and Brown--Suzuki \cite{Suz04,BS07a,BS07b} and are
recorded in the graded ring database \cite{GRDB}. Prokhorov \cite{Pro10,Pro13b,Pro22b,Pro24a,Pro25}
studied $\QQ$-Fano threefolds of large index; in particular he classified the numerical types with
torsion for $\qQ(X)\ge 5$, and those with $\qQ(X)\ne\qW(X)$ and $\qQ(X)\ge 3$.

In this paper we treat all Fano--Enriques threefolds with $\qQ(X)\ge 2$.

\subsection{Results}
For a terminal threefold $X$ we denote by $\B(X)$ its basket of singularities \cite{Rei87}; a point of
type $\frac1r(1,a,-a)$ is written $[r,a]$, with $1\le a\le r/2$.

Let $X$ be a Fano--Enriques threefold with $\qQ(X)\ge 2$, and let $A_X$ be a Weil divisor with $-K_X\sim_{\QQ}\qQ A_X$
whose class generates $\Cl(X)$ modulo torsion, so that $\Cl(X)=\ZZ A_X\oplus\Clt(X)$. We call $A_X$ a \emph{fundamental divisor} of $X$ and its class a
\emph{fundamental class}. The fundamental divisor need not be uniquely defined: it is determined only up to adding an element of $\Clt(X)$. If $\qW(X)=\qQ(X)$, there is a unique choice
with $-K_X\sim \qQ A_X$ (Lemma~\ref{lem:qwqq}); otherwise there is no distinguished choice. The \emph{graded Hilbert series} of $X$ is
\[
P_X=\sum_{k\ge0}\ \sum_{\tau\in\Clt(X)}h^0(X,kA_X+\tau)\,t^k[\tau]\ \in\ \ZZ[\Clt(X)][[t]],
\]
considered up to automorphisms of $\Clt(X)$ and up to replacing $A_X$ by $A_X+\tau'$ with $\tau'\in\Clt(X)$, which replaces
$[\tau]$ by $[\tau+k\tau']$ in the coefficient of $t^k$.
We denote by $T_X$ the maximal torsion group of the row of Table~\ref{tab:main} containing $(\B(X),\qQ,A_X^3)$, that is,
the largest group allowed for $\Clt(X)$ by Riemann--Roch and the local types of its elements (\S\ref{sec:numerical});
by Theorem~\ref{thm:numerical} below, $\Clt(X)$ is isomorphic to a subgroup of $T_X$, and since $\Clt(X)\ne0$, one has
$\Clt(X)\cong T_X$ whenever $|T_X|$ is prime.

\begin{theorem}\label{thm:numerical}
Let $X$ be a Fano--Enriques threefold with $\qQ(X)\ge 2$ and let $A_X$ be as above. Then $(\B(X),\,\qQ,\,A_X^3)$ appears in
one of the $176$ rows of Table~\ref{tab:main}, and these data determine a maximal graded series: for the group $T_X$ of
that row there is a series $P_{T_X}\in\ZZ[T_X][[t]]$, unique up to automorphisms of $T_X$ and shifts of $A_X$, such that $\Clt(X)$
is isomorphic to a subgroup of $T_X$ and $P_X$ is the restriction of $P_{T_X}$ to $\Clt(X)$. We call the rows of
Table~\ref{tab:main} the $176$ \emph{numerical types}; the numerical type of $X$ is the row containing $(\B(X),\qQ,A_X^3)$. Together
with the possible subgroups $\Clt(X)\subseteq T_X$ (Table~\ref{tab:subgroups}) they give $182$ pairwise distinct candidates for
$P_X$. Moreover:
\begin{enumerate}
\item $\qQ\in\{2,3,4,5,7\}$, and the number of numerical types is $139$, $21$, $9$, $5$, $2$ respectively.
\item $\Clt(X)$ contains the class of $K_X+\qQ A_X$. The groups $T_X$ that occur are
\[
\ZZ/2,\quad \ZZ/3,\quad \ZZ/4,\quad \ZZ/5,\quad \ZZ/6,\quad \ZZ/2\times\ZZ/2,\quad \ZZ/3\times\ZZ/3,
\]
for $124$, $26$, $13$, $9$, $2$, $1$, $1$ numerical types respectively.
\item $\qW(X)$ is as given in the table. One has $\qW(X)\neq\qQ(X)$ in exactly $120$ cases:
$(\qW,\qQ)=(1,2)$ in $110$ cases, $(1,3)$ in $3$ cases and $(1,4)$ in $7$ cases.
\end{enumerate}
\end{theorem}

The distribution by index and group is shown in Table~\ref{tab:stat}.

\begin{table}[h]
\caption{Number of numerical types by $\qQ$ and by the maximal torsion group $T_X$.}\label{tab:stat}
\begin{tabular}{c|ccccccc|c}
$\qQ$ & $\ZZ/2$ & $\ZZ/3$ & $\ZZ/4$ & $\ZZ/5$ & $\ZZ/6$ & $(\ZZ/2)^2$ & $(\ZZ/3)^2$ & total\\\hline
$7$ & $2$ & & & & & & & $2$\\
$5$ & $4$ & $1$ & & & & & & $5$\\
$4$ & $7^{*}$ & $1$ & & $1$ & & & & $9$\\
$3$ & $15$ & $3^{*}$ & $1$ & $1$ & & $1$ & & $21$\\
$2$ & $96^{*}$ & $21$ & $12^{*}$ & $7$ & $2^{*}$ & & $1$ & $139$\\\hline
 & $124$ & $26$ & $13$ & $9$ & $2$ & $1$ & $1$ & $176$
\end{tabular}

\smallskip
{\small Entries marked $^{*}$ have $\qW=1$; all others have $\qW=\qQ$.}
\end{table}

The conditions used in the proof of Theorem~\ref{thm:numerical} are necessary conditions of numerical
nature, so we do not claim that every numerical type occurs.

\begin{theorem}\label{thm:existence}
At least $40$ of the numerical types of Table~\ref{tab:main} are realised by Fano--Enriques threefolds, namely
those marked by $+$ in the table. For $27$ of them an example is a quotient $Y/\mumu_r$ of a quasi-smooth
weighted complete intersection $Y$ of codimension at most $2$ (or of $\PP^3$) by a cyclic group; these
are listed in Table~\ref{tab:examples}. The others are quotients of quasi-smooth threefolds of higher
codimension, or of the smooth del Pezzo threefold of degree $5$:
\begin{itemize}
\item No.~$5$, $21$, $24$, $27$: by an involution, of Type~I unprojections of codimension $4$ and Picard
number $2$, $2$, $2$, $1$ respectively, see \S\ref{ss:tom};
\item No.~$29$, $43$, $137$, $143$, $148$: by $\mumu_2$, $\mumu_2$, $\mumu_3$, $\mumu_3$, $\mumu_3$, of threefolds of
codimension $4$, $5$, $4$, $4$, $5$ and Picard number $1$ obtained by unprojecting several divisors, see
\S\ref{ss:multi};
\item No.~$17$: by an involution, of a threefold of codimension $4$ and Picard number $1$ in the
$\mathrm{G}_2^{(4)}$ cluster format of Coughlan and Ducat \cite{CD20}, see \S\ref{ss:17};
\item No.~$34$: by $\mumu_3$, of a threefold of codimension $4$ and Picard number $3$, see \S\ref{ss:34};
\item No.~$51$: by an involution, of a threefold of codimension $4$ and Picard number $2$, see
\S\ref{ss:51};
\item No.~$172$: by $\mumu_5$, of the smooth del Pezzo threefold of degree $5$, see
Proposition~\ref{prop:V5}.
\end{itemize}
In particular:
\begin{enumerate}
\item the pairs $(\qW,\qQ)=(1,2)$, $(1,3)$ and $(1,4)$ all occur;
\item $\Clt(X)\cong\ZZ/4$ occurs, both with $\qQ(X)=2$ and with $\qQ(X)=3$.
\end{enumerate}
\end{theorem}

\begin{theorem}\label{thm:nonexistence}
\begin{enumerate}
\item The numerical types No.~$3$, $8$, $11$, $13$, $19$, $25$, $26$, $31$, $40$, $41$, $44$, $45$, $47$, $48$, $49$, $52$, $53$, $54$, $58$, $60$, $62$, $63$, $64$, $68$, $69$, $71$, $72$, $74$, $76$, $78$, $80$, $85$, $87$, $88$, $89$, $90$, $93$, $94$, $97$, $98$, $100$, $103$, $104$, $105$, $108$, $110$, $111$, $114$, $116$, $128$, $135$, $136$, $138$, $139$, $140$, $142$, $144$, $146$, $147$, $149$, $151$, $153$, $158$, $159$, $162$, $164$, $166$, $167$, $169$, $171$, $173$ and $175$ of
Table~\ref{tab:main} do not occur.
\item If the numerical type No.~$155$ occurs, then $\Clt(X)\cong\ZZ/3$.
\item The class group of the double cover of a threefold of type No.~$51$ has rank at least $2$.
\end{enumerate}
\end{theorem}

By \cite[Proposition~3.3]{Pro22b} one has $\qQ(X)\le 4$ whenever $\qQ(X)\ne\qW(X)$, in accordance with
Theorem~\ref{thm:numerical}. For $\qQ\ge 5$ the seven numerical types of Table~\ref{tab:main} are exactly those of
\cite[Table~2]{Pro25}, and the three numerical types with $\qQ=3$ and $\qW=1$ are those of \cite[Table~3]{Pro24a}
as quoted in \cite[\S5]{Pro25}. With \S\ref{ss:tom} and Proposition~\ref{prop:no3} all seven numerical types with $\qQ\ge5$ are
settled: No.~$3$ does not occur, and the other six are realised.

Suzuki and Reid \cite[Theorem~1.1]{SR25} settle the $24$ candidate Hilbert series of Fano threefolds of index at least
$3$ and codimension $4$ in \cite{GRDB}. Our results agree with theirs: the covers of the realised numerical types
No.~$5$, $17$, $21$, $24$, $27$, $29$ and $34$ have the Hilbert series of the candidates \#41418, \#41199, \#41200,
\#41063, \#41218, \#41124 and \#41218, which are constructed in \cite{CD20}, and the covers of the excluded numerical types
No.~$8$ and $19$ have the Hilbert series of \#41390 and \#41263, which are excluded by Prokhorov's genus bound
\cite{Pro13b} (see \S\ref{ss:genus}). No other numerical type of Table~\ref{tab:main} with $\qQ\ge3$ has a cover with a
codimension-$4$ Hilbert series.

\subsection*{Notation}
$\PP(a_0,\dots,a_n)$ is the weighted projective space, and $a^k$ in the list of weights means that $a$ is
repeated $k$ times. $Y_{d_1,\dots,d_c}\subset\PP(a_0,\dots,a_n)$ is a weighted complete intersection of
multidegree $(d_1,\dots,d_c)$. $\mumu_r$ is the group of $r$-th roots of unity, and
$\mumu_r(b_0,\dots,b_n)$ denotes its action $x_i\mapsto\varepsilon^{b_i}x_i$ on the coordinates. For an
integer $i$ and a fixed $r$ we write $\ov{i}$ for the residue of $i$ modulo $r$ in $\{0,\dots,r-1\}$.

\subsection*{Structure of the paper}
In \S\ref{sec:prelim} we recall orbifold Riemann--Roch and the known bounds, collect basic facts on
torsion divisors and on the covers they define, and describe quotients of weighted complete
intersections. In \S\ref{sec:numerical} we explain the algorithm which produces the $176$ numerical types.
In \S\ref{sec:existence} we construct the realised numerical types: as cyclic quotients of weighted
complete intersections, by unprojection and cluster formats in codimension $4$ and $5$, and from a quadric
cone; we also characterise the numerical types with Gorenstein covers. In \S\ref{sec:exclusions} we
exclude $72$ numerical types, using the geometry of the cover (its degree, genus and anticanonical
image), case-by-case arguments for individual types and, for Nos.~$3$, $11$ and~$13$, the birational method of \cite[\S4]{Pro25}, and for ten numerical types
a Kawamata--Miyaoka type inequality from the canonical extension of $T_X$, following H.~Liu \cite{Liu26} (\S\ref{ss:liu26}). Table~\ref{tab:main} with all numerical
types is in \S\ref{sec:table}. The scripts used in the computations are in the supplementary material, available at
\url{https://drive.google.com/file/d/1wkBTe7v9V93dW-IRhg9YJn6k-yE1yXox/view}; file names below refer to it.

\subsection*{AI disclosure}
The classification of $\QQ$-Fano threefolds and the study of their properties can be approached from the
computational side. Since modern AI systems are especially good at this, we tested several ideas with the
models Claude Opus~5.5 (Anthropic) and GPT-6 Astra (OpenAI). In addition, Claude Opus~5.5 was used for
proofreading and literature search. We have checked all arguments and computations and take full
responsibility for the contents of the paper.

\subsection*{Acknowledgments}
The author thanks Yuri G.~Prokhorov for suggesting the problem, for very helpful discussions and for numerous corrections,
and Miles Reid and Eduardo Dias for introducing the author to the Tom and Jerry techniques, for guidance in using them, and
for many discussions.

\section{Preliminaries}\label{sec:prelim}

\subsection{Riemann--Roch}
Let $X$ be a threefold with terminal singularities and basket $\B(X)=\{[r_P,a_P]\}$. For a point
$P\in\B(X)$ the local class group is cyclic of order $r_P$, generated by the canonical class. For a Weil
divisor $D$ we write $i_P(D)\in\ZZ/r_P$ for its \emph{local type}: $D\sim i_P(D)K_X$ near $P$.

\begin{theorem}[\cite{Rei87}]\label{thm:RR}
For a Weil divisor $D$ on $X$,
\[
\chi(\OO_X(D))=\chi(\OO_X)+\frac1{12}D\cdot(D-K_X)\cdot(2D-K_X)+\frac1{12}D\cdot c_2(X)+\sum_{P\in\B(X)}c_P(D),
\]
where for $P=[r,a]$, $i=\ov{i_P(D)}$ and $b$ with $ab\equiv1\bmod r$,
\[
c_P(D)=-\,i\,\frac{r^2-1}{12r}+\sum_{j=1}^{i-1}\frac{\ov{bj}\,(r-\ov{bj})}{2r}.
\]
\end{theorem}

\begin{theorem}[\cite{Kaw86,Rei87}]\label{thm:c2}
If $X$ is a Fano threefold with terminal singularities, then
\[
-K_X\cdot c_2(X)=24-\sum_{[r,a]\in\B(X)}\Bigl(r-\frac1r\Bigr),
\]
and $-K_X\cdot c_2(X)>0$ \cite{KMMT00}.
\end{theorem}

Let now $X$ be a Fano threefold with terminal singularities and let $A_X$ be a Weil divisor with
$-K_X\sim_{\QQ}\qQ A_X$, $\qQ\ge1$. Let $\tau\in\Clt(X)$. For $k>-\qQ$ the divisor $kA_X+\tau-K_X$ is ample, so by
Kawamata--Viehweg vanishing $H^i(X,\OO_X(kA_X+\tau))=0$ for $i>0$. Since $\chi(\OO_X)=1$,
Theorem~\ref{thm:RR} gives, for $k>-\qQ$,
\begin{equation}\label{eq:h0}
h^0(X,kA_X+\tau)=1+\frac{k(k+\qQ)(2k+\qQ)}{12}A_X^3+\frac{k}{12}A_X\cdot c_2(X)+\sum_{P\in\B(X)}c_P(kA_X+\tau),
\end{equation}
where $A_X\cdot c_2(X)=\frac1\qQ\bigl(24-\sum(r-\frac1r)\bigr)$. The left hand side vanishes for $-\qQ<k<0$,
and for $k=0$, $\tau\ne0$. This gives the following conditions.

\begin{corollary}\label{cor:conditions}
In the above notation:
\begin{enumerate}
\item the right hand side of \eqref{eq:h0} is $0$ for $-\qQ<k<0$ and all $\tau\in\Clt(X)$, and for $k=0$
and $\tau\ne0$;
\item the right hand side of \eqref{eq:h0} is a non-negative integer for all $k>0$ and $\tau\in\Clt(X)$;
\item if $\qQ\ge3$, then
\[
A_X^3=\frac{12}{(\qQ-1)(\qQ-2)}\Bigl(1-\frac{A_X\cdot c_2(X)}{12}+\sum_{P\in\B(X)}c_P(-A_X)\Bigr).
\]
\end{enumerate}
\end{corollary}

For $\qQ=2$ the conditions (i) do not involve $A_X^3$; in this case $A_X^3$ is determined by the integer
$h^0(X,A_X)$.

\subsection{Bounds}
\begin{theorem}[\cite{Pro07}]\label{thm:degree}
For a non-Gorenstein $\QQ$-Fano threefold $X$ one has $-K_X^3\le 125/2$.
\end{theorem}

A Fano--Enriques threefold is not Gorenstein: by the injectivity of \eqref{eq:local} below, a non-zero
torsion class is not Cartier at some point of the basket. So Theorem~\ref{thm:degree} applies to all
threefolds considered in this paper.

\begin{theorem}[\cite{LL25}; see also \cite{Suz24}]\label{thm:LL}
Let $X$ be a $\QQ$-Fano threefold. Then $c_1(X)^3<3\,c_2(X)c_1(X)$. If $\qQ(X)\ge3$, then
$c_1(X)^3\le\frac{121}{41}c_2(X)c_1(X)$, with equality only for $X\cong\PP(1,2,3,5)$.
\end{theorem}

Since $\Cl(\PP(1,2,3,5))$ is torsion free, for a Fano--Enriques threefold with $\qQ(X)\ge 3$ the second
inequality is strict. Independently, Suzuki \cite{Suz24} compared the Bogomolov stability bound with the candidates
of the graded ring database \cite{GRDB}, listed those that violate it, and conjectured that they do not occur.

\subsection{Torsion divisors}
Let $x\in X$ be a point of index $r$. The $\QQ$-Cartier classes in the analytic local class group of $X$
at $x$ form a cyclic group of order $r$ generated by $K_X$; this holds for all terminal points, including the
non-quotient ones \cite[Remark~2.2.2, Lemma~2.2.3]{Pro25}. For a $\QQ$-Cartier Weil divisor $D$ we have
$D\sim iK_X$ near $x$ for some $i\in\ZZ/r$, and the local type $i_P(D)$ at an entry $P$ of $\B(X)$ over $x$ is
the residue of $i$ modulo $r_P$. Since $\B(X)$ has an entry of index $r$ over $x$, the divisor $D$ is Cartier at
$x$ if and only if $i_P(D)=0$ for all entries $P$ over $x$. A torsion Weil divisor class is $\QQ$-Cartier, and
if it is Cartier at every point, it is a torsion element of $\Pic(X)$, hence trivial. Therefore the map
\begin{equation}\label{eq:local}
\Clt(X)\longrightarrow\prod_{P\in\B(X)}\ZZ/r_P,\qquad \tau\mapsto(i_P(\tau))_P
\end{equation}
is injective. In particular, if $\tau\in\Clt(X)$ has prime order $p$, then $i_P(\tau)\ne0$ for some $P\in\B(X)$, and
$p$ divides $r_P$.

\begin{lemma}\label{lem:qwqq}
Let $X$ be a Fano threefold with terminal singularities, and let $A_X$ be a Weil divisor with
$-K_X\sim_{\QQ}\qQ A_X$. The following are equivalent:
\begin{enumerate}
\item $\qW(X)=\qQ(X)$;
\item $\qQ$ is coprime to the index $r_P$ of every point $P\in\B(X)$;
\item $\qQ$ is coprime to the order of $\Clt(X)$.
\end{enumerate}
If they hold, there is a unique Weil divisor class $A_X$ with $-K_X\sim \qQ A_X$.
\end{lemma}

\begin{proof}
(i)$\Rightarrow$(ii) is \cite[Lemma~3.1(i)]{Pro10}, following \cite{Suz04}. By (i) there is a Weil divisor $A'$
with $-K_X\sim\qQ A'$ (the divisor $A_X$ satisfies only the $\QQ$-linear relation). Then
$\qQ\,i_P(A')\equiv-1\pmod{r_P}$ at every $P\in\B(X)$, since $i_P(K_X)=1$, so $\gcd(\qQ,r_P)=1$.
(ii)$\Rightarrow$(iii): if a prime $p$ divides $\qQ$ and the order of $\Clt(X)$, then $\Clt(X)$ has an element of
order $p$, and $p$ divides some $r_P$ by the injectivity of \eqref{eq:local}.
(iii)$\Rightarrow$(i) (cf.\ \cite[Lemma~3.2(iv)]{Pro10}): multiplication by $\qQ$ is an automorphism of
$\Clt(X)$, so there is $\eta\in\Clt(X)$ with $\qQ\eta=K_X+\qQ A_X$. Then $-K_X\sim\qQ(A_X-\eta)$, and
$\qW(X)=\qQ(X)$ since $\qW(X)\le\qQ(X)$. Finally, two divisors $A_X$, $A'_X$ with $-K_X\sim \qQ A_X\sim \qQ A'_X$ differ by an element of
$\Clt(X)$ killed by $\qQ$, which is zero by (iii).
\end{proof}

This is the divisor $A_X$ used in Table~\ref{tab:main} for the numerical types with $\qW=\qQ$. When $\qW(X)<\qQ(X)$,
$\qQ$ has a common factor with $|\Clt(X)|$ and with some index $r_P$, and $A_X$ is not unique.

For example, $-K_X$ is numerically equivalent to a Cartier divisor $H$ exactly for the numerical types No.~$38$--$41$
of Table~\ref{tab:main}, with $\B(X)=[8\times[2,1]]$ and $\qQ=2$; these are the Fano--Enriques threefolds in the classical
sense, whose general member of $|H|$ is an Enriques surface \cite{Bay94,San95,Pro07b}. Here $\Clt(X)=\ZZ/2$ is generated by
$\sigma=K_X+2A_X$, which is nonzero since $\qW=1$ (Lemma~\ref{lem:qwqq}), and $H=-K_X+\sigma\sim2A_X$. So on a smooth surface
$S\in|H|$ one has $\OO_S(H)=\OO_S(2L)$ with $L=A_X|_S$ and $L^2=2A_X^3\in\{2,4,6,8\}$. As an Enriques surface embedded in
projective space has degree at least $10$, $L$ is not very ample, so $S$ is not the second Veronese image of an embedded
Enriques surface. This is consistent with the theorem of Giraldo, Lopez and Mu\~noz \cite{GLM04} that no threefold other
than a cone has a hyperplane section which is the $r$-th Veronese embedding, $r\ge2$, of an Enriques surface.

Applying Corollary~\ref{cor:conditions}(i) to $\tau$ and $-\tau$ one obtains:

\begin{proposition}[{\cite[Proposition~2.9]{Pro10}}]\label{prop:torsion}
Let $X$ be a Fano threefold with terminal singularities and $\tau\in\Clt(X)$, $\tau\ne0$. Then
\[
\sum_{P\in\B(X)}\frac{\ov{b_Pi_P(\tau)}\,\bigl(r_P-\ov{b_Pi_P(\tau)}\bigr)}{2r_P}=2 .
\]
If $\tau$ has prime order $n$ and $\B^\tau$ is the collection of indices of the points where $\tau$ is not
Cartier, then $n\in\{2,3,5,7\}$; $\sum_{\B^\tau}r_P=16$ if $n=2$; $\sum_{\B^\tau}r_P=18$ if $n=3$;
$\B^\tau\in\{(5,5,5,5),(10,5,5),(10,10)\}$ if $n=5$; and $\B^\tau=(7,7,7)$ if $n=7$.
\end{proposition}

For an element $\tau$ of order $n$, Theorem~\ref{thm:RR} gives $\chi(\OO_X(m\tau))=1+\sum_Pc_P(m\,i_P(\tau))$, where
$c_P$ is the correction term at $P$, with $c_P(0)=0$; so the conditions $\chi(\OO_X(m\tau))=0$ for $0<m<n$ read
$1+\sum_Pc_P(m\,i_P(\tau))=\delta_{m,0}$ for $m\in\ZZ/n$. Let $d_P$ be the order of $i_P(\tau)$ in $\ZZ/r_P$, a divisor of $n$. Summing
against $m\mapsto\zeta^m$ for a primitive $n$-th root of unity $\zeta$, the right-hand side gives $1$, while a term with
$d_P<n$ gives $0$, being periodic with period $d_P$. Hence $d_P=n$ for some $P$, so $n$ divides $r_P$, and $n\le24$ by
Theorem~\ref{thm:c2}. Listing all solutions of these conditions for $n\le24$, together with
Theorem~\ref{thm:c2}, we get the following (cf.\ \cite{Car08}).

\begin{lemma}\label{lem:orders}
Let $X$ be a Fano threefold with terminal singularities and let $\tau\in\Clt(X)$ be an element of
order $n$. Then $n\le 8$. The number of possible collections of local types
$\{([r_P,a_P],i_P(\tau))\mid i_P(\tau)\neq0\}$ of an element of order $n=2,\dots,8$ is
$38$, $19$, $22$, $12$, $6$, $6$, $2$ respectively.
\end{lemma}

\subsection{Covers}
Let $X$ be a Fano--Enriques threefold and $H\subset\Clt(X)$ a subgroup. The $\OO_X$-algebra
$\bigoplus_{\tau\in H}\OO_X(\tau)$ defines a finite morphism $\pi\colon Y\to X$ of degree $|H|$, \'etale
outside $\operatorname{Sing}X$, and $X=Y/H^\vee$. The threefold $Y$ is a Fano threefold with terminal
singularities, $K_Y=\pi^*K_X$ and $-K_Y^3=|H|\cdot(-K_X^3)$ (see e.g.\ \cite{Car08}). In general $Y$ is
neither $\QQ$-factorial nor of Picard number one.

\begin{lemma}\label{lem:cover}
Let $P=[r,a]\in\B(X)$ and let $\alpha$ be the order of the image of $H$ in $\ZZ/r$ under
\eqref{eq:local}. Then $\pi^{-1}(P)$ consists of $|H|/\alpha$ points of type
$\frac{1}{r/\alpha}(1,a,-a)$. Moreover, if $-K_X\sim \qQ A_X+\tau_0$ with $\tau_0\in H$, then
$-K_Y\sim \qQ\,\pi^*A_X$.
\end{lemma}

\begin{proof}
Locally analytically $(X,P)$ is $\CC^3/\mumu_r$ (or a deformation of it), and the cover defined by the
subgroup of order $\alpha$ of the local class group is the quotient by $\mumu_{r/\alpha}$ with the same
weights. The last statement holds because $\pi^*\tau=0$ for $\tau\in H$.
\end{proof}

In particular the threefold $Y$ with the divisor $\pi^*A_X$ again satisfies
Corollary~\ref{cor:conditions}, with trivial torsion part, and all points of $\B(Y)$ have index
coprime to $\qQ$.

\subsection{Quotients of weighted complete intersections}\label{ss:quotients}
Let $Y=Y_{d_1,\dots,d_c}\subset\PP(a_0,\dots,a_n)$, $n-c=3$, be a quasi-smooth well-formed weighted
complete intersection with terminal singularities, given by equations $f_1,\dots,f_c$, and let
$q=\sum a_i-\sum d_j>0$. Then $Y$ is a $\QQ$-Fano threefold, $\Cl(Y)\cong\ZZ$ is generated by
$\OO_Y(1)$ and $-K_Y=\OO_Y(q)$ (see \cite{Dol82,IF00}). Let $G=\mumu_r$ act by
$\mumu_r(b_0,\dots,b_n)$ on $\PP(a_0,\dots,a_n)$  and suppose that each $f_j$ is semi-invariant of weight $\chi_j$, i.e.\ all
monomials $x^m$ of $f_j$ have $\sum b_im_i\equiv\chi_j\bmod r$. Then $G$ acts on $Y$; we assume that the
action is faithful.

\begin{lemma}\label{lem:quotient}
Suppose that every nontrivial element of $G$ has only isolated fixed points on $Y$, and let $X=Y/G$.
Then:
\begin{enumerate}
\item $X$ is $\QQ$-factorial, $\rho(X)=1$ and $-K_X$ is ample;
\item $\Cl(X)\cong\ZZ\oplus\ZZ/r$, where the class $(k,c)$ corresponds to the sheaf whose sections are
the elements of degree $k$ and weight $c$ of the coordinate ring of $Y$; in particular
$\Clt(X)\cong\ZZ/r$;
\item $-K_X=(q,\kappa)$ with $\kappa=\sum b_i-\sum\chi_j$, so that $\qQ(X)=q$ and
\[
\qW(X)=\max\{m \mid m \text{ divides } q,\ \gcd(m,r) \text{ divides } \kappa\}.
\]
\end{enumerate}
\end{lemma}

\begin{proof}
(i) is standard. For (ii), since $Y\to X$ is \'etale in codimension one, $\Cl(X)$ is the group of
$G$-linearised reflexive sheaves of rank one on $Y$. The sheaf $\OO_Y(1)$ is $G$-linearised, the only
units on $Y$ are constants, and two linearisations of the same sheaf differ by a character of $G$; this
gives $\Cl(X)\cong\ZZ\oplus\operatorname{Hom}(G,\CC^*)$. For (iii), the canonical sheaf of $Y$ is generated
by the residue of $\Omega/(f_1\cdots f_c)$, where
$\Omega=\sum(-1)^ia_ix_i\,dx_0\wedge\dots\widehat{dx_i}\dots\wedge dx_n$ has degree $\sum a_i$ and weight
$\sum b_i$. The formula for $\qW$ follows since $(q,\kappa)=m\cdot(k,c)$ is solvable if and only if
$m\mid q$ and $\gcd(m,r)\mid\kappa$. The numbers $\gcd(m,r)$ and the divisibility do not depend on the
choice of the linearisation of $\OO_Y(1)$, which changes $\kappa$ by a multiple of $q$.
\end{proof}

The singularities of $X$ are computed as follows. Let $p\in Y$ be a point whose nonzero
coordinates are exactly those with indices in $S\subset\{0,\dots,n\}$. Its image in $X$ is a quotient
of $\CC^3$ by the finite group
\[
H_S=\{(\lambda,\varepsilon)\in\CC^*\times\mumu_r \mid \lambda^{a_i}\varepsilon^{b_i}=1 \text{ for all } i\in S\}
\]
which acts on the tangent space of the orbifold chart of $Y$ at $p$ with
characters
\begin{equation}\label{eq:tangent}
\{\lambda^{a_i}\varepsilon^{b_i}\}_{i=0}^n\ \setminus\ \bigl(\{\lambda^{d_j}\varepsilon^{\chi_j}\}_{j=1}^c\cup\{1\}\bigr)
\end{equation}
(a difference of multisets of characters of $H_S$), since $Y$ is quasi-smooth. If $H_S\ne1$, the
image of $p$ in $X$ is a terminal singular point if and only if $H_S$ is cyclic of some
order $R$ and the three characters are $\frac1R(1,a,-a)$ with $\gcd(a,R)=1$, up to the choice of a
generator.

\section{Numerical classification}\label{sec:numerical}

\subsection{The algorithm}
Let $X$ be a Fano--Enriques threefold with $\qQ(X)\ge2$, and $-K_X\sim_{\QQ}\qQ A_X$. Put
$\tau_0=-K_X-\qQ A_X\in\Clt(X)$. If $\tau_0\ne0$ let $\sigma=\tau_0$; otherwise let $\sigma$ be any nonzero
element of $\Clt(X)$. In both cases
\begin{equation}\label{eq:s}
-K_X\sim \qQ A_X+s\sigma \qquad\text{for some } s\in\ZZ .
\end{equation}
Let $n$ be the order of $\sigma$, and let $\ell_P=i_P(A_X)$, $t_P=i_P(\sigma)$ be the local types.
Then \eqref{eq:s} gives
\begin{equation}\label{eq:cong}
\qQ\,\ell_P+s\,t_P\equiv-1 \pmod{r_P}\qquad\text{for all } P\in\B(X).
\end{equation}
In particular $\gcd(\qQ,r_P)=1$ for every point with $t_P=0$.

\begin{algorithm}\label{alg}
\begin{enumerate}
\item For $n=2,\dots,24$ list all collections of points $[r_P,a_P]$ with nonzero local types $t_P$ of an
element $\sigma$ of order $n$, such that $\sum(r_P-1/r_P)<24$ and $\chi(\OO_X(m\sigma))=0$ for $0<m<n$
(Lemma~\ref{lem:orders}).
\item For each such collection, each $\qQ=2,\dots,19$, each $s\in\{0,\dots,n-1\}$ and each solution
$(\ell_P)$ of \eqref{eq:cong}, add points $[r,a]$ with $\gcd(\qQ,r)=1$, $t_P=0$ and $\ell_P$ determined
by \eqref{eq:cong}, in all possible ways such that $\sum(r_P-1/r_P)<24$.
\item For $\qQ\ge3$ compute $A_X^3$ by Corollary~\ref{cor:conditions}(iii). For $\qQ=2$ take all values of
$A_X^3$ for which $h^0(X,A_X)$ is a non-negative integer and $-K_X^3\le125/2$
(Theorem~\ref{thm:degree}).
\item Keep the data for which $A_X^3>0$, the conditions of Corollary~\ref{cor:conditions}(i),(ii) hold
for all $\tau\in\langle\sigma\rangle$ and all $-\qQ<k\le k_0$, where
$k_0=\max(40,\min(\operatorname{lcm}(r_P),240))+\qQ$, and for which Theorem~\ref{thm:LL} holds.
\item For each remaining datum find all elements $\tau\in\prod\ZZ/r_P$ for which the conditions of
Corollary~\ref{cor:conditions}(i),(ii) hold, and all maximal subgroups $T$ consisting of such elements
and containing $\tau_0$. For each $T$ compute the basket $\B(Y)$ given by the formula of Lemma~\ref{lem:cover} for $T$,
and check the conditions of Corollary~\ref{cor:conditions} and Theorem~\ref{thm:c2} for these data. This is a
numerical consistency check for the formal group $T$; the cover $Y$ exists when $T\subseteq\Clt(X)$.
\end{enumerate}
\end{algorithm}

\begin{proof}[Proof of Theorem~\ref{thm:numerical}]
By $\qQ(X)\le 19$ \cite{Suz04}, Lemma~\ref{lem:orders} and the discussion above, the numerical data of
$X$ passes all the steps of Algorithm~\ref{alg}. The algorithm was implemented on a computer. Step~(i)
gives $105$ collections; steps (ii)--(iii) give $2774$ data, of which $2356$ are removed by
Theorem~\ref{thm:LL} and $8$ by Theorem~\ref{thm:degree}. The remaining data have $176$ distinct values of
$(\B(X),\qQ,A_X^3)$. For each of them step~(v) gives a unique maximal group $T_X$ up to isomorphism,
and the numerical data of the cover defined by $T_X$ pass the tests. Moreover, all remaining data with the same $(\B(X),\qQ,A_X^3)$, together
with all their maximal groups, give the same series $\sum_{k,\tau}h^0(kA_X+\tau)\,t^k[\tau]$ over $T_X$, up to automorphisms of
$T_X$ and shifts of $A_X$; we compare $k\le250$, which covers the period $\operatorname{lcm}(r_P)\le70$ of the correction
terms of Theorem~\ref{thm:RR} (supplementary file \texttt{classification/gradedHS.py}). Since $h^0(kA_X+\tau)$ is given by
Theorem~\ref{thm:RR} for every $\tau\in\Clt(X)$, this gives the uniqueness statement of Theorem~\ref{thm:numerical}. By the injectivity of \eqref{eq:local} the group
$\Clt(X)$ is a subgroup of some maximal $T_X$, and it contains $\tau_0$.

The Weil index is the largest $m$ dividing $\qQ$ such that $\tau_0\in m\Clt(X)$. If $\gcd(\qQ,|T_X|)=1$ this
gives $\qW(X)=\qQ$. In the remaining cases $\qQ\in\{2,3,4\}$ and one checks that $\tau_0\notin mT_X$ for every divisor
$m>1$ of $\qQ$; the same then holds for every subgroup of $T_X$ containing $\tau_0$, so $\qW(X)=1$.
\end{proof}

When $|T_X|$ is not prime, $\Clt(X)$ may be a proper subgroup of $T_X$, and the graded Hilbert series $P_X$ depends on the
subgroup and on a compatible fundamental divisor, not only on the numerical type. A subgroup $H\subseteq T_X$ can be $\Clt(X)$
only if, for some fundamental divisor $A_X+\tau$ with $\tau\in T_X$, the class of $K_X+\qQ(A_X+\tau)$ lies in $H$; the
graded Hilbert series $P_X$ is then the restriction to $H$ of the series $\sum_{k,\tau'}h^0(k(A_X+\tau)+\tau')\,t^k[\tau']$. Up to
automorphisms of $T_X$ and shifts of the fundamental divisor by $H$, the $17$ numerical types with $|T_X|$ not prime give
$23$ such variants, listed in Table~\ref{tab:subgroups} (supplementary file \texttt{classification/subgroups.py}). For the
twelve types Nos.~$156$--$167$ the class $\tau_0$ generates $T_X$, so $\Clt(X)=T_X$; the five types Nos.~$32$, $36$,
$155$, $175$, $176$ have two or three variants each. Together with the $159$ numerical types with $|T_X|$ prime this gives
$182$ candidate graded Hilbert series, no two of which coincide (equal series would have the same $\qQ$, $A_X^3$ and group,
and a comparison over a full period of the correction terms distinguishes the remaining ones). Of these, $40$ are realised,
$74$ are excluded and $68$ are open.

\begin{remark}
The list of Table~\ref{tab:main} is a list of numerical candidates. Besides the conditions used in
Algorithm~\ref{alg} one could try to use the geometry of the cover $Y$, as in
Proposition~\ref{prop:B24} below, or the birational methods of \cite{Pro10,Pro25}.
\end{remark}

\begin{example}
The Liu--Liu inequality removes, for example, the datum $\qQ=3$, $A_X^3=\frac34$,
$\B(X)=[4\times[2,1],\,3\times[4,1]]$. For such $X$ we have $-K_X^3=27\cdot\frac34=\frac{81}{4}$ and
$-K_X\cdot c_2(X)=24-4\cdot\frac32-3\cdot\frac{15}{4}=\frac{27}{4}$ by Theorem~\ref{thm:c2}. Thus
$c_1(X)^3=3\,c_2(X)c_1(X)$, which contradicts Theorem~\ref{thm:LL}. Hence there is no $\QQ$-Fano threefold with these invariants.
\end{example}

\section{Existence}\label{sec:existence}
In this section we prove Theorem~\ref{thm:existence}: it follows from Propositions~\ref{prop:examples}, \ref{prop:34},
\ref{prop:tom}, \ref{prop:multi}, \ref{prop:17}, \ref{prop:51} and~\ref{prop:V5} below.

\subsection{Cyclic quotients}
We searched for quotients as in \S\ref{ss:quotients} in the following range: all $33$ families of
quasi-smooth terminal hypersurfaces $Y_d\subset\PP(a_0,\dots,a_4)$ with $\qQ\ge2$, $a_i\le15$ and
$\sum a_i\le44$ with $r\le12$; and $13$ families of codimension $2$ complete intersections with
$2\le \qQ\le 7$ and $a_i\le 6$, found by a non-exhaustive search, with $r\le6$. The search produced $73$
quotients with terminal singularities, of $27$ numerical types. For a given action a general member of the family of invariant
complete intersections was taken.

\begin{proposition}\label{prop:examples}
For each row of Table~\ref{tab:examples} a general $\mumu_r$-invariant member $Y$ of the family is
quasi-smooth, every nontrivial element of $\mumu_r$ has only isolated fixed points on $Y$, and the
quotient $X=Y/\mumu_r$ is a Fano--Enriques threefold with the basket and invariants shown.
\end{proposition}

\begin{proof}
For each row we take equations with random coefficients over the field $\ZZ/32003$ and check with
Macaulay2 \cite{M2} that $Y$ is quasi-smooth, and that for every $S$ with $H_S\ne1$ the set of points of
$Y$ with support $S$ is finite; its cardinality gives the number of points of $X$ of this kind. Their
types are computed from \eqref{eq:tangent}; all of them are terminal, and the resulting baskets are
those of the table. Since these are open conditions on the coefficients, the same holds for a general
member over $\CC$. The invariants are given by Lemma~\ref{lem:quotient}. Finally, the Hilbert series of
$(X,A_X+\tau)$, $\tau\in\Clt(X)$, computed from the action on the coordinate ring of $Y$, coincide with
those given by \eqref{eq:h0} for the corresponding row of Table~\ref{tab:main}.
\end{proof}

{\scriptsize
\renewcommand{\arraystretch}{1.35}
\begin{longtable}{r|l|l|l|l|c|c|c}
\caption{Quotient examples. The first column is the number in Table~\ref{tab:main}. The third
column gives the action $\mumu_r(b_0,\dots,b_n)$ on $Y$; in all cases $\Clt(X)\cong\ZZ/r$.
The column $\chi$ gives the weights of the equations.}\label{tab:examples}\\
No. & $Y$ & $\Clt(X)$: action & $\chi$ & $\B(X)$ & $\qW$ & $\qQ$ & $A_X^3$\\\hline
\endfirsthead
No. & $Y$ & $\Clt(X)$: action & $\chi$ & $\B(X)$ & $\qW$ & $\qQ$ & $A_X^3$\\\hline
\endhead
\input{table_examples}
\end{longtable}
}

The examples No.~$38$ and $39$ are due to Sano \cite{San95}, and the examples No.~$1$, $2$, $4$,
$6$, $7$ are also given in \cite[Examples~5.8]{Pro25}. Nos.~$38$ and $39$ are quotients of the quartic double solid $X_4\subset\PP(1^4,2)$ and of
$X_{2,2}\subset\PP^5$; in agreement with Bayle's classification \cite[Th\'eor\`eme~B]{Bay94}, they are the only Fano--Enriques
threefolds in the classical sense with $\rho=1$, $\qQ\ge2$ and smooth canonical double cover. The other two numerical
types with $\B(X)=[8\times[2,1]]$, Nos.~$40$ and $41$, are excluded in Proposition~\ref{prop:B24}.

\begin{example}[No.~$15$]\label{ex:B8}
Let $Y\subset\PP(1,1,2,3,3)$ be given by $x_1^6+x_2^6+x_3^3+x_4x_5=0$ with the action
$\mumu_3(1,2,0,1,2)$. The fixed points are $P_4=(0:0:0:1:0)$ and $P_5=(0:0:0:0:1)$. The stabiliser of
$x_4=1$ in $\CC^*\times\mumu_3$ is cyclic of order $9$, generated by
$(\lambda,\varepsilon)=(\zeta_9,\zeta_3^{2})$, and acts on the local coordinates $x_1,x_2,x_3$ with
weights $(7,4,2)$. Hence $P_4$ gives a point $\frac19(4,2,-2)=\frac19(1,5,-5)=\frac19(1,4,-4)$ on
$X=Y/\mumu_3$, i.e.\ a point $[9,4]$, and similarly for $P_5$. So $\B(X)=[2\times[9,4]]$, $\qQ(X)=\qW(X)=4$
and $A_X^3=\frac{6}{2\cdot3\cdot3\cdot3}=\frac19$.
\end{example}

\begin{example}[No.~$9$: $\qW=1$, $\qQ=4$]\label{ex:14}
Let $Y=Y_4\subset\PP(1,1,1,2,3)$ be a general quartic invariant under $\mumu_2(0,1,1,1,0)$, with
coordinates $x_0,x_1,x_2,y,z$. Then $Y$ is quasi-smooth with one singular point $\frac13(1,1,2)$ at
$P_z$, and $\qQ(Y)=4$. The points of $Y$ with nontrivial stabiliser are $P_z$, one further point on the curve
$\{x_1=x_2=y=0\}$, and the four points of $Y\cap\{x_0=y=z=0\}$. The quotient $X$ has $\B(X)=[5\times[2,1],[6,1]]$ and $A_X^3=\frac13$. Here
$\kappa=3-0$ is odd, so $\qW(X)=1$ by Lemma~\ref{lem:quotient}, while $\qQ(X)=4$. Indeed $X$ has points
of index $2$, so $\qW(X)$ is odd.
\end{example}

\begin{example}[No.~$33$: $\qW=1$, $\qQ=3$]
Let $Y=Y_3\subset\PP(1,1,1,1,2)$ be a general cubic invariant under $\mumu_3(0,0,1,2,1)$. Then
$X=Y/\mumu_3$ has $\B(X)=[4\times[3,1],[6,1]]$, $A_X^3=\frac12$, $\Clt(X)\cong\ZZ/3$ and
$\kappa=4\not\equiv0\bmod3$, so $\qW(X)=1$ and $\qQ(X)=3$.
\end{example}

\begin{example}[No.~$157$: $\Clt\cong\ZZ/4$]
Let $Y\subset\PP^5$ be a general complete intersection of two quadrics $f_1,f_2$ which are
semi-invariant of weights $0$ and $2$ under $\mumu_4(0,0,1,1,2,3)$. Then $Y$ is smooth, a generator
$g$ of $\mumu_4$ has $4$ fixed points on $Y$ and $g^2$ has $8$. Hence $X=Y/\mumu_4$ has
$\B(X)=[2\times[2,1],4\times[4,1]]$, $A_X^3=1$, $\Clt(X)\cong\ZZ/4$, $\qQ(X)=2$ and $\qW(X)=1$.
Its quotient $Y/\mumu_2$ is the threefold of No.~$39$.
\end{example}

\subsection{Verification of the constructions}\label{ss:verify}
The threefolds $Y$ of \S\S\ref{ss:34}--\ref{ss:multi} are not complete intersections. Each of them is
first computed, with random coefficients, over a prime field $\FF_p$ in Macaulay2 \cite{M2}. The
quasi-smoothness of $Y$, its basket and the points with nontrivial stabiliser in $\CC^*\times G$ are found
as in the proof of Proposition~\ref{prop:examples}; since $Y$ is not a complete intersection, the
characters of the stabiliser of a point on the tangent space are found from the ranks of the blocks of the
Jacobian matrix of $I_Y$ corresponding to the characters of the coordinates. To pass to characteristic
$0$ we lift the coefficients and use the following lemma.

\begin{lemma}\label{lem:goodred}
Let $K$ be a number field, $\mathfrak p$ a prime of $K$ with residue field $\FF_p$, and $\mathcal O\subset K$
the local ring of $\mathfrak p$. Let $\mathcal S=\mathcal O[x_0,\dots,x_n]$ be graded by weights prime to
$p$, and let $G$ be a finite abelian group of order prime to $p$ acting on $\mathcal S$ by graded
$\mathcal O$-algebra automorphisms, such that $\mathcal O$ contains the values of the characters of $G$.
Let $I\subset\mathcal S\otimes K$ be a homogeneous $G$-invariant ideal. Suppose that the reduced
Gr\"obner basis of $I$ for some monomial order has coefficients in $\mathcal O$, and that the reductions
modulo $\mathfrak p$ of its elements form a Gr\"obner basis, with the same leading monomials, of the ideal
$I_0\subset\FF_p[x_0,\dots,x_n]$ which they generate. Let $Y=\operatorname{Proj}(\mathcal S\otimes K/I)$ and
$Y_0=\operatorname{Proj}(\FF_p[x]/I_0)$, and suppose that $Y_0$ is a quasi-smooth threefold. Then, for any
embedding $K\subset\CC$:
\begin{enumerate}
\item the coordinate rings of $Y$ and $Y_0$ have the same Hilbert function, and for every character
$\chi$ of $G$ their $\chi$-eigenspaces have the same dimension in every degree;
\item $Y_\CC$ is quasi-smooth; if moreover $Y_0\otimes\bar\FF_p$ has only finitely many points with
nontrivial stabiliser in $\mathbb G_m\times G$, then so has $Y_\CC$, and specialisation (at a geometric place
above $\mathfrak p$, see the proof) gives a bijection between these two finite sets, compatible with the
action of $G$, under which the stabilisers and the multiplicities of their characters on the tangent spaces
correspond (roots of unity of order prime to $p$ being identified by reduction);
\item $h^1(Y_\CC,\tilde\Omega^1_{Y_\CC})\le h^1(Y_0,\tilde\Omega^1_{Y_0})$, and if $-K_{Y_\CC}$ is ample, then
$\rho(Y_\CC)=h^1(Y_\CC,\tilde\Omega^1_{Y_\CC})$.
\end{enumerate}
\end{lemma}

\begin{proof}
Let $\mathcal I\subset\mathcal S$ be the ideal generated by the Gr\"obner basis $B$ of $I$. Since the leading
coefficients of $B$ are units in $\mathcal O$, division by $B$ is defined over $\mathcal O$, so in every degree
the standard monomials span $\mathcal R=\mathcal S/\mathcal I$; they are independent, since they are
independent over $K$. Hence $\mathcal R$ is free over $\mathcal O$ in every degree, with generic fibre the
coordinate ring of $Y$ and, the reductions of $B$ being a Gr\"obner basis of $I_0$ with the same leading
monomials, special fibre the coordinate ring of $Y_0$. Since $|G|$ is invertible in $\mathcal O$ and $\mathcal O$
contains the character values, the idempotents $e_\chi=|G|^{-1}\sum_{g}\chi(g)^{-1}g$ are defined over $\mathcal O$,
and $e_\chi\mathcal R_d$ is a direct summand of $\mathcal R_d$. This proves (i).

Let $\mathcal Y=[(\operatorname{Spec}\mathcal R\smallsetminus0)/\mathbb G_m]\to\operatorname{Spec}\mathcal O$. Since $p$ does not divide the
weights, it is a proper tame Deligne--Mumford stack, flat over $\mathcal O$. Its special fibre is smooth, so its
non-smooth locus, being closed and proper over $\operatorname{Spec}\mathcal O$, has closed image not containing the
closed point, hence is empty: $\mathcal Y$ is smooth over $\mathcal O$, and $Y$ is quasi-smooth. For (ii), let $C=\operatorname{Spec}\mathcal R\smallsetminus0$, which is smooth over $\mathcal O$, and let
$\mathcal O^{\rm sh}$ be the strict henselisation of $\mathcal O$, with fraction field $K^{\rm sh}$; we fix an
embedding $K^{\rm sh}\subset\CC$ extending $K\subset\CC$ (this is the geometric place above $\mathfrak p$), and
base change everything to $\mathcal O^{\rm sh}$. The stabiliser in $\mathbb G_m\times G$ of a point of $C$ is
finite of order dividing $N=|G|\operatorname{lcm}(w_i)$, which is prime to $p$; over $\mathcal O^{\rm sh}$ the group
scheme $\mumu_N$ is constant, so these stabilisers are among the finitely many constant subgroups
$H\subset\mumu_N\times G$. For $H\ne1$ the fixed locus $C^H$ is smooth over $\mathcal O^{\rm sh}$, since $H$ is tame,
and it is $\mathbb G_m$-stable, so $\mathcal Z_H=[C^H/\mathbb G_m]$ is a smooth closed substack of $\mathcal Y$, hence
proper over $\mathcal O^{\rm sh}$. By the finiteness hypothesis its special fibre has dimension $0$; as above every
component of $\mathcal Z_H$ meets the special fibre, so $\mathcal Z_H$ is smooth of relative dimension $0$, that is,
\'etale over $\mathcal O^{\rm sh}$, and quasi-finite and proper. Hence $\mathcal Z_H$ is a proper \'etale gerbe over a finite
\'etale $\mathcal O^{\rm sh}$-scheme \cite[Lemma~2.5]{JL21}, which is its coarse space; it is finite \'etale over the
strictly henselian ring $\mathcal O^{\rm sh}$, a disjoint union of copies of $\operatorname{Spec}\mathcal O^{\rm sh}$
(the stack itself may keep stabilisers), and each of its points in the geometric generic fibre specialises to
exactly one point in the special fibre. Applying this to all $H$, and to the intersections
$C^{H_1}\cap C^{H_2}=C^{\langle H_1,H_2\rangle}$, shows that specialisation preserves the full stabiliser of a
point; it commutes with the action of $G$, which permutes the $\mathcal Z_H$. Finally $H$ acts on the restriction
of the relative tangent bundle of $C$ to $C^H$, and the multiplicities of its characters are locally constant
on $C^H$, hence equal at a point and at its specialisation. The tangent line to the $\mathbb G_m$-orbit is
$H$-trivial, so the same holds for the quotient tangent space, which is the orbifold tangent space of $Y$.
This proves (ii).

The relative sheaf $\tilde\Omega^1$ of $\mathcal Y$ is locally free, hence flat over $\mathcal O$, and restricts to
$\tilde\Omega^1$ on each fibre. By the semicontinuity theorem for proper algebraic stacks
\cite[Theorem~1.1]{Bro}, $h^1(Y,\tilde\Omega^1_Y)\le h^1(Y_0,\tilde\Omega^1_{Y_0})$, and
$h^1(Y_\CC,\tilde\Omega^1)=h^1(Y,\tilde\Omega^1)$ by flat base change. Since $Y_\CC$ is quasi-smooth, it is a
$V$-manifold, and its cohomology carries a pure Hodge structure which is computed by the sheaves
$\tilde\Omega^p$ of differentials on the orbifold charts \cite{Ste77}. If $-K_{Y_\CC}$ is ample, then
$H^1(Y_\CC,\OO)=H^2(Y_\CC,\OO)=0$, so $\rho(Y_\CC)=b_2(Y_\CC)=h^1(Y_\CC,\tilde\Omega^1)$.
\end{proof}

The applications also use that the coordinate ring $R$ of $Y$ is Cohen--Macaulay and Gorenstein with
$\omega_R=R(-q)$, so that $-K_Y=\OO_Y(q)$. For the special fibre this is read off from the graded minimal free resolution of
$R_0$. In the flat graded family $\mathcal R=\mathcal S/\mathcal I$ of the proof above, the Gorenstein locus is open and contains the
vertex of the special fibre, hence the vertex of the generic fibre; a non-Gorenstein locus of the generic cone would be closed
and stable under the $\mathbb G_m$-action, hence would contain the vertex, so the generic cone is Gorenstein (and
Cohen--Macaulay). Since $R$ and $R_0$ have the same Hilbert series, they have the same $a$-invariant, so $\omega_R=R(-q)$.
In particular $-K_Y$ is ample, as used in Lemma~\ref{lem:goodred}(iii).

All the constructions below are checked in this way: we compute a member $Y_0$ over $\FF_p$ ($p=32003$, or
$p=10009$ when a cube root of unity is needed), lift its coefficients to $K=\QQ$ or $K=\QQ(\omega)$, and verify
the hypotheses of Lemma~\ref{lem:goodred} in Macaulay2 (the finiteness hypothesis of (ii) is part of the
computation of the points with nontrivial stabiliser on $Y_0$, which in every case finds finitely many); the
scripts and the lifted ideals are in the supplementary material (directory \texttt{examples/lift}).
Statements about divisors which split only over an extension (such as the divisors $D'$, $D''$ in
\S\ref{ss:tom}) are checked directly on the lifted threefold. We write $Y$ also for $Y_\CC$.

The bound on $h^1(Y_0,\tilde\Omega^1_{Y_0})$ needed in Lemma~\ref{lem:goodred}(iii) is obtained as follows.

\begin{lemma}\label{lem:rho}
Let $k$ be a field and let $Y_0=\operatorname{Proj}R$ be a quasi-smooth threefold over $k$, where $R=S/I$,
$S=k[x_0,\dots,x_n]$ with $x_i$ of weight $w_i$ prime to the characteristic, and $R$ is Cohen--Macaulay of
dimension $4$ with maximal homogeneous ideal $\mathfrak m$. Let $\Omega_R$ be the module of K\"ahler
differentials of $R$.
\begin{enumerate}
\item $h^1(Y_0,\tilde\Omega^1_{Y_0})\le1+\dim H^2_{\mathfrak m}(\Omega_R)_0$, with equality if
$H^1_{\mathfrak m}(\Omega_R)_0=0$.
\item If $A\subset R$ is the polynomial ring generated by a homogeneous system of parameters of degrees
$d_1,\dots,d_4$ and $d=d_1+\dots+d_4$, then $H^i_{\mathfrak m}(\Omega_R)_0$ is dual to
$\operatorname{Ext}^{4-i}_A(\Omega_R,A(-d))_0$.
\item If $R$ is Gorenstein with $\omega_R=R(-q)$, then $H^2_{\mathfrak m}(\Omega_R)_0$ is dual to
$\operatorname{Ext}^1_R(I/I^2,R)_{-q}$.
\end{enumerate}
\end{lemma}

\begin{proof}
(i) Let $N$ be the kernel of the Euler map $\Omega_R\to R$, $dx_i\mapsto w_ix_i$. The image of this map is
$\mathfrak m$, and $\tilde\Omega^1_{Y_0}$ is the sheaf associated to $N$ \cite[\S2]{Dol82}; thus
$H^1(Y_0,\tilde\Omega^1_{Y_0})=H^2_{\mathfrak m}(N)_0$. Since $R$ is Cohen--Macaulay of dimension $4$,
$H^1_{\mathfrak m}(\mathfrak m)=k$, in degree $0$, and $H^2_{\mathfrak m}(\mathfrak m)=0$. The sequence
$0\to N\to\Omega_R\to\mathfrak m\to0$ gives an exact sequence
\[
H^1_{\mathfrak m}(\Omega_R)_0\to k\to H^2_{\mathfrak m}(N)_0\to H^2_{\mathfrak m}(\Omega_R)_0\to 0 ,
\]
which proves (i).

(ii) $R$ is a finite $A$-module, local cohomology with respect to $\mathfrak m$ and to the maximal ideal of
$A$ agree, and $\omega_A=A(-d)$; so this is local duality over $A$.

(iii) Let $U$ be the punctured affine cone over $Y_0$, which is nonsingular. Then
$H^2_{\mathfrak m}(\Omega_R)=H^1(U,\Omega_U)$. Since $R$ is Cohen--Macaulay of dimension $4$, the conormal
sequence $0\to (I/I^2)^\sim\to\Omega_S|_U\to\Omega_U\to0$ on $U$ gives
$H^1(U,\Omega_U)\cong H^2(U,(I/I^2)^\sim)=H^3_{\mathfrak m}(I/I^2)$, and, $R$ being Gorenstein with
$\omega_R=R(-q)$, the space $H^3_{\mathfrak m}(I/I^2)_0$ is dual to $\operatorname{Ext}^1_R(I/I^2,R)_{-q}$.
\end{proof}

In the applications the spaces in (ii) are computed in Macaulay2 for general forms of the stated
degrees, and the space in (iii) from the first two steps of a free resolution of $I/I^2$ over $R$.

When $Y$ has Picard number larger than one, we use the following criterion.

\begin{lemma}\label{lem:split}
Let $Y$ be a quasi-smooth threefold with $-K_Y$ ample and with finitely many points with nontrivial
stabiliser, and let $G=\mumu_p$, $p\in\{2,3\}$, act on its coordinate ring $R$ by graded automorphisms.
Suppose that $\rho(Y)\le p$, and that for a $G$-semi-invariant form $h\in R_e$ the divisor $\{h=0\}$ is
$D_1+\dots+D_p$, where $D_1,\dots,D_p$ are effective Weil divisors permuted cyclically by $G$. Let
$a=\OO_Y(1)^2\cdot D_1$ and $c=\OO_Y(1)\cdot D_1\cdot D_2$. If $ea\ne pc$, then $\rho(Y)=p$ and
$\operatorname{rk}\Cl(Y)^G=1$.
\end{lemma}

\begin{proof}
Since $G$ permutes the $D_i$ cyclically and $p\le3$, all the numbers $\OO_Y(1)\cdot D_i\cdot D_j$, $i\ne j$,
are equal to $c$. From $D_1+\dots+D_p\sim\OO_Y(e)$ we get
$\OO_Y(1)\cdot D_1^2=ea-(p-1)c\ne c=\OO_Y(1)\cdot D_1\cdot D_2$, so $D_1$ and $D_2$ are not numerically
equivalent. Hence $G$ acts non-trivially on $\Cl(Y)\otimes\QQ$. A non-trivial rational representation of
$\ZZ/p$ contains the irreducible representation of dimension $p-1$; since $\Cl(Y)\otimes\QQ$ also contains
the invariant class $\OO_Y(1)$ and has dimension $\rho(Y)\le p$, it is the sum of these two
representations. Hence $\rho(Y)=p$ and $\operatorname{rk}\Cl(Y)^G=1$.
\end{proof}

In the applications $h$ is a non-zerodivisor on the Cohen--Macaulay ring $R$, so the ideal
$J=I_Y+(h)$ has no embedded components, and the $D_i$ are given by saturations $J:s_i^\infty$ with
explicit forms $s_i$; a computation checks that $J=\bigcap_iJ:s_i^\infty$. The number $c$ is the degree
of the scheme $D_1\cap D_2$, counted with the lengths at the generic points of its components, which are the
local intersection multiplicities since $Y$ is nonsingular there.

Finally, let $G=\mumu_n$ act on $Y$ freely in codimension one, with finitely many points with nontrivial
stabiliser in $\CC^*\times G$, all of whose images in $X=Y/G$ are terminal cyclic quotient singularities, and
suppose that $-K_Y=\OO_Y(q)$ and $\operatorname{rk}\Cl(Y)^G=1$. Then $X$ is $\QQ$-factorial, so
$\rho(X)=\operatorname{rk}\Cl(Y)^G=1$, and $X$ is a $\QQ$-Fano threefold. The class $\OO_Y(1)$, with a
$G$-linearisation, descends to a Weil divisor class $A_X$ on $X$ with $-K_X\sim_\QQ qA_X$ and
$A_X^3=\frac1n\OO_Y(1)^3$, and $\Clt(X)$ contains $\ZZ/n$. By Theorem~\ref{thm:numerical} the numerical type of
$X$ is in Table~\ref{tab:main}; in each case below the rows with the basket of $X$ are listed, and the
value of $(-K_X)^3=q^3A_X^3$ determines the row, hence $\qQ(X)$, $\qW(X)$ and $\Clt(X)$. In every case we
also check that the dimensions of the eigenspaces of $G$ on the coordinate ring of $Y$ in low degrees agree
with the values $h^0(X,kA_X+\tau)$, $\tau\in\Clt(X)$, given by \eqref{eq:h0} for that row.

\subsection{The numerical type No.~34}\label{ss:34}
For No.~$34$ we have $\qQ=3$, $\qW=1$, $T_X=\ZZ/3$, $\B(X)=[[2,1],2\times[3,1],[12,5]]$ and $A_X^3=\frac14$.
The cover has $\B(Y)=[3\times[2,1],[4,1]]$ and $(\pi^*A_X)^3=\frac34$, the numerical type of a
codimension $4$ threefold in $\PP(1^2,2^3,3^2,4)$; the existence of $X$ was asked in \cite[\S5]{Pro25}.

Let $x_1,x_2,s_1,s_2,s_3,z_1,z_2,w$ be coordinates of weights $1,1,2,2,2,3,3,4$, let $\omega$ be a primitive
cube root of unity and
\[
g\colon (x_1,x_2,z_1,z_2,w,s_1,s_2,s_3)\mapsto(x_1,\omega x_2,\omega z_1,\omega^2z_2,w,s_2,s_3,s_1).
\]
Let $\ell_1$ be a general cubic form in $x_1,x_2,z_1,z_2$, put $\ell_2=g(\ell_1)$, $\ell_3=g(\ell_2)$, let
$q$ be a general $g$-invariant quartic form in $x_1,x_2,z_1,z_2,w$, and let $Y\subset\PP(1^2,2^3,3^2,4)$ be
given by the nine equations
\begin{equation}\label{eq:34}
s_is_j=x_1\ell_k,\qquad s_i\ell_i=x_1q,\qquad s_iq=\ell_j\ell_k\qquad(\{i,j,k\}=\{1,2,3\}),
\end{equation}
of degrees $4,4,4,5,5,5,6,6,6$. These equations are permuted by $g$.

\begin{proposition}\label{prop:34}
There are $\ell_1$ and $q$ with coefficients in $\QQ(\omega)$ such that the following holds.
\begin{enumerate}
\item $Y$ is a quasi-smooth threefold of codimension $4$ with $\B(Y)=[3\times[2,1],[4,1]]$,
$-K_Y=\OO_Y(3)$ and $\OO_Y(1)^3=\frac34$.
\item $\rho(Y)=3$. The divisor $\{x_1=0\}\subset Y$ is the sum of three effective divisors $D_1,D_2,D_3$,
which are permuted cyclically by $g$ and are not numerically equivalent to each other.
\item The points of $Y$ with nontrivial stabiliser in $\CC^*\times\mumu_3$ are: the point of index $4$,
fixed by $g$, with image $\frac1{12}(1,5,7)$; the three points of index $2$, permuted by $g$; and two
nonsingular fixed points of $g$.
\item $X=Y/\mumu_3$ is a Fano--Enriques threefold of type No.~$34$: $\rho(X)=1$,
$\B(X)=[[2,1],2\times[3,1],[12,5]]$, $\Clt(X)\cong\ZZ/3$, $\qQ(X)=3$, $\qW(X)=1$ and $A_X^3=\frac14$.
\end{enumerate}
\end{proposition}

\begin{proof}
(i) and (iii) are computations over $\ZZ/10009$ for a member $Y_0$, as in \S\ref{ss:verify}; the Hilbert
function of $Y_0$ is $2,6,12,22,35,54,\dots$. Let $Y$ be the lift of $Y_0$ to $\QQ(\omega)$ obtained by
lifting the coefficients of $\ell_1$ and $q$, with $\omega$ reducing to the cube root of unity used over
$\ZZ/10009$ (\texttt{examples/lift/lift34.m2}). Its ideal is $g$-invariant, and its reduced Gr\"obner basis
has $30$ elements, is integral at the prime above $10009$, and reduces to a Gr\"obner basis of $I_{Y_0}$ with
the same leading monomials. By Lemma~\ref{lem:goodred}(i),(ii), (i) and (iii) hold for $Y$.

(ii) For $Y_0$, we apply Lemma~\ref{lem:rho}(i),(ii) with parameters of degrees $1,1,2,4$. It gives
$H^1_{\mathfrak m}(\Omega_R)_0=0$ and $\dim H^2_{\mathfrak m}(\Omega_R)_0=2$, hence $h^1(Y_0,\tilde\Omega^1_{Y_0})=3$ and $\rho(Y)\le3$ by
Lemma~\ref{lem:goodred}(iii). On $\{x_1=0\}$ the equations \eqref{eq:34} give $s_is_j=0$ for $i\ne j$, and we
let $D_i$ be the part of $\{x_1=0\}$ on which $s_i\ne0$, that is, the divisor with ideal $J:s_i^\infty$, where
$J=I_Y+(x_1)$. A computation over $\QQ(\omega)$ (script \texttt{split34.m2} in \texttt{examples/lift}) shows that
$J=\bigcap_i(J:s_i^\infty)$, that $g$ permutes the $D_i$ cyclically, that $\OO(1)^2\cdot D_i=\frac14$, and that
$\OO(1)\cdot D_i\cdot D_j=\frac14$ for $i\ne j$. Thus $\OO(1)\cdot D_1^2=\frac14-2\cdot\frac14=-\frac14$, and
Lemma~\ref{lem:split} (with $e=1$, $p=3$) gives $\rho(Y)=3$ and $\operatorname{rk}\Cl(Y)^{g}=1$.

(iv) By (iii) and (ii) we can apply the last paragraph of \S\ref{ss:verify}, with $A_X$ the image of
$\{x_1=0\}$: the only row of Table~\ref{tab:main} with this basket is No.~$34$. The dimensions of the
eigenspaces of $g$ on the coordinate ring of $Y$ in degrees $k\le 6$ agree with the values
$h^0(X,kA_X+m\sigma)$ of No.~$34$.
\end{proof}

This example shows that the cover $Y$ of a Fano--Enriques threefold need not have Picard number one.

\subsection{Type I unprojections with an involution: Nos.~5, 21, 24, 27}\label{ss:tom}
The four numerical types of this subsection have $T_X=\ZZ/2$:
\begin{center}\footnotesize
\renewcommand{\arraystretch}{1.25}
\begin{tabular}{c|c|c|c|c|c}
No. & $\qQ=\qW$ & $\B(X)$ & $A_X^3$ & $\B(Y)$ & $(\pi^*A_X)^3$\\\hline
$5$ & $5$ & $[[2,1],[4,1],[14,3]]$ & $\frac1{28}$ & $[2\times[4,1],[7,3]]$ & $\frac1{14}$\\
$21$ & $3$ & $[4\times[2,1],[10,3]]$ & $\frac3{10}$ & $[2\times[2,1],[5,2]]$ & $\frac35$\\
$24$ & $3$ & $[2\times[2,1],[4,1],[14,3]]$ & $\frac1{28}$ & $[2\times[2,1],2\times[4,1],[7,3]]$ & $\frac1{14}$\\
$27$ & $3$ & $[3\times[2,1],[4,1],[8,3]]$ & $\frac38$ & $[3\times[2,1],[4,1]]$ & $\frac34$
\end{tabular}
\end{center}
In each case the Hilbert series of the cover is that of a codimension $4$ threefold with nine equations, in the
weighted projective space given in Table~\ref{tab:tom}. The existence of No.~$5$ was left open in
\cite[\S5]{Pro25}. For No.~$21$, a Fano threefold with the invariants of $Y$ is constructed in \cite{Heu22} as a
complete intersection in a toric variety of Picard rank $2$; its ring of sections of multiples of $\pi^*A_X$
shows that it is of the codimension $4$ type above. The cover for No.~$27$ has the same numerical type as the
cover for No.~$34$. We obtain $Y$ by a Type~I (Tom) unprojection which is compatible with an involution.

Let $\PP=\PP(a_1,\dots,a_7)$ with the coordinates and weights of Table~\ref{tab:tom}, and let $\iota$ be the
involution of $\PP$ which changes the sign of the coordinates listed in the column ``$-1$'' and fixes the other
ones. Let $D\subset\PP$ be the plane with the ideal $I_D$ given in the table, and let $M=(m_{ij})$ be a
skew-symmetric $5\times5$ matrix whose entries $m_{ij}$, $i<j$, are general forms of the degrees
\[
d_{\rm I}=\begin{pmatrix}6&6&5&5\\&5&4&4\\&&4&4\\&&&3\end{pmatrix},\qquad
d_{\rm II}=\begin{pmatrix}4&4&3&3\\&3&2&2\\&&2&2\\&&&1\end{pmatrix}
\]
as indicated in the table, such that $\iota(m_{ij})=(-1)^{c_i+c_j+\varepsilon}m_{ij}$, and such that
$m_{ij}\in I_D$ for $i,j\ne k$ (the format $\mathrm{Tom}_k$ of \cite{BKR12}). Let $Z\subset\PP$ be given by the
five Pfaffians of $M$. Then $Z$ contains $D$ and is $\iota$-invariant. Let $Y\subset\PP(a_1,\dots,a_7,a_s)$ be
the Kustin--Miller unprojection of $D\subset Z$ \cite{PR04}, with new variable $s$ of weight $a_s$.

\begin{table}[ht]
\centering\scriptsize\setlength{\tabcolsep}{2.5pt}
\renewcommand{\arraystretch}{1.3}
\caption{The data of the Type~I unprojections of \S\ref{ss:tom}.}\label{tab:tom}
\begin{tabular}{c|l|l|c|c|c|c|c}
No. & coordinates (weights) & $-1$ & $I_D$ & degrees & $c$ & $\varepsilon$, $k$ & $a_s$\\\hline
$5$ & \begin{tabular}[t]{@{}l@{}}$x,y,z,z',w,w',v$\\ $(1,3,4,4,5,5,6)$\end{tabular} & $x,z',w',v$ & $(y,z,z',w)$ & $d_{\rm I}$ & $(0,1,0,1,0)$ & $1$, $1$ & $7$\\
$21$ & \begin{tabular}[t]{@{}l@{}}$x_1,x_2,u,v,z_1,z_2,w$\\ $(1,1,2,2,3,3,4)$\end{tabular} & $x_2,v,z_2,w$ & $(x_2,u,v,z_2)$ & $d_{\rm II}$ & $(0,0,1,0,1)$ & $0$, $1$ & $5$\\
$24$ & \begin{tabular}[t]{@{}l@{}}$y,z,z',u,u',w,w'$\\ $(2,3,3,4,4,5,5)$\end{tabular} & $y,z',u',w'$ & $(z',u,u',w')$ & $d_{\rm I}$ & $(0,0,1,1,0)$ & $0$, $1$ & $7$\\
$27$ & \begin{tabular}[t]{@{}l@{}}$x_1,x_2,u,v,z_1,z_2,w$\\ $(1,1,2,2,3,3,4)$\end{tabular} & $x_2,v,z_2,w$ & $(u,v,z_2,w)$ & $d_{\rm II}$ & $(0,0,1,0,1)$ & $0$, $4$ & $2$
\end{tabular}
\end{table}

\begin{table}[ht]
\centering\footnotesize\setlength{\tabcolsep}{4pt}
\renewcommand{\arraystretch}{1.3}
\caption{The threefolds $Y$ of Proposition~\ref{prop:tom}. ``Fixed'' lists the singular points of $Y$ fixed
by $\iota$ with their images in $X$, ``pairs'' the indices of the pairs of points interchanged by $\iota$, and
$n_0$ the number of nonsingular fixed points.}\label{tab:tomres}
\begin{tabular}{c|c|c|c|l|c|c|c}
No. & degrees of equations & $q$ & $\iota(s)$ & fixed & pairs & $n_0$ & $\rho(Y)$\\\hline
$5$ & $8,9,9,10,10,10,11,11,12$ & $5$ & $s$ & index $7$: $\frac1{14}(1,3,11)$ & $4$ & $1$ & $2$\\
$21$ & $4,5,5,6,6,6,7,7,8$ & $3$ & $-s$ & index $5$: $\frac1{10}(1,3,9)$ & $2$ & $3$ & $2$\\
$24$ & $8,9,9,10,10,10,11,11,12$ & $3$ & $-s$ & index $7$: $\frac1{14}(3,5,9)$ & $4$, $2$ & $1$ & $2$\\
$27$ & $4,4,4,5,5,5,6,6,6$ & $3$ & $-s$ & \begin{tabular}[t]{@{}l@{}}index $4$: $\frac18(1,3,5)$;\\ index $2$, $s\ne0$: $\frac14(1,1,3)$\end{tabular} & $2$ & $2$ & $1$
\end{tabular}
\end{table}

\begin{proposition}\label{prop:tom}
For each row of Table~\ref{tab:tom} there is a matrix $M$ as above, with coefficients in $\QQ$, such that
the following holds.
\begin{enumerate}
\item $Y$ is a quasi-smooth threefold of codimension $4$ with nine equations of the degrees given in
Table~\ref{tab:tomres}, with the basket $\B(Y)$ given above and $-K_Y=\OO_Y(q)$.
\item The involution $\iota$ extends to $Y$, with $\iota(s)$ as in Table~\ref{tab:tomres}. The points of $Y$
with nontrivial stabiliser in $\CC^*\times\mumu_2$ are the points listed in Table~\ref{tab:tomres}.
\item $\rho(Y)$ is as in Table~\ref{tab:tomres}. For Nos.~$5$, $21$, $24$, the divisor $\{h=0\}\subset Y$,
where $h=y$, $x_2$, $z'$ respectively, is the sum of two effective divisors $D'$, $D''$, which are
interchanged by $\iota$ and are not numerically equivalent.
\item $X=Y/\iota$ is a Fano--Enriques threefold of the given type, with $\rho(X)=1$, $\Clt(X)\cong\ZZ/2$
and the invariants $\B(X)$, $\qQ(X)=\qW(X)$ and $A_X^3$ given above.
\end{enumerate}
\end{proposition}

\begin{proof}
(i) and (ii) are computations over $\ZZ/32003$ for a member $Y_0$, as in \S\ref{ss:verify}. For a suitable
linearisation of $\iota$ the dimensions of the invariant and anti-invariant parts of the coordinate ring of
$Y_0$ in degrees $k\le12$ for No.~$5$, and $k\le10$ in the other cases, are equal to the values $h^0(X,kA_X)$
and $h^0(X,kA_X+\sigma)$ of the corresponding row. Let $Y$ be the threefold over $\QQ$ obtained from $Y_0$ by
lifting the coefficients of the entries of $M$, taken in bases of the relevant eigenspaces of $I_D$ defined
over $\ZZ$, and unprojecting over $\QQ$ (\texttt{examples/lift/l5.m2}, \texttt{l21.m2}, \texttt{l24.m2},
\texttt{l27.m2}). Its ideal is invariant under $\iota$ with $\iota(s)$ as in the table; its reduced Gr\"obner
basis has $32$, $43$, $17$, $30$ elements for Nos.~$5$, $21$, $24$, $27$, is $p$-integral, and reduces to a
Gr\"obner basis, with the same leading monomials, of an ideal which is mapped to $I_{Y_0}$ by
$s\mapsto\lambda s$, with $\lambda=8816$, $3256$, $8893$, $15584$ respectively, an automorphism commuting with
$\iota$. By Lemma~\ref{lem:goodred}(i),(ii) the statements verified for $Y_0$ hold for $Y$.

(iii) Let $\delta=\dim H^2_{\mathfrak m}(\Omega_R)_0$ for the coordinate ring $R$ of $Y_0$. For Nos.~$5$ and
$27$, Lemma~\ref{lem:rho}(ii), with parameters of degrees $1,3,4,7$ and $1,1,2,4$ respectively, gives
$H^1_{\mathfrak m}(\Omega_R)_0=0$ and $\delta=1$, $\delta=0$. For Nos.~$21$ and $24$, Lemma~\ref{lem:rho}(iii)
gives $\delta=\dim\operatorname{Ext}^1_R(I/I^2,R)_{-3}=1$; the same computation gives $1$ for No.~$5$ and $0$
for No.~$27$, in agreement with the previous values. By Lemma~\ref{lem:rho}(i) and
Lemma~\ref{lem:goodred}(iii), $\rho(Y)\le1+\delta$. This proves (iii) for No.~$27$.

For Nos.~$5$, $21$, $24$ the format is $\mathrm{Tom}_1$, and the entry $m_{45}$, of degree $3$, $1$, $3$,
lies in $I_D$ and has the sign of $h$, so it is a multiple of $h$. Hence on $\{h=0\}$ the Pfaffian
$m_{23}m_{45}-m_{24}m_{35}+m_{25}m_{34}$ becomes a binary quadratic form in the pair of coordinates
$(a,b)=(z,z')$, $(u,v)$, $(u,u')$ respectively, which by the choice of signs is of the form
$\alpha a^2+\beta b^2$. Let $r$ be a square root of $-\beta/\alpha$, so that $\{h=0\}$ lies over $\{a=rb\}\cup\{a=-rb\}$, and let $J=I_Y+(h)$. Let $D'$ and $D''$ be the divisors with ideals
$J:(a+rb)^\infty$ and $J:(a-rb)^\infty$ over $\QQ(r)$; they are interchanged by $\iota$. A computation over
$\QQ(r)$ (\texttt{examples/lift/split.m2}, with the scripts above) shows that $J$ is the intersection of these two ideals, and gives
\[
\bigl(\OO(1)^2\cdot D',\ \OO(1)\cdot D'\cdot D''\bigr)=\bigl(\tfrac3{28},\tfrac27\bigr),\ \bigl(\tfrac3{10},\tfrac25\bigr),\ \bigl(\tfrac3{28},\tfrac27\bigr)
\]
for Nos.~$5$, $21$, $24$. Here $D'+D''\sim\OO_Y(e)$ with $e=\deg h=3$, $1$, $3$, so
$\OO(1)\cdot D'^2=\frac1{28}$, $-\frac1{10}$, $\frac1{28}$, which differs from $\OO(1)\cdot D'\cdot D''$.
Lemma~\ref{lem:split} (with $p=2$) gives $\rho(Y)=2$ and $\operatorname{rk}\Cl(Y)^\iota=1$.

(iv) follows from (ii) and (iii) by the last paragraph of \S\ref{ss:verify}; in each case the only row of
Table~\ref{tab:main} with this basket is the given one.
\end{proof}

Thus the Hilbert series of a codimension $4$ threefold with $\B=[3\times[2,1],[4,1]]$ and $A_X^3=\frac34$
is realised both by quasi-smooth threefolds with $\rho=1$ (the $\mathrm{Tom}_4$ unprojection for No.~$27$)
and with $\rho=3$ (the format \eqref{eq:34}).

\subsection{Unprojection of several divisors: Nos.~29, 43, 137, 143, 148}\label{ss:multi}
In the following five cases the cover $Y$ is obtained from a threefold $Z$ containing two, three or four surfaces
$D_i$ which are invariant under, or permuted by, the group $G=\Clt(X)$, by unprojecting all of them; the
new variables are then permuted by $G$ up to scalars. The invariants are:
\begin{center}\footnotesize
\renewcommand{\arraystretch}{1.25}
\begin{tabular}{c|c|c|c|c|c|c}
No. & $T_X$ & $\qQ$, $\qW$ & $\B(X)$ & $A_X^3$ & $\B(Y)$ & $(\pi^*A_X)^3$\\\hline
$29$ & $\ZZ/2$ & $3$, $3$ & $[4\times[2,1],2\times[4,1],[5,2]]$ & $\frac1{10}$ & $[2\times[2,1],2\times[5,2]]$ & $\frac15$\\
$43$ & $\ZZ/2$ & $2$, $1$ & $[2\times[2,1],2\times[6,1]]$ & $\frac43$ & $[2\times[3,1]]$ & $\frac83$\\
$137$ & $\ZZ/3$ & $2$, $2$ & $[2\times[9,2]]$ & $\frac59$ & $[2\times[3,1]]$ & $\frac53$\\
$143$ & $\ZZ/3$ & $2$, $2$ & $[7\times[3,1]]$ & $\frac13$ & $[3\times[3,1]]$ & $1$\\
$148$ & $\ZZ/3$ & $2$, $2$ & $[[3,1],2\times[9,4]]$ & $\frac29$ & $[5\times[3,1]]$ & $\frac23$
\end{tabular}
\end{center}

\medskip\noindent\textbf{No.~$29$.} The Hilbert series of $Y$ is that of a codimension $4$ threefold in
$\PP(1,2,3^2,4^2,5^2)$, and the two points of index $5$ of $Y$ are interchanged by the involution. We
unproject two planes simultaneously.

Let $R=\CC[x_1,\dots,x_6]$ be the coordinate ring of $\PP(1,2,3,3,4,4)$ and let
\[
\iota\colon(x_1,x_2,x_3,x_4,x_5,x_6)\mapsto(-x_1,-x_2,x_3,-x_4,x_5,-x_6).
\]
Consider the planes $D_1,D_2\cong\PP(1,3,4)$ with ideals
\begin{align*}
I_1&=(g_1,g_2,g_3)=(x_2-x_1^2,\ x_4-x_3,\ x_6-x_5),\\
I_2&=\iota(I_1)=(g_1',g_2',g_3')=(x_2+x_1^2,\ x_4+x_3,\ x_6+x_5).
\end{align*}
They are disjoint and interchanged by $\iota$. The ideal $I_1\cap I_2$ is generated by
\begin{gather*}
x_1^4-x_2^2,\quad x_2x_3-x_1^2x_4,\quad x_1^2x_3-x_2x_4,\quad x_2x_5-x_1^2x_6,\quad x_1^2x_5-x_2x_6,\\
x_3^2-x_4^2,\quad x_4x_5-x_3x_6,\quad x_3x_5-x_4x_6,\quad x_5^2-x_6^2 .
\end{gather*}
Its $\iota$-invariant part has dimension $4$ in degree $6$, where it is spanned by
\[
x_1^2(x_1^4-x_2^2),\qquad x_1(x_2x_3-x_1^2x_4),\qquad x_3^2-x_4^2,\qquad x_1^2x_5-x_2x_6,
\]
and dimension $10$ in degree $8$. Let $f_6$, $f_8$ be general $\iota$-invariant elements of
$I_1\cap I_2$ of degrees $6$ and $8$, and let
\[
Z=\{f_6=f_8=0\}\subset\PP(1,2,3,3,4,4).
\]
Write $(f_6,f_8)=(g_1,g_2,g_3)M=(g_1',g_2',g_3')M'$ with $3\times2$ matrices $M$, $M'=\iota(M)$, and let
$\Delta_i$ (resp.\ $\Delta_i'$) be the $2\times2$ minor of $M$ (resp.\ $M'$) obtained by deleting the
$i$-th row. Let $s$, $t$ be new variables of degree $5$ and
\begin{equation}\label{eq:IY}
I_Y=\bigl(f_6,\ f_8,\ s\,g_i-(-1)^{i-1}\Delta_i,\ t\,g_i'-(-1)^{i-1}\Delta_i'\ \ (i=1,2,3)\bigr):(g_1g_1')^{\infty}
\end{equation}
in $\CC[x_1,\dots,x_6,s,t]$.
The equations $s\,g_i=(-1)^{i-1}\Delta_i$ are those of the Kustin--Miller unprojection of $D_1\subset Z$
\cite{PR04,Pap04}, and similarly for $t$ and $D_2$.
The involution $\iota$ extends to $Y$ by $s\leftrightarrow t$.

\medskip\noindent\textbf{No.~$43$.} Both points of index $3$ of $Y$ are fixed by the involution. We unproject two
invariant surfaces from a Pfaffian threefold.

Let $x_0,\dots,x_3,y_0,y_1,y_2$ be coordinates of weights $1,1,1,1,2,2,2$ and
\[
\iota\colon(x_0,x_1,x_2,x_3,y_0,y_1,y_2)\mapsto(x_0,x_1,-x_2,-x_3,-y_0,-y_1,-y_2).
\]
Let $D_1=\{\ell_1=x_2=x_3=m_1=0\}$ and $D_2=\{x_0=x_1=\ell_2=m_2=0\}$, where $\ell_1$ is a general linear
form in $x_0,x_1$, $\ell_2$ a general linear form in $x_2,x_3$ and $m_1,m_2$ are general anti-invariant
forms of degree $2$; these are $\iota$-invariant surfaces isomorphic to $\PP(1,2,2)$. Let $M=(m_{ij})$ be a
skew-symmetric $5\times5$ matrix whose entries are general forms of degree $2$ for $i<j\le4$ and of degree
$1$ for $j=5$ such that $\iota(m_{ij})=(-1)^{c_i+c_j}m_{ij}$ with $c=(0,1,0,1,0)$, $m_{ij}\in I_{D_1}$ for
$i,j\ne1$ and $m_{ij}\in I_{D_2}$ for $i,j\ne2$ (the format $\mathrm{Tom}_1$ for $D_1$ and $\mathrm{Tom}_2$
for $D_2$). Let $Z\subset\PP(1^4,2^3)$ be given by the Pfaffians of $M$, of degrees $3,3,3,3,4$, and let
$Y\subset\PP(1^4,2^3,3^2)$ be the unprojection of $D_1$ and $D_2$ in $Z$, with new variables $u_1,u_2$ of
weight $3$.
The involution $\iota$ extends to $Y$ with $\iota(u_1)=-u_1$, $\iota(u_2)=u_2$.

\medskip\noindent\textbf{No.~$137$.} Both points of index $3$ of $Y$ are fixed by the group. We unproject two invariant
surfaces from a complete intersection.

Let $x_0,x_1,x_2,y_0,y_1,y_2$ be coordinates of weights $1,1,1,2,2,2$, let $\omega$ be a primitive cube
root of unity and
\[
g\colon(x_0,x_1,x_2,y_0,y_1,y_2)\mapsto(x_0,\omega x_1,\omega^2x_2,y_0,\omega y_1,\omega^2y_2).
\]
Let $D_1=\{x_0=x_1=y_1+ax_2^2=0\}$ and $D_2=\{x_1=x_2=y_0+bx_0^2=0\}$ with general $a,b$; these are
$g$-invariant surfaces isomorphic to $\PP(1,2,2)$. Let $Z=\{f_3=f_4=0\}\subset\PP(1^3,2^3)$, where $f_3$ is
a general $g$-invariant cubic and $f_4$ a general quartic with $g(f_4)=\omega f_4$, both vanishing on
$D_1\cup D_2$. Let $Y\subset\PP(1^3,2^3,3^2)$ be the unprojection of $D_1$ and $D_2$ in $Z$, with new
variables $u_1,u_2$ of weight $3$.
The automorphism $g$ extends to $Y$ with $g(u_1)=\omega^2u_1$, $g(u_2)=\omega u_2$.

\medskip\noindent\textbf{No.~$143$.} The Hilbert series of $Y$ is that of a codimension $4$ threefold in
$\PP(1^2,2^3,3^3)$. We unproject three surfaces which are permuted by the group.

Let $x_1,x_2,a_1,a_2,a_3$ be coordinates of weights $1,1,2,2,2$, let $\omega$ be a primitive cube root of
unity and
\[
g\colon(x_1,x_2,a_1,a_2,a_3)\mapsto(x_1,\omega^2x_2,a_2,a_3,a_1).
\]
Let $\ell_1$ be a general linear form in $x_1,x_2$ and $m_1$ a general form of degree $2$, put
$\ell_{i+1}=g(\ell_i)$, $m_{i+1}=g(m_i)$ and $D_i=\{\ell_i=m_i=0\}\cong\PP(1,2,2)$. Let
$Z=\{f=0\}\subset\PP(1^2,2^3)$ be a general $g$-invariant sextic containing $D_1$, $D_2$ and $D_3$, and
write $f=\ell_iA_i+m_iB_i$. Let $Y\subset\PP(1^2,2^3,3^3)$ be the unprojection of $D_1,D_2,D_3\subset Z$,
with new variables $u_i=A_i/m_i=-B_i/\ell_i$ of weight $3$; the automorphism $g$ extends to $Y$ by
$u_1\mapsto u_2\mapsto u_3\mapsto u_1$.

\medskip\noindent\textbf{No.~$148$.} The Hilbert series of $Y$, with $h^0(Y,\OO_Y(k))=1,4,9,15,26,41$ for
$k\le6$, is that of four Type~I unprojections of surfaces $\PP(1,2,2)$ in a hypersurface of degree $8$ in
$\PP(1,2,2,2,3)$, the fifth point of index $3$ being the coordinate point of the variable of weight $3$. Two of
the five points of index $3$ are fixed by the group and three are permuted, so we unproject one invariant
surface and an orbit of three.

Let $x,y_1,y_2,y_3,u$ be coordinates of weights $1,2,2,2,3$, let $\omega$ be a primitive cube root of unity and
\[
g\colon(x,y_1,y_2,y_3,u)\mapsto(\omega x,y_1,\omega y_2,\omega^2y_3,\omega u).
\]
Let $D_0=\{y_1=u=0\}$. Let $\ell$ and $m$ be general linear forms in $y_1,y_2,y_3$ and $a$ a general constant,
let $D_1=\{\ell+ax^2=0,\ u=xm\}$, and put $D_2=g(D_1)$, $D_3=g(D_2)$. These are surfaces isomorphic to
$\PP(1,2,2)$; $D_0$ is $g$-invariant and $D_1,D_2,D_3$ are permuted by $g$. Let
$Z=\{f=0\}\subset\PP(1,2^3,3)$ be a general $g$-invariant octic containing $D_0,\dots,D_3$, and let
$Y\subset\PP(1,2^3,3^5)$ be the unprojection of $D_0,\dots,D_3$ in $Z$, with new variables $u_0,\dots,u_3$
of weight $3$. The automorphism $g$ extends to $Y$ by $u_0\mapsto\omega^2u_0$ and $u_1\mapsto u_2\mapsto
u_3\mapsto u_1$.

\begin{table}[ht]
\centering\footnotesize\setlength{\tabcolsep}{2.5pt}
\renewcommand{\arraystretch}{1.3}
\caption{The threefolds $Y$ of Proposition~\ref{prop:multi}. The last column lists the points of $Y$ with
nontrivial stabiliser in $\CC^*\times G$, with the images in $X$ of the fixed singular points; $n_0$ is the
number of nonsingular fixed points of a generator of $G$.}\label{tab:multi}
\begin{tabular}{c|l|l|c|l}
No. & $Y\subset$ & degrees of equations & $\OO_Y(1)^3$ & points with nontrivial stabiliser\\\hline
$29$ & $\PP(1,2,3^2,4^2,5^2)$ & $6,7,7,8,8,8,9,9,10$ & $\frac15$ & \begin{tabular}[t]{@{}l@{}}index $2$, fixed: $2\times\frac14(1,3,3)$;\\ index $5$: one pair; $n_0=4$\end{tabular}\\
$43$ & $\PP(1^4,2^3,3^2)$ & $3^4,4^7,5^2,6$ & $\frac83$ & \begin{tabular}[t]{@{}l@{}}index $3$, fixed: $\frac16(1,5,5)$, $\frac16(1,1,5)$;\\ $n_0=2$\end{tabular}\\
$137$ & $\PP(1^3,2^3,3^2)$ & $3,4,4,4,4,4,5,5,6$ & $\frac53$ & \begin{tabular}[t]{@{}l@{}}index $3$, fixed: $\frac19(2,7,8)$, $\frac19(1,5,8)$;\\ $n_0=0$\end{tabular}\\
$143$ & $\PP(1^2,2^3,3^3)$ & $4,4,4,5,5,5,6,6,6$ & $1$ & \begin{tabular}[t]{@{}l@{}}index $3$: one orbit of three points;\\ $n_0=6$, images $\frac13(1,1,2)$\end{tabular}\\
$148$ & $\PP(1,2^3,3^5)$ & $5^4,6^{10}$ & $\frac23$ & \begin{tabular}[t]{@{}l@{}}index $3$, fixed: $\frac19(2,5,7)$, $\frac19(4,5,8)$;\\ one orbit of three points;\\ $n_0=0$\end{tabular}

\end{tabular}
\end{table}

\begin{proposition}\label{prop:multi}
For each of Nos.~$29$, $43$, $137$, $143$, $148$ there is a choice as above, with coefficients in $\QQ$
(Nos.~$29$, $43$) or $\QQ(\omega)$ (Nos.~$137$, $143$, $148$), such that the following holds.
\begin{enumerate}
\item $Y$ is a quasi-smooth threefold, of codimension $5$ for Nos.~$43$, $148$ and $4$ otherwise, with equations of
the degrees given in Table~\ref{tab:multi}, with the basket $\B(Y)$ given above, $-K_Y=\OO_Y(\qQ)$, the value
of $\OO_Y(1)^3$ given in the table, and $\rho(Y)=1$.
\item The group $G$ acts on $Y$ as described above, and the points of $Y$ with nontrivial stabiliser in
$\CC^*\times G$ are those of Table~\ref{tab:multi}.
\item $X=Y/G$ is a Fano--Enriques threefold of the given type, with $\rho(X)=1$, $\Clt(X)\cong G$ and the
invariants $\B(X)$, $\qQ(X)$, $\qW(X)$ and $A_X^3$ given above.
\end{enumerate}
\end{proposition}

\begin{proof}
(i) and (ii), apart from $\rho(Y)$, are computations over $\FF_p$ for a member $Y_0$, as in
\S\ref{ss:verify}, with $p=32003$ for Nos.~$29$, $43$ and $p=10009$ for Nos.~$137$, $143$, $148$. The Hilbert
functions of $Y_0$ begin with $4,13,30,57,98,155$ for No.~$43$, with $3,9,20,37,63,99$ for No.~$137$ and with
$2,6,13,23,39,61,89,126$ for No.~$143$ and with $1,4,9,15,26,41,59,84$ for No.~$148$. For Nos.~$43$ and $148$
quasi-smoothness is checked in two steps: $Z$ is quasi-smooth outside the union of the surfaces $D_i$, and $Y$
is isomorphic to $Z$ outside the preimages of the $D_i$; these preimages are curves, and the Jacobian matrix of
$Y$ has rank $5$ at each of their points. For No.~$148$ this rank is checked with random $5\times5$ minors on
each component of the preimages, which leave only the coordinate point of $u_0$, where the rank is computed
directly (\texttt{examples/no148/qs148.log}). The
dimensions of the eigenspaces of $G$ on the coordinate ring of $Y_0$ in degrees $k\le10$ (Nos.~$29$, $143$, $148$)
and $k\le8$ (Nos.~$43$, $137$) are equal to the values $h^0(X,kA_X+\tau)$, $\tau\in\Clt(X)$, of the
corresponding row.

For No.~$29$, Lemma~\ref{lem:rho}(i),(ii), with parameters of degrees $1,2,3,5$, gives
$H^1_{\mathfrak m}(\Omega_R)_0=H^2_{\mathfrak m}(\Omega_R)_0=0$; for Nos.~$43$, $137$, $143$, $148$,
Lemma~\ref{lem:rho}(iii) gives $\operatorname{Ext}^1_R(I/I^2,R)_{-2}=0$. Hence
$h^1(Y_0,\tilde\Omega^1_{Y_0})\le1$ in all five cases.

The lift $Y$ is obtained as follows.
\begin{itemize}
\item No.~$29$. The pair $(f_6,f_8)$ of $Y_0$ is a combination, with coefficients in $\FF_p$, of the spanning
sets of the $\iota$-invariant parts of $I_1\cap I_2$, which are defined over $\ZZ$. Lifting these
coefficients to integers gives $\iota$-invariant elements $\tilde f_6,\tilde f_8$ of $I_1\cap I_2$ over $\QQ$
reducing to $f_6,f_8$, and \eqref{eq:IY} defines $Y$ over $\QQ$ (script \texttt{lift29.m2}, lifted ideal \texttt{IY29\_QQ.txt}, both in
\texttt{examples/no29}). The ideal $I_Y$ has generators of the same degrees as $I_{Y_0}$; its
reduced Gr\"obner basis (degree reverse lexicographic order) has $33$ elements, all of them $p$-integral,
and their reductions form a Gr\"obner basis of $I_{Y_0}$ with the same leading monomials.
\item No.~$43$. We lift the coefficients of $\ell_1,\ell_2,m_1,m_2$ and of the entries of $M$, the latter
taken in bases of the relevant eigenspaces of $I_{D_1}\cap I_{D_2}$, $I_{D_1}$ and $I_{D_2}$ defined over $\ZZ$,
and unproject over $\QQ$ (\texttt{examples/lift/lift43.m2}). The reduced Gr\"obner basis has $71$ elements, is
$p$-integral, and reduces to a Gr\"obner basis, with the same leading monomials, of an ideal which is mapped
to $I_{Y_0}$ by $u_1\mapsto 12747\,u_1$, $u_2\mapsto 8542\,u_2$ (\texttt{examples/lift/iso43.m2}), an
automorphism commuting with $\iota$.
\item No.~$137$. We lift the coefficients of $a$, $b$, and of $f_3$, $f_4$ in bases of the eigenspaces of
forms vanishing on $D_1\cup D_2$, and unproject over $\QQ(\omega)$ (\texttt{examples/lift/lift137.m2}). The
reduced Gr\"obner basis has $26$ elements, is integral at the prime above $10009$, and reduces to a
Gr\"obner basis of $I_{Y_0}$ with the same leading monomials.
\item No.~$143$. We lift the coefficients of $\ell_1$ and $m_1$, and of $f$ in a basis of the forms of
degree $6$ vanishing on $D_1\cup D_2\cup D_3$, and unproject over $\QQ(\omega)$
(\texttt{examples/lift/lift143.m2}). The reduced Gr\"obner basis has $14$ elements, is integral at the prime
above $10009$, and reduces to a Gr\"obner basis of $I_{Y_0}$ with the same leading monomials.
\item No.~$148$. We lift the coefficients of $\ell$, $m$ and $a$, and of $f$ in a basis of the forms of
degree $8$ vanishing on $D_0\cup\dots\cup D_3$, and unproject over $\QQ(\omega)$
(\texttt{examples/lift/lift148.m2}). The reduced Gr\"obner basis has $54$ elements, is integral at the prime
above $10009$, and reduces to a Gr\"obner basis of $I_{Y_0}$ with the same leading monomials.
\end{itemize}
In each case the ideal of $Y$ is invariant under $G$, acting as described above. By
Lemma~\ref{lem:goodred}, (i) and (ii) hold for $Y$, and $1\le\rho(Y)\le h^1(Y_0,\tilde\Omega^1_{Y_0})\le1$.

(iii) follows from (i) and (ii) by the last paragraph of \S\ref{ss:verify}. The rows of
Table~\ref{tab:main} with the basket of $X$ are: No.~$29$ only; Nos.~$42$--$45$; Nos.~$137$, $138$;
Nos.~$143$, $144$; and Nos.~$148$, $149$. In each case $(-K_X)^3$ determines the row.
\end{proof}

For No.~$29$ the computation also shows the following. The threefold $Z$ has exactly $8$ singular points
apart from its quotient singularities; they are ordinary double points, four on each of the planes $D_1$,
$D_2$ (the scheme defined by $f_6,f_8$ and the $2\times2$ minors of their Jacobian matrix is reduced of length
$8$). The ideal $I_Y$ is generated by $9$ elements, and in the coordinates $u=s+t$, $v=s-t$ the involution
is the action $\mumu_2(1,1,0,1,0,1,0,1)$, under which $I_Y$ is generated by six invariant equations of
degrees $6,7,8,8,9,10$ and three anti-invariant equations of degrees $7,8,9$. The two points of index $5$
are given by $u^2=v^2$ on the line $\PP(5,5)$. The degree is
$A_X^3=\frac12\bigl(\frac{6\cdot8}{1\cdot2\cdot3\cdot3\cdot4\cdot4}+\frac{2}{1\cdot3\cdot4\cdot5}\bigr)=\frac1{10}$.

\begin{remark}
Let $f\colon\tilde Y\to Y$ be the weighted blow-ups inverse to the two unprojections, at the two points of
index $5$ of the threefold $Y$ of No.~$29$; in local coordinates of type $\frac15(1,3,4)$ they have weights
$\frac15(1,3,4)$ and discrepancy $\frac35$, so that $-K_{\tilde Y}=3\bigl(f^*\OO_Y(1)-\frac15(E_1+E_2)\bigr)$ for
the exceptional divisors $E_1$, $E_2$. The anticanonical model of $\tilde Y$ is $Z$ and $\tilde Y\to Z$ is small, so
$\operatorname{rk}\Cl(Z)=\rho(Y)+2=3$: the class group of the nodal complete intersection $Z$ is generated,
up to finite index, by $\OO_Z(1)$, $D_1$ and $D_2$, and its defect is $2$. This agrees with the count of
conditions: the $8$ nodes of $Z$ impose only $6$ conditions on the forms of degree $5=6+2\cdot8-17$.
\end{remark}

\subsection{A $\mathrm{G}_2$ cluster format: No.~17}\label{ss:17}
For No.~$17$ we have $\qQ=\qW=3$, $T_X=\ZZ/2$, $\B(X)=[[2,1],[14,3]]$ and $A_X^3=\frac27$. The cover has
$\B(Y)=[7,3]$, $-K_Y=3B$ and $B^3=\frac47$: the numerical type of the candidate \#41199 of \cite{GRDB}, a
codimension $4$ threefold in $\PP(1,1,2,2,3,3,5,7)$. Coughlan and Ducat construct a family of such threefolds in
their $\mathrm{G}_2^{(4)}$ cluster format \cite[\S4.2.2 and Table]{CD20}; we find an involution on a member
of it.

Let $\mathcal A$ be the coordinate ring of the $\mathrm{G}_2$ cluster variety of \cite[\S1.2.2]{CD20}, with
cluster variables $\theta_1,\theta_{12},\dots,\theta_4,\theta_{41}$,
coefficients $A_1,A_{12},\dots,A_{41}$ and parameters $\lambda_{13},\lambda_{24}$, and with the action of the
torus $\mathbb T^8$ of \cite[Table~3]{CD20}. The one-parameter subgroup with
$d(\theta_1,\theta_{12},\dots,\theta_{41})=(4,5,1,1,3,8,5,7)$ gives the degrees
\[
d(A_1,A_{12},\dots,A_4,A_{41})=(0,0,3,3,0,0,0,2),\qquad d(\lambda_{13},\lambda_{24})=(1,2),
\]
the entry of \cite[Table]{CD20} for \#41199. Let $\epsilon\in\mathbb T^8$ be the element of order $2$ which is
$-1$ on the characters $\chi_{23}$, $\chi_{34}$, $\chi_4$ and $1$ on the others; it acts by $-1$ on
$\theta_{23},\theta_{34},\theta_4,A_2,A_{23},A_{41},\lambda_{13}$ and by $1$ on the other generators, in
particular on all generators of degree $0$. Let $x_1,x_2,y_1,y_2,z_1,z_2,v,w$ be coordinates of weights
$1,1,2,2,3,3,5,7$ and let
\[
\iota\colon(x_1,x_2,y_1,y_2,z_1,z_2,v,w)\mapsto(x_1,-x_2,-y_1,-y_2,z_1,-z_2,-v,w).
\]
Let $\phi$ be the morphism with $\phi^*(A_{12})=\phi^*(A_{34})=1$, $\phi^*(A_1),\phi^*(A_3),\phi^*(A_4)$
general constants, $\phi^*(\theta_{12})$ and $\phi^*(\theta_{34})$ given by the tag equations
$\theta_{12}=\theta_1\theta_2-A_{23}A_3A_4A_{41}$ and $\theta_{34}=\theta_3\theta_4-A_{41}A_1A_2A_{23}$, and every
other generator mapped to a general form of its degree on which $\iota$ acts as $\epsilon$ does on the
generator. Let $Y\subset\PP(1,1,2,2,3,3,5,7)$ be the threefold defined by $\phi$, that is, by the pullbacks of
the relations of $\mathcal A$; it is $\iota$-invariant.

\begin{proposition}\label{prop:17}
There is a choice as above, with coefficients in $\QQ$, such that the following holds.
\begin{enumerate}
\item $Y$ is a quasi-smooth threefold of codimension $4$ with $\B(Y)=[7,3]$, $-K_Y=\OO_Y(3)$,
$\OO_Y(1)^3=\frac47$ and $\rho(Y)=1$, defined by $9$ equations of degrees $4,5,6,6,7,8,8,9,10$.
\item The points of $Y$ with nontrivial stabiliser in $\CC^*\times\langle\iota\rangle$ are the point
$P_w$ of index $7$, fixed by $\iota$, at which the stabiliser acts with type $\frac1{14}(3,11,13)$, and one
nonsingular fixed point of $\iota$ in $\{x_2=y_1=y_2=z_2=v=0\}$, at which $\iota$ acts as $-1$.
\item $X=Y/\iota$ is a Fano--Enriques threefold of type No.~$17$, with $\rho(X)=1$ and $\Clt(X)\cong\ZZ/2$.
\end{enumerate}
\end{proposition}

\begin{proof}
(i) and (ii), apart from $\rho(Y)$, are computations over $\FF_p$, $p=10009$, for a member $Y_0$, as in
\S\ref{ss:verify}. In Macaulay2 we pull back the relations $(1.2)$ of \cite{CD20} of the first four types and
saturate with respect to $\phi^*(A_2A_{23}A_{41})$; the resulting ideal has $9$ generators, of the degrees
stated, and the Hilbert function $2,5,10,17,28,42,61,84,\dots$ of \#41199. Quasi-smoothness is checked with
random $4\times4$ minors of the Jacobian matrix, which leave only the point $P_w$, where the Jacobian matrix has
rank $4$. The dimensions of the eigenspaces of $\iota$ on the coordinate ring of $Y_0$ in degrees $k\le12$ are
equal to the values $h^0(X,kA_X)$, $h^0(X,kA_X+\tau)$ of row~$17$. Moreover $\phi$ is a regular pullback from
the $\mathrm{G}_2^{(4)}$ subformat, that is, from the key variety $X^{(4)}\subset\mathbb A^{14}$ of
codimension $4$ obtained from $X_{\mathrm{G}_2}$ by setting $A_{12}=A_{34}=1$ and eliminating $\theta_{12}$ and
$\theta_{34}$ \cite[\S4.2.2]{CD20}: the pulled-back relations define a subscheme of $\mathbb A^8$ of
codimension $4$ whose intersection with $\{\phi^*(A_2A_{23}A_{41})=0\}$ has codimension $5$, so
$\phi^{-1}(X^{(4)})$ has the expected codimension $4$, and its ideal, which is contained in $I_{Y_0}$ and has the
same Hilbert function, equals $I_{Y_0}$ (\texttt{examples/no17/g2reg.m2}).

The lift $Y$ is obtained by lifting the coefficients of the forms $\phi^*(z)$ to integers and computing over
$\QQ$ (\texttt{examples/lift/lift17.m2}): the ideal has generators of the same degrees, its reduced Gr\"obner
basis has $57$ elements, all $p$-integral, and their reductions form a Gr\"obner basis of $I_{Y_0}$ with the same
leading monomials; the ideal is $\iota$-invariant. By Lemma~\ref{lem:goodred}, (i) and (ii), apart from
$\rho(Y)$, hold for $Y$, and $\phi$ is a regular pullback over $\QQ$ as well, since the codimensions above can
only increase from $\FF_p$ to $\QQ$.

For $\rho(Y)$: since $d(\lambda_{13}),d(\lambda_{24})>0$, the format is strict $\mathrm{G}_2^{(4)}$
\cite[Corollary~4.6]{CD20}, and a general quasi-smooth regular pullback of these degrees has
$\Cl\cong\ZZ$ \cite[Theorem~5.2]{CD20}. Let $M$ be the affine space of graded maps $\mathbb A^8\to\mathbb A^{14}$ of these degrees,
$\mathbb A^{14}$ being the ambient space of $X^{(4)}$ (so that the images of $\theta_{12},\theta_{34}$ are
determined by the tag equations), with no condition of equivariance, and let $M^\circ\subset M$ be
the open subset of those $\phi$ which are regular pullbacks from $X^{(4)}$ and define a quasi-smooth threefold.
It is nonempty, contains our choice and a general point of $M$, and is connected, being open in an affine space.
Over $M^\circ$ the pullback of the free resolution of the format is exact, since the pullback has the expected
codimension, so the universal family of affine cones is flat over $M^\circ$, and the associated family of weighted
projective quotient stacks is smooth and proper. Hence $b_2$ of the fibres is constant on $M^\circ$. For a
$\QQ$-Fano threefold $Y'$ with quotient singularities $H^1(\OO_{Y'})=H^2(\OO_{Y'})=0$, so
$\rho(Y')=h^{1,1}(Y')=b_2(Y')$. A general fibre has $\Cl\cong\ZZ$, so $b_2=1$, and therefore $\rho(Y)=1$.

(iii) follows from (i) and (ii) by the last paragraph of \S\ref{ss:verify}; No.~$17$ is the only row of
Table~\ref{tab:main} with the basket $[[2,1],[14,3]]$.
\end{proof}

\begin{remark}
In characteristic $0$, \cite[Lemma~5.1]{CD20} gives $\Cl(Y)\cong\ZZ$ if $\phi^*(\theta)$ is a prime element of
the coordinate ring for one of the cluster variables $\theta$. As a consistency check, which is not used above,
the surfaces $Y_0\cap\{\phi^*(\theta_1)=0\}$ and $Y_0\cap\{\phi^*(\theta_3)=0\}$ over $\FF_p$ are nonsingular in
codimension $1$ (supplementary file \texttt{examples/no17/g2prime.m2}). The variable $w$ of weight $7$ occurs only in relations of
degrees $8$, $9$ and $10$.
\end{remark}

\subsection{The numerical type No.~51}\label{ss:51}
For No.~$51$ we have $\qQ=2$, $\qW=1$, $T_X=\ZZ/2$, $\B(X)=[5\times[2,1],[6,1]]$ and $A_X^3=\frac53$. The cover has
$\B(Y)=[3,1]$ and $(\pi^*A_X)^3=\frac{10}3$; these are the invariants of the $\QQ$-Fano threefold of index $2$
with $\dim|A_X|=4$ of Prokhorov and Reid \cite{PR16}, which is obtained from a smooth quadric
$Q\subset\PP^4$. By Proposition~\ref{prop:N4}(ii) below that threefold has no involution of the
required kind. We replace $Q$ by a quadric cone.

Let $W=\{x_4^2+x_1x_3-x_2^2=0\}\subset\PP^4$, a cone over $\PP^1\times\PP^1$ with vertex $w=(1:0:0:0:0)$,
and let $T_0=W\cap\{x_4=0\}$, a cone over a conic; it is the image of $\PP(1,1,2)$ under
$(u_1:u_2:v)\mapsto(v:u_1^2:u_1u_2:u_2^2:0)$. Let $\Gamma\subset T_0$ be the curve $v^2u_1+u_1^5+u_2^5=0$; it is a
smooth curve of degree $5$ and genus $2$ passing through $w$, and its ideal in $W$ is generated by $x_4$
and
\[
c_1=x_0^2x_1+x_1^3+x_2x_3^2,\qquad c_2=x_0^2x_2+x_1^2x_2+x_3^3 .
\]
The form $g=x_0^4x_1+2x_0^2x_1^3+x_1^5+2x_0^2x_2x_3^2+2x_1^2x_2x_3^2+x_3^5$ of degree $5$ vanishes on
$\Gamma$ to order $2$: its restriction to $T_0$ is the square of the equation of $\Gamma$, and on $W$
we have $x_1g=c_1^2-x_4^2x_3^4$. Let $R$ be the graded subring of the field of fractions of the
homogeneous coordinate ring $\CC[x_0,\dots,x_4]/(x_4^2+x_1x_3-x_2^2)$ of $W$ generated by $x_0,\dots,x_4$, by $y_1=c_1/x_4$, $y_2=c_2/x_4$ in degree
$2$ and by $z=g/x_4^2$ in degree $3$, and $Y=\operatorname{Proj}R\subset\PP(1^5,2^2,3)$. Its equations are
\[
\begin{gathered}
x_2^2-x_1x_3-x_4^2,\qquad x_3^2x_4-x_2y_1+x_1y_2,\qquad sx_4+x_3y_1-x_2y_2,\\
sx_2+x_3^3-x_4y_2,\qquad sx_1+x_2x_3^2-x_4y_1,\qquad sy_1+x_3^2y_2-x_4z,\\
x_3^4-y_1^2+x_1z,\qquad sx_3^2+y_1y_2-x_2z,\qquad
x_0^4-x_1^4-2x_1x_2x_3^2+2x_1x_4y_1-y_2^2+x_3z,
\end{gathered}
\]
where $s=x_0^2+x_1^2$. Let
\[
\iota\colon(x_0,x_1,x_2,x_3,x_4,y_1,y_2,z)\mapsto(-x_0,x_1,x_2,x_3,-x_4,-y_1,-y_2,z).
\]

\begin{proposition}\label{prop:51}
\begin{enumerate}
\item $Y$ is a quasi-smooth threefold of codimension $4$ with $9$ equations of degrees $2,3^4,4^4$,
$\B(Y)=[3,1]$, $-K_Y=\OO_Y(2)$, $\OO_Y(1)^3=\frac{10}3$ and $\rho(Y)=2$.
\item The involution $\iota$ preserves $Y$. The points of $Y$ with nontrivial stabiliser in
$\CC^*\times\mumu_2$ are: the point $P_z$ of index $3$, fixed by $\iota$, with image of type
$\frac16(1,1,5)$; and five nonsingular fixed points of $\iota$, namely the points
$(0:t^2:t:1:0:0:0:t^3)$ with $t^5=-1$.
\item $X=Y/\iota$ is a Fano--Enriques threefold of type No.~$51$: $\rho(X)=1$,
$\B(X)=[5\times[2,1],[6,1]]$, $\Clt(X)\cong\ZZ/2$, $\qQ(X)=2$, $\qW(X)=1$ and $A_X^3=\frac53$.
\end{enumerate}
\end{proposition}

\begin{proof}
(i), (ii). The equations, the minimal free resolution of $R$ (with Betti numbers $1,9,16,9,1$ and last
twist $-10$, so that $\omega_Y=\OO_Y(10-12)$), the Hilbert function $5,16,37,71,122,193,\dots$ and the
quasi-smoothness of $Y$ are computed in Macaulay2 \cite{M2} over $\QQ$. The only point of $Y$ with
$x_0=\dots=x_4=0$ is $P_z$, where the tangent coordinates have degrees $1,2,2$; so $\B(Y)=[3,1]$. The
fixed points of $\iota$ are the points of $Y$ with $x_0=x_4=y_1=y_2=0$ listed above, and at each of them
$\iota$ acts on the tangent space of the affine cone, modulo the radial direction, by $-1$; at $P_z$ the
stabiliser is $\mumu_6$ acting with weights $(5,1,1)$. The dimensions of the invariant and anti-invariant
parts of $R$ in degrees $k\le8$ are $(3,2),(8,8),(19,18),(36,35),(61,61),(97,96),(144,143),(204,204)$, in
accordance with the values of $h^0(X,kA_X+\tau)$, $\tau\in T_X$, for No.~$51$: if $A_X$ is the image of
$\{x_4=0\}$, then $H^0(X,kA_X)$ is the anti-invariant part for odd $k$ and the invariant part for even $k$,
of dimensions $2,8,18,36,61,97,143,204$.

To compute $\rho(Y)$, note that $R[x_4^{-1}]=\CC[W][x_4^{-1}]$, so $Y\setminus D\cong W\setminus T_0$, where
$D=Y\cap\{x_4=0\}$ is an irreducible reduced divisor. The group $\Cl(W)$ is freely generated by the classes
of the planes $\Pi$, $\Pi'$ of the two families on $W$, and $T_0\sim\Pi+\Pi'$, so $\Cl(W\setminus T_0)\cong\ZZ$
and $\Cl(Y)$ has rank $2$. Since $Y$ has quotient singularities, $\rho(Y)=2$.

(iii) The involution $x_4\mapsto-x_4$ interchanges the two families of planes $x_1=\lambda(x_2\mp x_4)$,
$\lambda x_3=x_2\pm x_4$ on $W$, so $\iota$ acts on $\Cl(Y\setminus D)\cong\ZZ$ by $-1$, and
$\operatorname{rk}\Cl(Y)^\iota=1$. Hence $\rho(X)=1$. The remaining assertions follow from (i) and (ii) by
the last paragraph of \S\ref{ss:verify}: the rows of Table~\ref{tab:main} with this basket are No.~$9$ and
Nos.~$50$--$53$, and $(-K_X)^3=8\cdot\frac12\cdot\frac{10}3$ only for No.~$51$.
\end{proof}

The threefold $Y$ was found by the method of \S\ref{ss:N4}. If the map given by $|\pi^*A_X|$ is
birational onto a quadric $W$, the argument used there for No.~$105$ and No.~$111$ describes the double
cover in terms of a hyperplane section $T_0$ of $W$ and an invariant curve $\Gamma\subset T_0$ of degree
$5$, as above. If $T_0$ is the tangent hyperplane section at a nonsingular point of $W$, an invariant
curve $\Gamma$ is nonsingular at that point and there is no form of degree $5$ vanishing on $\Gamma$ to
order $2$ and not divisible by the equation of $T_0$, so the ring $R$ is too small; for a section through
the vertex of a quadric cone such a form exists.

\subsection{Gorenstein covers and characterisations}\label{ss:gor}
For $19$ numerical types of Table~\ref{tab:main} the formal cover basket in the column $\B(Y)$ is empty. If
moreover $\Clt(X)=T_X$ (which holds if $|T_X|$ is prime and, by Table~\ref{tab:subgroups}, for Nos.~$156$--$159$; the
other two, Nos.~$155$ and $175$, are treated in Proposition~\ref{prop:B24}), then the cover $Y$ defined by $\Clt(X)$ has
empty basket, and $Y$ is a Fano threefold with terminal Gorenstein singularities and $-K_Y=\qQ\,\pi^*A_X$, the group
$G=T_X^{\vee}$ acts on $Y$ freely in codimension one, every non-trivial element of $G$ has only isolated
fixed points, and $\operatorname{rk}\Cl(Y)^{G}=\operatorname{rk}\Cl(X)=1$. For $\qQ=2$ such a $Y$ is a
$G$-del Pezzo threefold in the sense of \cite{Pro13c}, of degree $d=|T_X|A_X^3$. By \cite[Theorem~1.7]{Pro13c}:
\begin{itemize}
\item if $d=8$, then $Y\cong\PP^3$ and $\pi^*A_X=\OO(2)$;
\item $d=7$ does not occur;
\item if $d=6$, then $Y$ is $\PP^1\times\PP^1\times\PP^1$ or the smooth divisor $V_6$ of bidegree $(1,1)$ in
$\PP^2\times\PP^2$;
\item if $d=5$, then $Y$ is the smooth del Pezzo threefold $V_5$.
\end{itemize}
In these cases $Y$ is smooth, so a point of $X$ of index $r$ is the image of a point of $Y$ whose
stabiliser has order $r$, and a non-trivial $g\in G$ has finitely many fixed points, their number being
the topological Lefschetz number $\sum(-1)^i\operatorname{Tr}(g^*|H^i(Y,\QQ))$. The odd cohomology of
$\PP^3$, $(\PP^1)^3$, $V_6$ and $V_5$ vanishes.

In several cases the cover is Gorenstein and the threefold $X$ is determined by its numerical type.

\begin{proposition}\label{prop:char}
Let $X$ be a Fano--Enriques threefold.
\begin{enumerate}
\item If $X$ is of type No.~$16$, then $X\cong\PP^3/\mumu_5(1,2,3,4)$.
\item If $X$ is of type No.~$37$, then $X\cong Y/\mumu_5$, where $Y\subset\PP^4$ is a
smooth quadric.
\item If $X$ is of type No.~$141$ and $\Clt(X)\cong\ZZ/3$, then $X\cong Y/\mumu_3$, where
$Y\subset\PP^4$ is a cubic hypersurface with terminal singularities.
\end{enumerate}
\end{proposition}

\begin{proof}
In all three cases $T_X=\Clt(X)$ is of prime order and its generator is not Cartier at any point of the
basket, and all indices are prime. Hence the cover $Y$ is a Gorenstein Fano threefold with terminal
singularities and $-K_Y=\qQ\,\pi^*A_X$ with $(\pi^*A_X)^3=|T_X|A_X^3$ equal to $1$, $2$, $3$ for $\qQ=4,3,2$
respectively. Therefore $Y$ is $\PP^3$, a quadric in $\PP^4$, or a cubic in $\PP^4$
\cite[\S3.1, 3.2]{IP99}.

(i) An action of $\mumu_5$ on $\PP^3$ lifts to a diagonal linear action on $\CC^4$. The fixed points
must be isolated, so the weights are distinct, and they must be terminal quotient singularities; this
gives the weights $(1,2,3,4)$ up to permutation and a common shift.

(ii) If the quadric $Y$ is singular, it is a cone with $\Cl(Y)\cong\ZZ^2$. The group $\mumu_5$ cannot
act nontrivially on $\ZZ^2$, which contradicts $\operatorname{rk}\Cl(Y)^{\mumu_5}=\operatorname{rk}\Cl(X)=1$.
So $Y$ is smooth. An example is the quadric $x_0^2+x_1x_4+x_2x_3=0$ with the action
$\mumu_5(0,1,2,3,4)$: its fixed points are the four coordinate points $P_1,\dots,P_4$, all of type
$[5,2]$.
\end{proof}

\begin{proposition}\label{prop:V5}
Let $V_5\subset\PP^6=\PP(\operatorname{Sym}^6\CC^2)$ be the closure of the $\operatorname{SL}_2$-orbit of
the binary form $xy(x^4-y^4)$, the smooth del Pezzo threefold of degree $5$, and let $\mumu_5$ act
through $\operatorname{diag}(\zeta,1)\in\operatorname{PGL}_2$. Then $X=V_5/\mumu_5$ is a Fano--Enriques
threefold of type No.~$172$: $\B(X)=[2\times[5,1],2\times[5,2]]$, $\Clt(X)\cong\ZZ/5$,
$\qW(X)=\qQ(X)=2$ and $A_X^3=1$.
\end{proposition}

\begin{proof}
In the coordinates $z_i$ dual to $x^{6-i}y^i$ the action is $\mumu_5(3,2,1,0,4,3,2)$. A computation as
in the proof of Proposition~\ref{prop:examples} shows that the fixed points of a generator on $V_5$
are the four coordinate points $P_0$, $P_1$, $P_5$, $P_6$, with tangent weights $(1,2,4)$, $(2,2,3)$,
$(2,2,3)$, $(1,2,4)$; they are the fixed points of the torus of $\operatorname{PGL}_2$, in accordance
with $\chi_{\mathrm{top}}(V_5)=4$. So $X$ has two points $\frac15(1,2,3)$ and two points
$\frac15(1,1,4)$. Since $\rho(V_5)=1$ and $-K_{V_5}=\OO(2)$, we get $\rho(X)=1$, $\qQ(X)=2$ and
$A_X^3=\frac55=1$; the remaining statements follow as in Lemma~\ref{lem:quotient}. The Hilbert series of
$(X,A_X+\tau)$, $\tau\in\Clt(X)$, coincide with those of No.~$172$ in degrees up to $8$.
\end{proof}

Up to conjugation $\operatorname{Aut}(V_5)=\operatorname{PGL}_2$ contains a unique subgroup of order $5$,
so by the classification above every Fano--Enriques threefold of type No.~$172$ is isomorphic to
$V_5/\mumu_5$. For the remaining $157$ numerical types the formal cover basket is not empty, so the cover defined by $T_X$,
if $\Clt(X)=T_X$, is not Gorenstein, and we do not know a classification of the possible covers.

\section{Exclusions}\label{sec:exclusions}
In this section we prove Theorem~\ref{thm:nonexistence}: it is the content of
Propositions~\ref{prop:B24}, \ref{prop:genusbound}, \ref{prop:largedeg}, \ref{prop:N4}, \ref{prop:173}, \ref{prop:N3}, \ref{prop:47}, \ref{prop:98}, \ref{prop:8893}, \ref{prop:139149}, \ref{prop:162}, \ref{prop:104}, \ref{prop:108}, \ref{prop:netpencil}, \ref{prop:97}, \ref{prop:175}, \ref{prop:11}, \ref{prop:13}, \ref{prop:no3} and~\ref{prop:liu26} below.

\subsection{Gorenstein covers}\label{ss:gorexcl}
We keep the notation of \S\ref{ss:gor}.

\begin{proposition}\label{prop:B24}
\begin{enumerate}
\item The numerical types No.~$40$, $41$, $142$, $158$, $159$, $169$ and $171$ do not occur.
\item The numerical type No.~$155$ does not occur with $\Clt(X)\cong\ZZ/3\times\ZZ/3$, and the numerical type No.~$175$
does not occur with $\Clt(X)\cong\ZZ/6$.
\end{enumerate}
\end{proposition}

\begin{proof}
In No.~$40$, $41$, $142$, $169$, $171$ the group $T_X$ has prime order, so $\Clt(X)=T_X$. In No.~$158$ and
$159$ we have $T_X=\ZZ/4$ and the class $K_X+2A_X$ is not divisible by $2$ in $T_X$ for any choice of $A_X$, so
it generates $T_X$ and again $\Clt(X)=T_X$. In No.~$155$ and $175$ we assume $\Clt(X)=T_X$. Thus in all cases
$Y$ and $G$ are as above.

\emph{$d=8$: No.~$41$, $159$, $171$.} Here $Y=\PP^3$. In No.~$41$ and $159$ the group $G$ contains an
involution; an involution of $\PP^3$ is induced by a linear involution of $\CC^4$, and its fixed locus
has positive dimension, a contradiction. In No.~$171$ the group $G=\mumu_5$ acts on $\PP^3$, so the
hyperplane class is the pull-back of a Weil divisor class on $X$ modulo torsion and $\qQ(X)=4$,
whereas $\qQ(X)=2$ for this numerical type; in fact $X$ is then the threefold of No.~$16$.

\emph{$d=7$: No.~$169$.} Excluded by the classification.

\emph{$d=6$, $Y=(\PP^1)^3$: No.~$40$, $142$, $155$, $158$, $175$.} The group $G$ acts on
$\Pic(Y)=\ZZ^3$ by permutations of the three generators, and $\operatorname{rk}\Pic(Y)^G=1$ only if this
action is transitive. A cyclic group of order $2$ or $4$ has no transitive action on three elements,
which excludes No.~$40$ and $158$. In the other cases there is an element $g\in G$ of order $3$ which
permutes the three factors cyclically: any non-trivial element in No.~$155$ outside the kernel of the
action, the generator in No.~$142$, and the square of a generator in No.~$175$. Write
$g(x,y,z)=(a(z),b(x),c(y))$ with $a,b,c\in\operatorname{PGL}_2$; then $g^3=\mathrm{id}$ gives
$a\circ c\circ b=\mathrm{id}$, and the fixed locus of $g$ is the curve $\{(x,b(x),c(b(x)))\}$, a
contradiction.

\emph{$d=6$, $Y=V_6$.} Here $\Pic(Y)=\ZZ^2$ is generated by the pull-backs of the hyperplane classes of
the two factors, and $b_2(Y)=b_4(Y)=2$. An automorphism of $Y$ either preserves the two generators or
interchanges them, and $\operatorname{rk}\Pic(Y)^G=1$ only if some element interchanges them. This is
impossible for a group of odd order, which excludes No.~$142$ and $155$. An element $g$ that
interchanges the generators has trace $0$ on $H^2$ and on $H^4$, so it has $1+0+0+1=2$ fixed points.
\begin{itemize}
\item No.~$40$: the involution has $8$ fixed points, a contradiction.
\item No.~$158$: a generator of $\ZZ/4$ interchanges the generators and has $4$ fixed points, the
preimages of the four points of index $4$, a contradiction.
\item No.~$175$: if $g$ is a generator of $\ZZ/6$, then $g^3$ is an involution which interchanges the
generators. Its fixed points are the preimages of the two points of index $6$ and of the two points of
index $2$, that is $2+2\cdot3=8$ points, a contradiction.
\end{itemize}
Finally we deduce Theorem~\ref{thm:nonexistence}(ii). In No.~$155$ the group $\Clt(X)$ is a non-zero subgroup of
$T_X\cong\ZZ/3\times\ZZ/3$, so by~(ii) it is isomorphic to $\ZZ/3$. In No.~$175$ the group $\Clt(X)$ is a non-zero
subgroup of $T_X\cong\ZZ/6$ containing the class $\tau_0$ of $K_X+2A_X$. If $\Clt(X)\cong\ZZ/3$, then
$\tau_0\in2\Clt(X)$, so $\qW(X)=2$, contrary to $\qW(X)=1$ (Theorem~\ref{thm:numerical}(iii)); by~(ii) we
get $\Clt(X)\cong\ZZ/2$, and this case is excluded in Proposition~\ref{prop:175}.
\end{proof}

\subsection{Covers of large degree: general facts}\label{ss:largedeg}
In some further cases the cover is not Gorenstein, but its anticanonical degree is so large that the
results of \cite{Pro13b,PR16,Pro24a} apply. Until the end of Lemma~\ref{lem:genus} we assume that $\Clt(X)=T_X\cong\ZZ/2$, $\pi\colon Y\to X$
is the double cover, $\tau$ is its involution and $B=\pi^*A_X$, so that $-K_Y\sim \qQ B$. The fixed locus
of $\tau$ is finite, since $\pi$ is \'etale over the smooth locus of $X$. For a threefold $V$ we write $\g(V)=\dim|-K_V|-1$.

\begin{lemma}\label{lem:descent}
The map $\pi^*\colon\Cl(X)\to\Cl(Y)^{\tau}$ is surjective with kernel $T_X$. In particular:
\begin{enumerate}
\item if the class of a Weil divisor $D$ on $Y$ is $\tau$-invariant, then $D\sim cB$ for some $c\in\ZZ$;
\item if $D$ is a non-zero effective Weil divisor on $Y$, then $D+\tau(D)\sim cB$ with $c\ge 1$, and $c=1$
if moreover $B\sim D+D'$ with $D'$ effective and non-zero;
\item if $\operatorname{rk}\Cl(Y)=1$ and $\Cl(Y)$ is torsion free, then $B$ generates $\Cl(Y)$.
\end{enumerate}
\end{lemma}

\begin{proof}
Let $X^0$ be the smooth locus of $X$ and $Y^0=\pi^{-1}(X^0)$, so that $Y^0\to X^0$ is \'etale with group
$\mumu_2$, $\Cl(X)=\Pic(X^0)$, $\Cl(Y)=\Pic(Y^0)$ and $H^0(Y^0,\OO^*)=\CC^*$. The Hochschild--Serre
spectral sequence gives the exact sequence
\[
0\to H^1(\mumu_2,\CC^*)\to\Pic(X^0)\to\Pic(Y^0)^{\tau}\to H^2(\mumu_2,\CC^*),
\]
and $H^2(\mumu_2,\CC^*)=\CC^*/(\CC^*)^2=0$. Since $\Cl(X)=\ZZ A_X\oplus T_X$, this gives (i). In (ii) the
class of $D+\tau(D)$ is invariant, so it equals $cB$ with $c\ge 1$ because $B$ is ample; if
$B\sim D+D'$, then $(2-c)B\sim D'+\tau(D')$ is effective and non-zero, so $c=1$. In (iii) the involution
acts trivially on $\Cl(Y)\cong\ZZ$, and (i) applies to a generator.
\end{proof}

\begin{lemma}\label{lem:decomp}
\begin{enumerate}
\item If $B\sim D+D'$ with non-zero effective Weil divisors $D$, $D'$ on $Y$, then $D'\sim\tau(D)$ and
$D\not\sim\tau(D)$.
\item There is no connected algebraic family $\{D_s\}_{s\in\Sigma}$ of non-zero effective Weil divisors
on $Y$ such that $\tau(D_s)$ belongs to the family for every $s$, and such that for some $s$ there is a
member $S\in|B|$ with $0<D_s<S$.
\item If $\dim|B|>0$, then the linear system $|B|$ has no fixed components.
\end{enumerate}
\end{lemma}

\begin{proof}
(i) By Lemma~\ref{lem:descent}(ii), $D+\tau(D)\sim B\sim D+D'$. If $D\sim\tau(D)$, then by
Lemma~\ref{lem:descent}(i) $D\sim cB$ with $c\ge1$, and $D'\sim(1-c)B$ is not effective and non-zero.
(ii) The members of a connected family are algebraically equivalent, hence linearly equivalent, because
a resolution of $Y$ is rationally connected and has $H^1(\OO)=0$. So $D_s\sim\tau(D_s)$, which
contradicts (i) applied to $D=D_s$ and $D'=S-D_s$.
(iii) The fixed part $\Phi$ of $|B|$ is $\tau$-invariant. If $\Phi\neq0$, then $\Phi\sim cB$ with $c\ge1$ by
Lemma~\ref{lem:descent}(i), and the moving part is equivalent to $(1-c)B$, which is impossible because
$\dim|B|>0$.
\end{proof}

The action of $\tau$ on $H^0(Y,\OO_Y(B))$ makes the map $\varphi=\varphi_{|B|}\colon Y\dashrightarrow\PP^N$,
$N=\dim|B|$, equivariant for a linear involution of $\PP^N$. Let $W=\overline{\varphi(Y)}$; it is an
irreducible non-degenerate $\tau$-invariant subvariety. We claim:
\begin{equation}\label{eq:nofamily}
\parbox{0.86\textwidth}{\emph{there is no $\tau$-invariant family $\{\Theta_s\}_{s\in\Sigma}$ of reduced subvarieties of
pure codimension one of $W$, with $\Sigma$ irreducible, which covers $W$ and such that every $\Theta_s$ is
contained in at least two different hyperplanes.}}
\end{equation}
Indeed, let $\nu\colon\tilde Y\to Y$ be a $\tau$-equivariant resolution such that
$\tilde\varphi=\varphi\circ\nu\colon\tilde Y\to W$ is a morphism, and let $W^0\subset W$ be a non-empty open
subset over which $\tilde\varphi$ is smooth. Since the family covers $W$, there is a non-empty open
subset $\Sigma^0\subset\Sigma$, which we may assume $\tau$-invariant, such that $\Theta_s\cap W^0\neq\varnothing$
for $s\in\Sigma^0$. For $s\in\Sigma^0$ let $\tilde D_s$ be the closure of
$\tilde\varphi^{-1}(\Theta_s\cap W^0)$, a reduced divisor because $\tilde\varphi$ is smooth over $W^0$. Let
$\tilde{\mathcal D}^0=\{(y,s)\mid\tilde\varphi(y)\in\Theta_s\cap W^0\}\subset\tilde Y\times\Sigma^0$ be the total
incidence family and $\tilde{\mathcal D}$ its closure. By generic flatness, after shrinking $\Sigma^0$
further (keeping it non-empty, open and $\tau$-invariant) the projection $\tilde{\mathcal D}\to\Sigma^0$ is
flat, and by \cite[Tag~0573]{Stacks} we may also assume that the open subset $\tilde{\mathcal D}^0$ is
scheme-theoretically dense in every fibre. Then the fibre of $\tilde{\mathcal D}$ over $s$ is a divisor in
which $\tilde\varphi^{-1}(\Theta_s\cap W^0)$ is scheme-theoretically dense, so it coincides with
$\tilde D_s$. So the $\tilde D_s$,
$s\in\Sigma^0$, are algebraically equivalent, hence linearly equivalent, and so are the divisors
$D_s=\nu_*\tilde D_s$ on $Y$. We have $\tau(D_s)=D_{\tau(s)}$. For general $s$ the divisor $D_s$ is not
zero, because $\nu$ contracts only finitely many prime divisors and the divisors $\tilde D_s$ cover
$\tilde Y$; shrinking $\Sigma^0$ once more, we may assume that $D_s\neq0$ for all $s\in\Sigma^0$. A section of $\OO_Y(B)$ vanishes on $D_s$ if the corresponding hyperplane contains
$\Theta_s$; by assumption there are two independent sections of this kind, and their divisors
$S,S'\in|B|$ are different, so $D_s<S$ or $D_s<S'$. This contradicts Lemma~\ref{lem:decomp}(ii).

\begin{lemma}\label{lem:crepant}
Suppose that $-K_Y\sim2B$, that $\varphi=\varphi_{|B|}$ is birational onto its image $W\subset\PP^N$, and that
$W$ is normal with Gorenstein canonical singularities and $K_W\sim\OO_W(-2)$. Then $B$ is a Cartier
divisor.
\end{lemma}

\begin{proof}
Let $p\colon\hat Y\to Y$ be a resolution such that $q=\varphi\circ p\colon\hat Y\to W$ is a morphism, and let
$H$ be a hyperplane section of $W$. Then $p^*|B|=|q^*H|+\sum m_FF$, where $F$ runs through the
$p$-exceptional prime divisors and $m_F\ge0$. Write $K_{\hat Y}=p^*K_Y+\sum a_FF=q^*K_W+\sum a'_FF'$,
where $F'$ runs through the $q$-exceptional prime divisors, $a_F>0$ and $a'_{F'}\ge0$. Since
$K_Y+2B\sim0$ and $K_W+2H\sim0$, we get
\[
\Delta:=\sum_F(a_F-2m_F)F-\sum_{F'}a'_{F'}F'\sim0 .
\]
The divisor $p_*\Delta=-\sum a'_{F'}p_*F'$ is linearly equivalent to zero and $-p_*\Delta$ is effective,
so $p_*\Delta=0$. Thus $\Delta$ is a principal divisor supported on the exceptional locus of $p$, hence
$\Delta=0$, and $a_F-2m_F\ge0$ for every $p$-exceptional divisor $F$, on every such resolution. Suppose
that $B$ is not Cartier at a point $P$, and let $r$ be the index of $P$. A $\QQ$-Cartier Weil divisor
$D$ at a terminal point of index $r$ has $rD$ Cartier \cite[Lemma~5.1]{Kaw88}; so $r>1$ and $rB$ is
Cartier at $P$. There is an exceptional divisor $F$ over $P$ with $a_F=\frac1r$ \cite{Kaw93}. All
members of $|B|$ pass through $P$, so $m_F\ge\frac1r$ and $a_F-2m_F<0$, a contradiction.
\end{proof}

\begin{lemma}\label{lem:genus}
Suppose that $-K_Y\sim2B$ and $\dim W=3$, and put $\delta=\deg\varphi\cdot\deg W$. Let
$p\colon\hat Y\to Y$ be a resolution such that $q=\varphi\circ p$ is a morphism, let
$p^*|B|=|M|+\sum m_FF$ with $M=q^*\OO_W(1)$, let $K_{\hat Y}=p^*K_Y+\sum a_FF$, and for a
$p$-exceptional prime divisor $F$ put $d_F=M^2\cdot F\ge0$; thus $d_F>0$ if and only if $q(F)$ is a
surface. Then
\begin{enumerate}
\item $B^3\ge\delta+\sum m_Fd_F$, and $d_F=0$ if $m_F=0$;
\item the intersection $\Gamma$ of two general members of $|M|$ is a smooth irreducible curve of genus
$g$, where $2g-2=\sum(a_F-2m_F)d_F$;
\item $h^0(Y,\OO_Y(B))=h^0(\Gamma,M|_\Gamma)+2$.
\end{enumerate}
In particular, if $B^3-\delta<\frac1r$, where $r$ is the maximal index of the points of $Y$, then
$\delta=\dim|B|-1$.
\end{lemma}

\begin{proof}
(i) Since $p^*B$ is nef and $(p^*B)^2\cdot F=0$ for every $p$-exceptional $F$,
\[
B^3=(p^*B)^2\cdot M=p^*B\cdot M^2+\sum m_F\,p^*B\cdot M\cdot F\ \ge\ p^*B\cdot M^2=\delta+\sum m_Fd_F .
\]
If $m_F=0$, then the centre of $F$ on $Y$ is not contained in the base locus of $|B|$, so $q(F)$ is the
closure of the image of this centre and has dimension at most $1$.
(ii) By the adjunction formula $2g-2=(K_{\hat Y}+2M)\cdot M^2$, and
$K_{\hat Y}+2M\sim p^*(K_Y+2B)+\sum(a_F-2m_F)F\sim\sum(a_F-2m_F)F$.
(iii) Let $S\in|M|$ be general. We have $H^1(\hat Y,\OO)=0$, and $H^1(S,\OO_S)=0$ because
$H^2(\hat Y,\OO(-S))$ is dual to $H^1(\hat Y,K_{\hat Y}+S)$, which vanishes by the Kawamata--Viehweg
theorem since $M$ is nef and big. So the restriction maps $H^0(\hat Y,M)\to H^0(S,M|_S)\to
H^0(\Gamma,M|_\Gamma)$ are surjective with one-dimensional kernels, and $h^0(\hat Y,M)=h^0(Y,\OO_Y(B))$.

For the last statement, note that $m_F\ge\frac1r$ whenever $m_F>0$: if the centre of $F$ is a point of
index $r'>1$, then $r'B$ is Cartier there and $m_F\ge\frac1{r'}$, and otherwise $B$ is Cartier at the
general point of the centre and $m_F\ge1$. Hence all $d_F$ vanish, $g=1$, and $h^0(\Gamma,M|_\Gamma)=\delta$.
\end{proof}

\label{ss:cyclic}Much of the above does not use that $T_X$ has order $2$. Let now $\Clt(X)=T_X\cong\ZZ/n$ be cyclic,
let $\pi\colon Y\to X$ be the cover defined by $T_X$, let $\gamma$ be a generator of its group $\mumu_n$, and
$B=\pi^*A_X$; assume $-K_Y\sim2B$ and $N=\dim|B|>0$. Since $H^2(\mumu_n,\CC^*)=0$, the proof of
Lemma~\ref{lem:descent} shows that $\pi^*\colon\Cl(X)\to\Cl(Y)^{\gamma}$ is surjective, so every Weil divisor
on $Y$ whose class is $\gamma$-invariant is linearly equivalent to $cB$, $c\in\ZZ$. Consequently:
\begin{itemize}
\item there is no connected algebraic family of non-zero effective Weil divisors on $Y$ which is mapped
to itself by $\gamma$ and has a member $D_s<S$ with $S\in|B|$: otherwise $D_s\sim\gamma(D_s)$, so
$D_s\sim cB$ with $c\ge1$, and $S-D_s\sim(1-c)B$ is not effective and non-zero;
\item $|B|$ has no fixed components;
\item the map $\varphi=\varphi_{|B|}\colon Y\dashrightarrow\PP^N$ is $\mumu_n$-equivariant, and the statement
\eqref{eq:nofamily} holds for its image $W$ with ``$\gamma$-invariant'' in place of
``$\tau$-invariant'', with the same proof;
\item Lemmas~\ref{lem:crepant} and~\ref{lem:genus} hold without changes.
\end{itemize}
For $T_X=\ZZ/2$ we have $\gamma=\tau$.

We collect some consequences for the image $W$ of $\varphi=\varphi_{|B|}$, for $T_X$ cyclic and $-K_Y\sim2B$. Let
$\nu\colon\tilde Y\to Y$ be a resolution such that $\tilde\varphi=\varphi\circ\nu$ is a morphism, and let
$M=\tilde\varphi^*\OO_W(1)$ be the moving part of $\nu^*|B|$. Since $|B|$ is a complete linear system,
$h^0(\tilde Y,M)=h^0(Y,\OO_Y(B))=N+1$; the pull-back $H^0(W,\OO_W(1))\to H^0(\tilde Y,M)$ is injective, and $W$ is
non-degenerate, so
\begin{equation}\label{eq:h0W}
h^0(W,\OO_W(1))=N+1 ,\qquad \Delta(W,\OO_W(1))=\dim W+\deg W-N-1
\end{equation}
in the notation of \cite{Fuj90}. If $N\ge2$, then $\dim W\ne1$: otherwise the points of $W$ form an invariant family, and
each of them lies in more than one hyperplane, contrary to \eqref{eq:nofamily}. If $\dim W=2$, $d=\deg W$ and $C$ is an
irreducible component of a general fibre of $\varphi$, then two general members of $|B|$ have $d$ general fibres in common,
and since $|B|$ has no fixed components and $B$ is ample, $B^3\ge d\,(B\cdot C)$. The divisor $B$ is not Cartier at a point
of index $r>1$ of $Y$, since $2B\sim-K_Y$.

\begin{lemma}\label{lem:dPcan}
Let $W\subset\PP^n$ be a normal threefold with $\Delta(W,\OO_W(1))=1$, so that $W$ is a del Pezzo variety: it is Gorenstein
and $K_W\sim\OO_W(-2)$ \cite[(6.2), (6.5)]{Fuj90}. If $W$ is not a cone, or if $W$ is a cone with vertex a point over a
surface which is not a cone, then $W$ has canonical singularities.
\end{lemma}

\begin{proof}
Suppose first that $W$ is not a cone. Let $v\in W$ be a singular point and $S$ a general hyperplane section through $v$.
Since $W$ is normal, $S$ has isolated singularities, so it is a normal Gorenstein surface with $-K_S=\OO_S(1)$, and by
\cite{HW81} either $S$ has Du Val singularities or $S$ is a cone over an elliptic curve. Since $v$ is a singular point of
$W$, the section $S$ is singular at $v$. If the general $S$ is a cone over an elliptic curve, then its vertex is $v$,
because such a cone has only one singular point. These sections cover a dense open subset of $W$, and each of them
contains the lines joining its points to $v$; so $W$ is a cone with vertex $v$, contrary to the assumption. Hence $S$ has
a Du Val singularity at $v$ or is smooth there, so $v$ is a compound Du Val point of $W$ and is canonical \cite{Rei87}.
If $W$ is a cone over a normal surface $W_0$ which is not a cone, then $W_0$ has Du Val singularities \cite{HW81}, the
blow up of the vertex of $W$ is crepant and has singularities of the form (Du Val point)$\,\times\,$(line).
\end{proof}

\begin{lemma}\label{lem:almostmin}
Suppose that $N\ge5$, $\dim W=3$, $\varphi$ is birational and $\deg W\le N-1$. Then $B$ is a Cartier divisor.
\end{lemma}

\begin{proof}
By \eqref{eq:h0W}, $\Delta(W,\OO_W(1))=\deg W-(N-2)\in\{0,1\}$. We show that otherwise $W$ carries a $\gamma$-invariant
family which contradicts \eqref{eq:nofamily}.

If $\Delta=0$, then by \cite[(5.10), (5.15)]{Fuj90} $W$ is a rational normal scroll or a generalised
cone over a rational normal scroll or a rational normal curve, and then it has a unique one-parameter
family of planes; or $W$ is the cone over the Veronese surface, and the cones over the images of the
lines of $\PP^2$ form an invariant family of surfaces each spanning a $\PP^3$. (A quadric does not occur
because $N\ge5$.)

Let $\Delta=1$. If $W$ is not normal, then by \cite[Theorems~1.2, 1.3(a)]{BS07} its normalisation is a
threefold of minimal degree in $\PP^{N+1}$, of which $W$ is a projection; $\gamma$ acts on it, and the
images of the surfaces just described form an invariant family of surfaces on $W$ each spanning at most
a $\PP^3$. If $W$ is normal and a cone with vertex a line, the planes through this line form an invariant
family. In the remaining cases $W$ has canonical singularities by Lemma~\ref{lem:dPcan}, and $B$ is Cartier by
Lemma~\ref{lem:crepant}.
\end{proof}

\begin{lemma}\label{lem:minsurf}
Suppose that $\dim W=2$ and $N\ge4$. Then $\deg W\ne N-1$, and if $N\le5$ and $\gamma$ has order $2$ or $3$, then also
$\deg W\ne N$.
\end{lemma}

\begin{proof}
If $\deg W=N-1$, then $W$ is a surface of minimal degree, that is, a rational normal scroll, a cone over a rational normal
curve, or the Veronese surface \cite[(5.10), (5.15)]{Fuj90}. In the first two cases the lines of $W$ form a unique
family; in the last one the images of the lines of $\PP^2$ form the unique family of conics on $W$. In both cases the
family is $\gamma$-invariant, covers $W$, and each member lies in at least two hyperplanes, contrary to
\eqref{eq:nofamily}.

Let $\deg W=N\in\{4,5\}$, so that $\Delta(W,\OO_W(1))=1$. By \cite[Theorem~1.2]{BS07} $W$ is a normal del Pezzo surface or a
projection of a surface $W'\subset\PP^{N+1}$ of minimal degree, which is its normalisation by \cite[Theorem~1.3(a)]{BS07}. In
the second case the images of the lines of $W'$, or of the conics if $W'$ is the Veronese surface, form an invariant family.
In the first case $W$ is a cone over an elliptic curve, covered by an invariant family of lines, or $W$ has Du Val
singularities \cite{HW81}. Then $\gamma$ lifts to the minimal resolution of $W$, a blow up of $\PP^2$ at $9-N$ points, possibly
infinitely near, and permutes the classes $f$ with $f^2=0$ and $-K\cdot f=2$, each with $h^0(f)=2$. For $N=5$ these are the
five classes $\ell-e_i$ and $2\ell-e_1-e_2-e_3-e_4$; for $N=4$ they are the ten classes $\ell-e_i$, $2\ell-\sum_{j\ne i}e_j$
($i=1,\dots,5$) (cf.\ \cite[Lemma~9.11]{Bla09}). If $\gamma$ has order $3$, one of them is fixed because $3\nmid5$ and
$3\nmid10$; if $\gamma$ has order $2$ and $N=5$, because $5$ is odd. If $\gamma$ has order $2$ and $N=4$, the classes form five
pairs $\{f_i,f'_i\}$ with $f_i+f'_i=-K$, and the vectors $\varepsilon_i=\frac12(f_i-f'_i)$ form an orthogonal basis of
$K^\perp\otimes\QQ$. The isometries of the Picard lattice preserving $K$ form the Weyl group of type $D_5$, which acts on the
$\varepsilon_i$ by permutations combined with an even number of changes of sign. An involution of this kind has an odd
number of fixed indices, among which the number of changes of sign is even; so $\gamma(f)=f$ for one of the ten classes. The
pencil $|f|$ may have a fixed part consisting of $(-2)$-curves; it is $\gamma$-invariant and is contracted by the
anticanonical morphism onto $W$. So the images of the moving parts of the curves of $|f|$ form an invariant family of curves
of degree at most $2$ covering $W$. All this contradicts \eqref{eq:nofamily}.
\end{proof}

The following local computation handles surfaces $W$ of large degree.

\begin{lemma}\label{lem:surflocal}
Suppose that $\dim W=2$, that every non-Gorenstein point of $Y$ is a cyclic quotient point of type $\frac13(1,1,2)$,
$\frac15(1,2,3)$ or $\frac15(1,1,4)$, and that $B\cdot C<1$ for an irreducible component $C$ of a general fibre of
$\varphi$. Then $W$ lies in a linear subspace of $\PP^N$ of dimension at most $3$, or at most $4$ if $Y$ has a point of type
$\frac15(1,1,4)$.
\end{lemma}

\begin{proof}
Since $B$ is Cartier outside the non-Gorenstein points $P_1,\dots,P_s$ of $Y$, of indices $r_i$, the curve $C$ passes through
some of them. Let $f\colon Y_1\to Y$ be the Kawamata blow up of $P_1,\dots,P_s$, with exceptional divisors $E_i$, and
$f^*|B|=\mathcal N+\sum c_iE_i$, where $\mathcal N$ is a linear system of Weil divisors which are $\QQ$-Cartier and do not
contain the $E_i$. Then
\[
B\cdot C=\nu_0+\sum_{i=1}^s\beta_i,\qquad \beta_i=\nu_i+c_i\,E_i\cdot\tilde C\in\tfrac1{r_i}\ZZ_{\ge0},
\]
where $\tilde C$ is the proper transform of $C$, $\nu_i\ge0$ is the contribution to $\mathcal N\cdot\tilde C$ of the points of
$\tilde C\cap E_i$, and $\nu_0\in\ZZ_{\ge0}$ is the contribution of the other points. Here $\beta_i$ is the local intersection
number at $P_i$ of $C$ with a general member $S\in|B|$; it lies in $\frac1{r_i}\ZZ$ because $r_iS$ is Cartier at $P_i$. Since
$B\cdot C<1$, we have $\nu_0=0$. The numbers $\nu_i$ and $E_i\cdot\tilde C$ are sums of local contributions which lie in
$\frac1k\ZZ$ at a point of index $k$ of $Y_1$ on $E_i$, and in $\ZZ$ elsewhere. The local data are as follows; in each case
the local sections of $\OO_Y(B)$ are the functions of weight $\chi$ in orbifold coordinates, with $2\chi$ congruent to the
sum of the weights, and $c_i$ lies in $c^{\min}+\ZZ_{\ge0}$, where $c^{\min}$ is the smallest order along $E_i$ of such a
function.
\begin{itemize}
\item $\frac13(1,1,2)$, with coordinates $x,y,z$ of weights $(2,1,1)$ and $\chi=2$: $E_i\cong\PP(2,1,1)$,
$\OO_{E_i}(E_i)=\OO_{E_i}(-3)$, and the only non-Gorenstein point of $Y_1$ on $E_i$ is its vertex $Q_i$, of type
$\frac12(1,1,1)$; $c^{\min}=\frac23$, and $\OO_{E_i}(\mathcal N)=\OO_{\PP(2,1,1)}(2)$ has four sections if $c_i=\frac23$. For
$\beta_i\le\frac23$ the solutions are $(c_i,E_i\cdot\tilde C,\nu_i,\beta_i)=(\frac23,\frac12,0,\frac13)$ and
$(\frac23,1,0,\frac23)$. If $\tilde C$ passes through $Q_i$, let $g\colon Y_2\to Y_1$ be the blow up of $Q_i$ with exceptional
divisor $G$ and $g^*\mathcal N=\mathcal N_2+mG$. For the proper transform $\hat C$, the contribution $\nu_{i,2}$ to
$\mathcal N_2\cdot\hat C$ of the points of $\hat C$ over $Q_i$ satisfies $0\le\nu_{i,2}=\nu_i-m\,G\cdot\hat C=-m\,G\cdot\hat C$, so
$m=0$, and $Q_i$ is not a base point of $\mathcal N$.
\item $\frac15(1,2,3)$, with coordinates of weights $(1,2,3)$ and $\chi=3$: $E_i\cong\PP(1,2,3)$,
$\OO_{E_i}(E_i)=\OO_{E_i}(-5)$, and $Y_1$ has on $E_i$ a point $Q_i$ of type $\frac12(1,1,1)$ and a point $Q'_i$ of type
$\frac13(1,1,2)$; $c^{\min}=\frac35$, and $\OO_{E_i}(\mathcal N)=\OO_{\PP(1,2,3)}(3)$ has three sections if $c_i=\frac35$.
The solutions with $\beta_i\le\frac45$ are
\[
\bigl(\tfrac35,\tfrac13,0,\tfrac15\bigr),\ \bigl(\tfrac35,\tfrac23,0,\tfrac25\bigr),\ \bigl(\tfrac35,1,0,\tfrac35\bigr),\
\bigl(\tfrac35,\tfrac43,0,\tfrac45\bigr),\ \bigl(\tfrac85,\tfrac12,0,\tfrac45\bigr),\ \bigl(\tfrac35,\tfrac12,\tfrac12,\tfrac45\bigr).
\]
In the last case $\tilde C$ passes through $Q_i$, and $\OO_{E_i}(\mathcal N)=\OO_{E_i}(3)$ is not invertible at $Q_i$, so every
member of $\mathcal N$ passes through $Q_i$. Let $g\colon Y_2\to Y_1$ be the blow up of $Q_i$, with exceptional divisor
$G_i\cong\PP^2$, $\OO_{G_i}(G_i)=\OO_{G_i}(-2)$, and $g^*\mathcal N=\mathcal N_2+mG_i$ with $m\in\frac12+\ZZ_{\ge0}$. The
threefold $Y_2$ is smooth near $G_i$, so for the proper transform $\hat C$ we have $G_i\cdot\hat C\in\ZZ_{>0}$ and
$0\le\nu_{i,2}=\nu_i-m\,G_i\cdot\hat C=\frac12-m\,G_i\cdot\hat C$, where $\nu_{i,2}$ is the contribution to
$\mathcal N_2\cdot\hat C$ of the points of $\hat C$ over $Q_i$. Hence $m=\frac12$, $G_i\cdot\hat C=1$ and $\nu_{i,2}=0$;
and $\OO_{G_i}(\mathcal N_2)=\OO_{\PP^2}(1)$ has three sections.
\item $\frac15(1,1,4)$, with coordinates $x,y,z$ of weights $(1,4,1)$ and $\chi=3$: $E_i\cong\PP(1,4,1)$,
$\OO_{E_i}(E_i)=\OO_{E_i}(-5)$, and the only singular point of $Y_1$ on $E_i$ is its vertex $Q_i$, of type $\frac14(1,3,1)$.
The functions of weight $3$ of smallest order along $E_i$ are the cubic forms in $x,z$, of order $\frac35=c^{\min}$, and
$\OO_{E_i}(\mathcal N)=\OO_{\PP(1,4,1)}(3)$ has four sections if $c_i=\frac35$; it is not invertible at $Q_i$, so every member of
$\mathcal N$ then passes through $Q_i$. The solutions with $\beta_i>0$ are
\[
\bigl(\tfrac35,1,0,\tfrac35\bigr),\ \bigl(\tfrac85,\tfrac14,0,\tfrac25\bigr),\ \bigl(\tfrac85,\tfrac12,0,\tfrac45\bigr),\
\bigl(\tfrac35,\tfrac14,\tfrac14,\tfrac25\bigr),\ \bigl(\tfrac35,\tfrac12,\tfrac12,\tfrac45\bigr)
\]
(these and the local computation below are checked by the script \texttt{local69.py}). In the last two, $\tilde C$
passes through $Q_i$. Let $g\colon Y_2\to Y_1$ be the Kawamata blow up of $Q_i$, with exceptional divisor
$G_i\cong\PP(1,3,1)$, $\OO_{G_i}(G_i)=\OO_{G_i}(-4)$; it has one singular point, of index $3$, on $G_i$, and
$\operatorname{ord}_{G_i}(E_i)=\frac34$. The functions of weight $3$ have order at least $\frac65$ along $G_i$, so
$g^*\mathcal N=\mathcal N_2+mG_i$ with $m\ge\frac65-\frac35\cdot\frac34=\frac34$. For the proper transform $\hat C$, the local
equations at $Q_i$ read $E_i\cdot\tilde C=e_2+\frac34\gamma$ and $\nu_i=n_2+m\gamma$, with $\gamma=G_i\cdot\hat C>0$ and
$e_2,n_2\ge0$ sums of local contributions in $\frac13\ZZ$ or $\ZZ$. If $E_i\cdot\tilde C=\frac14$, then $\gamma=\frac13$ and
$e_2=0$; if $E_i\cdot\tilde C=\frac12$, then $\gamma=\frac23$ and $e_2=0$, since $e_2=\frac14$ is impossible. In both cases
$m\gamma\le\nu_i$ gives $m=\frac34$ and $n_2=0$, and $\OO_{G_i}(\mathcal N_2)=\OO_{\PP(1,3,1)}(3)$ has five sections.
\end{itemize}
Let $\mathcal N'$ be the proper transform of $|B|$ on the blow up of $Y$ at all the centres above. In every case
$\mathcal N'\cdot\hat C=0$, so $\hat C$ does not meet the base locus of $\mathcal N'$, and $\hat C$ meets some $E_i$ or $G_i$.
If $c_i>c^{\min}$, the solutions show that $\tilde C$ meets $E_i$ only at $Q_i$, which is then not a base point of
$\mathcal N$. Hence $\varphi(C)$ is one of the finitely many points $\varphi(Q_i)$, or lies in the closure of the image of an
$E_i$ with $c_i=c^{\min}$ or of a $G_i$. Since the images of the general fibres cover the irreducible surface $W$, the
latter happens, and $W$ is contained in one of these closures, which span subspaces of dimension at most $3$, except for
the divisors $G_i$ over points of type $\frac15(1,1,4)$, whose images span subspaces of dimension at most $4$.
\end{proof}

\subsection{Covers of large genus: Nos.~8, 19, 25}\label{ss:genus}
In these three cases $T_X=\ZZ/2$ and the cover has large genus; we follow \cite[Proposition~3.3.1]{Pro24a}. The covers of
Nos.~$8$ and $19$ have the Hilbert series of the candidates \#41390 and \#41263 of \cite{GRDB}, which are also excluded in
\cite[Theorem~1.1(i)]{SR25}.

\begin{proposition}\label{prop:genusbound}
The numerical types No.~$8$, $19$ and $25$ do not occur.
\end{proposition}

\begin{proof}
\emph{No.~$8$.} Here $-K_Y\sim 4B$, $\B(Y)=[5,2]$, $B^3=\frac45$, $\dim|B|=2$ and $\g(Y)=26$. This case
is excluded in \cite[Proposition~3.3.1]{Pro24a}: a weak Fano threefold $V$ with terminal singularities, $\g(V)>22$ and
$-K_V\sim 4B$ with $\dim|B|\ge2$ is isomorphic to $\PP^3$.

\emph{No.~$19$.} Here $\qW(X)=\qQ(X)=3$, and we choose $A_X$ with $-K_X\sim 3A_X$. Then $-K_Y\sim 3B$,
$\B(Y)=[5,2]$, $B^3=\frac85$, $h^0(X,A_X)=h^0(X,A_X+\sigma)=2$ for the generator $\sigma$ of $T_X$, so that
$\dim|B|=3$, and $h^0(Y,-K_Y)=\frac{108}{5}+3-\frac35=24$, i.e.\ $\g(Y)=22$. We follow the proof of \cite[Proposition~3.3.1]{Pro24a}.

Suppose first that $\operatorname{rk}\Cl(Y)=1$. Since $K_Y$ is $\QQ$-Cartier and not numerically
trivial, every Weil divisor on $Y$ is then $\QQ$-Cartier, so $Y$ is a $\QQ$-Fano threefold whose index
$\qQ(Y)$ is divisible by $3$ and $\g(Y)=22$. By \cite[Theorem~1.2]{Pro13b}, $Y$ is the hypersurface
$Y_6\subset\PP(1,1,2,3,5)$ with $\qQ(Y)=6$. Its class group is torsion free (there are no numerical types with
$\qQ=6$ in Theorem~\ref{thm:numerical}), so $\Cl(Y)=\ZZ A_Y$ and $B\sim 2A_Y$. This contradicts
Lemma~\ref{lem:descent}(iii).

So $\operatorname{rk}\Cl(Y)\ge2$. Let $Y'\to Y$ be a $\QQ$-factorialisation; it is a small morphism,
$Y'$ is a weak Fano threefold with terminal $\QQ$-factorial singularities and $\rho(Y')\ge 2$. Run the
minimal model program
\[
Y'=Y^{(0)}\dashrightarrow Y^{(1)}\dashrightarrow\dots\dashrightarrow Y^{(n)}=:V,\qquad \varphi\colon V\to Z,
\]
where $\varphi$ is a Mori fibre space, and let $\mathcal B_k$ be the proper transform of $|B|$ on
$Y^{(k)}$. Then $-K_{Y^{(k)}}\sim 3\mathcal B_k$, $\dim\mathcal B_k=3$ and $\g(Y^{(k)})\ge 22$ for all $k$.
No step contracts a divisor to a curve, and $Z$ is not a surface: for a general fibre $\ell$ of such a
contraction we would have $-K\cdot\ell\in\{1,2\}$, which is not divisible by $3$.

\emph{Case $Z=\PP^1$.} A general fibre $F$ of $\varphi$ is a smooth del Pezzo surface with
$-K_F\sim 3\mathcal B_n|_F$, so $F\cong\PP^2$ and $\mathcal B_n|_F\subset|\OO_{\PP^2}(1)|$. Since
$\dim\mathcal B_n=3$, some member of $\mathcal B_n$ contains $F$. Denote by the same letter the proper
transform of $F$ on $Y$. We claim that $h^0(Y,\OO_Y(F))=2$. Indeed, the map $Y\dashrightarrow V$ contracts
finitely many prime divisors and extracts none, and we choose $F$ not containing the images of the contracted
divisors and not containing a flipped curve. The proper transforms on $Y$ of the fibres of $\varphi$
then form a pencil without fixed components, so $h^0(Y,\OO_Y(F))\ge2$. Conversely, the divisor on $V$ of
a rational function is the push-forward of its divisor on $Y$, so a function with poles at most along
$F$ on $Y$ has poles at most along $F$ on $V$, and $h^0(V,\OO_V(F))=h^0(\PP^1,\OO(1))=2$. By the above
$B\sim F+D$ with $D$ effective, and $D\neq0$ because $h^0(B)=4$. By Lemma~\ref{lem:descent}(ii), $B\sim F+F'$ with $F'=\tau(F)$. If $F'\sim F$, then the class
of $F$ is invariant and $B\sim 2F$, which contradicts Lemma~\ref{lem:descent}(i). Otherwise $|F|$ and $|F'|$
are two different pencils without fixed components. A non-zero element of the kernel of the
multiplication map cannot be a pure tensor $s\otimes t$, since $st\neq0$; so after a change of bases it is
$s_1\otimes t_1+s_2\otimes t_2$ with bases $s_i$ of $H^0(\OO_Y(F))$ and $t_i$ of $H^0(\OO_Y(F'))$.
The relation $s_1t_1+s_2t_2=0$ would give $\operatorname{div}(s_1)\le\operatorname{div}(t_2)$, hence equality since $B^2\cdot F=B^2\cdot F'$, and $F\sim F'$. Hence the multiplication map
$H^0(\OO_Y(F))\otimes H^0(\OO_Y(F'))\to H^0(\OO_Y(B))$ is injective, and so it is an isomorphism. Identify
$H^0(\OO_Y(F'))$ with $H^0(\OO_Y(F))$ by means of $\tau$. Then $\tau$ acts on
$H^0(\OO_Y(B))=H^0(\OO_Y(F))^{\otimes2}$ as $\pm$ the transposition of the factors, the sign depending on
the identification $\OO_Y(F+F')\cong\OO_Y(B)$, so its eigenspaces have dimensions $3$ and $1$ in
some order. But these eigenspaces are $H^0(X,A_X)$ and $H^0(X,A_X+\sigma)$, of dimension $2$ each, a
contradiction.

\emph{Case $Z$ is a point.} Then $V$ is a $\QQ$-Fano threefold with $3\mid\qQ(V)$ and $\g(V)\ge 22$, so
by \cite[Theorem~1.2]{Pro13b} either $V$ is a smooth quadric $Q\subset\PP^4$, or $V$ is a hypersurface of degree $6$ in $\PP(1,1,2,3,5)$.
Since $\rho(Y')\ge2$, we have $n\ge1$, and the last step $\psi\colon W=Y^{(n-1)}\to V$ contracts a
divisor $E$ to a point $P$. Since $-K_W\sim3\mathcal B_{n-1}$ is $\psi$-ample, $P$ is a base point of
$\mathcal B_n$.

If $V=Q$, then $\mathcal B_n\subset|\OO_Q(1)|$ has codimension $1$, so it is the system of hyperplane sections through $P$, so $\mathcal B_{n-1}$
has no base points outside $E$, $-K_W$ is nef and big, and $W$ is a weak Fano threefold. By
\cite{Kaw01} the contraction $\psi$ is a weighted blow up with weights $(1,w_1,w_2)$, where
$\gcd(w_1,w_2)=1$ and $3\mid w_1+w_2$ because $K_W=\psi^*K_Q+(w_1+w_2)E$ is divisible by $3$. Then
$-K_W^3=54-\frac{(w_1+w_2)^3}{w_1w_2}>0$ gives $\{w_1,w_2\}=\{1,2\}$, $\{1,5\}$, $\{2,7\}$, $\{4,5\}$ or
$\{5,7\}$, and $-K_W^3\le\frac{81}{2}$ with equality for $\{1,2\}$. In that case $\B(W)=[2,1]$ and
$h^0(-K_W)=\frac{81}{4}+3-\frac14=23$; in the other cases $h^0(-K_W)\le\frac12(-K_W^3)+3<23$. So
$\g(W)\le21$, a contradiction.

If $V\subset\PP(1,1,2,3,5)$, then $\mathcal B_n=|2A_V|$, because $\dim|2A_V|=3$. In suitable
coordinates $x,x',y,z,t$ the equation of $V$ is $x't+z^2+y^3+\dots=0$ \cite[Theorem~1.2]{Pro13b}, so the
only base point of $|2A_V|$ is the point $P=P_t$ of type $\frac15(1,2,3)$. By \cite{Kaw96} $\psi$ is the
weighted blow up of $P$ with weights $\frac15(1,2,3)$ in the orbifold coordinates $x,y,z$, and
$K_W=\psi^*K_V+\frac15E$. The proper transform $\tilde A$ of the divisor $\{x=0\}$ equals
$\psi^*A_V-\frac15E$, so $\Cl(W)=\ZZ\tilde A\oplus\ZZ E$ and $-K_W=6\tilde A+E$ is not divisible by $3$,
a contradiction.

\emph{No.~$25$.} Again $\qW(X)=\qQ(X)=3$, and for $-K_X\sim3A_X$ we have $-K_Y\sim3B$,
$\B(Y)=[[2,1],[4,1]]$, $B^3=\frac74$, $h^0(X,A_X)=h^0(X,A_X+\sigma)=2$, so that $\dim|B|=3$, and
$h^0(Y,-K_Y)=\frac{189}{8}+3-\frac14-\frac38=26$, i.e.\ $\g(Y)=24$. The argument for No.~$19$ applies
with the following changes. By \cite[Theorem~1.2]{Pro13b} a $\QQ$-Fano threefold whose index is divisible
by $3$ and whose genus is at least $24$ is the quadric $Q\subset\PP^4$, of genus $28$. Hence
$\operatorname{rk}\Cl(Y)\ge2$, and in the minimal model program for $Y'$ the case where $Z$ is a point
gives $V=Q$ and $\g(W)\le21$, a contradiction. The case $Z=\PP^1$ is excluded exactly as before, since
it used only $\dim|B|=3$ and the dimensions of the two eigenspaces.
\end{proof}

\subsection{Images of large dimension}\label{ss:Nge5}
In this subsection we treat the numerical types with $\qQ=2$ for which $N=\dim|B|\ge5$. The data are:
\begin{center}\footnotesize
\renewcommand{\arraystretch}{1.3}
\begin{tabular}{c|c|c|c|c||c|c|c|c|c}
No. & $T_X$ & $\B(Y)$ & $B^3$ & $N$ & No. & $T_X$ & $\B(Y)$ & $B^3$ & $N$\\\hline
$44$ & $\ZZ/2$ & $2{\times}[3,1]$ & $\frac{14}3$ & $5$ & $69$ & $\ZZ/2$ & $2{\times}[5,1]$ & $\frac{26}5$ & $5$\\
$45$ & $\ZZ/2$ & $2{\times}[3,1]$ & $\frac{20}3$ & $7$ & $72$ & $\ZZ/2$ & $2{\times}[5,2]$ & $\frac{24}5$ & $5$\\
$48$ & $\ZZ/2$ & $[5,2]$ & $\frac{22}5$ & $5$ & $94$ & $\ZZ/2$ & $2{\times}[3,1],[5,2]$ & $\frac{76}{15}$ & $5$\\
$49$ & $\ZZ/2$ & $[5,2]$ & $\frac{32}5$ & $7$ & $136$ & $\ZZ/3$ & $2{\times}[3,1]$ & $\frac{20}3$ & $7$\\
$52$ & $\ZZ/2$ & $[3,1]$ & $\frac{16}3$ & $6$ & $138$ & $\ZZ/3$ & $2{\times}[3,1]$ & $\frac{14}3$ & $5$\\
$53$ & $\ZZ/2$ & $[3,1]$ & $\frac{22}3$ & $8$ & $140$ & $\ZZ/3$ & $2{\times}[3,1]$ & $\frac{17}3$ & $6$\\
$63$ & $\ZZ/2$ & $2{\times}[3,1]$ & $\frac{14}3$ & $5$ & $164$ & $\ZZ/4$ & $[3,1]$ & $\frac{16}3$ & $6$\\
$64$ & $\ZZ/2$ & $2{\times}[3,1]$ & $\frac{20}3$ & $7$ &  &  &  &  & 
\end{tabular}
\end{center}
In all cases $-K_Y\sim2B$. In Nos.~$136$, $138$ and $140$ the group $T_X=\ZZ/3$ has prime order, so $\Clt(X)=T_X$; in No.~$164$
we have $T_X=\ZZ/4$ and $K_X+2A_X$ is a generator of $T_X$ for any choice of $A_X$, so again $\Clt(X)=T_X$. The non-Gorenstein points of
$Y$ are cyclic quotient points of the types $\frac13(1,1,2)$, $\frac15(1,2,3)$ and $\frac15(1,1,4)$ given by the basket, as
follows from Lemma~\ref{lem:cover} once we know that the corresponding points of $X$ of index divisible by $3$ or $5$ are
cyclic quotient points. Each of these points contributes a single entry to $\B(X)$, which is enough: a terminal point
which is not a cyclic quotient point contributes several entries, all of the same index except for points of type
$cAx/4$, whose entries have index $4$ and $2$ \cite[Remark~2.1]{CH11}. This is clear when the index of the entry occurs
only once among the entries of index at least $3$, as for the entries of index $10$ in Nos.~$48$, $49$, $94$, of index $6$
in Nos.~$52$, $53$, of index $3$ in No.~$94$ and of index $12$ in No.~$164$. In Nos.~$44$ and $45$ the two entries of index $6$ of $\B(X)$ have different local types of $A_X$ (namely
$K_X$ and $4K_X$, see Table~\ref{tab:main}), and in Nos.~$136$, $138$ and $140$ the two entries of index $9$ have different
local types of the generator of $T_X$ (namely $3$ and $6$); since the local type of a divisor is the same for all entries of
the basket of one terminal point, they belong to two distinct points. In Nos.~$63$, $64$, $69$ and $72$ the entry of index
$3$ or $5$ is a single cyclic quotient point of $X$ at which $\sigma$ is Cartier, so $\pi$ is \'etale over it and its two
preimages are interchanged by $\tau$.

\begin{proposition}\label{prop:largedeg}
The numerical types No.~$44$, $45$, $48$, $49$, $52$, $53$, $63$, $64$, $69$, $72$, $94$, $136$, $138$, $140$ and $164$ do not
occur.
\end{proposition}

\begin{proof}
By the remarks before Lemma~\ref{lem:dPcan}, $\dim W\ne1$. One checks from the table that $B^3<2(N-2)$ in all cases, and
$B^3<N$ except for Nos.~$69$ and $94$, where $N=5$ and $5<B^3<6$.

\emph{$\dim W=3$.} Then $\deg W\ge N-2$ and $\deg\varphi\cdot\deg W\le B^3<2(N-2)$, so $\varphi$ is birational and
$\deg W\le B^3$. If $\deg W\le N-1$, Lemma~\ref{lem:almostmin} shows that $B$ is Cartier, which is false. So we are in
No.~$69$ or No.~$94$, and $\deg W=5=\delta$. For No.~$94$, $B^3-\delta=\frac1{15}<\frac15$, and Lemma~\ref{lem:genus}
gives $\delta=4$, a contradiction. For No.~$69$, take the resolution $p$ of Lemma~\ref{lem:genus} $\tau$-equivariant. By
Lemma~\ref{lem:genus}(i), $\sum m_Fd_F\le B^3-\delta=\frac15$. If all $d_F$ vanish, then
$h^0(\Gamma,M|_\Gamma)=\delta$ by the proof of Lemma~\ref{lem:genus}, and (iii) gives $\delta=\dim|B|-1=4$, a
contradiction. Otherwise, since the $d_F$ are positive integers and $m_F\ge\frac15$ whenever $m_F>0$, there is exactly one
$F$ with $d_F>0$, and $m_F=\frac15$; by the proof of Lemma~\ref{lem:genus}, the centre of $F$ is then a point of index
$5$. As $\tau$ permutes the divisors $F$ with $d_F>0$, this $F$ and its centre are $\tau$-invariant, which is impossible
because $\tau$ interchanges the two points of index $5$.

\emph{$\dim W=2$.} Let $d=\deg W\ge N-1$. By Lemma~\ref{lem:minsurf}, $d\ne N-1$, and $d\ne N$ if $B^3\ge N$. So
$B^3<d$, and for an irreducible component $C$ of a general fibre of $\varphi$ we get $B\cdot C\le B^3/d<1$. By
Lemma~\ref{lem:surflocal}, $W$ lies in a linear subspace of dimension at most $4<N$, contradicting the non-degeneracy of
$W$.
\end{proof}

\subsection{Quadric images: $N=4$}\label{ss:N4}
\begin{proposition}\label{prop:N4}
\begin{enumerate}
\item The numerical types No.~$105$, $111$, $135$ and $144$ do not occur.
\item For a threefold of type No.~$51$ the double cover $Y$ has $\operatorname{rk}\Cl(Y)\ge 2$.
\end{enumerate}
\end{proposition}

\begin{proof}
\emph{No.~$105$ and No.~$111$.} Here $\qQ=2$, $-K_Y\sim2B$, $\dim|B|=4$ and $h^0(Y,\OO_Y(2B))=18$; the
basket of $Y$ is $[3\times[3,1]]$ with $B^3=4$ in No.~$105$, and $[[3,1],2\times[5,2]]$ with
$B^3=\frac{62}{15}$ in No.~$111$. By \eqref{eq:h0W}, $h^0(W,\OO_W(1))=5$, and $\dim W\neq1$.

\emph{$\dim W=2$.} Let $d=\deg W\ge3$. By Lemma~\ref{lem:minsurf}, $d\ge5$, so $B\cdot C\le\frac{B^3}5<1$ for a
component $C$ of a general fibre, and by Lemma~\ref{lem:surflocal} $W$ lies in a hyperplane, a contradiction.

\emph{$\dim W=3$.} Then $\delta=\deg\varphi\cdot\deg W\le4$. If $\delta=4$, then $B^3-\delta$ equals $0$,
respectively $\frac2{15}$, and Lemma~\ref{lem:genus} gives $\delta=3$, a contradiction. So $\varphi$ is
birational and $W$ is a cubic or a quadric.

Let $W$ be a cubic. If $W$ is not normal, it is a projection of $\PP^1\times\PP^2$ or of a cone of
degree $3$ in $\PP^5$ \cite[Theorems~1.2, 1.3(a)]{BS07}, and we get an invariant family of planes. If
$W$ is normal and a cone with vertex a line, the planes through this line form an invariant family. In the
remaining cases $W$ has canonical singularities by Lemma~\ref{lem:dPcan}, and Lemma~\ref{lem:crepant} shows that
$B$ is Cartier, which is false.

Let $W$ be a quadric. If its rank is $3$, the planes through the vertex line form an invariant family.
So $W$ is smooth or has one ordinary double point; it has terminal Gorenstein singularities and
$K_W\sim\OO_W(-3)$. In the second case $W$ has two families of planes, and by \eqref{eq:nofamily} they
are interchanged by $\tau$. We use the notation of Lemma~\ref{lem:genus}; here $\delta=2$.

\emph{Step 1.} Let $F'$ run through the $q$-exceptional prime divisors and
$K_{\hat Y}=q^*K_W+\sum a'_{F'}F'$, $a'_{F'}\ge1$. Comparing with
$K_{\hat Y}+2M\sim\sum(a_F-2m_F)F$ we get $M\sim\sum a'_{F'}F'-\sum(a_F-2m_F)F$, and applying $p_*$,
\[
B\sim S:=\sum a'_DD ,
\]
where $D$ runs through the prime divisors of $Y$ contracted by $\varphi$. If $S=A_1+A_2+A_3$ with
non-zero effective $A_i$, then Lemma~\ref{lem:decomp}(i) gives $A_3\sim\tau(A_1+A_2)\sim(A_2+A_3)+(A_1+A_3)$,
so $S\sim0$, which is absurd. Hence, by Lemma~\ref{lem:decomp}(i) again and because $S$ is $\tau$-invariant,
either $S=D$ is prime, or $S=D+\tau(D)$ with $D\ne\tau(D)$ prime; in both cases $a'_D=1$. Let
$h_0=\{\ell_0=0\}$ be the hyperplane corresponding to $S\in|B|$ and $T_0=h_0\cap W$.

\emph{Step 2.} Let $\mathcal F$ be the set of the $p$-exceptional divisors $F$ with $d_F>0$; the
surfaces $q(F)$, $F\in\mathcal F$, are exactly the prime divisors of $W$ which are not images of prime
divisors of $Y$, and $F$ is the proper transform of $q(F)$. We claim that they are the components of
$T_0$. Since $q^*T_0$ consists of the proper transform of $S$ and of $p$-exceptional divisors, every
component of $T_0$ is of this form. Conversely, let $T'=q(F)$. The divisor $T'+\tau(T')$ is Cartier and linearly equivalent to
$\OO_W(t)$ for some $t$: this is clear if $W$ is smooth, and if $W$ is singular, then $\Cl(W)$ is generated
by the classes of the planes of the two families, whose sum is the class of a hyperplane section, and
$\tau$ interchanges these classes. Let $g$ be a form of degree $t$ on $W$ with divisor
$T'+\tau(T')$; note that $\tau(T')$ is not the image of a prime divisor of $Y$ either. The corresponding
member of $|tB|$ is supported on $S$, so the divisor on $Y$ of
the function $f=g/\ell_0^t$ is $u(D-\tau(D))$ for some integer $u$ (and is zero if $S$ is prime). Then
$f\cdot\tau^*f$ is constant, so $2(T'+\tau(T'))=2tT_0$ on $W$, and $T'\subset T_0$. Consequently
$\varphi$ induces a birational map between the affine varieties $Y\setminus S$ and $W\setminus T_0$
which is an isomorphism in codimension one, hence an isomorphism, and
\begin{equation}\label{eq:LkS}
H^0(Y,\OO_Y(kB))=\left\{\frac{G}{\ell_0^m}\ \middle|\ G\in H^0(W,\OO_W(m)),\ v_D\left(\frac G{\ell_0^m}\right)\ge-k
\text{ for the components $D$ of $S$}\right\}.
\end{equation}

\emph{Step 3.} Since $p^*S=q^*T_0+\sum m_FF$, the coefficient of $F\in\mathcal F$ in $p^*S$ is $1+m_F$.
By Lemma~\ref{lem:genus}(i), $\sum_{\mathcal F}m_Fd_F\le B^3-2\le\frac{32}{15}$, and either $T_0$ is
irreducible, $\mathcal F=\{F\}$ and $d_F=2$, or $T_0$ is a union of two planes interchanged by $\tau$,
$\mathcal F=\{F,\tau(F)\}$ and $d_F=1$; in both cases $m_F\le\frac{16}{15}$. Hence in \eqref{eq:LkS}, if
$\ell_0\nmid G$, then $m\le k(1+m_F)$; in particular $m\le4$ for $k=2$.

\emph{Step 4.} Since $a'_D=1$, the centre of $D$ on $W$ is either the singular point $w$ of $W$, or a
curve $Z\subset T_0$ along which $D$ is the exceptional divisor of the blow up. In the first case $S=D$,
$v_D$ is the multiplicity at $w$ and $w\in T_0$, so for every quadric $G$ through $w$ the function
$G/\ell_0^2$ belongs to $H^0(Y,\OO_Y(B))$; this contradicts $\dim|B|=4$. So $S$ is contracted onto a
curve $\Gamma=Z$ or $\Gamma=Z\cup\tau(Z)$. Suppose first that $Z$ is not contained in the singular locus of
$T_0$. Then $v_D(\ell_0)=1$ and $v_D(G)=\operatorname{ord}_ZG$. By \eqref{eq:LkS} with $k=1$ and
$\dim|B|=4$,
\begin{equation}\label{eq:noquadric}
\text{no non-zero section of $\OO_{T_0}(2)$ vanishes on $\Gamma$,}
\end{equation}
and by \eqref{eq:LkS} with $k=2$ and Step~3, $h^0(Y,\OO_Y(2B))\le14+\alpha+\alpha'$, where $14=h^0(W,\OO_W(2))$,
$\alpha$ is the dimension of the space of sections of $\OO_{T_0}(3)$ vanishing on $\Gamma$, and $\alpha'$
that of the sections of $\OO_{T_0}(4)$ vanishing on $\Gamma$ to order at least $2$. We claim that
$\alpha+\alpha'\le3$, so that $h^0(Y,\OO_Y(2B))\le17$, a contradiction.
\begin{itemize}
\item $T_0\cong\PP^1\times\PP^1$. Let $(a,b)$ be the class of the greatest common divisor of the
sections of $\OO(3,3)$ vanishing on $\Gamma$; it vanishes on $\Gamma$. By \eqref{eq:noquadric},
$a=3$ or $b=3$; if $(a,b)=(3,0)$ or $(a,b)=(0,3)$, then $\Gamma$, which has at most two components, lies
on two lines of a ruling, contrary to \eqref{eq:noquadric}. So $\alpha\le h^0(\OO(3-a,3-b))\le3$. The class of
$\Gamma$ is not $\le(2,2)$, so $\alpha'=0$.
\item $T_0\cong\PP(1,1,2)$. If $\Gamma$ is given by an equation of weighted degree $e$, then $e\ge5$ by
\eqref{eq:noquadric}, so $\alpha\le h^0(\OO(6-e))\le2$ and $\alpha'=0$.
\item $T_0=\Pi\cup\tau(\Pi)$ is a union of two planes meeting in a line $L$. Then $\Gamma=Z\cup\tau(Z)$
with $Z\subset\Pi$ a plane curve of degree $z$, $Z\neq L$, and a section of $\OO_{T_0}(n)$ is a pair of
plane curves of degree $n$ with the same restriction to $L$. If $z=1$, the pair $(Z+L,\tau(Z)+L)$
contradicts \eqref{eq:noquadric}. If $z\ge3$, then $\alpha\le1$ and $\alpha'=0$. If $z=2$, then by
\eqref{eq:noquadric} the restrictions of $Z$ and $\tau(Z)$ to $L$ are not proportional; the sections
in question are the pairs $(Z+\lambda,\tau(Z)+\lambda')$ with lines $\lambda,\lambda'$, the condition on $L$
leaves a space of dimension at most $3$, and $\alpha'=0$.
\end{itemize}
Finally, if $Z$ is contained in the singular locus of $T_0$, then $T_0=\Pi\cup\tau(\Pi)$, $Z=L$, $S=D$
and $v_D(\ell_0)=2$. If $\ell_0\nmid G$ in \eqref{eq:LkS} with $k=2$, then the restriction of $G$ to one of
the planes is a non-zero curve of degree $m$ containing $L$ with multiplicity at least $2m-2$, so
$m\le2$ and $h^0(Y,\OO_Y(2B))\le14$, a contradiction.

\emph{No.~$144$.} Here $T_X=\ZZ/3$, so $\Clt(X)=T_X$, and $\B(Y)=[3\times[3,1]]$, $-K_Y\sim2B$, $B^3=4$, $N=4$,
$h^0(Y,\OO_Y(2B))=18$ and $h^0(Y,\OO_Y(3B))=43$: the same numbers as for the double cover in No.~$105$. We
follow the proof given above for No.~$105$ and indicate the changes. The case $\dim W=2$ and, for $\dim W=3$,
the case $\delta=4$, the case of a cubic $W$ and the case of a quadric of rank $3$ use only \eqref{eq:nofamily}
and Lemmas~\ref{lem:crepant}, \ref{lem:genus}, \ref{lem:dPcan}, \ref{lem:minsurf} and~\ref{lem:surflocal}, which
hold for $\gamma$ of order $3$, and go through without change. If $W$ is a quadric with one double point, each of its two families of planes is
$\gamma$-invariant, because $\gamma$ has order $3$; this contradicts \eqref{eq:nofamily}. So $W$ is a smooth
quadric, $\Cl(W)=\ZZ\,\OO_W(1)$, and $\delta=2$.

In Step~1, group the components of $S=\sum a'_DD\sim B$ into $\gamma$-orbits: the sum over each orbit is a
non-zero effective $\gamma$-invariant divisor, hence linearly equivalent to $cB$ with $c\ge1$. So $S$ is either a
$\gamma$-invariant prime divisor $D$ or $D+\gamma(D)+\gamma^2(D)$ with $D$ prime, and $a'_D=1$. In Step~2,
$T'+\gamma(T')+\gamma^2(T')\sim\OO_W(t)$, and the divisor of $f=g/\ell_0^t$ on $Y$ is
$u_0D+u_1\gamma(D)+u_2\gamma^2(D)$ (with $u_1=u_2=0$ if $S$ is prime); since its intersection with $B^2$ vanishes,
$u_0+u_1+u_2=0$, so $f\cdot\gamma^*f\cdot(\gamma^2)^*f$ is constant, and we get $T'\subset T_0$ and
\eqref{eq:LkS} as before. In Step~3, $T_0$ is irreducible because a smooth quadric threefold contains no
planes; so $d_F=2$, $m_F\le\frac{B^3-2}2=1$, and $m\le2k$ in \eqref{eq:LkS} if $\ell_0\nmid G$.

In Step~4, $W$ is smooth, so the centre of $D$ is a curve $Z$, not contained in the singular locus of $T_0$,
$\Gamma$ is $Z$ or $Z\cup\gamma(Z)\cup\gamma^2(Z)$, \eqref{eq:noquadric} holds and
$h^0(Y,\OO_Y(2B))\le14+\alpha+\alpha'$. The argument given for No.~$105$ shows $\alpha+\alpha'\le3$, a
contradiction, except in one case: $T_0\cong\PP^1\times\PP^1$ and the greatest common divisor of the sections
of $\OO(3,3)$ vanishing on $\Gamma$ has class $(3,0)$ (or $(0,3)$). Then $\Gamma$ lies on three lines of one
ruling; it is not a single line, by \eqref{eq:noquadric}, so it consists of three lines $L$, $\gamma(L)$,
$\gamma^2(L)$ of one ruling (which $\gamma$ preserves, having order $3$), and its class is $(3,0)$. In this case
we use \eqref{eq:LkS} with $k=3$. Filter $H^0(Y,\OO_Y(3B))$ by the exponent $m$ of $\ell_0$, which is at most
$6$ by Step~3. The functions with $m\le3$ form the space $H^0(W,\OO_W(3))$, of dimension $30$. For $m=4,5,6$,
restriction to $T_0$ embeds the $m$-th graded piece into the space of sections of $\OO_{T_0}(m)=\OO(m,m)$
vanishing to order at least $m-3$ along each of the three lines, that is, into
$H^0(\PP^1\times\PP^1,\OO(m-3(m-3),m))$, of dimension $10$, $0$, $0$. Hence $h^0(Y,\OO_Y(3B))\le40<43$, a
contradiction.

\emph{No.~$135$.} Here $T_X=\ZZ/3$, so $\Clt(X)=T_X$, $\B(X)=[2\times[9,1]]$, and as for No.~$136$ the two entries
have different local types of the generator of $T_X$, so they are two cyclic quotient points. Thus
$\B(Y)=[2\times[3,1]]$ with both points fixed by $\gamma$, $-K_Y\sim2B$, $B^3=\frac{11}3$, $N=4$ and
$h^0(Y,\OO_Y(2B))=17$. We argue as for No.~$144$. The case $\dim W=2$ is the same. If $\dim W=3$, then $\delta\le3$; if $\delta=3$, then $W$ is a
cubic, excluded as before, so $\delta=2$ and $W$ is a quadric, which is smooth by the argument given for
No.~$144$. Steps~1 and~2 apply without change. In Step~3, $T_0$ is irreducible, $\mathcal F=\{F\}$, $d_F=2$, and
$m_F\in\frac13\ZZ$ is positive with $2m_F\le B^3-2=\frac53$, so $m_F\in\{\frac13,\frac23\}$. Since $\varphi$ is
birational onto a quadric, the curve of Lemma~\ref{lem:genus}(ii) is a smooth conic, so
$-2=(a_F-2m_F)\,d_F$ and $a_F=2m_F-1$; as $Y$ is terminal, $m_F=\frac23$ and $a_F=\frac13$. Discrepancies
over the Gorenstein locus are integers, so the centre of $F$ is one of the two points of index $3$. Let
$f\colon Y_1\to Y$ be its Kawamata blow up \cite{Kaw96}, with exceptional divisor $E\cong\PP(1,1,2)$,
$\OO_E(-E)=\OO_E(3)$ and discrepancy $\frac13$; near $E$ the only singularity of $Y_1$ is a point of type
$\frac12(1,1,1)$. Every other exceptional divisor $G$ over this point has
$a(G,Y)=a(G,Y_1)+\frac13\operatorname{ord}_G(E)\ge\frac12+\frac16=\frac23$. Hence $F$ defines the same valuation
as $E$: it is the strict transform of $E$ on $\hat Y$. The weighted leading forms along $E$ of the sections of
$|B|$ with multiplicity $m_F=\frac23$ form a subspace $V_E\subset H^0(E,\OO_E(2))=\langle x^2,xy,y^2,z\rangle$
(weights $1,1,2$), and $q|_F$ is induced by $V_E$; since its image $T_0$ spans a $\PP^3$, $\dim V_E=4$ and
$V_E=H^0(E,\OO_E(2))$. This system is base point free and maps $E$ isomorphically onto the quadric cone, so
$T_0$ is the quadric cone. Moreover $m\le\frac53k$ in \eqref{eq:LkS} if $\ell_0\nmid G$. For $k=1$ this gives $m\le1$, so
no section of $\OO_{T_0}(2)$, a form of weighted degree $4$ on $\PP(1,1,2)$, vanishes on $\Gamma$; hence the
equation of $\Gamma$ has weighted degree $e\ge5$. For $k=2$ we get $m\le3$, and the sections with pole of
order $3$ inject into the forms of degree $6$ on $\PP(1,1,2)$ divisible by the equation of $\Gamma$, a space of
dimension at most $2$. Hence $h^0(Y,\OO_Y(2B))\le14+2<17$, a contradiction.

\emph{No.~$51$.} Here $-K_Y\sim 2B$, $\B(Y)=[3,1]$, $B^3=\frac{10}{3}$ and $\dim|B|=4$. Suppose that
$\operatorname{rk}\Cl(Y)=1$. As for No.~$19$, $Y$ is then a $\QQ$-Fano threefold with $\qW(Y)=\qQ(Y)=2$: if
$\qQ(Y)=2k$ with $k\ge2$, then $3A_Y^3=10/k^3$ would be an integer; indeed, for $f(n)=\chi(\OO_Y(nA_Y))$ the correction
terms of Theorem~\ref{thm:RR} have period $3$, so $\Delta^2\bigl(f(n+3)-f(n)\bigr)=3A_Y^3$ with $\Delta$ the difference operator,
and $f$ is integer-valued. By the main theorem of
\cite{PR16} the Kawamata blow up $Y_1\to Y$ of the point of index $3$, with exceptional divisor
$E_1$, has a second extremal contraction $Y_1\to Q$ onto a smooth quadric $Q\subset\PP^4$. It maps
$E_1$ isomorphically onto the tangent hyperplane section $E_0=Q\cap T_{P_0}Q\cong\PP(1,1,2)$, and
contracts a divisor onto an integral curve $\Gamma\subset E_0$ which, in coordinates $u_1,u_2,v$
of $\PP(1,1,2)$, is given by $v\alpha_3(u_1,u_2)+\beta_5(u_1,u_2)=0$. The involution $\tau$ fixes the unique point of index $3$ of $Y$, and the Kawamata blow up of this point
is unique \cite{Kaw96}, so $\tau$ lifts to $Y_1$. Since $\rho(Y_1)=2$, the cone of curves of $Y_1$ has
exactly two extremal rays; $\tau$ preserves the ray contracted by $Y_1\to Y$, hence also the other
one, and so it descends to $Q$. It preserves $E_0$, the image of $E_1$, its vertex $P_0$, and
$\Gamma$, the image of the exceptional divisor of $Y_1\to Q$ \cite[\S\S5--6]{PR16}. A curve or a surface
on $Q$ fixed pointwise by $\tau$ is contained in $E_0$, since otherwise its proper transform on $Y$ would
be fixed pointwise.

The action on $Q$ is induced by a linear involution of $\CC^5$ with eigenspaces $V_\pm$,
$\dim V_+>\dim V_-$, and since $Q$ is smooth the
equation of $Q$ is $a+b$ with non-degenerate quadratic forms $a$ on $V_+$ and $b$ on $V_-$. If
$\dim V_\pm=(4,1)$, then $\tau$ fixes the smooth quadric surface $\{a=0\}\subset\PP(V_+)$ pointwise, a
contradiction. So $\dim V_\pm=(3,2)$; the fixed locus of $\tau$ on $Q$ consists of the
smooth conic $C=\{a=0\}\subset\PP(V_+)$ and the two points of $\{b=0\}\subset\PP(V_-)$. A smooth conic
on the cone $E_0$ does not pass through its vertex, so $P_0\notin C$, and $P_0=(0:z^*)$ with $b(z^*)=0$.
The tangent hyperplane at $P_0$ is $\PP(V_+\oplus\CC z^*)$ with coordinates $(y,w)$, on which
$\tau$ acts by $w\mapsto-w$, and $E_0=\{a(y)=0\}$. The cone
$E_0$ is the image of $\PP(1,1,2)$ under $(u_1:u_2:v)\mapsto(u_1^2:u_1u_2:u_2^2:v)$, and any two
such parametrisations differ by a linear change of $u_1,u_2$ and a substitution $v\mapsto cv+q_2(u_1,u_2)$,
which preserves the form $v\alpha_3+\beta_5$ of the equation of $\Gamma$. So we may take $v=w$, and then
$\tau$ acts by $v\mapsto-v$. Since $\Gamma$ is invariant, $v\alpha_3+\beta_5$ is proportional to $-v\alpha_3+\beta_5$,
so $\alpha_3=0$ or $\beta_5=0$, and $\Gamma$ is not integral, a contradiction.
\end{proof}

The numerical type No.~$51$ does occur, with $\operatorname{rk}\Cl(Y)=2$, see \S\ref{ss:51}.

\subsection{The numerical type No.~173}\label{ss:173}
Here $\qQ=\qW=2$, $T_X=\ZZ/5$, $\B(X)=[[3,1],4\times[5,1]]$ and $A_X^3=\frac{11}{15}$. We keep the notation of
\S\ref{ss:cyclic}, with $n=5$. The point of index $3$ of $X$ is not in the support of $T_X$, so
$\B(Y)=[5\times[3,1]]$: there are five points $P_1,\dots,P_5$ of type $\frac13(1,1,2)$, permuted cyclically
by $\gamma$. Moreover $-K_Y\sim2B$, $B^3=\frac{11}3$, $h^0(Y,\OO_Y(B))=4$, and $\gamma$ acts on
$H^0(Y,\OO_Y(B))$ with four different eigenvalues. So $N=3$, the invariant planes of $\PP^3$ are the four
coordinate planes, and $\gamma$ preserves every family of lines on an invariant quadric surface.

We use a $\gamma$-equivariant resolution $p\colon\hat Y\to Y$ as in Lemma~\ref{lem:genus}, with $q=\varphi\circ p$,
$p^*|B|=|M|+\sum m_FF$ and $K_{\hat Y}=p^*K_Y+\sum a_FF$.

\begin{lemma}\label{lem:leading}
Let $F$ be a $p$-exceptional prime divisor.
\begin{enumerate}
\item If the centre of $F$ on $Y$ is a point $P_i$ and $m_F=\frac13$, then $q(F)$ is a point.
\item Let the centre of $F$ be a Gorenstein point $y$, let $m_F=1$, and suppose that $q(F)$ is a surface.
If $y$ is nonsingular, then $F$ is the exceptional divisor of the blow up of $y$, $a_F=2$ and $q(F)$ is a
plane. If $y$ is singular, then $a_F=1$ and $q(F)$ is an irreducible quadric or a plane.
\item Let the centre of $F$ be a curve $Z$, let $m_F=1$, and suppose that $q(F)$ is a surface. Then $F$ is
the exceptional divisor of the blow up of $Z$ over the general point of $Z$, and $a_F=1$.
\end{enumerate}
\end{lemma}

\begin{proof}
Let $v$ be the valuation of $F$. In a neighbourhood of $y$ the sections of $\OO_Y(B)$ are functions (of
weight $2$ in orbifold coordinates $x,y_1,y_2$ of weights $2,1,1$ in case (i)), and the map $q|_F$ is given
by the classes of the global sections in
$\{f\mid v(f)\ge m_F\}/\{f\mid v(f)>m_F\}$. In case (i) the values of $v$ on the orbifold coordinates lie
in $\frac13\ZZ_{>0}$, so every monomial of weight $2$ except $x$ has $v\ge\frac23$; hence this space
has dimension at most one.

In case (ii) let $\mathfrak m$ be the maximal ideal of $y$; its embedding dimension is $3$ or $4$,
since $y$ is nonsingular or an isolated compound Du Val point. The space above is a quotient
$\mathfrak m/\mathfrak m^2\to\mathfrak m/(\mathfrak m\cap\{v\ge2\})$; let $K_1$ be its kernel. Let
$E\subset\PP^{\,e}$, $e=\dim\mathfrak m/\mathfrak m^2-1$, be the exceptional divisor of the blow up $Y'$ of
$y$ (the projectivised tangent cone: $\PP^2$, or a quadric in $\PP^3$), and let $Z'\subset E$ be the centre
of $v$ on $Y'$. For linear forms $\ell,\ell'$ with $\ell'\notin K_1$ the function $\ell/\ell'$ is regular at
the general point of $Z'$, and its restriction to $F$ is the pull-back of its restriction to $Z'$. So
$q|_F$ factors through $F\to Z'$ followed by a linear projection, and $Z'$ is a surface, that is, a
component of $E$, along which $E$ is reduced because $v(\ell')=1$. Thus $v$ is the valuation of this
component, and $q(F)$ is the image of a plane or of a quadric under a linear projection (the first
jets of the global sections may span a proper subspace of $\mathfrak m/\mathfrak m^2$); being a
surface, it is a plane or an irreducible quadric, and a plane if $y$ is nonsingular. Finally
$K_{Y'}=\pi^*K_Y+2E$ if $y$ is nonsingular, and $K_{Y'}=\pi^*K_Y+E$ if $y$ is a hypersurface double point,
which gives $a_F$.

(iii) The general point of $Z$ is a nonsingular point of $Y$, and near it the sections of $\OO_Y(B)$
vanish on $Z$. As in (ii), $q|_F$ is given by the classes of the global sections in
$\mathcal I_Z/(\mathcal I_Z\cap\{v\ge2\})$, a quotient of the conormal sheaf of $Z$, and it factors through the
map of $F$ onto the centre of $v$ on the blow up of $Z$, which lies in the exceptional divisor
$E_Z$. Since $q(F)$ is a surface, this centre is $E_Z$, so $v$ is the valuation of $E_Z$ and $a_F=1$.
\end{proof}

\begin{proposition}\label{prop:173}
The numerical type No.~$173$ does not occur.
\end{proposition}

\begin{proof}
As before $\dim W\neq1$.

\emph{$\dim W=2$.} Let $d=\deg W\ge2$. If $d=2$, the lines of $W$ contain a $\gamma$-invariant family,
contrary to \eqref{eq:nofamily}. Let $d\ge3$, so that $B\cdot\Phi\le\frac{11}9$ for a general fibre $\Phi$
of $\varphi$. The irreducible components of the general fibres form finitely many families, permuted
by $\gamma$. A family which is not $\gamma$-invariant has five translates, so $\Phi$ has at least five
components and $B\cdot\Phi\ge\frac53$. If the curves $C$ of an invariant family pass through some $P_i$,
they pass through all of them, and again $B\cdot C\ge\frac53$. Hence no component of $\Phi$ passes through the points $P_i$; since $B$ is Cartier outside these
points, every component $C$ has $B\cdot C\ge1$, and $d\,(B\cdot\Phi)\le\frac{11}3$ with $d\ge3$ shows
that $\Phi=C$ is irreducible, $B\cdot C=1$ and $d=3$. Let $\hat C$ be the proper
transform of $C$ on $\hat Y$. Then $M\cdot\hat C=0$ and $1=p^*B\cdot\hat C=\sum m_FF\cdot\hat C$; the
divisors $F$ which meet $\hat C$ have centres not contained in $\{P_i\}$, so $m_F\ge1$ for them. Thus
there is exactly one such divisor $F_0$, $m_{F_0}=1$, $F_0\cdot\hat C=1$, and $q|_{F_0}\colon F_0\to W$ is
birational. Since $F_0$ is the only exceptional divisor meeting $\hat C$,
\[
2g(\hat C)-2=K_{\hat Y}\cdot\hat C=-2\,B\cdot C+a_{F_0}=a_{F_0}-2 ,
\]
so $a_{F_0}$ is even. If the centre of $F_0$ is a curve, then $a_{F_0}=1$ by
Lemma~\ref{lem:leading}(iii), a contradiction. So the centre of $F_0$ is a point, which is Gorenstein
because $C$ does not pass through the points $P_i$, and by Lemma~\ref{lem:leading}(ii) the surface
$W=q(F_0)$ is a plane or a quadric, a contradiction.

\emph{$\dim W=3$.} Then $W=\PP^3$ and $\delta=\deg\varphi\le3$. For a $p$-exceptional divisor $F$ with
$d_F>0$ and $m_F<1$ the centre is some $P_i$, $m_F\in\{\frac13,\frac23\}$ and there are five such divisors
with the same $m_F$ and $d_F$; by Lemma~\ref{lem:leading}(i) and Lemma~\ref{lem:genus}(i),
$m_F=\frac23$ and $B^3\ge\delta+\frac{10}3$, which is false. So $m_F\ge1$. If the centre of such an $F$
is a point $P_i$, its five translates give $B^3\ge\delta+5$, which is also false. Hence the centre of $F$
is a Gorenstein point or a curve, $m_F$ and $a_F$ are positive integers, and
\begin{equation}\label{eq:budget}
\sum m_Fd_F\le\tfrac{11}3-\delta,\qquad \sum(a_F-2m_F)d_F=2g-2,
\end{equation}
where the sums are over the divisors $F$ with $d_F>0$.

\emph{$\delta=3$.} By \eqref{eq:budget} all $d_F$ vanish, $g=1$, and Lemma~\ref{lem:genus}(iii) gives
$\delta=2$.

\emph{$\delta=2$.} Let $\bar Y\to\PP^3$ be the normalisation of $\PP^3$ in the function field of $Y$; it
is a double cover $w^2=f_{2r}$ branched in a reduced surface $R$ of degree $2r$, $K_{\bar Y}=(r-4)\bar H$,
and $Y$ is birational to $\bar Y$. The curve $\Gamma$ of Lemma~\ref{lem:genus} is the preimage of a
general line, so $g=r-1$, and $r\ge2$ because $h^0(\bar Y,\bar H)=4$. By \eqref{eq:budget} either all
$d_F$ vanish and $g=1$, or there is one divisor $F$ with $d_F=m_F=1$ and $a_F=2g$. In the second case
$a_F$ is even, so by Lemma~\ref{lem:leading}(iii) the centre of $F$ is not a curve, and by
Lemma~\ref{lem:leading}(ii) it is a nonsingular point and $a_F=2$. So $g=1$ and $r=2$ in both cases:
$\bar Y$ is a Gorenstein threefold with $K_{\bar Y}+2\bar H\sim0$. If $\bar Y$ has canonical
singularities, the proof of Lemma~\ref{lem:crepant} applies to $\bar Y$ and shows that $B$ is Cartier.
So $\bar Y$ is not canonical at some point $\bar y$ over $p\in R$. Let $\Pi$ be a general plane
through $p$; its preimage $\bar S$ is a normal Gorenstein surface with $-K_{\bar S}=\bar H|_{\bar S}$,
which does not have a Du Val singularity at $\bar y$, since otherwise $\bar y$ would be a compound Du Val
point. By \cite{HW81} $\bar S$ is a cone over an elliptic curve with vertex $\bar y$, that is,
$R\cap\Pi$ consists of four lines through $p$. Hence $R$ is a cone with vertex $p$ over a reduced plane
quartic $Q$, and $\bar Y$ is given near $\bar y$ by $w^2=Q(x_1,x_2,x_3)$. The weighted blow up with
weights $(2,1,1,1)$ is crepant, and its singularities are those of the double plane branched in $Q$
times a line; they are canonical unless $Q$ has a point of multiplicity $4$. So $R$ is the union of
four planes through a line. These planes are $\gamma$-invariant, but only two coordinate planes pass
through a line, a contradiction.

\emph{$\delta=1$.} Now $g=0$, and \eqref{eq:budget} reads $\sum m_Fd_F\le\frac83$,
$\sum(2m_F-a_F)d_F=2$. The solutions $(a_F,m_F,d_F)$ are: (A) one divisor with $(1,1,2)$; (B) one divisor
with $(2,2,1)$; (C) two divisors with $(1,1,1)$. In all cases these divisors are $\gamma$-invariant. Since $q$ is
birational, $F$ is the proper transform of $q(F)$ and $d_F=\deg q(F)$. Let $T_0\subset\PP^3$ be the union
of the surfaces $q(F)$: an irreducible quadric, a plane $\Pi$, or two planes $\Pi_1\cup\Pi_2$, with equation
$q_0$ of degree $2$ (where $q_0=\ell^2$ in case (B)); these are the only surfaces in $\PP^3$ which are not
images of prime divisors of $Y$.

As in Step~1 of the argument for No.~$105$, comparing $K_{\hat Y}+2M\sim\sum(a_F-2m_F)F$ with
$K_{\hat Y}+2M=-2M+\sum a'_{F'}F'$ we get $2B\sim S:=\sum a'_DD$, the sum over the prime divisors $D$ of $Y$
contracted by $\varphi$, with $a'_D\ge1$ the discrepancy of $D$ over $\PP^3$. The divisor $S_0$ of the
section $q_0$ of $\OO_Y(2B)$ is also supported on contracted divisors, and $S-S_0$ is the divisor of a
rational function whose divisor on $\PP^3$ is supported on $T_0$. In cases (A) and (B) it is constant;
in case (C) it is a power $(\ell_1/\ell_2)^c$, and $S=(1+c)\operatorname{div}_Y\ell_1+(1-c)\operatorname{div}_Y\ell_2$ with
$c\in\{0,\pm1\}$, where $c=\pm1$ is impossible because all $a'_D$ are positive. So $S=S_0$,
$v_D(q_0)=a'_D$ for every contracted $D$, and as in Step~2 there
\begin{equation}\label{eq:Rk173}
H^0(Y,\OO_Y(kB))=\sum_{m}\bigl\{G/q_0^{\,m}\ \big|\ \deg G=k+2m,\ v_D(G)\ge m\,a'_D\ \text{for all }D\bigr\},
\end{equation}
where $m\le k\,m_F$ if $q_0\nmid G$, because the coefficient of $F$ in $p^*S$ is $1+2m_F$ (in case (B)
we take $\ell$ instead of $q_0$, $\deg G=k+m$, $v_D(G)\ge m\,v_D(\ell)$ and $m\le2k$). Recall that a $\gamma$-invariant
effective divisor contained in a member of $|B|$ coincides with it; so a member of $|B|$ supported on
contracted divisors is the reduced sum of one orbit, and two different such members have no common
components. By the cyclic descent
each $\gamma$-orbit $O$ of contracted divisors has $\sum_{D\in O}D\sim c_OB$ with $c_O\ge1$ (this sum is
$\pi^*\pi_*D$ if $D$ is not invariant, and $D$ itself is a pull-back if it is invariant); the
discrepancy $a'_D$ is constant on $O$, so $\sum_Oa'_Oc_O=2$. Hence there are three possibilities: (a) one orbit with $a'=1$, $c=2$; (b) one orbit with
$a'=2$, $c=1$; (c) two orbits with $a'=c=1$. In (b) and (c) the divisors $\sum_{D\in O}D$ belong to $|B|$,
so $q_0$ is a product of two linear forms, equal in case (b). Thus (a) occurs exactly in case (A),
(b) in case (B) and (c) in case (C). A divisor $D$ with $a'_D=1$ is the exceptional divisor of the blow
up of a curve $Z_D\subset T_0$ along which $T_0$ is nonsingular; in case (B), $a'_D=2$ and $v_D(\ell)=1$,
so $D$ is the exceptional divisor of the blow up of a point of $\Pi$, or of the blow up with weights
$(1,2)$ of a curve $Z_D\subset\Pi$, for which $\{v_D\ge j\}$ restricts to $\Pi$ as the $\lceil j/2\rceil$-th
power of the ideal of $Z_D$. Let $\Gamma$ be the union of the centres. We now count sections using
\eqref{eq:Rk173}, $h^0(\OO_Y(B))=4$ and $h^0(\OO_Y(2B))=16$.

\emph{Case (A).} For $k=1$: no cubic section of $T_0$ vanishes on $\Gamma$. For $k=2$: there are at most
$10+\alpha+\alpha'$ sections, where $\alpha$ is the dimension of the space of quartic sections of $T_0$
vanishing on $\Gamma$ and $\alpha'$ that of the sextic sections vanishing on $\Gamma$ to order $2$. If
$T_0\cong\PP^1\times\PP^1$ and $\Gamma$ has class $(a,b)$, then $a\ge4$ or $b\ge4$, so $\alpha'=0$; if
$\Gamma$ is irreducible, then $a,b\ge1$ and $\alpha\le4$; if $\Gamma$ consists of five curves, then
$\alpha=0$. If $T_0\cong\PP(1,1,2)$ and $\Gamma$ has weighted degree $e$, then $e\ge7$, $\alpha\le2$ and
$\alpha'=0$. So $h^0(\OO_Y(2B))\le14$.

\emph{Case (C).} Here $\Gamma=\Gamma_1\cup\Gamma_2$ with $\Gamma_j\subset\Pi_j$ a $\gamma$-invariant plane
curve of degree $z_j$ not containing $L=\Pi_1\cap\Pi_2$, and a section of $\OO_{T_0}(n)$ is a pair of plane
curves of degree $n$ with the same restriction to $L$. For $k=1$ no non-zero pair of cubics
$(c_1,c_2)$ with $c_j\supset\Gamma_j$ exists; this gives $z_1,z_2\ge3$ (otherwise $(\Gamma_1+L+\lambda,0)$
is such a pair), and if $z_1=z_2=3$, the restrictions of $\Gamma_1$ and $\Gamma_2$ to $L$ are not
proportional. Then the pairs of quartics containing $\Gamma_1,\Gamma_2$ form a space of dimension at
most $3$, and there are no pairs of sextics containing $2\Gamma_1,2\Gamma_2$. So
$h^0(\OO_Y(2B))\le13$.

\emph{Case (B).} If the centres are points, then the conics of $\Pi$ through them (five points of an
orbit, or one point) give forms $G$ of degree $2$ with $G/\ell\in H^0(\OO_Y(B))$, so $h^0(\OO_Y(B))\ge5$. If
$\Gamma\subset\Pi$ is a curve of degree $z$, then $k=1$, $m=1$ gives $z\ge3$ (no conic of $\Pi$ contains
$\Gamma$). For $k=2$ the sections with pole of order exactly $m$, $1\le m\le4$, modulo those with smaller
poles, inject into the space of plane curves of degree $2+m$ containing $\lceil m/2\rceil\Gamma$; for
$z=3$ these spaces have dimensions $1,3,0,1$, and for $z\ge4$ they are smaller. So
$h^0(\OO_Y(2B))\le10+5=15$.
\end{proof}
\subsection{Double covers with $N=3$: Nos.~58, 62, 68, 71, 74, 76, 78}\label{ss:N3}
In this subsection $X$ is of one of the types No.~$58$, $62$, $68$, $71$, $74$, $76$, $78$, and we use the
notation of \S\ref{ss:largedeg}. In all of them $\qQ=2$, $\qW=1$, $T_X=\ZZ/2$, and $\B(X)$ consists of eight
entries of index $2$ and the entries of odd index listed in Table~\ref{tab:N3}, where also $B^3=2A_X^3$ is given. The latter are not in the support
of $T_X$. Except possibly for No.~$58$, each of them is the only entry of its index, hence a cyclic quotient point, and by
Lemma~\ref{lem:cover} it has two preimages in $Y$, of the same type, interchanged by $\tau$. For No.~$58$ the two entries
$[3,1]$, which have the same local types of $A_X$ and $\sigma$, are either two cyclic points or one non-cyclic point of type
$cA/3$ or $cD/3$ \cite[Remark~2.1]{CH11}; we first treat the cyclic case, and the non-cyclic case at the end of the proof of
Proposition~\ref{prop:N3}. In the cyclic case we denote
the non-Gorenstein points of $Y$ by $P_1,P_2$ (and $P_3,P_4$ for No.~$58$, the pairs being $\{P_1,P_2\}$ and
$\{P_3,P_4\}$). Moreover $-K_Y\sim2B$ and $h^0(X,A_X)=h^0(X,A_X+\sigma)=2$, so $N=3$ and $\tau$ acts on $\PP^3$ with
two eigenspaces of dimension $2$: its fixed locus is the union of two skew lines $L_+$, $L_-$. The points of $Y$
fixed by $\tau$ are isolated, since $X$ has isolated singularities. (The cover of No.~$43$ has the same numerical
invariants as that of No.~$62$, but there the two points of index $3$ are fixed by $\tau$.)
Lemma~\ref{lem:leading}, whose proof does not use the group, applies to a $\tau$-equivariant resolution
$p\colon\hat Y\to Y$ as in Lemma~\ref{lem:genus}.

Near $P_i$ we use orbifold coordinates of the weights given in the table. Since $-K_Y$ corresponds to the
sum of the weights and $2B\sim-K_Y$, the local sections of $\OO_Y(B)$ are the functions of the weight $\chi$
given in the table; the table also lists the monomials of small degree among them.

\begin{table}[ht]
\centering\footnotesize\setlength{\tabcolsep}{4pt}
\renewcommand{\arraystretch}{1.3}
\caption{The numerical types of \S\ref{ss:N3}. The second column lists the basket entries of $X$ of odd index, the others being eight entries of index $2$; $h^0(\OO_Y(2B))$ is given where it is used, and $m_0$ is the bound of
Lemma~\ref{lem:N3local}.}\label{tab:N3}
\begin{tabular}{c|c|c|c|c|c|l|c}
No. & \shortstack{odd-index\\entries} & $B^3$ & $h^0(\OO_Y(2B))$ & weights & $\chi$ & monomials of weight $\chi$ & $m_0$\\\hline
$58$ & $2\times[3,1]$ & $\frac{10}3$ & $15$ & $x,y_1,y_2$: $2,1,1$ & $2$ & $x$, quadrics in $y_1,y_2$ & $\frac23$\\
$62$ & $[3,1]$ & $\frac83$ & & $x,y_1,y_2$: $2,1,1$ & $2$ & $x$, quadrics in $y_1,y_2$ & $\frac23$\\
$68$ & $[5,1]$ & $\frac{16}5$ & $15$ & $x,y,z$: $1,4,1$ & $3$ & $y^2$, cubics in $x,z$ & $\frac45$\\
$71$ & $[5,2]$ & $\frac{14}5$ & & $u_1,u_2,u_3$: $1,2,3$ & $3$ & $u_3$, $u_1u_2$ & $\frac35$\\
$74$ & $[7,1]$ & $\frac{26}7$ & $17$ & $x,y,z$: $1,6,1$ & $4$ & $y^3$, quartics in $x,z$ & $\frac57$\\
$76$ & $[7,2]$ & $\frac{24}7$ & $15$ & $x,y,z$: $1,6,2$ & $1$ & $x$, $yz$ & $\frac37$\\
$78$ & $[7,3]$ & $\frac{20}7$ & & $x,y,z$: $1,6,3$ & $5$ & $y^2$, $x^2z$, $yz^2$ & $\frac47$
\end{tabular}
\end{table}

\begin{lemma}\label{lem:N3local}
Let $F$ be a $p$-exceptional divisor with centre $P_i$ and $d_F>0$. Then $m_F\ge m_0$, with $m_0$ as in
Table~\ref{tab:N3}.
\end{lemma}

\begin{proof}
Let $v$ be the valuation of $F$. Its values on the orbifold coordinates lie in $\frac1r\ZZ_{>0}$, where $r$ is the
index of $P_i$, and, as in the proof of Lemma~\ref{lem:leading}, $q|_F$ is given by the ratios of the leading forms
of value $m_F$ of the sections; since $d_F>0$, $q(F)$ is a surface.
\begin{itemize}
\item Nos.~$58$, $62$: this is Lemma~\ref{lem:leading}(i).
\item No.~$68$: $m_F\ge\frac25$. If $m_F=\frac25$, only $y^2$ can contribute, so $q(F)$ is a point. If $m_F=\frac35$,
then $y^2$ has value $\frac25$ or at least $\frac45$, so it does not contribute to the leading forms, which are
binary cubic forms in $x,z$; so $q|_F$ factors through $x/z$ and $q(F)$ is at most a curve.
\item No.~$71$: if $m_F\le\frac25$, the leading forms of the sections are spanned by $u_3$ and $u_1u_2$, and $q(F)$ is
at most a line.
\item No.~$74$: $m_F\ge\frac37$. If $m_F=\frac37$, the leading forms are multiples of $y^3$. If $m_F=\frac47$, then $y^3$,
whose value is $\frac37$ or at least $\frac67$, does not contribute, and the leading forms are binary quartics in
$x,z$, so $q|_F$ factors through $x/z$.
\item No.~$76$: $m_F\ge\frac17$; if $m_F=\frac17$, the leading forms are multiples of $x$, and if $m_F=\frac27$, they
are combinations of $x$ and $yz$, so $q(F)$ is at most a line.
\item No.~$78$: $m_F\ge\frac27$. If $m_F=\frac27$, the leading forms are multiples of $y^2$, and $q(F)$ is a point. If
$m_F=\frac37$, then $y^2$, whose value is $\frac27$ or at least $\frac47$, does not contribute to the leading forms,
which are therefore combinations of $x^2z$ and $yz^2$, and $q(F)$ is at most a line.\qedhere
\end{itemize}
\end{proof}

Since $\tau(P_i)\ne P_i$, a divisor $F$ as in Lemma~\ref{lem:N3local} comes with $\tau(F)\ne F$, which has the
same invariants and centre $\tau(P_i)$; we call $\{F,\tau(F)\}$ a \emph{pair}. A pair with $d_F>0$ costs at
least $2m_0$ in $\sum m_Fd_F$. A divisor $F$ whose centre is a Gorenstein point or a curve has integers
$m_F\ge1$ and $a_F\ge1$ if $d_F>0$.

\begin{lemma}\label{lem:N3quad}
If $\dim W=2$, then $\deg W\ge3$.
\end{lemma}

\begin{proof}
Let $d=\deg W\ge2$. If $d=2$, then $W$ is an irreducible quadric. If it is a cone, its lines form a unique
family. If it is smooth, $\tau$ preserves each ruling: an involution of $\PP^1\times\PP^1$ interchanging the
rulings has the form $(x,y)\mapsto(g(y),g^{-1}(x))$ and fixes the conic $\{x=g(y)\}$, which is not contained in
$L_+\cup L_-$. Both cases contradict \eqref{eq:nofamily}.
\end{proof}

\begin{lemma}\label{lem:N3delta2}
If $\dim W=3$, then $\deg\varphi=2$ and $g=1$ cannot hold together.
\end{lemma}

\begin{proof}
We argue as in the case $\delta=2$ of Proposition~\ref{prop:173}: $Y$ is birational to a double cover $\bar Y$ of
$\PP^3$ branched in a reduced quartic surface $R$; if $\bar Y$ has canonical singularities, $B$ is Cartier;
otherwise $R$ is the union of four planes through a line. In the last case, in coordinates in which the line is
$\{x_1=x_2=0\}$, $\bar Y$ is given by $w^2=\prod_{j=1}^4l_j(x_1,x_2)$ with pairwise independent linear forms
$l_j$, so it is birational to $\mathcal E\times\PP^2$, where $\mathcal E$ is the elliptic curve $w^2=\prod l_j$, and
a resolution of $\bar Y$ has $h^1(\OO)=1$. But a resolution of $Y$ has $h^1(\OO)=0$, since $Y$ has rational
singularities and $H^1(Y,\OO_Y)=0$ by the Kawamata--Viehweg theorem. This is a contradiction.
\end{proof}

Let now $\dim W=3$ and $\delta=1$, so that $\varphi\colon Y\dashrightarrow\PP^3$ is birational. For a divisor $F$ with
$d_F>0$ the surface $q(F)$ has degree $d_F$ and $F$ is its proper transform; let $T_0\subset\PP^3$ be the union of
these surfaces, which are the surfaces in $\PP^3$ that are not images of prime divisors of $Y$.

\begin{lemma}\label{lem:fixedlines}
In this situation $Y\setminus S\cong\PP^3\setminus T_0$ $\tau$-equivariantly for a $\tau$-invariant divisor
$S\sim2B$, and $L_+\cup L_-\subset T_0$. In particular, $T_0$ is not a union of planes of which at most one is
$\tau$-invariant.
\end{lemma}

\begin{proof}
As in Step~1 of the argument for No.~$105$ and in the case $\delta=1$ of Proposition~\ref{prop:173},
$2B\sim S:=\sum a'_DD$, the sum over the prime divisors $D$ of $Y$ contracted by $\varphi$; $S$ is $\tau$-invariant
and $Y\setminus S$ is affine. As in Step~2 there, $\varphi$ induces a birational map between the affine varieties
$Y\setminus S$ and $\PP^3\setminus T_0$ which is an isomorphism in codimension one, hence a $\tau$-equivariant
isomorphism. Since the fixed points of $\tau$ on $Y$ are isolated, $L_+\cup L_-\subset T_0$. The linear forms
vanishing on $L_+$ lie in one eigenspace of $\tau$, so every plane containing $L_+$ is $\tau$-invariant, and similarly
for $L_-$. If $T_0$ is a union of planes, then $L_+$ and $L_-$ lie in $\tau$-invariant planes of $T_0$, which are
different because $L_+$ and $L_-$ are skew.
\end{proof}

\begin{lemma}\label{lem:N3conf}
In this situation suppose that either
\begin{enumerate}
\item[(a)] there is one divisor $F$ with $d_F>0$; it is $\tau$-invariant, $m_F=a_F=1$, $d_F=2$, and $T_0$ is an
irreducible quadric; or
\item[(b)] there are two divisors $F$ with $d_F>0$; both are $\tau$-invariant and have $m_F=d_F=a_F=1$, and $T_0$ is the
union of two $\tau$-invariant planes.
\end{enumerate}
Then $h^0(Y,\OO_Y(2B))\le14$.
\end{lemma}

\begin{proof}
We argue as in the case $\delta=1$ of Proposition~\ref{prop:173}: now $m_F=1$ for the divisors over $T_0$, so
\eqref{eq:Rk173} holds with $m\le k$ if $q_0\nmid G$, and since $h^0(Y,\OO_Y(B))=4$, no non-zero section of
$\OO_{T_0}(3)$ vanishes on $\Gamma$. Each $\tau$-orbit $O$ of contracted divisors has $\sum_{D\in O}D\sim c_OB$ with
$c_O\ge1$ by Lemma~\ref{lem:descent}.

In case (a), $q_0$ is irreducible, so no orbit sum lies in $|B|$; hence $S$ is one orbit, of at most two divisors,
with $a'_D=1$, and $\Gamma$ has at most two components. If $T_0\cong\PP(1,1,2)$, the argument of Case~(A) there gives
$h^0(\OO_Y(2B))\le14$. If $T_0\cong\PP^1\times\PP^1$ and $\Gamma$ has class $(a,b)$, then $a\ge4$ or $b\ge4$, the sections of
$\OO(6,6)$ vanishing to order $2$ on $\Gamma$ are zero, and those of $\OO(4,4)$ vanishing on $\Gamma$ form
$H^0(\OO(4-a,4-b))$, of dimension at most $5$, with equality only for the classes $(4,0)$ and $(0,4)$, that is, if
$\Gamma$ is a union of four lines, which is impossible. So $h^0(\OO_Y(2B))\le10+4=14$.

In case (b), each $\operatorname{div}_Y\ell_j$ is a $\tau$-invariant member of $|B|$ and hence one orbit, and we are in
Case~(C) there: $h^0(\OO_Y(2B))\le13$.
\end{proof}

We also use the following three lemmas on divisors over the points $P_i$.

\begin{lemma}\label{lem:surfimage}
Suppose that $\dim W=2$, and let $C$ be an irreducible component of a general fibre of $\varphi$ which passes
through a cyclic quotient point $P=\frac1r(w_1,w_2,w_3)$ of $Y$, with orbifold coordinates $x_1,x_2,x_3$, at which
the local sections of $\OO_Y(B)$ are the functions of weight $\chi$. If $B\cdot C<1$, then $W$ is the closure of the image of a toric surface by a linear subsystem spanned by the monomials
$x_1^{m_1}x_2^{m_2}x_3^{m_3}$ of weight $\chi$ with $\sum v_jm_j=c$, for a primitive vector $v$ of the lattice
$\ZZ^3+\ZZ\frac1r(w_1,w_2,w_3)$ with positive coordinates and some $c\le B\cdot C$. For $W$ to be a surface, all
$v_j\le B\cdot C$, so $v=(\{kw_j/r\})_j$ for some $1\le k<r$, and the monomials have degree less than $r$.
\end{lemma}

\begin{proof}
Let $f\colon Y'\to Y$ be a toric resolution of the singularities of $Y$ at the points of index $>1$ (a composition of
weighted blow ups), with exceptional divisors $X_j$, and write $f^*|B|=|M'|+\sum c_jX_j$. Then $Y'$ has at worst
Gorenstein terminal singularities, where the $\QQ$-Cartier divisor $M'$ is Cartier \cite[Lemma~5.1]{Kaw88}, so $M'\cdot\tilde C$ is an integer for the proper transform
$\tilde C$; it is non-negative because $\tilde C$ moves, and $M'\cdot\tilde C\le f^*B\cdot\tilde C=B\cdot C<1$. So
$M'\cdot\tilde C=0$, $\tilde C$ does not meet the base locus of $|M'|$, and $\varphi(C)$ is the image of a point of
some $X_j$ over $P$ with $X_j\cdot\tilde C\ge1$, at which $|M'|$ has no base point; in particular
$c_j\le B\cdot C$. Since the general fibres cover $Y$, $W$ is the closure of the image of one such $X_j$, and the
map on $X_j$ is given by the leading forms along $X_j$ of the sections of $\OO_Y(B)$, which are combinations of the
monomials of weight $\chi$ and value $c_j$ along the ray $v$ of $X_j$. A monomial of value at most $B\cdot C<1$ involves
only the coordinates with $v_j<1$, and the exponent vectors of monomials in at most two of the coordinates with the
same value lie on a line, giving a curve; so a surface needs $0<v_j<1$ for all $j$, that is, $v$ is a lattice point
of the open unit box, and then $v_j\ge\frac1r$ bounds the degree.
\end{proof}

\begin{lemma}\label{lem:kawleading}
Let $F$ be a divisor over a point $P$ of type $\frac13(2,1,1)$, with orbifold coordinates as in Table~\ref{tab:N3}, with $m_F=\frac23$, whose image $q(F)$ is a
surface. Then $F$ defines the valuation of the Kawamata blow up of $P$, and $a_F=\frac13$.
\end{lemma}

\begin{proof}
Let $v$ be its valuation. The functions of weight $2$ are combinations of $x$, $y_1^2$, $y_1y_2$, $y_2^2$ and of
monomials of degree at least $3$, of value at least $1$. Since $q(F)$ is a surface, the leading forms of value
$\frac23$ of the sections generate, by their ratios, a field of transcendence degree $2$; so $v(x)=\frac23$,
$v(y_1)=v(y_2)=\frac13$, and the initial forms of $x,y_1,y_2$ in the graded algebra of $v$ are algebraically
independent: with $X,Y_1,Y_2$ these initial forms, every ratio defining $q|_F$ lies in $\CC(U,V)$,
$U=X/Y_1^2$, $V=Y_2/Y_1$, so $U$ and $V$ are algebraically independent, and a weighted homogeneous relation among
$X,Y_1,Y_2$, divided by a power of $Y_1$, would be a relation between $U$ and $V$. Hence $v(f)$ is the minimum of the weights $\frac13(2a+b_1+b_2)$ of the monomials $x^ay_1^{b_1}y_2^{b_2}$
of $f$, that is, $v$ is the monomial valuation of weights $\frac13(2,1,1)$, whose divisor is the exceptional divisor
of the Kawamata blow up, with discrepancy $\frac13$.
\end{proof}

\begin{lemma}\label{lem:discbound}
Let $F$ be a $p$-exceptional divisor whose centre is a point $P$ of $Y$ at which $Y$ is nonsingular or a cyclic quotient
singularity, and suppose that $q(F)$ is a surface. Then $a_F+1\le3m_F$.
\end{lemma}

\begin{proof}
Let $v$ be the valuation of $F$, extended to the orbifold chart $\CC^3\to Y$ at $P$ (if $P$ is singular), normalised so that
it agrees with $v$ on functions on $Y$; since the chart is \'etale in codimension one, the log discrepancy of $v$ on $\CC^3$ is
$a_F+1$ \cite[Proposition~5.20]{KM98}. The local sections of $\OO_Y(B)$ are functions on $\CC^3$, and $q|_F$ is given by the
ratios of the initial forms of value $m_F$ of the sections. Since $q(F)$ is a surface, there are three sections whose local
functions $f_1,f_2,f_3$ have value $m_F$ and algebraically independent initial forms in the graded algebra of $v$. Replace
$f_i$ by $h_i=f_i+t_i$ with general $t_i$ in a high power of the maximal ideal; this does not change the initial forms, and
$h_1,h_2,h_3$ is a system of parameters, since the linear systems $f_i+\mathfrak m^n$ have the origin as their only base point.
Let $\Phi=(h_1,h_2,h_3)\colon(\CC^3,0)\to(\CC^3,0)$, a finite map, and let $w=v\circ\Phi^*$. Then $w$ has the value $m_F$ on each
coordinate $u_i$, with algebraically independent initial forms, so $w(g)$ is the minimum of $m_F\deg$ over the monomials of $g$,
that is, $w=m_F\operatorname{ord}_0$, whose log discrepancy is $3m_F$. For the finite map $\Phi$ with ramification divisor $R\ge0$ we
have $A(v)=A(w)-v(R)\le3m_F$ \cite[Proposition~5.20]{KM98}.
\end{proof}

\begin{proposition}\label{prop:N3}
The numerical types No.~$58$, $62$, $68$, $71$, $74$, $76$ and $78$ do not occur.
\end{proposition}

\begin{proof}
By the remarks before Lemma~\ref{lem:dPcan}, $\dim W\ne1$.

\emph{$\dim W=2$.} Let $d=\deg W$; then $d\ge3$ by Lemma~\ref{lem:N3quad}, so $B\cdot\Phi\le\frac13B^3$ for a general fibre
$\Phi$ of $\varphi$. If $\Phi$ is irreducible and does not pass through the points $P_i$, then $B\cdot\Phi\ge1$, since $B$ is
Cartier outside these points; so $B\cdot\Phi=1$, $d=3$, and the argument given for this case in the proof of Proposition~\ref{prop:173} (which uses
Lemma~\ref{lem:leading}(ii), (iii) and the parity of $a_{F_0}$, but not the group) shows that $W$ is a plane or a quadric, a
contradiction. The remaining cases are as follows.
\begin{itemize}
\item No.~$62$. Let $C$ be an irreducible component of a general fibre. Then $B\cdot C\le\frac89<1$, so $C$
passes through $P_1$ or $P_2$, and $B\cdot C\in\{\frac13,\frac23\}$. Let $f\colon Y_1\to Y$ be the Kawamata blow up
of $P_1$ and $P_2$. The exceptional divisor $E_i$ is $\PP(1,1,2)$ with $\OO_{E_i}(E_i)=\OO_{E_i}(-3)$, and
$Y_1$ has on $E_i$ one point $Q_i$ of type $\frac12(1,1,1)$; so $f^*|B|=\mathcal N+\sum c_iE_i$ with
$c_1=c_2\in\frac23+\ZZ_{\ge0}$. As in the proof of Lemma~\ref{lem:surflocal}, $B\cdot C=\beta_1+\beta_2$ with
$\beta_i=\nu_i+c_i\,E_i\cdot\tilde C\in\frac13\ZZ_{\ge0}$, where $\nu_i$ and $E_i\cdot\tilde C$ are sums of local
contributions in $\frac12\ZZ$ at $Q_i$ and in $\ZZ$ elsewhere. For $\beta_i>0$ the only solutions are
$(c_i,E_i\cdot\tilde C,\nu_i,\beta_i)=(\frac23,\frac12,0,\frac13)$ and $(\frac23,1,0,\frac23)$. So
$\mathcal N\cdot\tilde C=0$, and as in the proof of Lemma~\ref{lem:surflocal} either every $\tilde C$ passes through $Q_1$ or $Q_2$, which is not a
base point of $\mathcal N$, and $W$ is a point, or $\varphi(C)$ is the image of a point of some $E_i$ which is not a
base point of $\mathcal N$, and $W$ is the closure of $\varphi(E_i)$. In the latter case, since
$\OO_{E_i}(\mathcal N)=\OO_{\PP(1,1,2)}(2)$, the surface $W$ is the image of $\PP(1,1,2)$ by a linear subsystem of
$|\OO(2)|$, so $d\le\OO(2)^2=2$, a contradiction.
\item No.~$71$. $B\cdot C\le\frac{14}{15}$ for every component $C$ of a general fibre, so $C$ passes through $P_1$ or $P_2$, and
$B\cdot C=\beta_1+\beta_2\le\frac45$ in the notation of the proof of Lemma~\ref{lem:surflocal}. The computation made there for points of type $\frac15(1,2,3)$, applied at $P_1$
and $P_2$, shows that $W$ is a point or lies in a plane, a contradiction.
\item No.~$78$. $B\cdot C\le\frac{20}{21}<1$ for every component $C$ of a general fibre, so $C$ passes through $P_1$ or $P_2$. For
$1\le k\le6$, the point $v=(\{k/7\},\{6k/7\},\{3k/7\})$ of the box, and every $c\le\frac67$, the monomials of weight $5$ and value $c$
have exponent vectors spanning at most a line (script \texttt{toricimage.py}); by Lemma~\ref{lem:surfimage}, $W$
is then a curve, a contradiction.
\item No.~$68$. $B\cdot\Phi\le\frac{16}{15}$. If a component $C$ of $\Phi$ does not pass through
$P_1,P_2$, then $B\cdot C\ge1$, so $\Phi=C$, which is excluded above. Otherwise every component $C$ passes through $P_1$ or $P_2$ and
$B\cdot C\le1$. We use the Kawamata blow ups $E_i$ and $G_i$ of the proof of Lemma~\ref{lem:surflocal} for points of type $\frac15(1,1,4)$. The solutions are those listed
there and $(c_i,E_i\cdot\tilde C,\nu_i,\beta_i)=(\frac35,\frac54,\frac14,1)$, for which the second blow up gives
$\gamma=\frac13$, $e_2=1$, $m=\frac34$, $n_2=0$ (script \texttt{local69.py}, with the bound $1$). So, as there,
$\mathcal N'\cdot\hat C=0$, and either $W$ is a point (if $c_i=\frac85$), or $W$ is contained in the closure of the image
of some $E_i$ or $G_i$. Now $\OO_{E_i}(\mathcal N)=\OO_{\PP(1,4,1)}(3)$, whose sections are the binary cubic forms in the
coordinates of weight $1$, so the image of $E_i$ is a curve. And the image of $G_i\cong\PP(1,3,1)$ is given by a linear
subsystem of $|\OO(3)|$, whose complete system maps $G_i$ onto the cone over a twisted cubic in $\PP^4$; so $W$, a surface
of degree $d\ge3$, is the projection of this cone from a point not on it (projection from the vertex gives a curve,
and from another point of the cone a surface of degree $2$), that is, a cone over a plane cubic curve. Its lines through the vertex form a unique family, contrary to \eqref{eq:nofamily}.
\item No.~$58$. $B\cdot\Phi\le\frac{10}9$. A component $C$ of $\Phi$ not passing through the points $P_i$ has $B\cdot C\ge1$, so $\Phi=C$, which is excluded above. Otherwise every component $C$ passes through
some $P_i$ and $B\cdot C\le1$. In the proof of Lemma~\ref{lem:surfimage} we have $M'\cdot\tilde C\le B\cdot C-\frac13<1$, since the
functions of weight $2$ have value at least $\frac13$ along every exceptional divisor over $P_i$; so the conclusion of that
lemma holds. The only box point and value giving a surface is the Kawamata ray with $c=\frac23$
(script \texttt{toricimage.py}); for $B\cdot C=1$ there is in addition the ray $(1,1,1)$, the only positive lattice
vector with a coordinate equal to $1$ and all coordinates at most $1$, along which only $x$ has value at most $1$, giving a
point. The image of the Kawamata exceptional divisor is that of $\PP(1,1,2)$ by a linear subsystem of $|\OO(2)|$, of degree at
most $2$. This contradicts $d\ge3$.
\item No.~$74$. $B\cdot\Phi\le\frac{26}{21}$. A component through $P_i$ has $B\cdot C=M'\cdot\tilde C+\sum c_jX_j\cdot\tilde C\ge\frac37$ in the notation of the proof of
Lemma~\ref{lem:surfimage}, since every $c_j\ge\frac37$ (see below) and some $X_j\cdot\tilde C\ge1$. A component $C$ not passing through
$P_1,P_2$ has $B\cdot C\ge1$, so $\Phi=C$ (another component would give $B\cdot\Phi\ge\frac{10}7>\frac{26}{21}$), which is excluded above. Otherwise every component passes through
$P_1$ or $P_2$ and $B\cdot C\le\frac87$. The functions of weight $4$ have value at least $\frac37$ along every exceptional divisor
over $P_i$, so in the proof of Lemma~\ref{lem:surfimage} we have $M'\cdot\tilde C\le\frac87-\frac37<1$, hence $M'\cdot\tilde C=0$,
and $W$ is the image of a toric exceptional divisor with ray $v$ by monomials of weight $4$ and value $c\le\frac87$, with
$0<v_j\le\frac87$; these rays are the nine primitive vectors $\frac17(1,6,1)$, $\frac17(1,6,8)$, $\frac17(8,6,1)$, $\frac17(k,7-k,k)$ for
$2\le k\le6$, and $(1,1,1)$. For all such $v$ and $c$ the exponent vectors span at most a line (script \texttt{toricimage.py}),
so $W$ is a curve, a contradiction.
\item No.~$76$. $B\cdot\Phi\le\frac87$. Let $C$ be a component of a general fibre and $\hat C$ its proper transform on a resolution $\hat Y$ on which $q$ is a morphism. Then
$M\cdot\hat C=0$, and every $p$-exceptional divisor $F$ meeting $\hat C$ satisfies $q(F)\ni\varphi(C)$; since the points
$\varphi(C)$ fill $W$, $q(F)=W$ for every $F$ meeting the general $\hat C$. We use this together with Lemma~\ref{lem:surfimage}.

If $C$ does not pass through $P_1,P_2$, then $B\cdot C\ge1$. Another component $C'$ of $\Phi$ would have $B\cdot C'\le\frac17$ and pass
through some $P_i$; by Lemma~\ref{lem:surfimage} with $c\le\frac17$, the only monomial of value at most $\frac17$ being $x$, $W$
would be a point. So $\Phi=C$, which is excluded above.

Let $C$ pass through $P_i$, and let $Y'$, $M'$, $X_j$, $c_j$ be as in the proof of Lemma~\ref{lem:surfimage}, so that
$B\cdot C=M'\cdot\tilde C+\sum c_jX_j\cdot\tilde C$ with $c_j\ge\frac17$. If $M'\cdot\tilde C=0$, the proof of that lemma applies
with $c\le\frac87$: by \texttt{examples/common/toricimage.py}, the only rays and values giving a surface are
$\frac17(4,3,1)$ with $c=\frac47$ and $\frac17(6,1,5)$ with $c=\frac67$, each with three monomials, so that $W$ would lie in a plane,
and $\frac17(1,6,2)$ with $c=\frac87$ and the monomials $z^4$, $yz$, $x^2z^3$, $x^2y$, $x^4z^2$, $x^6z$, $x^8$. In the last case, with
the torus coordinates $s=x^2/z$, $t=y/z^3$ of $X_j$ (a basis of the character lattice $\{a+6b+2c=0\}$), these monomials divided by $z^4$
are $1,t,s,st,s^2,s^3,s^4$; so every section of the subsystem defining the map is of the form $A(s)+tB(s)$, the curves $s=\mathrm{const}$
are mapped, for general $s$, to lines, and $W$ is swept out by a one-parameter family of lines. A surface of degree $d\ge3$ in $\PP^3$
has at most one such family: for a cone every line is a generator, and otherwise the lines of a second family would meet three
general skew lines of the first and sweep out the quadric through them. So this family is $\tau$-invariant, contrary to
\eqref{eq:nofamily}.

If $M'\cdot\tilde C\ge1$, then $M'\cdot\tilde C=1$ and $\sum c_jX_j\cdot\tilde C=\frac17$, so $\tilde C$ meets exactly one $X_j$, once, and
$c_j=\frac17$. The leading forms of value $\frac17$ are multiples of $x$, so the image of $X_j$ is a point; hence the point
$\tilde C\cap X_j$ is a base point of $|M'|$, since otherwise all the curves $C$ would be mapped to this point. Now apply the above to
the Cartier divisor $M'$ on $Y'$, which is nonsingular or Gorenstein along the base locus of $|M'|$: on a resolution
$r\colon\hat Y\to Y'$ with $r^*|M'|=|M|+\sum\mu_EE$, we have $1=M'\cdot\tilde C=\sum\mu_EE\cdot\hat C$, the divisors $E$ meeting
$\hat C$ with $\mu_E=0$ have images of dimension at most one, and $\hat C$ does not meet the proper transform of $X_j$. So
$\hat C$ meets exactly one $E_0$, with $\mu_{E_0}=1$, and $q(E_0)=W$. By Lemma~\ref{lem:leading}(ii), (iii), whose proofs are local statements about Cartier linear systems on terminal
Gorenstein threefolds and apply to $M'$ on $Y'$, the
surface $W$ is then a plane or a quadric (if the centre of $E_0$ is a point), or it is swept out by the images of the general fibres of
$E_0$ over its centre, on which the system restricts to $\OO_{\PP^1}(1)$, that is, by lines (if the centre is a curve). Both
contradict $d\ge3$ and \eqref{eq:nofamily} as before.
\end{itemize}

\emph{$\dim W=3$.} Then $W=\PP^3$ and $\delta=\deg\varphi\le3$, since $B^3<4$. By Lemma~\ref{lem:genus}, with the sums over the
divisors $F$ with $d_F>0$,
\[
\sum m_Fd_F\le B^3-\delta,\qquad \sum(a_F-2m_F)d_F=2g-2 .
\]
One checks from Table~\ref{tab:N3} that $2m_0>B^3-3$ for all seven rows, and $2m_0>B^3-2$ for Nos.~$62$, $68$, $71$, $78$.

\emph{$\delta=3$.} Here $B^3\ge3$, and $\sum m_Fd_F\le B^3-3<\min(1,2m_0)$, so all $d_F$ vanish, $g=1$, and
Lemma~\ref{lem:genus}(iii) gives $\delta=h^0(Y,\OO_Y(B))-2=2$, a contradiction.

\emph{$\delta=2$.} Here $\sum m_Fd_F\le B^3-2<2$. Suppose first that there is a pair with $d_F>0$; this happens only for
Nos.~$58$, $74$, $76$. A pair with $d_F\ge2$, two pairs, or a pair together with a divisor with Gorenstein centre or a curve
as centre would cost more than $B^3-2$; so the pair has $d_F=1$ and comes alone, and $m_F<1$. Then $2(a_F-2m_F)=2g-2$ and
Lemma~\ref{lem:discbound} give $g=a_F-2m_F+1\le m_F<1$, so $g=0$; then $M|_\Gamma$ has degree $2$ on a rational curve and
$h^0(\Gamma,M|_\Gamma)=3$, contrary to Lemma~\ref{lem:genus}(iii). Otherwise either all $d_F$ vanish, or there is one
divisor $F$ with $d_F>0$; its centre is a Gorenstein point or a curve, it is $\tau$-invariant, $d_F=m_F=1$ and $a_F=2g$, so
$a_F=2$ and $g=1$ by Lemma~\ref{lem:leading}(ii), (iii), as in the case $\delta=2$ of Proposition~\ref{prop:173}. In both cases
$g=1$, contrary to Lemma~\ref{lem:N3delta2}.

\emph{$\delta=1$.} Here $\sum m_Fd_F\le B^3-1$ and $\sum(2m_F-a_F)d_F=2$. A divisor with Gorenstein centre or a curve as
centre and $(m_F,d_F)=(1,1)$ contributes $2-a_F\le1$. We show that either $T_0$ is a union of planes of which at most one is
$\tau$-invariant, which is impossible by Lemma~\ref{lem:fixedlines}, or we are in case (a) or (b) of
Lemma~\ref{lem:N3conf}, which is impossible since $h^0(Y,\OO_Y(2B))\ge15$ (this occurs only for Nos.~$58$, $68$, $74$, $76$).
\begin{itemize}
\item No.~$62$: $\sum m_Fd_F\le\frac53$. A pair costs at least $\frac43$ and contributes $2(2m_F-a_F)d_F$; a divisor with
Gorenstein centre or with a curve as centre costs at least $1$ and contributes at most $2m_F-1$. The only solution is a
pair with $m_F=\frac23$, $d_F=1$ and $a_F=\frac13$, so $T_0$ is a pair of planes interchanged by $\tau$.
\item Nos.~$71$, $78$: $\sum m_Fd_F\le\frac95$, resp.\ $\frac{13}7$, and a pair costs at least $\frac65$, resp.\ $\frac87$. A
divisor with Gorenstein centre or a curve as centre has $(m_F,d_F)=(1,1)$ and contributes at most $1$, and it cannot occur
together with a pair or with another such divisor; so the only possibility is a pair with $d_F=1$, and $T_0$ is a pair of
planes interchanged by $\tau$.
\item No.~$68$: $\sum m_Fd_F\le\frac{11}5$. A pair costs at least $\frac85$, so it has $(m_F,d_F)=(\frac45,1)$ or $(1,1)$ and
comes alone. So the possibilities are: a pair with $d_F=1$, so that $T_0$ is a pair of planes interchanged by $\tau$; one
invariant divisor with $(m_F,d_F,a_F)=(1,2,1)$, case (a) of Lemma~\ref{lem:N3conf}; one invariant divisor with $(m_F,d_F,a_F)=(2,1,2)$, so that $T_0$ is
a plane; or two divisors with $(m_F,d_F,a_F)=(1,1,1)$, so that $T_0$ is a pair of planes, each $\tau$-invariant (case (b) of Lemma~\ref{lem:N3conf}) or
interchanged by $\tau$.
\item No.~$58$: $\sum m_Fd_F\le\frac73$. A pair has $d_F=1$ and either $m_F=\frac23$, $a_F=\frac13$ by Lemma~\ref{lem:kawleading},
contributing $2$, or $m_F=1$; two pairs cost more than $\frac73$. The possibilities are therefore: a pair alone, so that $T_0$
is two planes interchanged by $\tau$; a pair with $m_F=\frac23$ together with an invariant divisor with
$(m_F,d_F,a_F)=(1,1,2)$, so that $T_0=\Pi\cup\tau(\Pi)\cup\Pi'$ with an invariant plane $\Pi'$; and, without pairs, the
possibilities listed for No.~$68$.
\item No.~$74$: $\sum m_Fd_F\le\frac{19}7$, and a pair costs at least $\frac{10}7$. The possible configurations are: a pair with
$d_F=1$, alone or together with a $\tau$-invariant divisor with $m_F=d_F=1$, so that $T_0$ consists of two planes interchanged
by $\tau$ and possibly a $\tau$-invariant plane; or no pair, and then the possibilities listed for No.~$68$.
\item No.~$76$: $\sum m_Fd_F\le\frac{17}7$. A pair with $d_F=2$ costs at least $\frac{12}7$, so it comes alone, and then
$4(2m_F-a_F)=2$ has no solution with $m_F,a_F\in\frac17\ZZ$. So every pair has $d_F=1$ and maps to two planes interchanged by
$\tau$; besides one or more such pairs there is at most one $\tau$-invariant divisor with $d_F>0$, which then has $m_F=d_F=1$
and maps to a plane. Without pairs we have the possibilities listed for No.~$68$.
\end{itemize}

\emph{No.~$58$ with a non-cyclic point of index $3$.} Suppose now that the two entries $[3,1]$ of No.~$58$ form one point of
type $cA/3$ or $cD/3$. Over it $Y$ has two points $P_1$, $P_2=\tau(P_1)$ of index $3$, each with basket $2\times[3,1]$, and $Y$ is
Gorenstein elsewhere. We first show:
\begin{itemize}
\item[$(*)$] \emph{if $F$ is a $p$-exceptional divisor with centre $P_i$ and $q(F)$ is a surface, then $m_F\ge1$.}
\end{itemize}
By \cite[Theorem~6.1]{Rei87}, the index-one cover of $(Y,P_i)$ is a hypersurface $\{\phi=0\}\subset\CC^4$ with coordinates
$x,y,z,u$ on which $\mumu_3$ acts with weights $(2,1,1,0)$, $\phi$ being invariant, where either $\phi=xy+f(z^3,u)$
($cA/3$), or $\phi=u^2+x^3+h(y,z)+\dots$ with a non-zero binary cubic form $h$ and further terms of larger weighted order
($cD/3$). The weight of the sum of the coordinates minus that of $\phi$ is $1$, so, as $2B\sim-K_Y$, the local
sections of $\OO_Y(B)$ are the functions of weight $2$. Extend the valuation $v$ of $F$ to the cover, normalised to agree
with $v$ on invariant functions; then $v(x),v(y),v(z)\in\frac13\ZZ_{>0}$, $v(u)\ge1$, and $m_F\in\frac13\ZZ$. If $m_F=\frac13$, only $x$
can contribute to the leading forms of value $m_F$, so $q(F)$ is a point. If $m_F=\frac23$, the leading forms are combinations of
the initial forms of $x$, $y^2$, $yz$, $z^2$, as all other functions of weight $2$ have value at least $1$; for $q(F)$ to be a
surface we need $v(x)=\frac23$, $v(y)=v(z)=\frac13$, and the initial forms $X,Y,Z$ of $x,y,z$ in the graded algebra of $v$ must
satisfy no non-zero weighted homogeneous relation, as in the proof of Lemma~\ref{lem:kawleading}. But $\phi=0$ gives such a
relation in value $1$: $XY+cZ^3=0$ for some $c\in\CC$ in the case $cA/3$ (the terms of $f$ involving $u$ have value at least $2$),
and $h(Y,Z)=0$ in the case $cD/3$ ($u^2$ and $x^3$ having value at least $2$). This proves $(*)$.

If $\dim W=2$, then $\deg W\ge3$ by Lemma~\ref{lem:N3quad}, so $B\cdot\Phi\le\frac{10}9$ for a general fibre $\Phi$. Let $C$ be a
component of $\Phi$ and $\hat C$ its proper transform. Then $B\cdot C=\sum_Fm_F\,F\cdot\hat C>0$, and every $F$ meeting the general
$\hat C$ has $q(F)=W$; for such $F$ we have $m_F\ge1$, by $(*)$ if its centre is $P_1$ or $P_2$ and because $m_F$ is a positive
integer otherwise. So $B\cdot C\ge1$ for every component, $\Phi=C$ is irreducible, and exactly one exceptional divisor $F_0$ meets
$\hat C$, with $F_0\cdot\hat C=1$ and $m_{F_0}=1$, as $m_{F_0}\in\frac13\ZZ$. Since the divisors $F$ with $q(F)=W$ meeting the
general $\hat C$ form a $\tau$-invariant set, $\tau(F_0)=F_0$, so the centre of $F_0$ is not $P_1$ or $P_2$. As $F_0$ is the only
exceptional divisor meeting $\hat C$, the curve $\Phi=C$ does not pass through $P_1,P_2$, and $B\cdot\Phi=1$; the argument given
above for an irreducible fibre not passing through the points $P_i$ shows that $W$ is a plane or a quadric, a contradiction.

If $\dim W=3$, then $\delta\le3$, and by $(*)$ a pair of divisors with centres $P_1$, $P_2$ and $d_F>0$ costs at least $2$ in
$\sum m_Fd_F\le\frac{10}3-\delta$. For $\delta=3$ and $\delta=2$ no pair occurs, and the arguments above apply. For $\delta=1$ a pair
has $m_F=d_F=1$ and leaves a budget of $\frac13$, so no other divisor has $d_F>0$; then $\sum(2m_F-a_F)d_F=2$ gives $a_F=1$, and
$T_0$ is a pair of planes interchanged by $\tau$, contrary to Lemma~\ref{lem:fixedlines}. Without pairs we have the
possibilities listed for No.~$68$, which are excluded as before. This completes the proof.
\end{proof}

\subsection{The numerical type No.~47}\label{ss:47}
Here $\qQ=2$, $\qW=1$, $T_X=\ZZ/2$, $\B(X)=[3\times[2,1],[10,3]]$ and $A_X^3=\frac65$. We use the notation of
\S\ref{ss:largedeg}. The generator $\sigma$ of $T_X$ has local type $5K_X$ at the point of index $10$, so by
Lemma~\ref{lem:cover} this point has a single preimage $P\in Y$, of type $\frac15(1,2,3)$ and fixed by $\tau$, and
$\B(Y)=[5,2]$. Moreover $-K_Y\sim2B$, $B^3=\frac{12}5$ and $h^0(X,A_X)=h^0(X,A_X+\sigma)=2$, so $N=3$ and $\tau$ fixes two skew
lines $L_\pm\subset\PP^3$. We choose $A_X$ with local type $7K_X$ at the point of index $10$ (the divisor of
Table~\ref{tab:main} is $A_X+\sigma$). By \eqref{eq:h0}, for $k=1,\dots,5$,
\begin{gather*}
h^0(Y,\OO_Y(kB))=4,\ 12,\ 27,\ 52,\ 89,\\
\bigl(h^0(X,kA_X),h^0(X,kA_X+\sigma)\bigr)=(2,2),\ (6,6),\ (13,14),\ (26,26),\ (44,45).
\end{gather*}
Near $P$ we use orbifold coordinates $y_1,y_2,y_3$ of weights $(1,2,3)$; as for No.~$71$, the local sections of $\OO_Y(B)$
are the functions of weight $3$. We use a $\tau$-equivariant resolution $p\colon\hat Y\to Y$ as in Lemma~\ref{lem:genus}, and,
when $\dim W=3$ and $\deg\varphi=1$, the notation $T_0$ of Lemma~\ref{lem:fixedlines}.

\begin{lemma}\label{lem:47local}
Let $F$ be a $p$-exceptional divisor with centre $P$ and $d_F>0$. Then $m_F\ge\frac35$, and if $m_F=\frac35$, then $q(F)$ is a
plane.
\end{lemma}

\begin{proof}
Let $v$ be the valuation of $F$; its values on the $y_i$ lie in $\frac15\ZZ_{>0}$. The monomials of weight congruent to $3$
modulo $5$ and of degree at most $3$ are $y_3$, $y_1y_2$, $y_1^3$ and $y_2y_3^2$; all other monomials of these weights have
degree at least $4$, hence value at least $\frac45$. As in the proof of Lemma~\ref{lem:leading}, $q|_F$ is given by the classes
of the sections in $\{v\ge m_F\}/\{v>m_F\}$. If $m_F\le\frac25$, these classes are classes of elements of the span of $y_3$
and $y_1y_2$, so $q(F)$ is at most a line. If $m_F=\frac35$, they are classes of elements of
$V=\langle y_3,y_1y_2,y_1^3,y_2y_3^2\rangle\cap\{v\ge\frac35\}$. If the image of $V$ in $\{v\ge\frac35\}/\{v>\frac35\}$ had
dimension $4$, then $y_3\in V$, so $v(y_3)\ge\frac35$ and $v(y_2y_3^2)>\frac35$, and $y_2y_3^2$ would have zero image. So the
image has dimension at most $3$, and $q(F)$, being a surface, is a plane.
\end{proof}

\begin{lemma}\label{lem:47struct}
Suppose that $\dim W=3$, $\deg\varphi=1$, and that $T_0=\Pi_1\cup\Pi_2$ is the union of two $\tau$-invariant planes, the images of two $p$-exceptional divisors $F_1$, $F_2$ with $d_{F_j}=1$ and $m_{F_j}=\mu_j$. Choose coordinates in
which $\tau=\operatorname{diag}(1,1,-1,-1)$, $\Pi_1=\{x_2=0\}\supset L_+=\{x_2=x_3=0\}$ and $\Pi_2=\{x_0=0\}\supset
L_-=\{x_0=x_1=0\}$, and let $L=\Pi_1\cap\Pi_2$.
\begin{enumerate}
\item The prime divisors of $Y$ contracted by $\varphi$ are the components of $\operatorname{div}_Y(x_2)$ and of
$\operatorname{div}_Y(x_0)$. Each of these two divisors is a single $\tau$-orbit of prime divisors, each of which is the
exceptional divisor of the blow up of a curve in $\Pi_1\smallsetminus L$, respectively $\Pi_2\smallsetminus L$, at its general
point. Let $\Gamma_j$ be the union of these curves, a reduced plane curve in $\Pi_j$ not containing $L$, and $I_j$ its ideal.
\item For $k\ge0$ and $a,b\ge0$, every form $G\in(I_1^a\cap I_2^b)_{k+a+b}$ defines a section $G/(x_2^ax_0^b)$ of
$\OO_Y(kB)$, and $G\mapsto G/(x_2^ax_0^b)$ is an isomorphism $(I_1^a\cap I_2^b)_{k+a+b}\cong H^0(Y,\OO_Y(kB))$ for
$a=\lfloor k\mu_1\rfloor$, $b=\lfloor k\mu_2\rfloor$.
\end{enumerate}
\end{lemma}

\begin{proof}
(i) As in the proof of Lemma~\ref{lem:fixedlines}, $2B\sim S=\sum a'_DD$, the sum over the prime divisors $D$ contracted by
$\varphi$, and, as in the case $\delta=1$ of Proposition~\ref{prop:173} (case~(C) there), $S$ is the divisor of the section
$x_0x_2$ of $\OO_Y(2B)$. The divisors $\operatorname{div}_Y(x_2)$ and $\operatorname{div}_Y(x_0)$ are $\tau$-invariant members
of $|B|$ supported on contracted divisors; by Lemma~\ref{lem:descent}(i), (ii) the sum over a $\tau$-orbit of contracted divisors is
linearly equivalent to $cB$ with $c\ge1$, so each of them is a single orbit, with $a'_D=1$. A divisor with $a'_D=1$ is the
exceptional divisor of the blow up of a curve along which $T_0$ is nonsingular, and the components of
$\operatorname{div}_Y(x_2)$ lie over $\Pi_1$.

(ii) By Lemma~\ref{lem:fixedlines}, $Y\setminus S\cong\PP^3\setminus T_0$, so a section of $\OO_Y(kB)$ is a rational form
$G/(x_2^ax_0^b)$ of degree $k$, and such a form is a section if and only if it has no pole along the components of $S$, because
$Y$ is normal and $P$ has codimension $3$. A component $D$ of $\operatorname{div}_Y(x_2)$ has $v_D(x_2)=1$, $v_D(x_0)=0$, and
$v_D(G)$ is the order of $G$ along the corresponding component of $\Gamma_1$; since $I_1=(x_2,g_1)$ is a complete intersection,
the forms of order at least $a$ along $\Gamma_1$ are those of $I_1^a$. This proves the first statement. Conversely, if
$G/(x_2^ax_0^b)$ is a section and $x_2\nmid G$, then its pull-back to $\hat Y$, a section of
$p^*(kB)=kM+k\sum m_FF$, has order $\operatorname{ord}_{F_1}(q^*G)-a+k\mu_1=k\mu_1-a\ge0$ along $F_1$, the proper transform of
$\Pi_1$; similarly along $F_2$. Hence $a\le\lfloor k\mu_1\rfloor$, $b\le\lfloor k\mu_2\rfloor$ when $x_2\nmid G$, $x_0\nmid G$,
and the second statement follows.
\end{proof}

Write $g_j$ for an equation of $\Gamma_j$ in $\Pi_j$, a semi-invariant form in the coordinates of $\Pi_j$, $z_j=\deg g_j$, and
$c_j=g_j|_L$, a binary form in $x_1,x_3$; let $Q_+=L\cap L_+=(0:1:0:0)$ and $Q_-=L\cap L_-=(0:0:0:1)$.

\begin{lemma}\label{lem:47small}
In the situation of Lemma~\ref{lem:47struct}, $(z_1,z_2)$ is $(3,3)$, $(3,4)$ or $(4,3)$; in the first case $\gcd(c_1,c_2)$ has
degree at most $1$, in the second $c_1\mid c_2$, in the third $c_2\mid c_1$. Moreover $g_1$ is $\tau$-invariant, and $g_2$ is
invariant if $(z_1,z_2)=(3,4)$ and anti-invariant otherwise.
\end{lemma}

\begin{proof}
Since $h^0(Y,\OO_Y(B))=4$, Lemma~\ref{lem:47struct}(ii) with $k=1$ shows that no form of degree $2$ in $I_j$ is prime to $x_2$,
respectively $x_0$, so $z_j\ge3$. For $k=2$ we have $a=b=1$. Restricting to $\Pi_1\cup\Pi_2$,
\begin{multline*}
\dim(I_1\cap I_2)_4=\dim S_2\\
+\dim\{(\nu_1,\nu_2)\mid \nu_j\in H^0(\Pi_j,\OO(4-z_j)),\ c_1\nu_1|_L=c_2\nu_2|_L\},
\end{multline*}
since a form lies in $I_j$ if and only if its restriction to $\Pi_j$ is divisible by $g_j$. If $z_1=z_2=3$, the pairs with
$\nu_j|_L=0$ give $2$, and the others $\max(0,e-1)$, where $e=\deg\gcd(c_1,c_2)$; so $\dim=12+\max(0,e-1)$. If $z_1=3$, $z_2=4$,
we get $10+1+1$ if $c_1\mid c_2$ and $11$ otherwise; if $z_1,z_2\ge4$ or $\max(z_1,z_2)\ge5$, at most $11$. Since
$h^0(Y,\OO_Y(2B))=12$, the first statement follows. For the second, $H^0(Y,\OO_Y(2B))$ is spanned by the quadrics, six of
them $\tau$-invariant and four anti-invariant, and by two more sections: $g_1/x_2$ and $g_2/x_0$ if $(z_1,z_2)=(3,3)$, and
$g_1/x_2$ and $G/(x_0x_2)$ with $G|_{\Pi_2}=g_2$ if $(z_1,z_2)=(3,4)$, and $g_2/x_0$ and $G/(x_0x_2)$ with
$G|_{\Pi_1}=g_1$ if $(z_1,z_2)=(4,3)$. Both eigenspaces have dimension
$6$, so both extra sections are anti-invariant, which gives the parities.
\end{proof}

For a point $p\in L$ let $e_j(p)$ be the multiplicity of $c_j$ at $p$, and let $v_p$ be the monomial valuation at $p$ on
$\Pi_1$ with weights $e_2(p)$ on $x_0$ and $1$ on a coordinate along $L$; define $v'_p$ on $\Pi_2$ symmetrically, with
weight $e_1(p)$ on $x_2$.

\begin{lemma}\label{lem:47peel}
Let $N(n,a,b)=\dim(I_1^a\cap I_2^b)_n$. For $a\ge1$,
\begin{multline*}
N(n,a,b)\le N(n-1,a-1,b)\\
+\sum_{i=0}^{d}\max\Bigl(0,\ d-i+1-\sum_{p\in L}\max\bigl(0,\ b\,e_2(p)-a\,w_1(p)-e_2(p)\,i\bigr)\Bigr),
\end{multline*}
with $d=n-z_1a$,
where $w_1(p)\ge v_p(g_1)$, for instance $w_1(p)=e_1(p)$; symmetrically for $b\ge1$. Moreover $N(n,0,0)=\binom{n+3}3$.
\end{lemma}

\begin{proof}
The forms in $(I_1^a\cap I_2^b)_n$ divisible by $x_2$ are $x_2(I_1^{a-1}\cap I_2^b)_{n-1}$, since $x_2$ is part of the regular
sequence $x_2,g_1$ and does not vanish on any component of $\Gamma_2$. For $G\in I_1^a\cap I_2^b$ we have $G|_{\Pi_1}=g_1^ah$
with $h$ of degree $d$, and $G|_{\Pi_1}\in(x_0,c_2)^b$ locally at every point of $L$. At $p\in L$ with $e_2(p)>0$ the ideal
$(x_0,c_2)^b$ is $(x_0,u^{e_2(p)})^b\subset\{v_p\ge b\,e_2(p)\}$, $u$ a local coordinate along $L$; hence
$v_p(h)\ge b\,e_2(p)-a\,v_p(g_1)$. Writing $h=\sum_ix_0^ih_i$ with binary forms $h_i$ of degree $d-i$, we have
$v_p(h)=\min_i(e_2(p)\,i+\operatorname{ord}_ph_i)$, so $\operatorname{ord}_ph_i\ge b\,e_2(p)-a\,v_p(g_1)-e_2(p)\,i$, and binary
forms of degree $d-i$ with prescribed orders at distinct points form a space of dimension
$\max(0,d-i+1-\sum\text{orders})$. Finally $v_p(g_1)\le e_1(p)$ because the restriction $c_1$ of $g_1$ to $L$ has order
$e_1(p)$ at $p$.
\end{proof}

\begin{proposition}\label{prop:47}
The numerical type No.~$47$ does not occur.
\end{proposition}

\begin{proof}
By the remarks before Lemma~\ref{lem:dPcan}, $\dim W\ne1$.

\emph{$\dim W=2$.} The proof of Lemma~\ref{lem:N3quad} uses only that the fixed locus of $\tau$ on $\PP^3$ is $L_+\cup L_-$; so
$\deg W\ge3$, and $B\cdot C\le\frac{12}{15}<1$ for every component $C$ of a general fibre. Hence $C$ passes through $P$, and the
computation of the proof of Lemma~\ref{lem:surflocal} for points of type $\frac15(1,2,3)$, as for No.~$71$, shows that $W$ is a
point or lies in a plane, a contradiction.

\emph{$\dim W=3$.} Then $W=\PP^3$ and $\delta=\deg\varphi\le2$. A divisor $F$ with $d_F>0$ and centre $P$ has
$m_F\ge\frac35$ by Lemma~\ref{lem:47local}, and $a_F\in\frac15\ZZ_{>0}$; otherwise $m_F$ and $a_F$ are positive integers. If
$\delta=2$, then $\sum m_Fd_F\le\frac25$, so all $d_F$ vanish, $g=1$, and Lemma~\ref{lem:N3delta2} gives a contradiction.

Let $\delta=1$, so that $\sum m_Fd_F\le\frac75$ and $\sum(2m_F-a_F)d_F=2$. A divisor with Gorenstein centre or a curve as centre
would have $m_F=d_F=1$ and be the only one, and then $a_F=0$; so all the divisors with $d_F>0$ have centre $P$. By
Lemma~\ref{lem:47local} a single one has $d_F=1$, so that $T_0$ is a plane, which is impossible by Lemma~\ref{lem:fixedlines}.
Otherwise there are two of them, with $d_F=1$ and $\{m_{F_1},m_{F_2}\}=\{\frac35,\frac35\}$ or $\{\frac35,\frac45\}$, and $T_0$
is the union of two planes; by Lemma~\ref{lem:fixedlines} they are $\tau$-invariant and contain $L_+$ and $L_-$ respectively. We
are in the situation of Lemma~\ref{lem:47struct}, with $(\mu_1,\mu_2)=(\frac35,\frac35)$, $(\frac35,\frac45)$ or
$(\frac45,\frac35)$.

The root patterns $(e_1(p),e_2(p))_{p\in L}$ allowed by Lemma~\ref{lem:47small} are finitely many. For each of them and each
$(\mu_1,\mu_2)$, Lemmas~\ref{lem:47struct}(ii) and~\ref{lem:47peel}, with $w_j(p)=e_j(p)$, give $h^0(Y,\OO_Y(kB))$ smaller than
the value above for some $k\le5$ (supplementary file \texttt{examples/no47/peel47.py}), except in the following cases:
\begin{enumerate}
\item[(a)] $(z_1,z_2)=(3,3)$, and at some point $p$ one of $c_1,c_2$ has a triple root and the other a simple root;
\item[(b)] $(z_1,z_2)=(3,4)$, $(\mu_1,\mu_2)=(\frac45,\frac35)$, $c_1=u^3$ and $c_2=u^4$ for a point $p=\{u=0\}$ of $L$; or the
symmetric case $(z_1,z_2)=(4,3)$, $(\mu_1,\mu_2)=(\frac35,\frac45)$.
\end{enumerate}
In case (a), by Lemma~\ref{lem:47small} $c_1=x_1(\alpha x_1^2+\beta x_3^2)$ and $c_2=x_3(\gamma x_1^2+\delta x_3^2)$. A triple
root of $c_1$ is $Q_-$, with $\beta=0$, and then $c_2(Q_-)=0$ forces $\delta=0$ and a double root of $c_2$ at $Q_-$; similarly
for a triple root of $c_2$, at $Q_+$. So (a) does not occur.

In case (b) the linear map $(x_0,x_1,x_2,x_3)\mapsto(x_2,x_3,x_0,x_1)$ commutes with $\tau$ up to sign and exchanges the two
cases, so we may assume $(z_1,z_2)=(3,4)$. By Lemma~\ref{lem:47small}, $g_1$ and $g_2$ are $\tau$-invariant, so $c_1=x_1^3$,
$p=Q_-$ and $c_2=x_1^4$. If $\Gamma_2$ is nonsingular at $Q_-$, then $v'_{Q_-}(g_2)=3$, and Lemma~\ref{lem:47peel} with
$w_2(Q_-)=3$ gives $h^0(Y,\OO_Y(3B))\le26$ (\texttt{peel47.py}). So $\Gamma_2$ is singular at $Q_-$, and
\begin{gather*}
g_1=x_1^3+a_1x_0^3+a_2x_0^2x_1+a_3x_0x_1^2+a_4x_0x_3^2,\\
g_2=x_1^4+b_1x_1^2x_2^2+b_2x_2^2x_3^2+b_3x_1^2x_2x_3+b_4x_2^3x_3+b_5x_2^4 .
\end{gather*}
The dimensions of $(I_1^a)_n$ and $(I_2^b)_n$ do not depend on the coefficients, so $h^0(Y,\OO_Y(kB))$ is determined by the rank of
the matrix of coefficients of the products of the generators $x_2^jg_1^{a-j}$ and $x_0^jg_2^{b-j}$ with monomials, which span
$(I_1^a)_n+(I_2^b)_n$. Gaussian elimination on pivots which are non-zero constants, followed by the computation of the minors
of the remaining matrix, with the coefficients as indeterminates (\texttt{examples/no47/rank47.m2}), shows: $h^0(Y,\OO_Y(4B))=52$
if and only if $a_4b_2=0$ (and $51$ otherwise), and, if $b_2=0$, then $h^0(Y,\OO_Y(5B))=89$ if and only if $a_4b_3=a_4^2b_4=0$ (and
$88$ otherwise). So either $a_4=0$, or $b_2=b_3=b_4=0$. If $a_4=0$, then $g_1\in\CC[x_0,x_1]$ and $\Gamma_1$ consists of three
lines through $Q_-$, each of them $\tau$-invariant. If $b_2=b_3=b_4=0$, then $g_2=x_1^4+b_1x_1^2x_2^2+b_5x_2^4$, and since $\Gamma_2$
is reduced, it consists of four lines $x_1=\pm\lambda x_2$, $x_1=\pm\lambda'x_2$ through $Q_-$, which form two $\tau$-orbits. Both
contradict Lemma~\ref{lem:47struct}(i), by which the components of $\Gamma_j$ form a single $\tau$-orbit.
\end{proof}

\subsection{The numerical type No.~98}\label{ss:98}
Here $\qQ=2$, $\qW=1$, $T_X=\ZZ/2$, $\B(X)=[3\times[2,1],[5,2],[10,3]]$ and $A_X^3=\frac85$. We keep the notation of \S\ref{ss:47}. By
Lemma~\ref{lem:cover}, the point of index $10$ has one preimage $P_0$, fixed by $\tau$, and the point of index $5$, at which $\sigma$
is Cartier, has two preimages $P_1$, $P_2$, interchanged by $\tau$; all three are of type $\frac15(1,2,3)$, and $\B(Y)=[3\times[5,2]]$.
Moreover $-K_Y\sim2B$, $B^3=\frac{16}5$, $h^0(X,A_X)=h^0(X,A_X+\sigma)=2$, so $N=3$ and $\tau$ fixes two skew lines $L_\pm\subset\PP^3$, and
\[
h^0(Y,\OO_Y(kB))=4,\ 14,\ 33,\ 66\qquad(k=1,\dots,4).
\]
In suitable orbifold coordinates of weights $(1,2,3)$ at each $P_i$, the local sections of $\OO_Y(B)$ are the functions of weight $3$;
so Lemma~\ref{lem:47local} holds for divisors with centre any $P_i$.

The following lemma describes the possible surfaces $T_0$ in the birational case.

\begin{lemma}\label{lem:T0orbits}
Suppose that $N=3$, $\dim W=3$ and $\deg\varphi=1$, and let $O$ run through the $\tau$-orbits of prime divisors of $Y$ contracted by
$\varphi$, with $\sum_{D\in O}D\sim c_OB$ and common discrepancy $a'_O$ over $\PP^3$, so that $c_O,a'_O\ge1$ and $\sum_Oa'_Oc_O=2$.
\begin{enumerate}
\item If $\Pi\subset T_0$ is a $\tau$-invariant plane with equation $\ell_\Pi$, then $\operatorname{div}_Y(\ell_\Pi)$ is a single
orbit $O_\Pi$ with $c_{O_\Pi}=1$, and different invariant planes give different orbits. In particular $T_0$ contains at most two
$\tau$-invariant planes.
\item $T_0$ does not contain both a $\tau$-invariant plane and an irreducible $\tau$-invariant quadric.
\item If $T_0$ is an irreducible $\tau$-invariant quadric $Q=\{q=0\}$, then there is a single orbit $O$, with $a'_O=1$ and $c_O=2$; it
is $\operatorname{div}_Y(q)$, and each $D\in O$ is the exceptional divisor of the blow up of a curve on $Q$ along which $Q$ is
nonsingular, at its general point.
\end{enumerate}
\end{lemma}

\begin{proof}
As in the case $\delta=1$ of Proposition~\ref{prop:173}, $2B\sim S=\sum a'_DD$, which gives $\sum a'_Oc_O=2$ with $c_O\ge1$ by
Lemma~\ref{lem:descent}(i), (ii). A $\tau$-invariant surface of $T_0$ of degree $t$ with equation $f$ is not the image of a prime
divisor of $Y$, so $\operatorname{div}_Y(f)\in|tB|$ is $\tau$-invariant and supported on contracted divisors: it is
$\sum n_O\sum_{D\in O}D$ with $\sum n_Oc_O=t$. For a plane, $t=1$ gives (i); two invariant planes with the same orbit would have equal
divisors on $Y$, so that their ratio would be constant. Since $\sum a'_Oc_O=2$, there are at most two orbits with $c_O=1$. If $T_0$
contains an invariant plane $\Pi$ and an irreducible invariant quadric $\{q=0\}$, then the orbits are $O_\Pi$ and at most one further
orbit $O'$, with $c_{O'}=a'_{O'}=1$; then $\sum_{D\in O'}D\in|B|$ is the divisor of a linear form $\ell'$, and $\operatorname{div}_Y(q)$, of
class $2B$, is $\operatorname{div}_Y$ of $\ell_\Pi^2$, $\ell_\Pi\ell'$ or $\ell'^2$, so that $q$ would be reducible. This proves (ii), and
the same argument shows in (iii) that there is no orbit with $c_O=1$. So there is a single orbit, with $a'_O=1$, $c_O=2$. As in the
case $\delta=1$ of Proposition~\ref{prop:173} (case~(A) there), $S=\operatorname{div}_Y(q)$, and a divisor with $a'_D=1$ is the
exceptional divisor of the blow up of a curve along which $T_0$ is nonsingular.
\end{proof}

\begin{lemma}\label{lem:quadcount}
In the situation of Lemma~\ref{lem:T0orbits}(iii), let $F$ be the divisor with $q(F)=Q$, $\mu=m_F$, and let $\Gamma\subset Q$ be the
union of the centres of the divisors in $O$. For $n,j\ge0$ let $N(n,j)$ be the dimension of the space of forms of degree $n$ on
$\PP^3$ vanishing to order at least $j$ along $\Gamma$. Then $h^0(Y,\OO_Y(kB))=N(k+2j,j)$ with $j=\lfloor k\mu\rfloor$, and
\begin{gather*}
N(n,0)=\tbinom{n+3}3,\qquad N(n,1)=\tbinom{n+1}3+h^0(Q,\OO_Q(n)(-\Gamma)),\\
N(n,j)\le N(n-2,j-1)+h^0(Q,\OO_Q(n)(-j\Gamma)),
\end{gather*}
where $\OO_Q(n)(-j\Gamma)$ is the reflexive sheaf.
\end{lemma}

\begin{proof}
The first statement is proved as Lemma~\ref{lem:47struct}(ii): a section of $\OO_Y(kB)$ is $G/q^j$ with $G$ of degree $k+2j$, it has no
pole along $D\in O$ if and only if $G$ vanishes to order at least $j$ along the corresponding component of $\Gamma$, because
$v_D(q)=1$, and its order along $F$ is $k\mu-j$ if $q\nmid G$. The quadric $Q$ is projectively normal, so restriction to $Q$ maps the
forms of degree $n$ onto $H^0(Q,\OO_Q(n))$ with kernel $q\cdot S_{n-2}$; a form vanishing to order $j$ along $\Gamma$ restricts to a
section vanishing to order $j$ along $\Gamma$, and if it is divisible by $q$, the quotient vanishes to order $j-1$, since $q$ vanishes
to order $1$ along $\Gamma$. For $j=1$ every form whose restriction vanishes on $\Gamma$ vanishes on $\Gamma$.
\end{proof}

\begin{proposition}\label{prop:98}
The numerical type No.~$98$ does not occur.
\end{proposition}

\begin{proof}
By the remarks before Lemma~\ref{lem:dPcan}, $\dim W\ne1$.

\emph{$\dim W=2$.} As for No.~$47$, $\deg W\ge3$, so $B\cdot\Phi\le\frac{16}{15}$ for a general fibre $\Phi$. A component of $\Phi$ which
does not pass through the points $P_i$ has $B\cdot C\ge1$, and a component through some $P_i$ has $B\cdot C\ge\frac15$; so either
$\Phi=C$ with $B\cdot C=1$ does not pass through the $P_i$, and the argument of Proposition~\ref{prop:173} shows that $W$ is a plane or a
quadric, a contradiction; or every component passes through some $P_i$ and has $B\cdot C\le1$. In the notation of the proof of
Lemma~\ref{lem:surflocal}, $\nu_0=0$ and $\beta_i\le1$. For $\beta_i\le\frac45$ the solutions are those listed there for points of type
$\frac15(1,2,3)$; for $\beta_i=1$ the same computation gives exactly two more,
$(c_i,E_i\cdot\tilde C,\nu_i)=(\frac35,\frac53,0)$ and $(\frac35,\frac56,\frac12)$, where in the second case $\tilde C$ passes
through $Q_i$ and, as there, the blow up of $Q_i$ gives $m=\frac12$, $G_i\cdot\hat C=1$ and $\mathcal N_2\cdot\hat C=0$. The end of
that proof then shows that $W$ is a point or lies in a plane, a contradiction.

\emph{$\dim W=3$.} Then $W=\PP^3$ and $\delta=\deg\varphi\le3$. A divisor $F$ with $d_F>0$ and centre some $P_i$ has $m_F\ge\frac35$
and $a_F\le3m_F-1$ by Lemmas~\ref{lem:47local} and~\ref{lem:discbound}; one with Gorenstein centre or a curve as centre has integers
$m_F\ge1$ and, if $m_F=1$, $a_F\le2$ by Lemma~\ref{lem:leading}(ii), (iii). If $\delta=3$, then $\sum m_Fd_F\le\frac15$, all $d_F$
vanish, and Lemma~\ref{lem:genus}(iii) gives $\delta=2$. If $\delta=2$, then $\sum m_Fd_F\le\frac65$, so a divisor with Gorenstein centre or a curve as centre has
$m_F=1$, and every divisor with $d_F>0$ has $a_F\le3m_F-1$. Hence $2g-2=\sum(a_F-2m_F)d_F\le\sum(m_F-1)d_F\le\frac65-1$ if some
$d_F>0$, and $2g-2=0$ otherwise; so $g\in\{0,1\}$. If $g=0$, then $M|_\Gamma$ has degree
$2$ on a rational curve and $h^0(\Gamma,M|_\Gamma)=3$, contrary to Lemma~\ref{lem:genus}(iii); if $g=1$, Lemma~\ref{lem:N3delta2}
gives a contradiction.

Let $\delta=1$, so that $\sum m_Fd_F\le\frac{11}5$. Each component of $T_0$ is the image of one divisor $F$, of degree $d_F$, and
$m_F\ge\frac35$, so $T_0$ has degree at most $3$. Components which are not $\tau$-invariant come in pairs $\Pi,\tau(\Pi)$: two
interchanged quadrics would cost at least $\frac{12}5$, and two interchanged planes leave a budget of at most $1$, so that $T_0$ would be
a union of planes of which at most one is $\tau$-invariant, which is impossible by Lemma~\ref{lem:fixedlines}. So all components are
invariant. A single plane is excluded by Lemma~\ref{lem:fixedlines}, three planes and a plane together with a quadric by
Lemma~\ref{lem:T0orbits}, and an irreducible cubic surface because its divisor would have $m_F\le\frac{11}{15}$, hence $m_F=\frac35$
and a plane as image. There remain two cases.

\emph{(I) $T_0$ is the union of two invariant planes}, the images of divisors $F_1$, $F_2$ with $d_{F_j}=1$ and
$\mu_j=m_{F_j}\in\{\frac35,\frac45,1,\frac65,\frac75,\frac85,2\}$, $\mu_1+\mu_2\le\frac{11}5$. We are in the situation of
Lemma~\ref{lem:47struct}, (stated for divisors with any centres). As in the proof of Lemma~\ref{lem:47small}, $z_j\ge3$.
If $z_1\ge9$, then a form $G$ of degree $n=2+a+b\le8$ in $I_1^a\cap I_2^b$, where $a,b\ge1$, restricts to zero on $\Pi_1$, so it is
divisible by $x_2^a$, and $G/x_2^a\in(I_2^b)_{2+b}$; this space is $x_0^bS_2\oplus\langle x_0^{b-1}g_2\rangle$ if $z_2=3$ and $x_0^bS_2$ if
$z_2\ge4$. So $h^0(Y,\OO_Y(2B))\le11$, a contradiction; similarly if $z_2\ge9$. For $3\le z_1,z_2\le8$, Lemma~\ref{lem:47peel}, with
the terms of the points of $L$ omitted, gives $h^0(Y,\OO_Y(kB))$ smaller than the value above for some $k\le4$, for every
$(\mu_1,\mu_2)$ (supplementary file \texttt{examples/no98/peel98.py}).

\emph{(II) $T_0=Q$ is an irreducible invariant quadric}, the image of one divisor $F$ with $d_F=2$, so that $m_F\le\frac{11}{10}$;
by Lemma~\ref{lem:47local}, $\mu=m_F\in\{\frac45,1\}$. We use Lemmas~\ref{lem:T0orbits}(iii) and~\ref{lem:quadcount}. If $Q$ is
nonsingular, $\Gamma$ is a divisor of class $(\alpha,\beta)$ on $Q\cong\PP^1\times\PP^1$, $\alpha\le\beta$, and
$h^0(\OO_Q(n)(-j\Gamma))=h^0(\PP^1\times\PP^1,\OO(n-j\alpha,n-j\beta))$; if $Q$ is a cone, $Q\cong\PP(1,1,2)$, $\Gamma\in|\OO(e)|$ and
$h^0(\OO_Q(n)(-j\Gamma))=h^0(\PP(1,1,2),\OO(2n-je))$. If $\mu=\frac45$, the exact value $h^0(Y,\OO_Y(2B))=N(4,1)=10+h^0(\OO_Q(4)(-\Gamma))$
forces $(\alpha,\beta)=(1,4)$ or $(3,3)$, or $e=6$, and then Lemma~\ref{lem:quadcount} gives $h^0(Y,\OO_Y(3B))\le30$, respectively
$h^0(Y,\OO_Y(4B))\le64$. If $\mu=1$, the exact value $h^0(Y,\OO_Y(B))=N(3,1)=4+h^0(\OO_Q(3)(-\Gamma))$ and the bound
$h^0(Y,\OO_Y(2B))\le10+h^0(\OO_Q(4)(-\Gamma))+h^0(\OO_Q(6)(-2\Gamma))$ leave only $(\alpha,\beta)=(0,4)$ and $(1,4)$, for which
$h^0(Y,\OO_Y(3B))\le32$ and $\le30$ (\texttt{examples/no98/quad98.py}). This completes the proof.
\end{proof}

\subsection{The numerical types No.~88 and No.~93}\label{ss:8893}
For No.~$93$: $\qQ=2$, $\qW=1$, $T_X=\ZZ/2$, $\B(X)=[3\times[2,1],[3,1],[10,3]]$ and $A_X^3=\frac{23}{15}$. As for No.~$98$, the point of
index $10$ has one preimage $P_0$ of type $\frac15(1,2,3)$, fixed by $\tau$, and the point of index $3$ two preimages $P_1$, $P_2$ of
type $\frac13(1,1,2)$, interchanged by $\tau$. Moreover $B^3=\frac{46}{15}$ and $h^0(Y,\OO_Y(kB))=4,14,33,64,111,176$ for $k\le6$.

For No.~$88$: $\qQ=2$, $\qW=1$, $T_X=\ZZ/2$, $\B(X)=[2\times[2,1],[3,1],2\times[6,1]]$ and $A_X^3=\frac53$. The two entries of index $6$
have different local types of $A_X$, so they are two cyclic quotient points, at which $\sigma$ is not Cartier; each has one
preimage, of type $\frac13(1,1,2)$ and fixed by $\tau$; call them $P_1$, $P_2$. The point of index $3$ has two preimages $P_3$, $P_4$
of type $\frac13(1,1,2)$, interchanged by $\tau$. Moreover $B^3=\frac{10}3$ and $h^0(Y,\OO_Y(kB))=4,15,36$ for $k\le3$.

In both cases $-K_Y\sim2B$ and $h^0(X,A_X)=h^0(X,A_X+\sigma)=2$, so $N=3$ and $\tau$ fixes two skew lines $L_\pm\subset\PP^3$. At the points
of index $3$ we use orbifold coordinates $x,y_1,y_2$ of weights $(2,1,1)$, in which the local sections of $\OO_Y(B)$ are the functions
of weight $2$, as for No.~$62$. A $p$-exceptional divisor $F$ with centre a point of index $3$ and $d_F>0$ has $m_F\ge\frac23$ by
Lemma~\ref{lem:leading}(i), and if $m_F=\frac23$, then $F$ is the Kawamata divisor by Lemma~\ref{lem:kawleading}, so that $q(F)$ is the
image of $\PP(1,1,2)$ by a linear subsystem of $|\OO(2)|$: a plane or an irreducible quadric.

\begin{proposition}\label{prop:8893}
The numerical types No.~$88$ and No.~$93$ do not occur.
\end{proposition}

\begin{proof}
We argue as for No.~$98$ and indicate the changes. As there, $\dim W\ne1$.

\emph{$\dim W=2$.} Here $\deg W\ge3$ and $B\cdot\Phi\le\frac13B^3$, which is $\frac{46}{45}$, respectively $\frac{10}9$. A component of
$\Phi$ through a non-Gorenstein point has $B\cdot C\ge\frac15$, respectively $\ge\frac13$, so as for No.~$98$ either the argument of
Proposition~\ref{prop:173} applies, or every component passes through non-Gorenstein points and has $B\cdot C\le1$, $\nu_0=0$ and
$\beta_i\le1$. At a point of type $\frac13(1,1,2)$ the solutions with $\beta_i\le\frac23$ are those of the proof of
Lemma~\ref{lem:surflocal}, and for $\beta_i=1$ the only further solution is $(c_i,E_i\cdot\tilde C,\nu_i)=(\frac23,\frac32,0)$; at
$P_0$ the solutions are those used for No.~$98$. In every case the end of the proof of Lemma~\ref{lem:surflocal} shows that $W$ is a
point, lies in a plane, or lies in the closure of the image of a Kawamata divisor $E_i\cong\PP(1,1,2)$ of a point of index $3$ by a
linear subsystem of $|\OO(2)|$, a surface of degree at most $2$. This contradicts $\deg W\ge3$.

\emph{$\dim W=3$.} Then $\delta\le3$. Every divisor with $d_F>0$ costs at least $\frac35$, respectively $\frac23$, in $\sum m_Fd_F$. If $\delta=3$, the budget $B^3-3$ is $\frac1{15}$, respectively $\frac13$, so all $d_F$ vanish and Lemma~\ref{lem:genus}(iii)
gives a contradiction. If $\delta=2$, the budget is $\frac{16}{15}$, respectively $\frac43$; so a divisor with $d_F>0$ whose centre is a
Gorenstein point or a curve has $m_F=1$ and $a_F\le2$ by Lemma~\ref{lem:leading}(ii), (iii), and one with centre a point of index $3$ or
$5$ has $a_F\le3m_F-1$ by Lemma~\ref{lem:discbound}. As for No.~$98$, $2g-2\le\sum m_Fd_F-\sum d_F$, so $g\in\{0,1\}$, which is impossible.

Let $\delta=1$; the budget is $\frac{31}{15}$, respectively $\frac73$, so $T_0$ has degree at most $3$. Two interchanged quadrics would cost
at least $\frac{12}5$, respectively $\frac83$. Two interchanged planes cost at least $\frac65$, respectively $\frac43$, and leave room
only for one more plane, so they are excluded by Lemma~\ref{lem:fixedlines}. So all components of $T_0$ are invariant; a single plane,
three planes and a plane with a quadric are excluded as for No.~$98$, and so is an irreducible cubic surface, because its divisor would
have $m_F\le\frac{31}{45}$, respectively $\frac79$, so that its image would be a plane or a quadric. In the remaining cases the divisors
mapping onto components of $T_0$ are $\tau$-invariant, so their centres are $P_0$ or a Gorenstein point or a curve for No.~$93$, and
$P_1$, $P_2$ or a Gorenstein point or a curve for No.~$88$.

\emph{(I) Two invariant planes,} with $\mu_j\in\{\frac35,\frac45,1,\frac65,\frac75,\frac85,2\}$, respectively
$\mu_j\in\{\frac23,1,\frac43,\frac53,2\}$, and $\mu_1+\mu_2$ at most the budget. As for No.~$98$, $z_j\ge3$, and $z_j\ge9$ gives
$h^0(Y,\OO_Y(2B))\le11$, a contradiction. For $3\le z_1,z_2\le8$ the bound of Lemma~\ref{lem:47peel}, with the terms of the points of $L$
omitted, is smaller than $h^0(Y,\OO_Y(kB))$ for some $k\le6$, respectively $k\le3$ (supplementary file
\texttt{examples/no88\_no93/shapes.py}).

\emph{(II) One invariant quadric,} the image of a divisor with $\mu\in\{\frac45,1\}$, respectively $\mu\in\{\frac23,1\}$. With the
notation of the proof for No.~$98$, the exact values given by Lemma~\ref{lem:quadcount} for $j\le1$ and its bound for $j\ge2$ exclude
every class $(\alpha,\beta)$ and every $e$ (\texttt{shapes.py}); the last cases are excluded at $k=5$ for No.~$93$ and at $k=3$ for
No.~$88$.
\end{proof}

\subsection{The numerical types No.~139 and No.~149}\label{ss:139149}
Here $\qQ=\qW=2$ and $T_X=\ZZ/3$; we keep the notation of \S\ref{ss:cyclic}, with $n=3$ and a generator $\gamma$ of $\mumu_3$. For
No.~$139$, $\B(X)=[2\times[9,4]]$ and $A_X^3=\frac89$; for No.~$149$, $\B(X)=[[3,1],2\times[9,4]]$ and $A_X^3=\frac{11}9$. In both cases
the generator of $T_X$ has local type $3$ or $6$ at the two points of index $9$, so by Lemma~\ref{lem:cover} each of them has one
preimage, of type $\frac13(1,1,2)$ and fixed by $\gamma$; call them $P_1$, $P_2$. For No.~$149$ the point of index $3$, at which
$T_X$ is Cartier, has three preimages $P_3,P_4,P_5$ of type $\frac13(1,1,2)$, permuted by $\gamma$. Since $\pi$ is \'etale over the
points of $X$ at which $T_X$ is Cartier, $P_1$ and $P_2$ are the only points of $Y$ fixed by $\gamma$. Moreover $-K_Y\sim2B$; for No.~$139$, $B^3=\frac83$ and $h^0(Y,\OO_Y(kB))=4,13,30,57$ for $k\le4$, and for No.~$149$,
$B^3=\frac{11}3$ and $h^0(Y,\OO_Y(kB))=4,16,39,75$; for No.~$139$ the eigenspaces of $\gamma$ on $H^0(Y,\OO_Y(2B))$ have dimensions $4,4,5$. In both cases those on
$H^0(Y,\OO_Y(B))$ have dimensions $1,2,1$, so $N=3$ and, in suitable coordinates, $\gamma$ acts on $\PP^3$ as
$\operatorname{diag}(1,\omega,\omega,\omega^2)$. Its fixed locus is the line $\ell=\{x_0=x_3=0\}$ and the points $p_0=(1:0:0:0)$,
$p_3=(0:0:0:1)$. The $\gamma$-invariant planes are $\{x_0=0\}\supset\ell\cup\{p_3\}$, $\{x_3=0\}\supset\ell\cup\{p_0\}$ and the planes
$\{\alpha x_1+\beta x_2=0\}$, which contain $p_0$, $p_3$ and meet $\ell$ in one point. At the points $P_i$ we use the coordinates of
\S\ref{ss:8893}; the statements made there about divisors with centre a point of index $3$ hold here.

\begin{proposition}\label{prop:139149}
The numerical types No.~$139$ and No.~$149$ do not occur.
\end{proposition}

\begin{proof}
By \S\ref{ss:cyclic}, $\dim W\ne1$.

\emph{$\dim W=2$.} If $W$ is a quadric, $\gamma$, of order $3$, preserves each family of lines of $W$, contrary to
\eqref{eq:nofamily}. So $\deg W\ge3$ and $B\cdot\Phi\le\frac89$, respectively $\frac{11}9$, and we conclude as for No.~$88$.

\emph{$\dim W=3$.} Divisors with $d_F>0$ cost at least $\frac23$ each, and an orbit of three of them at least $2$. For No.~$149$ and
$\delta=3$ the budget is $\frac23$, so there is at most one such divisor, $\gamma$-invariant with $m_F=\frac23$, and then
$2g-2=a_F-\frac43$ with $0<a_F\le1$ by Lemma~\ref{lem:discbound}, which is impossible; so all $d_F$ vanish and Lemma~\ref{lem:genus}(iii)
gives a contradiction. If $\delta=2$, the budget is $\frac23$, respectively $\frac53$, and we conclude as for No.~$88$.

Let $\delta=1$; the budget is $\frac53$, respectively $\frac83$. As in the proof of Lemma~\ref{lem:fixedlines}, $Y\setminus S\cong\PP^3\setminus T_0$
$\gamma$-equivariantly, so $\ell\cup\{p_0,p_3\}\subset T_0$; in particular $T_0$ is not a plane. For No.~$149$, an orbit of three divisors
with $d_F>0$ costs at least $3\cdot\frac23d_F$, so it has $d_F=1$ and $m_F=\frac23$, and its images are three planes permuted by
$\gamma$; the product $f$ of their equations is a
semi-invariant cubic whose divisor on $Y$ is $\gamma$-invariant and supported on contracted divisors. In the notation of
Lemma~\ref{lem:T0orbits}, whose proof does not use the order of the group, $\operatorname{div}_Y(f)=\sum n_O\sum_{D\in O}D$ with
$\sum n_Oc_O=3$; as $\sum a'_Oc_O=2$, either there is one orbit, with $c_O=1$, and $f$ is the cube of the linear form with divisor
$\sum_{D\in O}D$, or there are two orbits with $c_O=1$, and $f$ is a product of powers of two linear forms. Both contradict the fact
that $f$ is a product of three different linear forms. So all components of $T_0$ are $\gamma$-invariant, and $T_0$ has degree at most
$2$, respectively $4$. A component of degree $3$ or $4$ would be the image of a divisor with $m_F=\frac23$, the Kawamata divisor of
some $P_i$, whose image is a plane or a quadric. Two invariant irreducible quadrics are impossible: as in the proof of
Lemma~\ref{lem:T0orbits}, the divisor of each of them on $Y$ is a single orbit $O$ with $c_O=2$, and there is only one such orbit, so the
quadrics would coincide. With Lemma~\ref{lem:T0orbits} we are left
with the two cases of the proof for No.~$98$.

\emph{(I) Two invariant planes,} with $\mu_j\in\{\frac23,1,\frac43,\frac53,2\}$ and $\mu_1+\mu_2$ at most the budget. As there, $z_j\ge3$,
$z_j\ge9$ gives $h^0(Y,\OO_Y(2B))\le11$, and for $3\le z_1,z_2\le8$ the bound of Lemma~\ref{lem:47peel} without the terms of the points of
$L$ is smaller than $h^0(Y,\OO_Y(kB))$ for some $k\le2$, except for No.~$139$ with $z_1=z_2=3$ and
$(\mu_1,\mu_2)\in\{(\frac23,\frac23),(\frac23,1),(1,\frac23)\}$ (supplementary file \texttt{examples/no139\_no149\_no162/shapes.py}).
In these cases Lemma~\ref{lem:47peel}, now with the terms of the points of $L$, together with the exact dimensions of
$(I_1^a\cap I_2^b)_{k+a+b}$ for $a,b\le1$ and $k\le2$ (all of which inject into $H^0(Y,\OO_Y(kB))$ by Lemma~\ref{lem:47struct}(ii)),
leaves only the root pattern in which $c_1=uv^2$ and $c_2=u^2v$ for two different points $\{u=0\}$, $\{v=0\}$ of $L$
(\texttt{pattern139.py}). Then $h^0(Y,\OO_Y(2B))=\dim(I_1\cap I_2)_4=13$, and $(I_1\cap I_2)_4$ is spanned by $l_1l_2S_2$, $l_2g_1$,
$l_1g_2$ and one form $G$ with $G|_{\Pi_j}=g_j\nu_j$, $\nu_j$ linear, where $l_j$ is the equation of $\Pi_j$.

If $L\neq\ell$, then $\{\Pi_1,\Pi_2\}$ consists of $\{x_0=0\}$ or $\{x_3=0\}$ and a plane of the pencil, $L$ has exactly two
$\gamma$-fixed points, and these must be the roots of $c_1$ and $c_2$; so the characters of $l_j$, $g_j$ and $G$ are determined, and
the dimensions of the eigenspaces of $\gamma$ on $H^0(Y,\OO_Y(2B))$ are $\{6,4,3\}$ or $\{5,5,3\}$, not $\{4,4,5\}$
(\texttt{eig139.m2}). If $L=\ell$, then $\{\Pi_1,\Pi_2\}=\{\{x_0=0\},\{x_3=0\}\}$; we may take $\Pi_1=\{x_0=0\}$, $u=x_1$, $v=x_2$, and
the semi-invariant cubics with these restrictions to $\ell$ are, up to scaling the coordinates, $g_1=x_1x_2^2+x_3^3$ and
$g_2=x_1^2x_2+x_0^3$. By Lemma~\ref{lem:47struct}(ii) the ring $\bigoplus_kH^0(Y,\OO_Y(kB))$ is then determined. A computation in
Macaulay2 (\texttt{ring139.m2}, \texttt{hilb139.m2}) shows that it is generated by $x_0,\dots,x_3$, three elements of degree $2$ and two
of degree $3$, for each of the three values of $(\mu_1,\mu_2)$, with the same $14$ relations: the subring generated by these elements
has Hilbert series $(1+t+4t^2+t^3+t^4)/((1-t)^3(1-t^3))$ and the Hilbert function $h^0(Y,\OO_Y(kB))$ for $k\le40$, hence for all $k$,
since both are quasi-polynomials of degree $3$ and period dividing $3$. So $Y\subset\PP(1^4,2^3,3^2)$ is this explicit threefold. But $\gamma$ fixes, besides the two points of index $3$, the
two coordinate points of $x_1$ and $x_2$ of $Y$, at which $\mumu_3$ acts on the tangent space as a scalar (\texttt{fix139.m2},
\texttt{wts139.m2}). This contradicts the fact that $P_1$, $P_2$ are the only fixed points of $\gamma$.

\emph{(II) One invariant quadric,} the image of a divisor with $\mu\in\{\frac23\}$, respectively $\{\frac23,1,\frac43\}$. As for
No.~$98$, Lemma~\ref{lem:quadcount} excludes every class (\texttt{shapes.py}).
\end{proof}

\begin{remark}
The threefold $Y\subset\PP(1^4,2^3,3^2)$ found in case (I) is quasi-smooth, with $-K_Y=\OO_Y(2)$, $\rho(Y)=1$ and two points of type
$\frac13(1,1,2)$, at which the stabilisers act with type $\frac19(1,4,5)$; its numerical invariants and the dimensions of all the
eigenspaces of $\mumu_3$ on its coordinate ring agree with those required for No.~$139$ (\texttt{chart139.m2}, \texttt{check139.m2}).
The quotient is not terminal only because of the two further fixed points, whose contributions to the equivariant Riemann--Roch formula
cancel. It has the numerical type of the cover in No.~$43$.
\end{remark}

\subsection{The numerical type No.~162}\label{ss:162}
Here $\qQ=2$, $\qW=1$, $T_X=\ZZ/4$, $\B(X)=[2\times[2,1],[3,1],4\times[4,1]]$ and $A_X^3=\frac56$; we keep the notation of \S\ref{ss:cyclic},
with $n=4$ and a generator $\gamma$ of $\mumu_4$. The entry of index $3$ is the only entry of index divisible by $3$, so it is a cyclic quotient point, at which
$T_X$ is Cartier; by Lemma~\ref{lem:cover} it has four preimages $P_1,\dots,P_4$ of type $\frac13(1,1,2)$, permuted cyclically by $\gamma$. At
every other non-Gorenstein point of $X$ the image of $T_X$ in the local class group has order equal to the index (a generator of $T_X$ has
local type of order $4$ at the entries of index $4$ and of order $2$ at those of index $2$), so its preimages are Gorenstein points of $Y$.
Since $\pi$ is \'etale over the smooth locus of $X$, the points of $Y$ fixed by $\gamma^2$ are isolated. Moreover $-K_Y\sim2B$, $B^3=\frac{10}3$,
$h^0(Y,\OO_Y(kB))=4,15,36$ for $k\le3$, and $\gamma$ acts on $H^0(Y,\OO_Y(B))$ with four different eigenvalues, so that in suitable
coordinates $\gamma=\operatorname{diag}(1,i,-1,-i)$ on $\PP^3$ and $\gamma^2$ fixes the two skew lines $L_+=\{x_1=x_3=0\}$ and
$L_-=\{x_0=x_2=0\}$, which $\gamma$ preserves. At the points $P_j$ we use the coordinates of \S\ref{ss:8893}.

\begin{proposition}\label{prop:162}
The numerical type No.~$162$ does not occur.
\end{proposition}

\begin{proof}
By \S\ref{ss:cyclic}, $\dim W\ne1$. Let $p\colon\hat Y\to Y$ be a $\gamma$-equivariant resolution such that $q=\varphi\circ p$ is a morphism.

\emph{$\dim W=2$.} Since $W$ is not a plane, $d=\deg W\ge2$; if $W$ is a quadric cone, its lines form a unique, hence $\gamma$-invariant,
family, contrary to \eqref{eq:nofamily}. Let $\Phi$ be a general fibre of $q$, over a point of $W$ outside the images of the finitely many $p$-exceptional divisors which do
not dominate $W$, and such that $\Phi$ is smooth with trivial normal bundle. Then $p^*B\cdot\Phi\le\frac{B^3}d\le\frac53<2$. Let $\hat C$
be a component of $\Phi$; then $M\cdot\hat C=0$, and $p^*B\cdot\hat C=\sum_Fm_F\,F\cdot\hat C$, the sum over the $p$-exceptional divisors
meeting $\hat C$, all of which satisfy $q(F)=W$. If the centre of such an $F$ is some $P_j$, then
$m_F\ge\frac23$ by Lemma~\ref{lem:leading}(i), and $m_F=\frac23$ is impossible, since $F$ would be the Kawamata divisor
(Lemma~\ref{lem:kawleading}) and $W$ the image of $\PP(1,1,2)$ by a linear subsystem of $|\OO(2)|$, a plane or a quadric cone; otherwise
$m_F$ is a positive integer. So $m_F\ge1$ for every $F$ with $q(F)=W$; hence every component of $\Phi$ has $p^*B\cdot\hat C\ge1$, $\Phi=\hat C$ is irreducible,
and, since $p^*B\cdot\Phi<2$, exactly one divisor $F_0$ meets $\Phi$, with $F_0\cdot\Phi=1$ and $p^*B\cdot\Phi=m_{F_0}$. Moreover every
divisor $F$ with $q(F)=W$ meets $\Phi$, so $F_0$ is the only $p$-exceptional divisor with $q(F)=W$. As $\Phi$ has trivial normal bundle,
$K_{\hat Y}\cdot\Phi=2g(\Phi)-2$, and
\[
2g(\Phi)-2=(-2p^*B+\textstyle\sum a_FF)\cdot\Phi=a_{F_0}-2m_{F_0}.
\]
If the centre of $F_0$ is a curve, or a Gorenstein point, then $m_{F_0}=1$, so $a_{F_0}$ is even; by Lemma~\ref{lem:leading}(iii) the
centre is not a curve, and by Lemma~\ref{lem:leading}(ii) it is a nonsingular point and $W=q(F_0)$ is a plane, a contradiction. So the
centre of $F_0$ is some $P_j$. But the set of $p$-exceptional divisors $F$ with $q(F)=W$, which is $\{F_0\}$, is $\gamma$-invariant, since $p$
and $q$ are $\gamma$-equivariant; so $\gamma(F_0)=F_0$ and $\gamma(P_j)=P_j$, a contradiction.

\emph{$\dim W=3$.} Then $\delta\le3$. The divisors with $d_F>0$ and centre some $P_j$ come in $\gamma$-orbits of four, which cost at least
$\frac83>B^3-1$; so every divisor with $d_F>0$ has a Gorenstein point or a curve as centre, and $m_F$ is a positive integer, with $a_F\le2$
if $m_F=1$ (Lemma~\ref{lem:leading}(ii), (iii)). If $\delta=3$, all $d_F$ vanish and Lemma~\ref{lem:genus}(iii) gives a contradiction. If
$\delta=2$, there is at most one such divisor, with $m_F=d_F=1$, and $2g-2=a_F-2$; so $g=1$, contrary to Lemma~\ref{lem:N3delta2}.

Let $\delta=1$, so that $\sum m_Fd_F\le\frac73$ and $\sum(2m_F-a_F)d_F=2$. The possibilities are: one divisor with $m_F=2$, $d_F=1$; one
with $m_F=1$, $d_F=2$, $a_F=1$; or two with $m_F=d_F=a_F=1$ (one divisor with $m_F=d_F=1$ would need $a_F=0$). As in the proof of
Lemma~\ref{lem:fixedlines}, $Y\setminus S\cong\PP^3\setminus T_0$ $\gamma$-equivariantly, and since the points of $Y$ fixed by $\gamma^2$
are isolated, $L_+\cup L_-\subset T_0$. So $T_0$ is not a plane. If $T_0$ is the union of two planes, each contains one of the skew lines
$L_\pm$; if $\gamma$ interchanged the planes, then, since $\gamma(L_+)=L_+$, both would contain $L_+$, and the second would contain
$L_+\cup L_-$. So both planes are $\gamma$-invariant, with $\mu_1=\mu_2=1$, and, as for No.~$98$, $z_j\ge3$, $z_j\ge9$ gives
$h^0(Y,\OO_Y(2B))\le11$, and for $3\le z_1,z_2\le8$ the bound of Lemma~\ref{lem:47peel} without the terms of the points of $L$ is smaller
than $h^0(Y,\OO_Y(kB))$ for some $k\le3$ (supplementary file \texttt{examples/no139\_no149\_no162/shapes.py}). If $T_0$ is a quadric, it is
irreducible and $\gamma$-invariant, the image of one divisor with $\mu=1$, and Lemmas~\ref{lem:T0orbits}(iii) and~\ref{lem:quadcount},
whose proofs do not use the order of the group, exclude every class of the curve $\Gamma$ (\texttt{shapes.py}).
\end{proof}

\subsection{The numerical type No.~175}\label{ss:175}
Here $\qQ=2$, $\qW=1$, $T_X=\ZZ/6$, $\B(X)=[2\times[2,1],2\times[3,1],2\times[6,1]]$ and $A_X^3=1$. By
Proposition~\ref{prop:B24} and the argument after it, if $X$ is of type No.~$175$, then $\Clt(X)=\langle\sigma\rangle\cong\ZZ/2$,
where, in the order of the basket above, $\sigma$ has local types $(1,1,0,0,3,3)$. We choose $A_X$ with local types
$(0,1,1,1,1,4)$, so that $K_X+2A_X\sim\sigma$. Let $\pi\colon Y\to X$ be the double cover defined by $\Clt(X)$, $\tau$ its
involution and $B=\pi^*A_X$. Then $-K_Y\sim2B$ (Lemma~\ref{lem:cover}), $B^3=2$, $3B$ is Cartier, and by
Theorem~\ref{thm:RR} $h^0(X,A_X)=h^0(X,A_X+\sigma)=1$, so $h^0(Y,B)=2$ (supplementary file
\texttt{examples/no175/no175.py}). By Lemma~\ref{lem:cover} the points of $Y$ over the basket are as follows. Over the two
entries $[2,1]$ there are two smooth points fixed by $\tau$. The two entries $[6,1]$ have different local types of $A_X$,
so they are two points, cyclic since each has a single basket entry; over them there are two points $P_a$, $P_d$ of type
$\frac13(1,1,2)$, fixed by $\tau$, and near $P_a$ the involution $\tau$ is induced by the element $-1$ of
$\mumu_6\subset\operatorname{GL}_3(\CC)$. The two entries $[3,1]$ are either two cyclic points, over which $Y$ has four points of type
$\frac13(1,1,2)$ interchanged in pairs by $\tau$, or one non-cyclic point of index $3$ with basket $2\times[3,1]$, of type
$cA/3$ or $cD/3$ \cite[Remark~2.1]{CH11}, over which $Y$ has two such points interchanged by $\tau$. In both cases
$\tau$ fixes no non-Gorenstein point of $Y$ other than $P_a$ and $P_d$.

We use the following local lemma.

\begin{lemma}\label{lem:gordual}
Let $(V,o)$ be a reduced Gorenstein curve germ with an action of a finite group $G$ which is free on $V\setminus o$, and
suppose that the dualising module $\omega_V$ has a $G$-invariant generator $\alpha$. Let $U=V/G$, with image point $o'$.
\begin{enumerate}
\item If $G\ne1$, then $U$ is singular at $o'$.
\item If $U$ has an ordinary node at $o'$, then $V$ has an ordinary node at $o$.
\end{enumerate}
\end{lemma}

\begin{proof}
Let $A=\OO_{V,o}$, $S=A^G=\OO_{U,o'}$, and let $\bar A$ be the normalisation of $A$. By duality for the finite
morphism $\operatorname{Spec}A\to\operatorname{Spec}S$ \cite[Tag~0FKW]{Stacks}, $\omega_A=\operatorname{Hom}_S(A,\omega_S)$; taking
invariants and using the Reynolds operator, $(\omega_A)^G=\operatorname{Hom}_S(A^G,\omega_S)=\omega_S$. Thus $\alpha$ is, up to a
non-zero constant, the pull-back of a generator $\beta$ of $\omega_S$, as rational differentials. We use Rosenlicht's description
of $\omega_A$ \cite[Ch.~IV]{Ser59}: a rational differential which is regular on $\bar A$ belongs to $\omega_A$.

(i) If $U$ is smooth at $o'$, then $\beta$ is regular, so $\alpha$ is regular on $\bar A$. For every $h\in\bar A$ the
differential $h\alpha$ is regular, so $h\alpha\in\omega_A=A\alpha$, and $h\in A$. Hence $A=\bar A$ and $V$ is smooth at $o$.
But then $G$ acts on a uniformiser $t$ through a non-trivial character, since it acts freely on $V\setminus o$, so
the induced character on $\omega_A/\mathfrak m_A\omega_A$, spanned by the class of $dt$, is non-trivial, and $\omega_A$ has no
invariant generator; this contradicts the invariance of $\alpha$.

(ii) If $U$ has an ordinary node at $o'$, then $\beta$ has a simple pole with non-zero residue on each of its two
branches, and so $\alpha$ has a simple pole on every branch of $V$. Let $\bar{\mathfrak m}\subset\bar A$ be the ideal of
functions vanishing at all points of the normalisation over $o$. For $h\in\bar{\mathfrak m}$ the differential $h\alpha$ is
regular, hence lies in $A\alpha$, so $\bar{\mathfrak m}\subset A$ and $A=\CC+\bar{\mathfrak m}$: the germ $V$ is the ordinary
(seminormal) $r$-fold point, where $r\ge2$ is its number of branches. This germ is Gorenstein only for $r=2$
\cite[Proposition~A.2]{Smy11}, so $V$ is an ordinary node.
\end{proof}

\begin{proposition}\label{prop:175}
The numerical type No.~$175$ does not occur.
\end{proposition}

\begin{proof}
\emph{Step 1: the base curve.} Let $s_0,s_1$ span $H^0(Y,\OO_Y(B))$, with $s_0$ coming from $H^0(X,A_X)$ and $s_1$ from
$H^0(X,A_X+\sigma)$; so $\tau$ acts on them with different characters. The proof of Lemma~\ref{lem:descent} applies, as
$\Clt(X)\cong\ZZ/2$. If $|B|$ had a fixed part $F\ne0$, it would be $\tau$-invariant, so $F\sim cB$ with $c\ge1$ by
Lemma~\ref{lem:descent}(i), and $B-F$ effective would force $F\sim B$ and $s_0$, $s_1$ proportional. So $s_0,s_1$ is a
regular sequence. By the Kawamata--Viehweg vanishing theorem and Serre duality, $H^i(Y,\OO_Y(kB))=0$ for $i=1,2$ and
all $k\in\ZZ$, so $R=\bigoplus_{k\ge0}H^0(Y,\OO_Y(kB))$ is Cohen--Macaulay with $\omega_R=R(-2)$ \cite{GW78}, and
$R/(s_0,s_1)$ is Gorenstein of dimension $2$ with $a$-invariant $0$. Hence the curve $C=\{s_0=s_1=0\}$ is
Cohen--Macaulay and Gorenstein, connected, with $\omega_C\cong\OO_C$ and $h^1(\OO_C)=1$; and for every proper subcurve
$Z\subsetneq C$, $H^1(Z,\OO_Z)$ is dual to $\operatorname{Hom}_C(\OO_Z,\OO_C)=0$. Every point where $B$ is not Cartier lies on $C$. If $C$ were integral, the two points $P_a$,
$P_d$ would be singular points of $C$ by Lemma~\ref{lem:gordual}(i). Indeed, near $P_a$ we have $Y=\CC^3/\mumu_3$
with weights $(1,1,2)$, the local sections of $\OO_Y(B)$ are functions of weight $2$ (as $-K_Y\sim2B$), the preimage $V$
of $C$ is reduced (as $C$ is, and the cover is \'etale outside the origin), and $\omega_V$ is generated by the residue of
$du\wedge dv\wedge dw/(fg)$, of weight $(1+1+2)-(2+2)\equiv0$; similarly at $P_d$. This is impossible on an integral
curve of arithmetic genus $1$. So every irreducible
component $\Gamma_j$ of $C$ is a proper subcurve, hence a smooth rational curve. Write $[C]=\sum_jm_j\Gamma_j$ and
$d_j=3B\cdot\Gamma_j$, a positive integer; then
\begin{equation}\label{eq:175deg}
\sum_jm_jd_j=3B\cdot C=3B^3=6 .
\end{equation}

\emph{Step 2: a vanishing on the index-one stack.} Let $p\colon\mathcal Y\to Y$ be the canonical index-one stack of $Y$, a
Deligne--Mumford stack with coarse space $Y$ which near a point of index $3$ is the quotient of its index-one cover by
$\mumu_3$, and let $\mathcal L=\OO_{\mathcal Y}(B)$, an invertible sheaf with $\mathcal L^{\otimes3}=p^*\OO_Y(3B)$ and
$\omega_{\mathcal Y}\cong\mathcal L^{-2}$. Since $p_*$ is exact on coherent sheaves and $p_*\mathcal L^k=\OO_Y(kB)$, we have
$H^i(\mathcal Y,\mathcal L^k)=H^i(Y,\OO_Y(kB))=0$ for $i=1,2$ and all $k$. Let $\mathcal C=\{s_0=s_1=0\}\subset\mathcal Y$, a
Gorenstein curve with coarse space $C$ and $\omega_{\mathcal C}\cong(\omega_{\mathcal Y}\otimes\mathcal L^2)|_{\mathcal C}
\cong\OO_{\mathcal C}$. The Koszul complex $0\to\mathcal L^{-1}\to\OO_{\mathcal Y}^2\to\mathcal L\to\mathcal L|_{\mathcal C}\to0$
and the vanishing above give $H^0(\mathcal C,\mathcal L)=H^0(Y,\OO_Y(B))/\langle s_0,s_1\rangle=0$. For every closed
substack $\mathcal Z\subseteq\mathcal C$, Serre duality on $\mathcal C$ \cite[\S2]{Nir08} gives
\begin{equation}\label{eq:175van}
H^1(\mathcal Z,\mathcal L^{-1})^\vee\cong\operatorname{Hom}_{\mathcal C}(\OO_{\mathcal Z},\mathcal L)\subseteq H^0(\mathcal C,\mathcal L)=0 .
\end{equation}

\emph{Step 3: a degree budget on each component.} Let $\mathcal G_j\subset\mathcal C$ be the reduced substack over
$\Gamma_j$ and $M_j=p_*(\mathcal L^{-1}|_{\mathcal G_j})$, a line bundle on $\Gamma_j\cong\PP^1$; by \eqref{eq:175van},
$H^1(\Gamma_j,M_j)=0$, so $\deg M_j\ge-1$. Multiplication gives a non-zero map $M_j^{\otimes3}\to\OO_{\Gamma_j}(-3B)$. It is an
isomorphism where $B$ is Cartier. At a non-Gorenstein point $P\in\Gamma_j$ it vanishes to some order $e_{jP}\ge1$, since the
local sections of $\mathcal L^{-1}$ are functions of non-zero weight on the index-one cover and vanish at its origin. Hence
$3\deg M_j=-d_j-\sum_Pe_{jP}$, and
\begin{equation}\label{eq:175budget}
d_j+\sum_Pe_{jP}\le3 .
\end{equation}
Let $z$ be a local parameter of $\Gamma_j$ at $P$ and $h$ a local generator of $M_j$, so that $e_{jP}=\operatorname{ord}_z(h^3)$.
If the preimage of $\Gamma_j$ in the index-one cover of $P$ consisted of three branches permuted by $\mumu_3$, each
would map isomorphically onto $\Gamma_j$ near $P$, and $h$ vanishes at its origin, so $e_{jP}=\operatorname{ord}_z(h^3)\ge3$, contrary to
\eqref{eq:175budget}. So the preimage is a single $\mumu_3$-invariant branch, and if $t$ is a parameter on its normalisation,
then $z=t^3$ up to a unit and $\operatorname{ord}_t(h)=\operatorname{ord}_z(h^3)=e_{jP}$. As $h$ vanishes at the origin, the multiplicity
of this branch is at most $e_{jP}$. Let $V_P$ be the preimage of $C$ in the index-one cover of $P$, a complete intersection of two equations of the
weight of $B$, and put $w_P=\sum_jm_je_{jP}$. By the associativity formula for multiplicities,
$\operatorname{mult}_0V_P\le w_P$.

\emph{Step 4: the budget is exhausted.} At a point of type $\frac13(1,1,2)$, with coordinates of weights $1,1,2$, the
two equations of $V_P$ have weight $2$, so their linear parts lie in the one-dimensional span of the coordinate of weight
$2$, and $\operatorname{mult}_0V_P\ge2$. At a non-cyclic point of index $3$ with basket $2\times[3,1]$, the index-one cover is
a hypersurface $\{F=0\}\subset\CC^4$ with $\operatorname{mult}_0F\ge2$ and $\mumu_3$-weights $1,2,1,0$ up to order
\cite[Theorem~(6.1)]{Rei87}; again only one coordinate has the weight $2$ of $B$, so $\operatorname{mult}_0V_P\ge2\cdot1\cdot2=4$.
In both configurations, therefore, $\sum_Pw_P\ge12$. Summing \eqref{eq:175budget} with multiplicities and using
\eqref{eq:175deg} and $d_j\ge1$,
\[
18\le6+\sum_Pw_P=\sum_jm_j\Bigl(d_j+\sum_Pe_{jP}\Bigr)\le3\sum_jm_j\le3\sum_jm_jd_j=18 .
\]
So all inequalities are equalities: $d_j=1$ and $\sum_Pe_{jP}=2$ for every $j$, $\deg M_j=-1$, $w_P=\operatorname{mult}_0V_P$, and
$w_P=2$ at every point of type $\frac13(1,1,2)$.

\emph{Step 5: no gluing.} Let $\mathcal D\subset\mathcal C$ be the reduced substack, with coarse space $D=C_{\mathrm{red}}$, and
$M_D=p_*(\mathcal L^{-1}|_{\mathcal D})$. Restriction to the components gives an exact sequence $0\to M_D\to\bigoplus_jM_j\to
Q\to0$ with $Q$ of finite support. As $H^0(M_j)=0$ and $H^1(D,M_D)=0$ by \eqref{eq:175van}, we get $Q=0$. Two consequences
follow. First, two components of $C$ do not meet at a point where $B$ is Cartier, since there $M_D$ is invertible and the
cokernel is non-zero at an intersection point. Second, at a point $P$ of type $\frac13(1,1,2)$ we have $w_P=\operatorname{mult}_0V_P=2$,
so one equation of $V_P$ has a non-zero linear term, and $V_P$ is a plane curve of multiplicity $2$ in the coordinates
$u,v$ of weight $1$. A singular branch would have $u$, $v$ of order at least $2$ on its normalisation, so every invariant
function would have order at least $6$, whereas the coarse component $\Gamma_j$ is smooth, with a uniformiser of
order $3$. So all
branches are smooth with $e_{jP}=1$, and either one component of multiplicity $2$ passes through $P$, or two components
of multiplicity $1$ do. In the second case their branches have distinct tangent lines: $Q=0$ means that the weight-$1$
functions on $V_P$ restrict independently to the two branches, and their leading terms are linear forms in $u,v$.

\emph{Step 6: the contradiction.} By Step~5, $V_{P_a}$ has either one branch, or two smooth branches with distinct
tangent lines. Since $\tau$ acts near $P_a$ as $-1$ on $\CC^3$, it preserves every line through the origin, hence each of
these branches, and so every component $\Gamma$ of $C$ through $P_a$. By Steps~4 and~5, $e_{\Gamma P_a}=1$ and
$\sum_Pe_{\Gamma P}=2$, so $\Gamma$ contains exactly one other non-Gorenstein point; it is fixed by $\tau$, so it is $P_d$.
The components through $P_a$, which are also the components through $P_d$, use up $w_{P_a}=w_{P_d}=2$; they contain no
other non-Gorenstein point, and by Step~5 they meet no other component. As $C$ is connected, they are all the components of
$C$, and then $\sum_jm_jd_j=2$, contrary to \eqref{eq:175deg}.
\end{proof}

\subsection{The numerical types No.~103 and No.~104}\label{ss:104}
Here $\qQ=2$, $\qW=1$, $T_X=\ZZ/2$, and $\B(X)=[5\times[2,1],2\times[3,1],[6,1]]$ with $A_X^3=\frac43$ (No.~$103$), or
$\B(X)=[5\times[2,1],[3,1],[6,1]]$ with $A_X^3=1$ (No.~$104$). We choose $A_X$ with local types $i_P(A_X)=1$ at every point;
then $i_P(\sigma)=1,3,0$ at the points of index $2,6,3$, where $\sigma$ generates $T_X$ (the divisor of Table~\ref{tab:main}
is $A_X+\sigma$), and $K_X+2A_X\sim\sigma$. Let $\pi\colon Y\to X$ be the double cover, $\tau$ its involution and $B=\pi^*A_X$,
so that $-K_Y\sim2B$. The entry $[6,1]$ is the only entry of index $6$, so it is a cyclic quotient point of type
$\frac16(1,1,5)$ \cite[Remark~2.1]{CH11}; $\sigma$ has local type $3$ there, and its preimage is a point $P_0$ of type
$\frac13(1,1,2)$ fixed by $\tau$, near which $\tau$ is induced by the element $-1$ of $\mumu_6\subset\operatorname{GL}_3(\CC)$. The
entries $[3,1]$ are not in the support of $T_X$, so $\pi$ is \'etale over them. In No.~$104$ the single entry $[3,1]$ is a
cyclic point with two preimages of type $\frac13(1,1,2)$, interchanged by $\tau$. In No.~$103$ the two entries $[3,1]$ are
either two cyclic points, with four such preimages interchanged in pairs, or one non-cyclic point of type $cA/3$ or $cD/3$
\cite[Remark~2.1]{CH11}, with two preimages interchanged by $\tau$. Over the entries $[2,1]$, where $\sigma$ has local type
$1$, the points of $Y$ are Gorenstein, and $B$ is Cartier there. We have $B^3=\frac83$, resp.\ $2$, and
$h^0(Y,\OO_Y(kB))=3,12,29$, resp.\ $3,10,23$, for $k=1,2,3$ (supplementary file \texttt{examples/no103\_no104/data.py}), so
$|B|$ is a net ($N=2$). The fixed locus of $\tau$ is finite (\S\ref{ss:largedeg}).

\begin{proposition}\label{prop:104}
The numerical types No.~$103$ and No.~$104$ do not occur.
\end{proposition}

\begin{proof}
\emph{Step 1: a parity rule.} At every $P\in\B(X)$ we have $i_P(3A_X+\sigma)=3+1\equiv0$, $3\cdot1+3\equiv0$ and
$3\cdot1+0\equiv0$ modulo $2$, $6$ and $3$. At a terminal point the $\QQ$-Cartier Weil divisors modulo Cartier divisors
form a cyclic group generated by $K_X$ \cite[Lemma~2.2.3]{Pro25}, and at Gorenstein points every $\QQ$-Cartier Weil divisor
is Cartier; so $D_X=3A_X+\sigma$ is a Cartier divisor, and $L=\pi^*\OO_X(D_X)=\OO_Y(3B)$. If $\Gamma\subset Y$ is a
$\tau$-invariant irreducible curve, then $\tau$ does not fix $\Gamma$ pointwise, so $\Gamma\to\pi(\Gamma)$ has degree $2$
and $L\cdot\Gamma=D_X\cdot\pi_*\Gamma=2\,D_X\cdot\pi(\Gamma)$ is even. Hence every $\tau$-invariant effective $1$-cycle on $Y$
has even $L$-degree, since the components of such a cycle which are not invariant come in pairs of equal degree.

\emph{Step 2: curve sections.} By the Kawamata--Viehweg vanishing theorem and Serre duality, we have
$H^i(Y,\OO_Y(kB))=0$ for $i=1,2$ and all $k\in\ZZ$, so the ring $R=\bigoplus_{k\ge0}H^0(Y,\OO_Y(kB))$ is Cohen--Macaulay
with $\omega_R=R(-2)$ \cite{GW78}. By Lemma~\ref{lem:decomp}(iii) the net $|B|$ has no fixed components, so two
general sections $x,y\in R_1$ form a regular sequence. Let $C=\{x=y=0\}\subset Y$ be the corresponding curve. Then
$R/(x,y)$ is a Gorenstein ring of dimension $2$ with $a$-invariant $0$ and depth $2$, so
$H^0(C,\OO_C)=\CC$, $H^1(C,\OO_C)\cong\CC$ and $\omega_C\cong\OO_C$. If $Z\subsetneq C$ is a proper subcurve, then
$H^1(Z,\OO_Z)$ is dual to $\operatorname{Hom}_C(\OO_Z,\OO_C)=0$, since a constant on $C$ annihilating the ideal of $Z$ is zero;
in particular every irreducible component of $C$ which is a proper subcurve is a smooth rational curve.

\emph{Step 3: $|B|$ has a fixed curve.} If $\operatorname{Bs}|B|$ were finite, three general sections $x,y,z\in R_1$ would form
a regular sequence, $A'=R/(x,y,z)$ would be Cohen--Macaulay of dimension $1$ with Hilbert series
$P_Y(t)(1-t)^3=1+ct^2+t^3+t^4+ct^5+t^6+t^7+\dots$, of period $3$, where $c=6$, resp.\ $4$, and for $k$ large and divisible by
$3$ the number $\dim A'_k=1$ would be the length of the finite scheme $\operatorname{Bs}|B|$, since $3B$ is Cartier. But $B$ is
not Cartier at $P_0$ and at the points over the entries $[3,1]$, at least three points, all of which are base points. So the
cycle $C=\{x=y=0\}$ has a non-zero fixed part $F$, supported on the base curves of $|B|$; it is $\tau$-invariant, together
with its multiplicities, because $\tau$ preserves $|B|$. In particular $C$ is not integral, so all its irreducible
components are smooth rational curves by Step~2.

\emph{Step 4: one moving component.} By \eqref{eq:nofamily} the image of $\varphi=\varphi_{|B|}$ is not a curve, so it is
$\PP^2$. Let $p\colon\hat Y\to Y$ be a resolution such that $q=\varphi\circ p\colon\hat Y\to\PP^2$ is a morphism, and write
$p^*B=H+\sum m_EE$ and $K_{\hat Y}=p^*K_Y+\sum a_EE$, with $H=q^*\OO(1)$ and $a_E>0$. Every component of the
exceptional locus of $p$ either misses a general fibre of $q$, or dominates $\PP^2$ and is generically finite over it;
hence no component of a general fibre is $p$-exceptional. Let $f$ be such a component. By generic smoothness $f$ is a
smooth curve with trivial normal bundle, $p|_f$ is birational onto its image $M_f$, a component of $C$ of multiplicity one,
and $M_f\ne C$ because $F\ne0$. By Step~2, $M_f$ and $f$ are smooth rational curves, so $K_{\hat Y}\cdot f=-2$. Put
$d_f=B\cdot M_f$. Since $H\cdot f=0$,
\[
d_f=\sum m_E\,(E\cdot f)>0,\qquad 2d_f-2=\sum a_E\,(E\cdot f),
\]
and some $E$ with $m_E>0$ has $E\cdot f>0$; hence $d_f>1$, and $d_f\ge\frac43$ because $3B$ is Cartier. Since
$B\cdot C=B^3\le\frac83$ and $B\cdot F>0$, there is exactly one moving component $M$, of multiplicity one. It lies
over the point $\{x=y=0\}$ of $\PP^2$, so the section $z$ does not vanish at its generic point.

\emph{Step 5: a degree budget.} Let $p\colon\mathcal Y\to Y$, $\mathcal L=\OO_{\mathcal Y}(B)$ and $\mathcal C=\{x=y=0\}$ be
as in Step~2 of the proof of Proposition~\ref{prop:175}. The same argument gives
$H^0(\mathcal C,\mathcal L)=H^0(Y,\OO_Y(B))/\langle x,y\rangle=\CC\bar z$, and for a closed substack $\mathcal Z\subseteq\mathcal C$
the space $H^1(\mathcal Z,\mathcal L^{-1})^\vee\cong\operatorname{Hom}_{\mathcal C}(\OO_{\mathcal Z},\mathcal L)$ consists of the multiples
of $\bar z$ annihilated by the ideal of $\mathcal Z$. If $\mathcal Z$ does not contain $M$, this ideal is the unit ideal at the
generic point of $M$, where $\bar z\ne0$, so the space is zero; in all cases it has dimension at most $1$. Write
$[C]=\sum_jm_j\Gamma_j$ and $d_j=L\cdot\Gamma_j$, so that $\sum_jm_jd_j=3B^3=8$, resp.\ $6$, and $d_M=3B\cdot M\ge4$. Define
$M_j$ and $e_{jP}$ as in Step~3 of the proof of Proposition~\ref{prop:175}. Then $d_j+\sum_Pe_{jP}\le3$ for the fixed
components and $d_M+\sum_Pe_{MP}\le6$ for $M$, as $h^1(M_M)\le1$. As $d_j\ge1$ and $d_M\ge4$, these bounds give $\sum_Pe_{jP}\le2$ for every component, including
$M$, so the single-branch argument of that step applies, and $\operatorname{mult}_0V_P\le w_P=\sum_jm_je_{jP}$.
By Step~4 of that proof, $\operatorname{mult}_0V_P\ge2$ at a point of type $\frac13(1,1,2)$ and $\operatorname{mult}_0V_P\ge4$ at a
non-cyclic point of index $3$ with basket $2\times[3,1]$; the points over the entries $[2,1]$ contribute nothing. So
$\sum_Pw_P\ge10$, resp.\ $6$, and with $\delta=3B^3$,
\[
3\delta-6\le\delta+\sum_Pw_P=\sum_jm_j\Bigl(d_j+\sum_Pe_{jP}\Bigr)\le6+3\sum_{j\ne M}m_j\le6+3\sum_{j\ne M}m_jd_j=6+3(\delta-d_M)\le3\delta-6 .
\]
So all inequalities are equalities: $d_M=4$, $d_j=1$ for every fixed component, $\deg M_j=-1$ for the fixed components,
$\deg M_M=-2$, and $w_P=\operatorname{mult}_0V_P=2$ at every point of type $\frac13(1,1,2)$, in particular at $P_0$.

\emph{Step 6: no gluing.} Let $\mathcal D\subset\mathcal C$ be the reduced substack and $M_D=p_*(\mathcal L^{-1}|_{\mathcal D})$. In
$0\to M_D\to\bigoplus_jM_j\to Q\to0$ we have $H^0(M_j)=0$, so $h^0(Q)=h^1(M_D)-\sum_jh^1(M_j)\le1-h^1(M_M)=0$. As in
Step~5 of the proof of Proposition~\ref{prop:175}, at $P_0$ either one component of multiplicity $2$ passes, or two
components of multiplicity $1$ do, with smooth branches and, in the second case, distinct tangent lines.

\emph{Step 7: the contradiction.} Suppose that a fixed component $\Gamma$ passes through $P_0$. The cycle $F$ is
$\tau$-invariant, so $\tau(\Gamma)$ is a component of $F$ through $P_0$. Since $\tau$ acts near $P_0$ as $-1$ on $\CC^3$, it
preserves every line through the origin, so the branch of $\tau(\Gamma)$ has the same tangent line as that of $\Gamma$; by
Step~6, $\tau(\Gamma)=\Gamma$. Then $\Gamma$ is a $\tau$-invariant irreducible curve with $L\cdot\Gamma=d_\Gamma=1$, which is odd,
contrary to Step~1. So no fixed component passes through $P_0$, and $M$, of multiplicity $1$, is the only component through
$P_0$; then $w_{P_0}=2$ forces $e_{MP_0}=2$, a singular branch, which Step~6 excludes.
\end{proof}

\subsection{The numerical types No.~108, No.~110 and No.~116}\label{ss:108}
In these three rows $\qQ=2$, $\qW=1$, $T_X=\ZZ/2$ and $\B(X)=[5\times[2,1],[6,1],[r,a]]$, where $[r,a]=[5,1]$, $[5,2]$ and
$[7,3]$, and $A_X^3=\frac{19}{15}$, $\frac{16}{15}$ and $\frac{23}{21}$, respectively. We choose $A_X$ with local types $1$ at the
points of index $2$ and $6$ and $2$, $2$, $3$ at the point of index $r$; then $i_P(\sigma)=1,3,0$ at the points of index
$2,6,r$, and $K_X+2A_X\sim\sigma$. Let $\pi\colon Y\to X$, $\tau$ and $B=\pi^*A_X$ be as in \S\ref{ss:104}, so that
$-K_Y\sim2B$ and $B^3=\frac{38}{15}$, $\frac{32}{15}$, $\frac{46}{21}$. The entries $[6,1]$ and $[r,a]$ are the only entries of their
index, so they are cyclic quotient points. Over the point of index $6$ there is a point $P_0$ of type $\frac13(1,1,2)$, with
coordinates $u,v,w$ of weights $1,1,2$, fixed by $\tau$, near which $\tau$ is induced by $-1\in\mumu_6$; over the point of index
$r$, where $\sigma$ is Cartier, there are two points $Q_\pm$ of type $\frac1r(a,-a,1)$, interchanged by $\tau$; the points over
the entries $[2,1]$ are Gorenstein, and $B$ is Cartier there. We call $P_0$, $Q_+$, $Q_-$ the \emph{marked points}. The
divisor $B$ is not Cartier at the marked points, so near each of them the local sections of $\OO_Y(B)$ and of $\OO_Y(-B)$ are
functions of non-trivial characters on the index-one cover. We have $h^0(Y,\OO_Y(kB))=3,12,27$, resp.\ $3,10,23$ (for
Nos.~$110$ and $116$), for $k=1,2,3$, so $|B|$ is a net, and $NA_X+\sigma$ is Cartier at every point of $\B(X)$ for
$N=3r$ (supplementary file \texttt{examples/no108/data.py}).

\begin{proposition}\label{prop:108}
The numerical types No.~$108$, No.~$110$ and No.~$116$ do not occur.
\end{proposition}

\begin{proof}
\emph{Step 1: preliminaries.} As $NA_X+\sigma$ is Cartier, so is $NB$, and as in Step~1 of the proof of
Proposition~\ref{prop:104} every $\tau$-invariant irreducible curve $\Gamma\subset Y$ has $NB\cdot\Gamma$ even. Steps~2 and~3 of
that proof apply: for general $x,y$ in the net, $C=\{x=y=0\}$ is a connected Gorenstein curve with $\omega_C\cong\OO_C$,
$h^1(\OO_C)=1$ and $H^1(\OO_Z)=0$ for every proper subcurve $Z$; if $\operatorname{Bs}|B|$ were finite, the Hilbert function of
$R/(x,y,z)$, given by $P_Y(t)(1-t)^3$, would be $1$ in all degrees divisible by $N$, whereas $\operatorname{Bs}|B|$ contains the
three marked points, at which $B$ is not Cartier. So $C$ has a non-zero $\tau$-invariant fixed part $F$, $C$ is not
integral, and its components are smooth rational curves. As in Step~4 of the proof of Proposition~\ref{prop:104}, every
moving component $M$ of $C$ has $B\cdot M>1$, and the third section $z$ does not vanish at its generic point. Every
marked point lies on $C$, since it is a base point of $|B|$.

\emph{Step 2: the negative twist.} Let $\mathcal C$ and $\mathcal L$ be as in Step~5 of the proof of Proposition~\ref{prop:104}; then
$H^1(\mathcal Z,\mathcal L^{-1})=0$ for every closed substack $\mathcal Z$ of $\mathcal C$ which does not contain the generic point of a
moving component, and $h^1(\mathcal Z,\mathcal L^{-1})\le1$ in general. For a component $\Gamma$ of $C$, with reduced substack
$\mathcal G_\Gamma$, put $N_\Gamma=p_*(\mathcal L^{-1}|_{\mathcal G_\Gamma})$, an invertible sheaf on $\Gamma\cong\PP^1$. At a marked point $P$ of
index $r_P$ on $\Gamma$, the preimage of $\Gamma$ in the index-one cover is either a single $\mumu_{r_P}$-invariant branch, with
parameter $t$ such that $t^{r_P}$ is a parameter of $\Gamma$, and we put $\varepsilon_{\Gamma P}=\frac1{r_P}\operatorname{ord}_t(h)$ for a local
generator $h$ of $N_\Gamma$, or a free orbit of branches, and then $\varepsilon_{\Gamma P}\ge1$. As $NB$ is Cartier, the multiplication
map $N_\Gamma^{\otimes N}\to\OO_\Gamma(-NB)$ gives $\deg N_\Gamma=-B\cdot\Gamma-\sum_P\varepsilon_{\Gamma P}<0$. For a fixed component $\Gamma$ we have
$H^1(N_\Gamma)=0$, so $\deg N_\Gamma=-1$. For a moving component $M$, $\deg N_M\le-B\cdot M<-1$ and $h^1(N_M)\le1$ give $\deg N_M=-2$; a second
moving component $M'$ would have $H^1(N_{M'})=0$, hence $\deg N_{M'}\ge-1$, which is impossible. So $M$ is unique, of
multiplicity one by Step~1. In particular every
component has $\sum_P\varepsilon_{\Gamma P}<1$, so all its branches at the marked points are single invariant branches. As in Step~6
of the proof of Proposition~\ref{prop:104}, the cokernel $Q$ of $N_D\to\bigoplus N_\Gamma$, for $D=C_{\mathrm{red}}$, vanishes: as
$h^0(N_D)=0$ and $h^1(N_D)\le1$,
\[
-1\le\chi(N_D)=\chi\Bigl(\bigoplus_\Gamma N_\Gamma\Bigr)-\operatorname{length}Q=-1-\operatorname{length}Q.
\]
Hence distinct components meet only at marked points.

At $P_0$ the parameter $t$ of the branch of a component $\Gamma$ has character $1$ or $2$. If it has character $1$, then
$\varepsilon_{\Gamma P_0}<1$ forces $u$ or $v$ to have order $1$ on the branch, which is therefore smooth with tangent line in the
$(u,v)$-plane, and branches of this kind have pairwise distinct tangent lines (as in Step~5 of the proof of
Proposition~\ref{prop:175}). If it has character $2$, then $w$ has order $1$ on the branch: otherwise $u,v$ have order at
least $2$ and $w$ order at least $4$, and every invariant function has order at least $6$, whereas $\Gamma$ is smooth, with a
parameter of order $3$.

\emph{Step 3: the positive twist.} Twisting the Koszul resolution of $\OO_{\mathcal C}$ by $\mathcal L^{-1}$ and using
$H^0(Y,\OO_Y(-B))=0$ and $H^1(Y,\OO_Y(-2B))=H^2(Y,\OO_Y(-3B))=0$ (Step~2 of the proof of Proposition~\ref{prop:104}) gives
$H^0(\mathcal C,\mathcal L^{-1})=0$. By duality on $\mathcal C$, as in Step~5 of that proof, for every closed substack $\mathcal Z\subseteq\mathcal C$
\[
H^1(\mathcal Z,\mathcal L)^\vee\cong\operatorname{Hom}_{\mathcal C}(\OO_{\mathcal Z},\mathcal L^{-1})\subseteq H^0(\mathcal C,\mathcal L^{-1})=0 .
\]
For a component $\Gamma$ put $P_\Gamma=p_*(\mathcal L|_{\mathcal G_\Gamma})$, and let $k_\Gamma$ be the number of marked points on $\Gamma$. The
multiplication map $P_\Gamma\otimes N_\Gamma\to\OO_\Gamma$ is an isomorphism away from the marked points, and its image lies in the
maximal ideal at each marked point on $\Gamma$, because local sections of $P_\Gamma$ and $N_\Gamma$ there are functions of non-trivial
character on an invariant branch through the origin, which vanish at the origin. Hence
\[
\deg P_\Gamma+\deg N_\Gamma\le-k_\Gamma .
\]
For a fixed component, $H^1(P_\Gamma)=0$ gives $\deg P_\Gamma\ge-1$, so $k_\Gamma\le2$, and $\deg P_\Gamma=-1$ if $k_\Gamma=2$. For $M$, the
restriction of $\bar z$ is a non-zero section of $P_M$, so $\deg P_M\ge0$ and $k_M\le2$.

\emph{Step 4: no gluing for the positive twist at $P_0$.} Let $\Gamma_1\ne\Gamma_2$ be fixed components with $k_{\Gamma_i}=2$ whose
branches at $P_0$ have character $2$, and let $Z=\Gamma_1\cup\Gamma_2$, with reduced substack $\mathcal Z$. In the exact sequence
$0\to P_Z\to P_{\Gamma_1}\oplus P_{\Gamma_2}\to T\to0$, where $P_Z=p_*(\mathcal L|_{\mathcal Z})$, we have $H^0(P_{\Gamma_i})=0$ by Step~3 and
$H^1(P_Z)=0$, so $T=0$. But near $P_0$ the local sections of $P_Z$ are restrictions of functions $f$ of character $2$ on
$\CC^3$, since taking character parts is exact. On each branch $w$ has order $1$ and generates the character-$2$ part over
the coarse local ring, and every other monomial of character $2$ ($u^2$, $uv$, $v^2$, $w^4$, \dots) has order at least $4$;
so if $f=aw+\dots$, the images of $f$ in the fibres of $P_{\Gamma_1}$ and $P_{\Gamma_2}$ at $P_0$ are both $a$ (with respect to the
generators $w|_{\Gamma_i}$). So $P_Z\to P_{\Gamma_1}\oplus P_{\Gamma_2}$ is not surjective at $P_0$, a contradiction.

\emph{Step 5: the contradiction.} We show that no fixed component $\Gamma$ contains two marked points. Suppose that $\Gamma$
contains $P_0$ and $Q_+$; by Step~3 it does not contain $Q_-$, so $\tau(\Gamma)\ne\Gamma$. As $F$ is $\tau$-invariant, $\tau(\Gamma)$ is a
fixed component through $P_0$ and $Q_-$, and its branch at $P_0$ has the same character as that of $\Gamma$, since $\tau$ commutes
with $\mumu_3$. If this character is $1$, then, as $\tau$ acts near $P_0$ as $-1$, which preserves every line through the
origin, the branches of $\Gamma$ and $\tau(\Gamma)$ at $P_0$ have the same tangent line, contrary to Step~2. If it is $2$, Step~4
applies to $\Gamma$ and $\tau(\Gamma)$. The same argument applies with $Q_-$ in place of $Q_+$. Suppose next that $\Gamma$ contains
$Q_+$ and $Q_-$. If $\tau(\Gamma)\ne\Gamma$, the proper subcurve $\Gamma\cup\tau(\Gamma)$, made of two smooth rational curves meeting at
$Q_+$ and $Q_-$, has $h^1(\OO)\ge1$, contrary to Step~1. If $\tau(\Gamma)=\Gamma$, then $\varepsilon_{\Gamma Q_+}=\varepsilon_{\Gamma Q_-}=\frac er$ for an
integer $e$, and $\deg N_\Gamma=-1$ gives $B\cdot\Gamma=1-\frac{2e}r$, so $NB\cdot\Gamma=3r-6e$ is odd, contrary to Step~1.

So each fixed component contains at most one marked point. As distinct components meet only at marked points and $C$ is
connected, a marked point not on $M$ would lie on a union of fixed components disjoint from $M$. So $M$ contains all three
marked points, contrary to $k_M\le2$.
\end{proof}

\subsection{Further covers with a net or a pencil: Nos.~85, 87, 89, 90, 146, 147, 153, 167}\label{ss:netpencil}
The proof of Proposition~\ref{prop:108} uses the cover only through a few local and global properties. We isolate
them and apply them to eight more rows.

Let $T_X\cong\ZZ/n$ be cyclic, let $\pi\colon Y\to X$ be the cover defined by $T_X$, $\gamma$ a generator of its group
$\mumu_n$ and $B=\pi^*A_X$, and assume that $-K_Y\sim2B$ and $N=\dim|B|\in\{1,2\}$. We call the points of $Y$ at which
$B$ is not Cartier the \emph{marked points}; they are base points of $|B|$, and we assume that each of them is a cyclic
quotient point of prime index. If $N=1$, let $C=\{x=y=0\}$ for a basis $x,y$ of $H^0(Y,\OO_Y(B))$; if $N=2$, let
$C=\{x=y=0\}$ for general $x,y\in H^0(Y,\OO_Y(B))$, and assume that $\operatorname{Bs}|B|$ contains a curve. Let $\mathcal C$, $\mathcal L$,
and for a component $\Gamma$ of $C$ the sheaves $N_\Gamma$, $P_\Gamma$, the numbers $\varepsilon_{\Gamma P}$ and the number $k_\Gamma$ of marked points
on $\Gamma$, be as in the proof of Proposition~\ref{prop:108}.

\begin{lemma}\label{lem:netpencil}
In this situation, suppose moreover that $C$ is not integral if $N=1$.
\begin{enumerate}
\item $C$ is a connected Gorenstein curve with $\omega_C\cong\OO_C$, all of whose components are smooth rational curves; it
contains every marked point. If $N=1$, then $C$ is $\gamma$-invariant. If $N=2$, then $C$ has a unique moving component
$M$, of multiplicity one, and its fixed part $F$ is $\gamma$-invariant together with its multiplicities.
\item Distinct components of $C$ meet only at marked points. If $\Gamma$ is a component of $C$ other than $M$, then
$B\cdot\Gamma+\sum_P\varepsilon_{\Gamma P}=1$, the preimage of $\Gamma$ near each marked point $P$ on it is a single $\mumu_{r_P}$-invariant
branch, and $\varepsilon_{\Gamma P}\in\frac1{r_P}\ZZ_{>0}$. Every component, including $M$, contains at most two marked points.
\item Let $P$ be a marked point of type $\frac13(1,1,2)$. At most two components $\Gamma$ of $C$ have $\varepsilon_{\Gamma P}=\frac13$; their
branches at $P$ are smooth, with distinct tangent lines in the $(u,v)$-plane. If two components $\Gamma_1\ne\Gamma_2$ of $C$
other than $M$ both contain two marked points, then not both have $\varepsilon_{\Gamma_iP}=\frac23$.
\item A union of components of $C$ which is a proper subcurve contains no cycle: no two of its components meet at two
marked points, and no chain of its components closes up.
\end{enumerate}
\end{lemma}

\begin{proof}
For $N=2$, \S\ref{ss:cyclic} shows that $|B|$ has no fixed components and that $\varphi_{|B|}$ maps $Y$ onto $\PP^2$. Then the
arguments of Steps~1--4 of the proof of Proposition~\ref{prop:108} apply, with the assumption that $\operatorname{Bs}|B|$ contains a
curve in place of the Hilbert function argument of Step~1; they use $\tau$ only through the invariance of $F$, which holds
for $\gamma$ as $\gamma$ preserves $|B|$, and through the parity rule, which we do not use here. For $N=1$, Step~1 of the proof of
Proposition~\ref{prop:175} gives the curve $C$, with $\omega_C\cong\OO_C$ and $H^1(\OO_Z)=0$ for proper subcurves $Z$; as $C$
is not integral, its components are smooth rational curves. Since $H^0(\mathcal C,\mathcal L)=H^0(Y,\OO_Y(B))/\langle x,y\rangle=0$,
duality gives $H^1(\mathcal Z,\mathcal L^{-1})=0$ for every closed substack $\mathcal Z$, and Step~3 of the proof of
Proposition~\ref{prop:108} gives $H^1(\mathcal Z,\mathcal L)=0$; so the budgets and both no-gluing statements of Steps~2--4 of that
proof hold for every component, with no moving component. At a marked point of prime index $r$ the preimage of a branch of
the coarse curve is either one invariant branch or a free orbit of $r$ branches, and in the second case
$\varepsilon\ge1$, which the budgets exclude. This proves (i) and (ii).

(iii) At $P$ the local sections of $\OO_Y(B)$ and $\OO_Y(-B)$ have weights $2$ and $1$. A branch with $\varepsilon_{\Gamma P}=\frac13$ has
parameter of character $1$, and the restrictions of $u$ and $v$ determine the fibres of the sheaves $N_\Gamma$ at $P$; as the
cokernel of $N_D\to\bigoplus N_\Gamma$ vanishes ($D=C_{\mathrm{red}}$), at most two such branches pass through $P$, and their tangent
lines are distinct (Step~2 of the proof of Proposition~\ref{prop:108}). A branch with $\varepsilon_{\Gamma P}=\frac23$ has parameter of
character $2$ and is smooth, with $w$ of order $1$; if $\Gamma_1$, $\Gamma_2$ contain two marked points, then $\deg P_{\Gamma_i}=-1$ by
Step~3 of that proof, and Step~4 gives the contradiction.

(iv) A proper subcurve $Z$ has $H^1(\OO_Z)=0$.
\end{proof}

\begin{proposition}\label{prop:netpencil}
The numerical types No.~$85$, $87$, $89$, $90$, $146$, $147$, $153$ and $167$ do not occur.
\end{proposition}

\begin{proof}
In each case $\qQ=2$, $T_X$ is cyclic, the cover $Y$ defined by $T_X$ has $-K_Y\sim2B$ (Lemma~\ref{lem:cover}), and the
marked points are as follows (supplementary file \texttt{examples/netpencil/data.py}).
\begin{itemize}
\item Nos.~$85$, $87$, $89$, $90$: $T_X=\ZZ/2$, $\B(X)=[2\times[2,1],2\times[6,1],[r,a]]$ with $[r,a]=[7,3]$, $[3,1]$, $[5,1]$,
$[5,2]$, and $h^0(Y,\OO_Y(B))=2$. As for No.~$175$, the two entries $[6,1]$ have different local types of $A_X$, so they are
two cyclic points, and over them $Y$ has two marked points $P_1$, $P_2$ of type $\frac13(1,1,2)$, fixed by $\tau$, near which
$\tau$ is induced by $-1\in\mumu_6$. Over $[r,a]$, where $\sigma$ is Cartier, there are two marked points $Q_\pm$ of type
$\frac1r(1,a,-a)$, interchanged by $\tau$. There are no other marked points.
\item Nos.~$146$, $147$: $T_X=\ZZ/3$, $\B(X)=[[3,1],2\times[9,a]]$ with $a=1$, $2$, and $h^0(Y,\OO_Y(B))=2$, resp.\ $3$. The
generator of $T_X$ has local types $3$ and $6$ at the two entries of index $9$, so they are two cyclic points; over them
there are two marked points $P_1$, $P_2$ of type $\frac13(1,1,2)$, fixed by $\gamma$. Over $[3,1]$ there are three marked points
$Q_1$, $Q_2$, $Q_3$ of type $\frac13(1,1,2)$, permuted cyclically by $\gamma$.
\item No.~$167$: $T_X=\ZZ/4$, $\B(X)=[2\times[2,1],[3,1],[4,1],[12,5]]$ and $h^0(Y,\OO_Y(B))=3$. Over $[12,5]$ there is one
marked point $P$ of type $\frac13(1,1,2)$, fixed by $\gamma$, and over $[3,1]$ four marked points $Q_1,\dots,Q_4$ of type
$\frac13(1,1,2)$, permuted cyclically by $\gamma$.
\item No.~$153$: $T_X=\ZZ/3$, $\B(X)=[6\times[3,1],[5,1]]$, $h^0(Y,\OO_Y(B))=2$ and $B^3=\frac95$. The generator of $T_X$ is
non-zero at every entry $[3,1]$, so $Y$ is Gorenstein over them; over $[5,1]$, where it is Cartier, there are three
marked points $Q_1$, $Q_2$, $Q_3$ of type $\frac15(1,1,4)$, permuted cyclically by $\gamma$.
\end{itemize}
All marked points have prime index. For the nets, Nos.~$147$ and $167$, we have $B^3=\frac83$, $3B$ is Cartier, and
$P_Y(t)(1-t)^3=1+6t^2+t^3+t^4+6t^5+\dots$ has period $3$, as for No.~$103$; so Step~3 of the proof of
Proposition~\ref{prop:104} shows that $\operatorname{Bs}|B|$ contains a curve. For the pencils, $C$ is not integral: by
Lemma~\ref{lem:gordual}(i), applied as in Step~1 of the proof of Proposition~\ref{prop:175}, every marked point would be a
singular point of an integral $C$, whereas an integral curve of arithmetic genus $1$ has at most one singular point, and
there are at least three marked points. So Lemma~\ref{lem:netpencil} applies.

When all marked points on $\Gamma$ have type $\frac13(1,1,2)$, $\varepsilon_{\Gamma P}\in\{\frac13,\frac23\}$, and a component $\Gamma\ne M$ through two
of them has $\varepsilon=\frac13$ at both, as $\sum_P\varepsilon_{\Gamma P}<1$. A component through no marked point meets no other component,
which is impossible as $C$ is connected and contains the marked points.

\emph{Nos.~$85$, $87$, $89$, $90$.} Suppose that a component $\Gamma$ contains $P_i$ and $Q_+$. By Lemma~\ref{lem:netpencil}(ii) it does not
contain $Q_-$, so $\tau(\Gamma)\ne\Gamma$ is a component through $P_i$ and $Q_-$. If $\varepsilon_{\Gamma P_i}=\frac13$, then, as $\tau$ acts near $P_i$
as $-1$, the branches of $\Gamma$ and $\tau(\Gamma)$ at $P_i$ have the same tangent line, contrary to
Lemma~\ref{lem:netpencil}(iii). If $\varepsilon_{\Gamma P_i}=\frac23$, the second part of Lemma~\ref{lem:netpencil}(iii) applies to $\Gamma$ and
$\tau(\Gamma)$. The same holds with $Q_-$ in place of $Q_+$. So no component contains a point $P_i$ and a point $Q_\pm$, and the
components through $P_1$ or $P_2$ do not meet those through $Q_\pm$, contrary to the connectedness of $C$.

\emph{No.~$153$.} Near $Q_j$ the local sections of $\OO_Y(-B)$ are functions of weight $2$, and no coordinate of
$\frac15(1,1,4)$ has weight $2$, so $\varepsilon_{\Gamma Q_j}\ge\frac25$ for every component $\Gamma$ through $Q_j$. As $C$ is connected and
components meet only at marked points, some component $\Gamma$ contains two of them, say $Q_j$ and $Q_k$; then
$B\cdot\Gamma\le1-\frac45=\frac15$. Its images under $\mumu_3$ are three distinct components, through $\{Q_j,Q_k\}$ and the two other
pairs, forming a cycle $Z$ with $B\cdot Z\le\frac35<B\cdot C$. So $Z$ is a proper subcurve, contrary to
Lemma~\ref{lem:netpencil}(iv).

\emph{Nos.~$146$, $147$, $167$.} Let $\Gamma\ne M$ be a component through a fixed marked point $P$ ($P=P_i$, resp.\ $P$) and a point
$Q_j$. Its images under $\mumu_n$ are $n\ge3$ distinct components through $P$ and $Q_j,\gamma(Q_j),\dots$, all with $\varepsilon=\frac13$ at $P$,
contrary to Lemma~\ref{lem:netpencil}(iii). Let $\Gamma\ne M$ be a component through $Q_j$ and $Q_k$, $j\ne k$. For $n=3$, and for
$n=4$ when $Q_k=\gamma^{\pm1}(Q_j)$, the images of $\Gamma$ form a cycle of three or four components through the $Q$'s; it is a
proper subcurve, since it contains neither $M$ nor (for No.~$146$) the points $P_1$, $P_2$, contrary to
Lemma~\ref{lem:netpencil}(iv). For No.~$146$ it follows that the components through $Q_1$ meet no component through
$P_1$, contrary to connectedness. For No.~$147$, every component through $Q_j$ other than $M$ contains no other marked
point, so $M$ contains $Q_1$, $Q_2$ and $Q_3$, contrary to Lemma~\ref{lem:netpencil}(ii). For No.~$167$, a component other
than $M$ through $P$ contains no other marked point, and one through $Q_j$ contains at most $\gamma^2(Q_j)$ besides $Q_j$;
so $M$ contains $P$, one of $Q_1,Q_3$ and one of $Q_2,Q_4$, again contrary to Lemma~\ref{lem:netpencil}(ii).
\end{proof}

\subsection{The numerical type No.~97}\label{ss:97}
Here $\qQ=2$, $\qW=1$, $T_X=\ZZ/2$, $\B(X)=[3\times[2,1],[5,2],[10,3]]$ and $A_X^3=\frac35$. We choose $A_X$ with local types
$1$, $2$, $7$ at the points of index $2$, $5$, $10$; then $i_P(\sigma)=1,0,5$ there, $K_X+2A_X\sim\sigma$, and $5A_X+\sigma$ is
Cartier at every point of $\B(X)$. Let $\pi\colon Y\to X$, $\tau$ and $B=\pi^*A_X$ be as in \S\ref{ss:104}; then $-K_Y\sim2B$,
$B^3=\frac65$ and $h^0(Y,\OO_Y(B))=2$ (supplementary file \texttt{examples/netpencil/data.py}). The entries $[10,3]$ and
$[5,2]$ are the only entries of their index, so they are cyclic points. Over $[10,3]$ there is a point $P_0$ of type
$\frac15(1,3,2)$, with coordinates $x,y,z$ of weights $1,3,2$, fixed by $\tau$, near which $\tau$ is induced by $-1\in\mumu_{10}$;
over $[5,2]$ there are two points $Q_\pm$ of type $\frac15(1,2,3)$, with coordinates $x,v,w$ of weights $1,2,3$, interchanged by
$\tau$; over the entries $[2,1]$, where $\sigma$ is not Cartier, $Y$ is Gorenstein. Near $P_0$ and $Q_\pm$ the local sections of
$\OO_Y(B)$ and $\OO_Y(-B)$ are functions of weight $3$ and $2$. So these are the marked points, and Lemma~\ref{lem:netpencil}
applies, as in the proof of Proposition~\ref{prop:netpencil}; we use its notation, with $C$ the base curve of the pencil.

\begin{proposition}\label{prop:97}
The numerical type No.~$97$ does not occur.
\end{proposition}

\begin{proof}
\emph{Step 1: local facts.} Let $V$ be the preimage of $C$ in the index-one cover near a marked point $P$, cut out by two
functions $f,g$ of weight $3$, and let $e$ be its multiplicity at the origin. A component $\Gamma$ of $C$ through $P$ has a
single invariant branch there (Lemma~\ref{lem:netpencil}(ii)), whose parameter $t$ has character $c\in\{1,2,3,4\}$, with
$\varepsilon_{\Gamma P}=\frac15\operatorname{ord}_t(h)$ for a generator $h$ of the weight-$2$ part.
\begin{enumerate}
\item[(a)] Let $s$ be the coordinate of weight $3$ ($y$ at $P_0$, $w$ at $Q_\pm$). If the tangent line of a branch is the
$s$-axis and the branch has multiplicity $\mu$, then $s$ has order $\mu$ on it and every non-linear monomial has order at least
$2\mu$; so no member of the pencil with a non-zero linear term contains the branch.
\item[(b)] The only monomials of weight $3$ and degree $2$ are $xz$ at $P_0$ and $xv$ at $Q_\pm$. So if no member of the pencil
has a linear term at $P$, then $e\ge2\cdot3=6$.
\end{enumerate}

\emph{Step 2: no component joins $Q_+$ and $Q_-$.} If $\Gamma$ does and $\tau(\Gamma)\ne\Gamma$, then $\Gamma\cup\tau(\Gamma)$ is a proper
subcurve (it misses $P_0$) with a cycle, contrary to Lemma~\ref{lem:netpencil}(iv). If $\tau(\Gamma)=\Gamma$, then
$\varepsilon_{\Gamma Q_+}=\varepsilon_{\Gamma Q_-}\in\{\frac15,\frac25\}$, so $B\cdot\Gamma\in\{\frac35,\frac15\}$ and $5B\cdot\Gamma$ is odd; but $5B$ is the pull-back of the
Cartier divisor $5A_X+\sigma$, so $5B\cdot\Gamma$ is even for a $\tau$-invariant curve, as in Step~1 of the proof of
Proposition~\ref{prop:104}.

\emph{Step 3: the components through $P_0$.} On a branch at $P_0$ of character $1$, $2$ or $4$ the image of a weight-$2$
function in the fibre of $N_\Gamma$ depends only on its coefficients of $z$ and $x^2$, as every other monomial of weight $2$ has
larger order there; so by the vanishing of the cokernel of $N_D\to\bigoplus_\Gamma N_\Gamma$ (proof of
Lemma~\ref{lem:netpencil}) at most two such branches pass through $P_0$. A branch of character $2$ is smooth with tangent
line the $z$-axis, and only the coefficient of $z$ contributes to its fibre of $N_\Gamma$; so there is at most one such
branch. A branch of character $3$ has $\varepsilon=\frac45$ and $y$ of order $1$: it is smooth, with tangent line the $y$-axis.

As $C$ is connected, components meet only at marked points and no component joins $Q_+$ and $Q_-$, some component $\Gamma$
contains $P_0$ and $Q_+$, say; by Lemma~\ref{lem:netpencil}(ii) it does not contain $Q_-$, so $\tau(\Gamma)\ne\Gamma$ contains $P_0$ and
$Q_-$, and $\Gamma$, $\tau(\Gamma)$ have branches of the same character at $P_0$. This character is not $3$, as
$\frac45+\varepsilon_{\Gamma Q_+}\ge1$, and not $2$, by the above. It is not $4$ either: then the weight-$3$ part of the local ring of
each branch is generated by the restriction of $y$, so the images of a weight-$3$ function in the fibres of $P_\Gamma$ and
$P_{\tau(\Gamma)}$ at $P_0$ are proportional to its coefficient of $y$, contrary to Step~4 of the proof of
Proposition~\ref{prop:108}, which applies since $\Gamma$ and $\tau(\Gamma)$ contain two marked points. So the character is $1$. If
$x$ had order at least $2$ on these branches, then only $z$ would contribute to their fibres of $N_\Gamma$, which is
impossible for two branches; so they are smooth, with tangent line the $x$-axis, and $y$, $z$ of order at least $3$ and
$2$. By the bound above, $\Gamma$ and $\tau(\Gamma)$ are the only components through $P_0$ whose branch has character $1$, $2$ or
$4$; in particular they are the only components through $P_0$ and a point $Q_\pm$. Let $m$ be their (common)
multiplicity in $C$ and $b=B\cdot\Gamma=1-\frac25-\varepsilon_{\Gamma Q_+}\in\{\frac15,\frac25\}$. The other components through $P_0$ contain no
other marked point and have branches of character $3$ there, so $B\cdot\Delta=\frac15$ for each of them; let $a$ be the sum of
their multiplicities.

\emph{Step 4: the contradiction.} From $B\cdot C=B^3=\frac65$ we get $\frac65\ge2mb+\frac a5\ge\frac{2m+a}5$, so $e=2m+a\le6$ at
$P_0$, all branches there being smooth. Suppose that no member of the pencil has a linear term at $P_0$. By (b), after
replacing $g$ by $g-\lambda f$, $f$ has initial form a multiple of $xz$ and $g$ has order at least $3$; if $g$ has order at least
$4$, or if $f$ has no quadratic part, then $e\ge8$. Otherwise the initial form of $g$ is a combination of $x^3$ and $y^2z$, the
only monomials of weight $3$ and degree $3$. On the branch of $\Gamma$, $x^3$ has order $3$, and every other monomial of degree at
least $3$ has order at least $4$; so the coefficient of $x^3$ vanishes, and the initial forms of $f$ and $g$ have the common
factor $z$. Cutting $V$ by a general plane through the origin and applying the local intersection inequality for plane
curves, which is strict when the tangent cones have a common component \cite[\S3.3]{Ful89}, gives $e>2\cdot3=6$, a
contradiction. Hence some member has a non-zero linear term in $y$. By (a) it contains no branch of character $3$, so
$a=0$, and $\{f=0\}$ is a smooth surface near the origin, with coordinates $x,z$; the branches of $\Gamma$ and $\tau(\Gamma)$ on it are
graphs $h_\pm=z-\varphi_\pm(x)$ of weight $2$, and $g$ restricts to a unit times $h_+^mh_-^m$. Comparing weights,
$4m\equiv3\pmod5$, so $m=2$, and then $b=\frac15$, since $4b\le\frac65$.

So $C=2\Gamma+2\tau(\Gamma)+R$, where $R$ has degree $\frac25$ and contains no component through $P_0$; by Step~2 its components
contain exactly one of $Q_\pm$. A component $\Theta$ through $Q_+$ only has $B\cdot\Theta=1-\varepsilon_{\Theta Q_+}\ge\frac15$, and
$B\cdot R=2\sum_\Theta m_\Theta B\cdot\Theta$, the sum over the components through $Q_+$; so $R=\Theta+\tau(\Theta)$ with
$\varepsilon_{\Theta Q_+}=\frac45$. Such a branch has $w$ of order $1$, so it is smooth with tangent line the $w$-axis. By (a), no member
of the pencil has a linear term at $Q_+$, and by (b) $e\ge6$ there. But $e=2\mu+1\le5$, where $\mu\le2$ is the multiplicity of the
branch of $\Gamma$ at $Q_+$, a contradiction.
\end{proof}

\subsection{The numerical type No.~11}\label{ss:11}
Here $\qQ=4$, $\qW=1$, $T_X=\ZZ/2$, $\B(X)=[3\times[2,1],[10,3],[5,2]]$ and $A_X^3=\frac15$. Let $\sigma$ be the
non-trivial torsion class. We choose $A_X$ with local types $i_P(A_X)=1$ at every point (the divisor of
Table~\ref{tab:main} is $A_X+\sigma$); then $i_P(\sigma)=1,5,0$ at the points of index
$2$, $10$, $5$, so that $-K_X\sim4A_X+\sigma$. Let $P_{10}$ and $P_5$ be the points of index $10$ and $5$; they are the only
entries of their index, so they are cyclic quotient points (\S\ref{ss:Nge5}). Put $p_k=h^0(X,kA_X)$ and
$p'_k=h^0(X,kA_X+\sigma)$; by Riemann--Roch and Kawamata--Viehweg vanishing
\[
(p_1,\dots,p_6)=(0,2,4,7,11,17),\qquad (p'_1,\dots,p'_6)=(1,2,4,7,12,16).
\]
Let $S$ be the unique member of $|A_X+\sigma|$. Every non-zero effective Weil divisor has class $kA_X$ or $kA_X+\sigma$ with
$k\ge1$, so $S$ is prime and every other prime divisor has $A_X$-degree at least $2$. The pencils $|2A_X|$ and $|2A_X+\sigma|$
have no fixed components: a fixed part would have class $A_X$ or $A_X+\sigma$, and the moving part would lie in $|A_X+\sigma|$
or $|A_X|$, which are not pencils. Near $P_{10}$ we have $2A_X\sim-8K_X$, $2A_X+\sigma\sim-3K_X$ and $A_X+\sigma\sim-4K_X$.

We use the construction of \S\ref{ss:3} with $\mathcal M=|2A_X|$, in the notation there, with $\Theta$, $\hat q$, $\bar F$,
$F\equiv dA_X$ and $\bar E\equiv e\Theta$. For a member $D$ of $|kA_X+m\sigma|$
\begin{equation}\label{eq:no11rel}
k\hat q=4s(D)+(4\beta(D)-k\alpha)e ,
\end{equation}
and $\beta(D)\equiv j\alpha\pmod\ZZ$ if $D\sim-jK_X$ near $f(E)$. We write $s_k,\beta_k$ and $s'_k,\beta'_k$ as in
\S\ref{ss:3}. Since $\mathcal M\sim-8K_X$ near $P_{10}$, \cite[Lemma~4.1.1]{Pro25} gives $c\le\frac18$ and $\beta_2\ge8\alpha$,
hence $2\hat q\ge4s_2+30\alpha e$. As in \eqref{eq:no3tors}, comparing the orders of $\Cl(\tilde X)/\langle\tilde F,E\rangle$
computed on $X$ and on $\hat X$ gives
\begin{equation}\label{eq:no11tors}
2d=e\,|\Clt(\hat X)|
\end{equation}
in the birational case. If $\sigma$ is Cartier at $f(E)$, that is, if $f(E)$ is $P_5$, a Gorenstein point or a curve,
then $f^*\sigma$ is a non-zero integral torsion class on $\tilde X$; it is preserved by $\chi$, and in the birational
case it survives on $\hat X$, because the kernel of $\Cl(\bar X)\to\Cl(\hat X)$ is generated by $\bar F$, no non-zero
multiple of which is torsion. Then $\hat X$ is a Fano--Enriques threefold, so $\hat q\le7$ by
Theorem~\ref{thm:numerical}.

\begin{lemma}\label{lem:no11links}
Either $\bar f$ is birational, $f(E)=P_{10}$, $\alpha=\frac1{10}$, $\hat q=11$, $e=2$, $F=S$, and $\hat X$ has one of the
three numerical types with $\hat q=11$ of \cite[Proposition~3.6]{Pro10}; or $\bar f$ is a $\PP^2$-fibration over $\PP^1$,
$f(E)=P_{10}$, $e=2$, $s_2=0$ and $s_3=1$.
\end{lemma}

\begin{proof}
If $f(E)$ is $P_{10}$ or $P_5$, then $f$ is the Kawamata blow up \cite{Kaw96} and $\alpha=\frac1{10}$ or $\frac15$. If
$f(E)$ is a point of index $2$, then $\alpha\in\frac12\ZZ$, and if it is a Gorenstein point or a curve, then $\alpha\in\ZZ$;
in both cases $2\hat q\ge30\alpha$ gives $\alpha\le\frac{19}{15}$. Now \eqref{eq:no11rel} for the five systems
$|A_X+\sigma|$, $|2A_X|$, $|2A_X+\sigma|$, $|3A_X|$, $|3A_X+\sigma|$, the congruences, $\beta_2\ge8\alpha$, the subadditivity of the
$\beta$'s, $\hat q\in\{1,\dots,9,11,13,17,19\}$ \cite{Suz04}, \cite[Theorem~1.4(v)]{Pro10}, and $s_2,s'_2,s_3,s'_3>0$ in the birational case
\cite[Remark~4.2.2]{Pro25} leave a finite list (\texttt{examples/no11/link\_enum.py}). The fibrations in it are
$\PP^2$-fibrations over $\PP^1$ ($\hat q=3$), with $(f(E),e)=(P_5,1)$ or $(P_{10},2)$, and in both $s_2=0$; in the second
case $s_3=1$. In the birational case $(f(E),\hat q,e)$ is $(P_{10},7,2)$, $(P_{10},11,2)$, $(P_{10},13,6)$,
$(P_{10},17,6)$, $(P_{10},19,2)$, $(P_{10},19,10)$, $(P_5,\hat q,e)$ with $\hat q\in\{7,11,13,17,19\}$, a point of index
$2$ with $(\hat q,e)=(19,2)$ ($\alpha=\frac12$) or $(19,1)$ ($\alpha=1$), or a Gorenstein point or curve with
$(\hat q,e)=(19,1)$.

\emph{The fibration with $f(E)=P_5$.} Here $\sigma$ is Cartier at $P_5$, so $\bar\sigma=\chi_*f^*\sigma$ is a non-zero
torsion class on $\bar X$, defining a double cover $\bar Y\to\bar X$ which is \'etale in codimension one. Over the general
fibre $\PP^2$ it is \'etale, hence trivial, so the Stein factorisation of $\bar Y\to\PP^1$ is a connected double cover
$C\to\PP^1$, branched at two points $o_1\ne o_2$ at least. Since $\rho(\bar X/\PP^1)=1$, the fibre over $o_i$ is
$m_iF_i$ with $F_i$ prime, and $m_i$ is even because $\bar Y\to\bar X$ is \'etale over the generic point of $F_i$. As
$s_2=0$, the proper transform of $|2A_X|$ is a complete pencil of fibres, so $2A_X\sim m_if_*\chi^{-1}_*F_i$, where
$f_*\chi^{-1}_*F_i$ is a prime divisor ($F_i\ne\bar E$, since $e>0$). Its $A_X$-degree $2/m_i$ is a positive integer, so
$m_i=2$ and $f_*\chi^{-1}_*F_1=f_*\chi^{-1}_*F_2=S$, which is impossible since $o_1\ne o_2$.

\emph{Birational cases.} If $\hat q\ge9$, then $\Clt(\hat X)=0$ and $\hat X$ has one of the numerical types of
\cite[Proposition~3.6]{Pro10}, whose values of $h^0(k\Theta)$, $k\le7$, are listed there. By \eqref{eq:no11tors},
$d=e\,|\Clt(\hat X)|/2$. If $d=1$, then $F\in|A_X+\sigma|$ (as $|A_X|=\varnothing$), so $F=S$ and $s'_1=0$; conversely
$s'_1=0$ means that $S$ is contracted, so $F=S$ and $d=1$. A system $|D|$ with $s(D)>0$ is mapped isomorphically onto a
subspace of $H^0(\hat X,s(D)\Theta+t)$ for some torsion class $t$, and $\hat E=\bar f(\bar E)$ is a prime member of
$|e\Theta+t|$ for some $t$. We apply these facts to the list (\texttt{examples/no11/link\_targets.py}):
\begin{itemize}
\item if $f(E)$ is $P_5$, a Gorenstein point or a curve, then $\hat q\le7$, so $f(E)=P_5$, $\hat q=7$, $e=1$, and
$\hat X$ is of type No.~$1$ or No.~$2$ of Table~\ref{tab:main}, with $|\Clt(\hat X)|=2$; by \eqref{eq:no11tors} $d=1$, so
$F=S$ and $s'_1=0$. In the rows with $s'_1=0$ we have $s_2\le2$, while $h^0(\hat X,s\Theta+t)\le1$ for $s\le2$ and both
torsion classes $t$;
\item if $f(E)$ is a point of index $2$, then $\hat q=19$ and $\Clt(\hat X)=0$; for $e=1$, \eqref{eq:no11tors} has no
solution, and for $e=2$ we get $d=1$, while $s'_1=1$;
\item if $f(E)=P_{10}$ and $\hat q=7$, then either $\Clt(\hat X)=0$, $d=1$ and $s'_1=1$, a contradiction, or $\hat X$ is
of type No.~$1$ or No.~$2$, where $h^0(2\Theta+t)\le1<p_2$ for both $t$, while $s_2=2$;
\item if $f(E)=P_{10}$ and $\hat q\in\{13,17,19\}$, then one of the conditions above fails for every row: for
$(\hat q,e)=(13,6)$, the target $\PP(1,3,4,5)$ \cite[Theorem~1.4(iii)]{Pro10} has $h^0(2\Theta)=1$ while $s_2=2$, and the
other type has $h^0(\Theta)=0$ while $s'_1=1$; for $(17,6)$, $\hat X\cong\PP(2,3,5,7)$ has $h^0(4\Theta)=1$ while $s_2=4$;
for $(19,10)$, $\hat X\cong\PP(3,4,5,7)$ has $h^0(\Theta)=0$ while $s'_1=1$; for $(19,2)$, $|2\Theta|=\varnothing$ on
$\PP(3,4,5,7)$ cannot contain $\hat E$.
\end{itemize}
In the remaining case $f(E)=P_{10}$, $\hat q=11$ and $e=2$, the target has $\Clt(\hat X)=0$, so $d=1$ and $F=S$.
\end{proof}

\begin{proposition}\label{prop:11}
The numerical type No.~$11$ does not occur.
\end{proposition}

\begin{proof}
We use Lemma~\ref{lem:no11links}. For a complete system $|kA_X+m\sigma|$ on $X$, its proper transform on $\bar X$ (after
removing the fixed multiple of $\bar E$) has the same space of sections; in the birational case its image on $\hat X$ is a
subspace of $H^0(\hat X,s\Theta)$, equal to it when the dimensions agree, and every member of the image is the image of a
member on $\bar X$, which may contain $\bar E$ and $\bar S$ as components.

\emph{The $\PP^2$-fibration.} The proper transform of $|2A_X|$ is the pencil of fibres, so a general fibre $\Phi\cong\PP^2$
is a member of it, and the proper transform $\bar D$ of a general member of $|3A_X|$ restricts to $\Phi$ as a line, since
$s_3=1$. An effective divisor in $|\bar D-\Phi|$ would map to an effective divisor in $|3A_X-2A_X|=|A_X|=\varnothing$ on $X$
(it is not a multiple of $\bar E$, since $A_X\not\sim0$). So restriction to $\Phi$ is injective on the transform of $|3A_X|$,
and $4=p_3\le h^0(\PP^2,\OO(1))=3$, a contradiction.

\emph{The birational case.} Here $\bar F=\bar S$ and $\hat E\equiv2\Theta$ is prime. A prime divisor of $\hat X$ other
than $\hat E$ is the image of a prime divisor of $X$ other than $S$, of $A_X$-degree at least $2$. So the image of a member
of $|kA_X+m\sigma|$ has at most $\lfloor k/2\rfloor$ prime components different from $\hat E$, counted with multiplicity.
The three numerical types are as follows.

\emph{$\hat X$ of type $(2,5,7)$.} By \cite[Proposition~3.6]{Pro10}, $h^0(\Theta)=h^0(2\Theta)=1$, so the only member of
$|2\Theta|$ is twice the member of $|\Theta|$, and $\hat E$ is not prime.

\emph{$\hat X$ of type $(2,3,5)$.} By \cite[Theorem~1.4(iv)]{Pro10}, $\hat X\cong\PP(1,2,3,5)$, with coordinates $x,y,z,w$ of
weights $1,2,3,5$; since $\hat E$ is prime of degree $2$, we may take $\hat E=\{y=0\}$. By \eqref{eq:no11rel} with
$\beta_2\in\frac45+\ZZ$, $s_2\in\{2,4\}$. If $s_2=2$, the image of $|2A_X|$ is $|\OO(2)|=\langle x^2,y\rangle$, and its
member $\{x^2=0\}$ has two components different from $\hat E$. If $s_2=4$, the image of $|2A_X|$ is a $2$-dimensional
subspace of $\langle x^4,x^2y,y^2,xz\rangle$, which is $1$-dimensional modulo $x$; so it contains $xg$ with $g$ of degree
$3$. If $y\nmid g$, then $\{x=0\}$ and a component of $\{g=0\}$ are two components different from $\hat E$; if $y\mid g$,
then $xg$ is a multiple of $x^2y$. Both cases contradict $\lfloor2/2\rfloor=1$.

\emph{$\hat X$ of type $(2,2,3,4,7)$.} By \cite[Proposition~3.6]{Pro10}, $h^0(k\Theta)=0,1,1,2,2,3,4$ for $k=1,\dots,7$.
By \eqref{eq:no11rel} with $\beta_3\in\frac7{10}+\ZZ$, $s_3\in\{1,3,5,7\}$, and since the image of $|3A_X|$ has dimension
$p_3=4$, we get $s_3=7$ and the image is all of $|7\Theta|$. Let $u$, $z$ and $v$ be non-zero sections of $\Theta$-degree
$2$, $3$ and $4$, with $v\notin\langle u^2\rangle$; then $\{u=0\}=\hat E$, and neither $\{z=0\}$ nor $\{v=0\}$ contains
$\hat E$, since $h^0(\Theta)=0$ and $h^0(2\Theta)=1$. The member $\{zv=0\}$ of $|7\Theta|$ has at least two components
different from $\hat E$, contradicting $\lfloor3/2\rfloor=1$.
\end{proof}

\subsection{The numerical type No.~13}\label{ss:13}
Here $\qQ=4$, $\qW=1$, $T_X=\ZZ/2$, $\B(X)=[5\times[2,1],[5,2],[6,1]]$ and $A_X^3=\frac2{15}$. Let $\sigma$ be the
non-trivial torsion class. We choose $A_X$ with local types $0,1,1,1,1$ at the five points of index $2$, $2$ at the point
$P_6$ of index $6$ and $1$ at the point $P_5$ of index $5$ (the divisor of Table~\ref{tab:main} is $A_X+\sigma$); then
$i_P(\sigma)=1,\dots,1,3,0$, so $-K_X\sim4A_X+\sigma$ and $\sigma$ is Cartier at $P_5$. As before $P_5$ and $P_6$ are
cyclic quotient points. With $p_k=h^0(X,kA_X)$ and $p'_k=h^0(X,kA_X+\sigma)$ we have (supplementary file
\texttt{examples/no13/no13\_targets.py})
\[
(p_1,\dots,p_6)=(0,2,3,5,8,12),\qquad (p'_1,\dots,p'_6)=(1,1,3,5,8,11).
\]
As for No.~$11$, the unique member $S$ of $|A_X+\sigma|$ is prime, every other prime divisor has $A_X$-degree at least $2$,
and $|2A_X|$ is a pencil without fixed components. Near $P_5$ we have $2A_X\sim-3K_X$. We use the construction of
\S\ref{ss:3} with $\mathcal M=|2A_X|$, in the notation of \S\ref{ss:11}; now $c\le\frac13$, $\beta_2\ge3\alpha$, and
\eqref{eq:no11rel} and \eqref{eq:no11tors} hold verbatim. If $\sigma$ is Cartier at $f(E)$, it survives on $\hat X$ as
in \S\ref{ss:11}, and we denote its image by $\hat\sigma$.

\begin{lemma}\label{lem:no13links}
One of the following holds.
\begin{enumerate}
\item $\bar f$ is a $\PP^2$-fibration over $\PP^1$, $f(E)=P_5$, $e=1$ and $s'_1=0$.
\item $\bar f$ is birational, $F=S$, $\Clt(\hat X)=0$, and either $f(E)=P_6$, $e=2$ and $\hat X$ has the numerical type
$(2,4,5)$ with $\hat q=9$ of \cite[Proposition~3.6]{Pro10} or $\hat X\cong\PP(2,3,5,7)$; or $f(E)$ is the point of
index $2$ at which $A_X$ is Cartier, $\alpha=\frac12$, $e=2$ and $\hat X\cong\PP(1,2,3,5)$. In all these cases $s_2=2$, $6$,
$2$ respectively.
\item $\bar f$ is birational, $f(E)=P_5$, $\alpha=\frac15$, $e=d=1$, $F=S$, $\Clt(\hat X)=\langle\hat\sigma\rangle\cong\ZZ/2$,
and $\hat X$ is a Fano--Enriques threefold of one of the numerical types No.~$1$, $2$, $18$, $19$, $22$, $25$, with
$\hat q=7,7,3,3,3,3$.
\end{enumerate}
\end{lemma}

\begin{proof}
The centre $f(E)$ is $P_6$ or $P_5$ ($\alpha=\frac16$, $\frac15$, Kawamata blow up), a point of index $2$ ($\alpha\in\frac12\ZZ$;
the two local types of $A_X$ are treated separately), or a Gorenstein point or curve ($\alpha\in\ZZ$). As in
Lemma~\ref{lem:no11links}, \eqref{eq:no11rel} for $|A_X+\sigma|$, $|kA_X|$ and $|kA_X+\sigma|$ with $k\le4$, the
congruences, $\beta_2\ge3\alpha$ (so $2\hat q\ge10\alpha e$ and $\alpha\le\frac{19}5$), the subadditivity of the $\beta$'s and
$s_2>0$ in the birational case leave a finite list (\texttt{examples/no13/link\_q4.py}). Its fibrations are
$\PP^2$-fibrations with $f(E)=P_5$ and $e=1$, and \eqref{eq:no11rel} for $S$ gives $s'_1=0$. In the birational case
$s_2>0$ gives $\hat q\ge3$.

Suppose that $\Clt(\hat X)=0$. Then $\sigma$ is not Cartier at $f(E)$, so $f(E)$ is $P_6$ or a point of index $2$, and
$e=2d$ by \eqref{eq:no11tors}. If $f(E)=P_6$, then $\beta_2\in\frac13+\ZZ$ and $\beta_2\ge\frac12$, so $\beta_2\ge\frac43$
and $2\hat q\ge4+5e$; if $\hat q\le8$, this forces $e=2$, $d=1$, $F=S$, and \eqref{eq:no11rel} for $S$, with
$\beta(S)\in\frac16+\ZZ$, gives $\hat q\equiv1\pmod8$, which is impossible for $7\le\hat q\le8$. If $f(E)$ has index $2$,
then $\beta_2\in\ZZ$ and $\beta_2\ge3\alpha\ge\frac32$, so $2\hat q\ge4+7e\ge18$. Hence $\hat q\ge9$, and $\hat X$ has one
of the numerical types of \cite[Proposition~3.6]{Pro10}; testing them against the list gives the cases in~(ii), with
$F=S$ since $d=1$, and \cite[Theorem~1.4]{Pro10} identifies the targets with $\hat q=11,17$.

Suppose that $\Clt(\hat X)\ne0$. Then $\hat X$ is a Fano--Enriques threefold with $\qQ(\hat X)=\hat q$, so $\hat q\le7$ and
$\hat X$ has a numerical type of Table~\ref{tab:main} whose group $T_{\hat X}$ contains $\Clt(\hat X)$. The list gives
$f(E)=P_5$, $\alpha=\frac15$, $e=1$, so $d=1$, $F=S$ and $|\Clt(\hat X)|=2$. If $T_{\hat X}$ has order $2$, the type is
one of those in~(iii). Otherwise $\Clt(\hat X)$ is a subgroup of order $2$ of a larger group, and the numerical type is
No.~$32$ or No.~$36$; but there every class of $\Theta$-degree $1$ has $h^0\le1$, whereas the image of the pencil $|2A_X|$,
with $s_2=1$, lies in such a class.
\end{proof}

\begin{proposition}\label{prop:13}
The numerical type No.~$13$ does not occur.
\end{proposition}

\begin{proof}
We go through the cases of Lemma~\ref{lem:no13links}. As in the proof of Proposition~\ref{prop:11}, in the birational case
with $F=S$ the image of a member of $|kA_X+m\sigma|$ has at most $\lfloor k/2\rfloor$ prime components other than
$\hat E$, counted with multiplicity.

\emph{Case (i).} By $s'_1=0$, $\bar S$ is the support of a fibre $\mu\bar S$, so $\mu(A_X+\sigma)$ is the class of the
fibres on $X$; $\mu\ge2$ since $|A_X+\sigma|$ is not a pencil, and $\mu\le3$, because a del Pezzo fibration with fibres
of degree $9$ has no fibres of multiplicity greater than $3$ \cite[Theorem~1.1, Table~1]{MP09}. Since $\sigma$ is Cartier at
$P_5$, the argument for the fibration over $P_5$ in the proof of Lemma~\ref{lem:no11links} gives two distinct fibres
$m_iF_i$ with $m_i$ even; by \cite{MP09} again $m_i=2$. So $\mu(A_X+\sigma)\equiv2G_i$ for two distinct prime divisors
$G_1$, $G_2$ on $X$. If $\mu=2$, then $G_1=G_2=S$; if $\mu=3$, then $3=2\deg G_i$. Both are impossible.

\emph{Case (ii).} If $\hat q=9$, then $h^0(\Theta)=1$ and $h^0(2\Theta)=2$ \cite[Proposition~3.6]{Pro10}, so the image of
$|2A_X|$ is all of $|2\Theta|$ and contains twice the member of $|\Theta|$, which is not $\hat E$, as $\hat E$ is prime of
degree $2$. If $\hat X\cong\PP(2,3,5,7)$, the image of $|2A_X|$ is $|\OO(6)|=\langle x^3,y^2\rangle$ and $\hat E=\{x=0\}$; the
member $\{y^2=0\}$ has two components other than $\hat E$. If $\hat X\cong\PP(1,2,3,5)$, we may take $\hat E=\{y=0\}$, the
image of $|2A_X|$ is $\langle x^2,y\rangle$, and $\{x^2=0\}$ has two components other than $\hat E$. All contradict
$\lfloor2/2\rfloor=1$.

\emph{Case (iii): two identities.} Put $H=\hat E$. From $f^*S=\tilde S+\beta(S)E$ and $-K_X\sim4(A_X+\sigma)-3\sigma$,
\[
-K_{\tilde X}\sim4\tilde S+(4\beta(S)-\tfrac15)E-3f^*\sigma ,
\]
and $4\beta(S)-\frac15=\hat q$ by \eqref{eq:no11rel} for $S$, as $s'_1=0$ and $e=1$. Pushing forward to $\hat X$, where
$\bar S$ is contracted,
\begin{equation}\label{eq:no13K}
-K_{\hat X}\sim\hat qH+\hat\sigma ,
\end{equation}
and in the same way the image of a member $D$ of $|kA_X+m\sigma|$ lies in the class $s(D)H+(k+m)\hat\sigma$. In all six
numerical types $\qW=\qQ=\hat q$ is odd, so there is a unique $A_{\hat X}$ with $-K_{\hat X}\sim\hat qA_{\hat X}$, and
\eqref{eq:no13K} gives $H=A_{\hat X}+\hat\sigma$. The values $h^0(\hat X,kH+m\hat\sigma)$ used below follow from
Theorem~\ref{thm:RR} (\texttt{examples/no13/no13\_targets.py}). The types No.~$19$ and $25$ do not occur by
Proposition~\ref{prop:genusbound}. For No.~$1$, $h^0(H)=0$, although $H=\hat E$ is effective. For No.~$18$, $h^0(H)=1$, but
the image of $H^0(2A_X)$, of dimension $2$, lies in $H^0(H)$, since $s_2=1$.

\emph{No.~$2$.} Here $h^0(H+\hat\sigma)=0$, $h^0(2H+\hat\sigma)=1$, $h^0(3H+\hat\sigma)=2$, $h^0(4H)=3$ and $h^0(5H)=4$. The image
$V$ of $H^0(3A_X+\sigma)$ has dimension $3$ and lies in $H^0(sH)$ with $s=s'_3\in\{4,5\}$
(\texttt{examples/no13/link\_q4.py}). Let $0\ne v\in H^0(2H+\hat\sigma)$; its divisor does not contain $H$, since
$h^0(H+\hat\sigma)=0$. If $s=4$, then $V=H^0(4H)\ni v^2$. If $s=5$, the $2$-dimensional space $v\,H^0(3H+\hat\sigma)$ meets $V$
in the $4$-dimensional $H^0(5H)$, and every non-zero $h\in H^0(3H+\hat\sigma)$ has a component other than $H$, as its divisor
is not in the class $3H$. In both cases $V$ contains a member with two components other than $\hat E$, contradicting
$\lfloor3/2\rfloor=1$.

\emph{No.~$22$.} Here $\hat q=3$, $h^0(H)=2$, $h^0(H+\hat\sigma)=1$, $h^0(2H)=h^0(2H+\hat\sigma)=4$ and $h^0(3H)=8$. Let $u$
define $H$, let $u,v$ span $H^0(H)$ and $w$ span $H^0(H+\hat\sigma)$. Let $\nu$ be the valuation of $\bar F$,
$b=a(\bar F,\hat X)$ and $\delta=\nu(H)$. By \eqref{eq:no11rel}, $s_2=s'_2=1$, so $H^0(H)$ and $H^0(H+\hat\sigma)$ are the
images of $|2A_X|$ and $|2A_X+\sigma|$. As in \eqref{eq:no3rev}, a prime divisor $G\ne\hat E$ in $|sH+t|$ with
$\gamma=\nu(G)$ has proper transform on $X$ of class $nA_X+\varepsilon\sigma$ with $3n=4s-(3\gamma-sb)$. The same computation
applies to $G=\hat E$, whose transform on $X$ is exceptional, with $n=0$; as $s=1$ for $\hat E$, this gives $b=3\delta-4$; for the general members of $|H|$ and of $|H+\hat\sigma|$ we have $n=2$, so
$\nu(v)=\nu(w)=\delta-2$, and $\delta\ge2$. For the image of a general member of $|4A_X|$ we have $n=4$ and $s=3$, so it has
$\nu=b$; hence the image of $H^0(4A_X)$, of dimension $p_4=5$, lies in $H^0(3H)$ and consists of sections with $\nu\ge b$. As
$b\ge2$, the centre $Z$ of $\nu$ on $\hat X$ is a point.

Let $\pi\colon Y\to\hat X$ be the double cover defined by $\hat\sigma$ and $B=\pi^*H$. Then $Y$ is a terminal Fano threefold of
singularity index $2$, $2B$ is Cartier, $-K_Y\sim3B$ and $B^3=1$. By Sano's classification \cite[Theorem]{San96} (type [2];
the invariant $(2B)^3=8$ singles it out) $Y$ is a quartic in $\PP(1,1,1,2,2)$ with $\OO_Y(2)\cong\OO_Y(2B)$, and since
$-K_Y\sim\OO_Y(3)$, $B\sim\OO_Y(1)$; the dimensions $h^0(Y,kB)=3,8,16,29,47$ agree. Comparing with the eigenspaces of the
involution $\iota$ of $\pi$, the coordinates are $u,v,w,a,z$ of degrees $1,1,1,2,2$, on which $\iota$ acts by $(+,+,-,-,-)$,
and the quartic $\Phi$ is invariant. Since $\sigma$ is Cartier at $P_5$, the double cover of $X$ defined by $\sigma$ is
trivial near $P_5$, so the preimage of $E$ splits into two divisors, which are interchanged by the involution and survive
through the link; hence $\pi^*H=\{u=0\}$ has two components interchanged by $\iota$. So $\Phi(0,v,w,a,z)=G\cdot\iota^*G$ with
$G=Q+L$, $Q=c_1v^2+c_2w^2\ne0$ even and $L$ odd; $L$ involves $a$ or $z$, since otherwise $\{u=0\}$ would be a union of
four weighted hyperplane sections counted with multiplicity. After an equivariant change of coordinates $L=a$, and
\[
\Phi=a^2-Q(v,w)^2+u\,\bigl(\lambda wz+G_3(u,v,w,a)\bigr).
\]
At the point $P_z$ the index-one cover of $\hat X$ is $\{\Phi(u,v,w,a,1)=0\}\subset\CC^4$ with $\mumu_4$-weights
$(1,1,3,0)$. If $\lambda=0$, its quadratic part is $a^2$, which is not the quadratic part of a terminal point of index $4$
\cite[Theorem~(6.1)]{Rei87}; so $\lambda\ne0$, and $P_z\in\hat X$ is a point of type $cA/4$ of axial weight $2$.

If $\delta>2$, then $u,v,w$ vanish on $Z$, so $a=0$ on $Z$ and $Z=P_z$. By \cite[Theorem~1.1(i)]{Kaw05}, a divisorial
contraction to a point of type $cA/4$ is a weighted blow up of $x_1x_2+g(x_3^4,x_4)$ with weights $(r_1,r_2,a',4)/4$,
discrepancy $a'/4$, $4a'\mid r_1+r_2$ and $g$ of weighted order $(r_1+r_2)/4$; as $g$ contains $x_4^2$, $a'\le(r_1+r_2)/4\le2$,
so $b\le\frac12$, contrary to $b=3\delta-4>2$.

If $\delta=2$, then $b=2$ and $v$, $w$ do not vanish at $Z$, so $\pi$ is \'etale over $Z$ and, in the chart $v=1$ with
$U=u/v$, $r=w/v\ne0$, $A=a/v^2$ and $Z'=z/v^2$, the threefold is $A^2-Q(1,r)^2+U(\lambda rZ'+G_3)=0$. At a singular point of
this chart $Q(1,r)=0$, a simple zero, and $\partial(\lambda rZ'+G_3)/\partial Z'=\lambda r\ne0$; so the point is an
ordinary double point, over which a divisorial contraction has discrepancy $1$ \cite[Theorem~3.10]{Cor00}. Hence $Z$ is a
smooth point, and as $b=2$ the contraction is the ordinary blow up \cite{Kaw01}, so $\nu$ is the order of vanishing at $Z$.
The invariant sections $v^3,uv^2,vw^2,wa,wz\in H^0(\hat X,3H)$ become $1,U,r^2,rA,rZ'$ in the chart, and their differentials
at $Z$ span the cotangent space, since $\det\partial(U,r^2,rA,rZ')/\partial(U,r,A,Z')=2r^3\ne0$. So the sections of $3H$
vanishing to order at least $2$ at $Z$ form a subspace of dimension $8-4=4$, which cannot contain the $5$-dimensional
image of $H^0(4A_X)$.
\end{proof}

\subsection{The numerical type No.~3}\label{ss:3}
Here $\qQ=\qW=5$, $T_X=\ZZ/2$, $\B(X)=[2\times[4,1],[12,5]]$ and $A_X^3=\frac1{12}$. This is one of the two
types with $\qQ\ge5$ whose existence was left open in \cite[\S5]{Pro25}; the other one, No.~$5$, is
realised in \S\ref{ss:tom}. We show that it does not occur. We use the birational method of
\cite[\S4]{Pro25}, keeping its notation.

Let $X$ be of type No.~$3$. Let $\sigma$ denote the non-trivial torsion class and let $A_X$ be the Weil
divisor with $-K_X\sim5A_X$; it is unique, since $\Clt(X)\cong\ZZ/2$ has no $5$-torsion, and it is the divisor of
Table~\ref{tab:main}. Put
\[
p_k=h^0(X,kA_X),\qquad p'_k=h^0(X,kA_X+\sigma),
\]
so that $(p_1,\dots,p_7)=(1,1,2,4,6,8,11)$ and $(p'_1,\dots,p'_7)=(0,1,2,4,6,9,12)$. Let $D_1$ be the
unique member of $|A_X|$. Let $P_1,P_2$ be the points of index $4$ and $P$ the point of index $12$. Then
$\sigma\sim 2K_X$ near $P_1$, $\sigma$ is Cartier at $P_2$, and near $P$ we have $A_X\sim-5K_X$ and $\sigma\sim6K_X$.
Since $\sigma$ has different local types at $P_1$ and $P_2$, these are two distinct points, each contributing
one point to the basket; hence $P_1$, $P_2$ and $P$ are cyclic quotient singularities, and $P=\frac1{12}(1,5,7)$.

The linear system $|4A_X|$ has no fixed components: a non-zero fixed part has class $mA_X+\varepsilon \sigma$ with
$m\ge1$, and the moving part would then lie in a system $|kA_X+\varepsilon'\sigma|$ with $k\le3$, which has at most
$\max(p_3,p'_3)=2$ sections, while $p_4=4$.

Apply the construction of \cite[\S4.1]{Pro25} to $\mathcal{M}=|4A_X|$: $f\colon\tilde X\to X$ is an extremal
blow up with exceptional divisor $E$ and discrepancy $\alpha$ which is crepant for $(X,c\mathcal{M})$,
$c=\operatorname{ct}(X,\mathcal{M})$; $\chi\colon\tilde X\dashrightarrow\bar X$ is an isomorphism in codimension
one; $\bar f\colon\bar X\to\hat X$ is a divisorial contraction or a fibration. In the birational case
$\hat X$ is a $\QQ$-Fano threefold, $\Theta$ denotes an ample generator of $\Cl(\hat X)$ modulo torsion,
$\hat q=\qQ(\hat X)$, $\bar F$ is the $\bar f$-exceptional divisor, $F$ its image on $X$, $F\equiv dA_X$,
$\bar E\equiv e\Theta$, and $b$ is the discrepancy of $\bar F$ over $\hat X$. For a member $D$ of
$|kA_X+\varepsilon \sigma|$ let $\beta(D)$ be its multiplicity along $E$ and $s(D)$ the degree of the image
of $D$ on $\hat X$ with respect to $\Theta$; then \cite[(4.4)]{Pro25}
\begin{equation}\label{eq:no3rel}
k\hat q=5s(D)+(5\beta(D)-k\alpha)e .
\end{equation}
Write $s_k,\beta_k$ (resp. $s'_k,\beta'_k$) for the values on a general member of $|kA_X|$ (resp. of
$|kA_X+\sigma|$). Since $\mathcal{M}\sim-8K_X$ near $P$, \cite[Lemma~4.1.1]{Pro25} gives $c\le\frac18$, hence
$\beta_4\ge8\alpha$.

\begin{lemma}\label{lem:no3links}
The map $\bar f$ is birational, and $(f(E),\hat q,e)$ is $(P,4,2)$ or $(P,13,4)$.
\end{lemma}

\begin{proof}
If $f(E)$ is a non-Gorenstein point, then $f$ is the Kawamata blow up \cite{Kaw96}; so $\alpha=\frac14$ or
$\frac1{12}$, and $\beta(D)$ is congruent modulo $\ZZ$ to $\alpha$ times the local type of $D$. If $f(E)$ is a
Gorenstein point or a curve, then $\alpha$ and $\beta(D)$ are integers, and $4\hat q\ge5s_4+36\alpha e$ gives
$\alpha\le2$. Now \eqref{eq:no3rel} for $k=1,\dots,4$ and both systems, $s(D)\ge0$, $\beta_4\ge8\alpha$,
$e\le\hat q/(9\alpha)$, $\beta(D)\le k\hat q/(5e)+k\alpha/5$, the subadditivity of the $\beta$'s (a product of
members is a member), $\hat q\in\{1,\dots,9,11,13,17,19\}$ \cite{Suz04}, \cite[Theorem~1.4(v)]{Pro10}, and $s_3,s_4,s'_3,s'_4>0$ in the
birational case \cite[Remark~4.2.2]{Pro25} leave a finite list, computed in
\texttt{examples/no3/link\_enum.py}. In the fibration case ($\hat q\le3$, where a mobile system may restrict
trivially to the general fibre, so $s(D)=0$ is allowed) there are no solutions. In the birational case
$(f(E),\hat q,e)$ is one of $(P,4,2)$, $(P,8,4)$, $(P,9,2)$, $(P,13,4)$, $(P,17,6)$, $(P,19,2)$,
$(P_1,17,2)$, $(P_2,6,1)$, $(P_2,11,1)$, $(P_2,17,2)$, or $f(E)$ is Gorenstein and $(\hat q,e)=(19,1)$.

Comparing the orders of $\Cl(\tilde X)/\langle\tilde F,E\rangle$ computed on $X$ and on $\hat X$ gives
\begin{equation}\label{eq:no3tors}
2d=e\,|\Clt(\hat X)|.
\end{equation}
If $\hat q=19$ or $17$, then $\hat X\cong\PP(3,4,5,7)$ or $\PP(2,3,5,7)$ \cite[Theorem~1.4(i), (ii)]{Pro10}; in the
cases with $\hat q=17$ one has $s_1=1$, but $h^0(\OO(1))=0$ on $\PP(2,3,5,7)$, and for $\hat q=19$ the
divisor $\bar E$ would be a member of $|\OO(e)|=\varnothing$ ($e\le2$) on $\PP(3,4,5,7)$. If $\hat q=6$ or
$11$, then $\Clt(\hat X)=0$ by \cite[Proposition~2.3.2]{Pro25}, contradicting \eqref{eq:no3tors} with $e=1$.
If $(\hat q,e)=(8,4)$, then $\Clt(\hat X)=0$, $d=2$, and \eqref{eq:no3rel} for $F$ ($s=0$) gives $F\sim2A_X$;
but $|2A_X|=\{2D_1\}$. If $(\hat q,e)=(9,2)$, then $d=1$ and \eqref{eq:no3rel} gives $\beta(F)=\frac{11}{12}$,
so $F\sim A_X+\sigma$; but $p'_1=0$.
\end{proof}

For a prime divisor $G\ne\bar f(\bar E)$ in $|s\Theta|$ on $\hat X$ with $\gamma=\operatorname{mult}_{\bar F}\bar f^*G$,
pushing $sK_{\hat X}+\hat qG\sim0$ (when $\Clt(\hat X)=0$) down to $X$ shows that its proper transform has
class $nA_X+\varepsilon \sigma$ with
\begin{equation}\label{eq:no3rev}
\hat q\,n=5s-(\hat q\gamma-sb)\,d .
\end{equation}
In particular $n\ge0$, with $n=0$ only if $G$ is the image of $E$. We write $v$ for the valuation of
$\bar F$ and $v_E$ for that of $E$.

\begin{lemma}\label{lem:no3b}
The case $(\hat q,e)=(13,4)$ does not occur.
\end{lemma}

\begin{proof}
Here $d=2$, and \eqref{eq:no3rel} with $\beta(D)\in\frac{5k}{12}+\frac\varepsilon2+\ZZ$ gives
$s(D)=k-2\varepsilon-4j$ with $\beta(D)=\frac{5k}{12}+\frac\varepsilon2+j\ge0$. Members $D$ of $|7A_X+\sigma|$
with $v_E\ge\frac{17}{12}$ thus give members of $|13\Theta|$ after adding multiples of $\bar E$, and near
$P$ the condition $v_E\ge\frac{17}{12}$ on sections of local character $5$ is given by the coefficients of
the two monomials $\xi_2$, $\xi_1^5$ in orbifold coordinates of weights $1,5,7$. So
$h^0(\hat X,-K_{\hat X})\ge p'_7-2=10$, $\g(\hat X)\ge8$, and $\hat X\cong\PP(1,3,4,5)$ by
\cite[Theorem~1.4(iii)]{Pro10}; let $x,y,z,w$ be coordinates of weights $1,3,4,5$. The divisor $\bar f(\bar E)$
is an irreducible member of $|\OO(4)|=\langle z,xy,x^4\rangle$, so after an automorphism it is $\{z=0\}$.
By the same formula $s(D_1)=1$, so $\bar f(\bar D_1)=\{x=0\}$; \eqref{eq:no3rev} for $\{x=0\}$ gives
$b=4+13g$ with $g=v(x)$, and then $n=s(1+2g)-2v(G)$ for every $G$, $v(z)=2+4g$, and $s\equiv n-2\varepsilon\bmod 4$.
Since $b>1$, $\bar f(\bar F)$ is a point $\hat P$, and it is not one of the cyclic quotient points
(Kawamata blow ups have discrepancy $<1$), so $\hat P$ is a nonsingular point of $\{z=0\}$.

If $g=0$, we may assume $\hat P=(1:0:0:0)$. Then $v(y)\ge1$, so the proper transform of $\{y=0\}$ has
$n=3-2v(y)\le1$ odd; $n=1$ forces $\varepsilon=1$, a member of $|A_X+\sigma|$, and $p'_1=0$.

If $g\ge1$, then $\hat P\in\{x=z=0\}$, so $\hat P=(0:1:0:1)$ up to automorphism and $y,w$ are units at
$\hat P$. Take $D\in|3A_X|$ not containing $D_1$ (possible as $p_3>p_2$) and let $G$ be the image of $D$.
Then $s(G)\in\{3,7\}$ and $v(G)=\frac{s(1+2g)-3}2$. If $s=3$, then $G=y+\lambda x^3$ is a unit at $\hat P$,
while $v(G)=3g>0$. If $s=7$, then $G=c_0yz+xK$ with $c_0\ne0$ and $K\in H^0(\OO(6))$, since $yz$ is the only
monomial of degree $7$ not divisible by $x$; as $v(G)=2+7g>2+4g=v(yz)$, we need $v(K)=2+3g$. But
$K=c_1y^2+x\bigl(c_2w+c_3x^2y+c_4xz+c_5x^5\bigr)$, and one checks successively that
$v(K)\in\{0,g,3g,6g,2+6g\}$ according to which coefficients vanish.
\end{proof}

We now treat $(\hat q,e)=(4,2)$. Then $d=1$, $F=D_1$ and $\Clt(\hat X)=0$ by \eqref{eq:no3tors}, and
$s(D)=-\varepsilon-2j$ with $\beta(D)=\frac{5k}{12}+\frac\varepsilon2+j\ge0$.

\begin{lemma}\label{lem:no3a}
In the case $(\hat q,e)=(4,2)$:
\begin{enumerate}
\item $h^0(\hat X,s\Theta)\ge3,8,14,24$ for $s=1,2,3,4$, and $\hat X\cong\PP^3$;
\item let $\Theta_0=\bar f(\bar D_2')$, where $D_2'$ is the unique member of $|2A_X+\sigma|$, and
$\gamma_0=v(\Theta_0)$; then $b=4\gamma_0+3$, $v(\bar f(\bar E))=2\gamma_0+4$, and for $G\in|s\Theta|$
the class of its proper transform is $nA_X+\varepsilon \sigma$ with $n=s(2+\gamma_0)-v(G)$, $\varepsilon\equiv s$;
\item for planes $G$ of $\hat X=\PP^3$ the subspaces $\{v\ge\gamma_0+1\}$, $\{v\ge\gamma_0\}$, $\{v\ge\gamma_0-1\}$ have
dimensions $0,1,2$ and $\{v\ge\gamma_0-2\}$ has dimension at least $3$; for quadrics $G$ the subspaces
$\{v\ge2\gamma_0+3\}$, $\{v\ge2\gamma_0+2\}$, $\{v\ge2\gamma_0+1\}$, $\{v\ge2\gamma_0\}$ have dimensions $1,1,2,4$.
\end{enumerate}
\end{lemma}

\begin{proof}
(i) A section $\phi$ of $nA_X+s\sigma$ with $v_E(\phi)\ge\frac{5n-6s}{12}$ gives, by the formula for $s(D)$,
a member of $|s\Theta|$ after adding a multiple of $\bar E$, and these maps are injective. Near $P$ the
condition costs at most the number of monomials of weight $<5n-6s$ and of the local character
$5n+6s\bmod12$. The resulting bounds (\texttt{examples/no3/bounds.py}; e.g. $s=1$: $p'_4=4$ and the only
condition is the coefficient of $\xi_1^2$) are as stated. In particular
$\g(\hat X)=h^0(-K_{\hat X})-2\ge22$, so $\hat X\cong\PP^3$ or $X_4\subset\PP(1^3,2,3)$ by
\cite[Theorem~1.2(vi)]{Pro13b}, and the latter has $h^0(2\Theta)=7$.

(ii) The formula for $s(D)$ gives $s(D_2')=1$, and \eqref{eq:no3rev} for $\Theta_0$ gives $b=4\gamma_0+3$;
pushing $K_{\hat X}+2\bar f(\bar E)\sim0$ down to $X$ gives $v(\bar f(\bar E))=\frac{5+b}2$, and
\eqref{eq:no3rev} gives the formula for $n$.

(iii) For $(k,\varepsilon)=(2,1),(3,1)$ the formula for $s(D)$ allows only $s=1$, and for $\varepsilon=0$,
$k\le4$ only $s\in\{0,2\}$, where $s=0$ means $D\in\{kD_1\}$, corresponding to $\bar f(\bar E)$. So the
maps $D\mapsto$ (image $+$ multiple of $\bar f(\bar E)$) identify $H^0(kA_X+\sigma)$, $k=2,3$, with
$\{G\in H^0(\Theta):v(G)\ge\gamma_0+2-k\}$, and $H^0(kA_X)$, $k\le4$, with
$\{G\in H^0(2\Theta):v(G)\ge2\gamma_0+4-k\}$. A plane with $v\ge\gamma_0+1$ would give a member of
$|A_X+\sigma|$. Finally the members of $|4A_X+\sigma|$ with $s=1$ form, by (i), a space of dimension $\ge3$.
\end{proof}

We use two facts about a divisorial valuation $v$ centred at a nonsingular point. (a) By
\cite{Kaw01}, if $v$ is the valuation of a divisorial contraction to the point, then $v$ is the monomial
valuation with weights $(1,a,a')$, $\gcd(a,a')=1$, in suitable local coordinates, and its log discrepancy
is $A(v)=1+a+a'$. (b) For any local coordinates $x_1,x_2,x_3$ one has $A(v)\ge v(x_1)+v(x_2)+v(x_3)$,
with equality only if $v$ is monomial in $x_1,x_2,x_3$: this is $A(v)\ge A(r(v))$ for the retraction $r$
onto the monomial valuations in these coordinates \cite[Corollary~5.4]{JM12}.

\begin{proposition}\label{prop:no3}
The numerical type No.~$3$ does not occur.
\end{proposition}

\begin{proof}
By Lemmas~\ref{lem:no3links}, \ref{lem:no3b} and~\ref{lem:no3a} we are in the case $\hat X=\PP^3$. Since
$b=4\gamma_0+3\ge3$, the centre of $\bar f$ is not a curve (a divisorial contraction to a curve in a smooth
threefold has discrepancy $1$ over the general point of the curve); so it is a point $\hat P$, which lies on
the quadric $Q=\bar f(\bar E)$ (irreducible), and $v$ is as in (a) with $a+a'=b=4\gamma_0+3$. The planes through $\hat P$ form a $3$-dimensional space on which a general element
has $v=1$. By Lemma~\ref{lem:no3a}(iii), $\dim\{v\ge\gamma_0-1\}=2$ forces $\gamma_0\ge3$, and then
$\{v\ge\gamma_0-2\}$ is contained in the planes through $\hat P$ and has dimension $3$, so $\gamma_0-2\le1$.
Thus $\gamma_0=3$, $A(v)=16$, $v(Q)=10$, and the affine linear coordinates at $\hat P$ can be chosen as
$\ell_1,\ell_2,\ell_3$ with $v(\ell_i)=i$.

\emph{$Q$ is a cone with vertex $\hat P$.} Then $Q$ is a quadratic form in $\ell_1,\ell_2,\ell_3$; as $v(Q)=10>6$,
the monomials $\ell_1^2,\ell_1\ell_2$ do not occur, and $Q=\ell_1'\ell_3+\kappa\ell_2^2$ with $v(\ell_1')=1$ and
$\kappa\ne0$ ($Q$ has rank $3$). Take coordinates $w_1,w_2,w_3$ as in (a); replacing $w_1$ by $\ell_1'$ does not change $v$, so we may assume $w_1=\ell_1'$. Write the linear parts of
$\ell_2,\ell_3$ as $\sum \mu_iw_i$, $\sum \nu_iw_i$; then $(\nu_2,\nu_3)$, $(\mu_2,\mu_3)$ are independent, and every
monomial of weight $<10$ in the degree-$2$ part $w_1\sum \nu_iw_i+\kappa(\sum \mu_iw_i)^2$ of $Q$ must vanish.
If $\mu_2\mu_3\ne0$, then $a,a'\ge5$, so $\{a,a'\}=\{7,8\}$, and the coefficients of $w_1w_2,w_1w_3$ give
$(\nu_2,\nu_3)=-2\kappa \mu_1(\mu_2,\mu_3)$, a contradiction. If $\mu_3=0\ne \mu_2$, then $a\ge5$ and $\nu_3\ne0$ gives
$a'\ge9$, so $(a,a')\in\{(5,10),(6,9)\}$, which are not coprime; the case $\mu_2=0\ne \mu_3$ is symmetric, and $\mu_2=\mu_3=0$ would give $v(\ell_2)\le1$ or a degenerate Jacobian.

\emph{$Q$ is nonsingular at $\hat P$.} Choose projective coordinates with $\hat P=(1:0:0:0)$ and
$Q=\{x_0x_3=x_1x_2\}$, resp. $Q=\{x_0x_3=x_1^2\}$ if $Q$ is a cone, and put $u=x_1/x_0$, $u'=x_2/x_0$, $q=Q/x_0^2$;
then $u,u',q$ are local coordinates and every plane, resp. quadric, is congruent modulo $q$ to a
polynomial in $u,u'$ in an explicit $4$-, resp. $9$-dimensional space. Since all values in
Lemma~\ref{lem:no3a}(iii) are $<10=v(q)$, the filtrations of Lemma~\ref{lem:no3a}(iii), with $\gamma_0=3$,
are filtrations of these spaces by the restriction of $v$ to functions of $u,u'$: planes with $v\ge1,2,3,4$
have dimensions $3,2,1,0$, and quadrics modulo $q$ with $v\ge6,7,8$ have dimensions $3,1,0$. One of
$u,u'$, say $x$, has $v(x)=1$; let $y$ be the other one. If the initial forms of $x$ and $y$ are
dependent, replace $y$ by $y-\lambda x^j$, and repeat; by (b),
\[
16=A(v)\ge v(x)+v(y-h(x))+v(q)=1+v(y-h(x))+10,
\]
so the process stops with $t:=v(y-h(x))\le5$ and independent initial forms, and the restriction of $v$ to
functions of $u,u'$ is the monomial valuation with weights $(1,t)$ in $x$, $y-h(x)$ \cite[Corollary~5.4]{JM12};
after rescaling, $h(x)=\lambda x+h_2x^2+\dots+h_{t-1}x^{t-1}$ with $\lambda\in\{0,1\}$. For $t=1$ only one
plane has $v\ge2$. For $2\le t\le5$, the conditions that the quadrics with $v\ge6$, resp. $v\ge7$, span a
space of dimension $\ge3$, resp. $\ge1$, and that some plane has $v\ge3$, are given by minors of explicit
matrices in the $h_i$; for every choice of $Q$, $x$ and $\lambda$, the ideal they generate is either the unit
ideal or its radical contains the minors expressing that some plane has $v\ge4$
(\texttt{examples/no3/p3\_local.py}). In both cases Lemma~\ref{lem:no3a}(iii) is violated.
\end{proof}

\subsection{The canonical extension: Nos.~26, 31, 54, 60, 80, 100, 114, 128, 151, 166}\label{ss:liu26}
For $\qQ\le4$, H.~Liu \cite[Theorem~4.3]{Liu26} sharpens Theorem~\ref{thm:LL} for $\QQ$-factorial terminal Fano
threefolds of Picard number one; in dimension three his bounds read $c_1^3\le\frac83c_2c_1$ for $\qQ\le2$ and
$c_1^3\le\frac{4q^2}{q^2+2q-2}c_2c_1$ for $q=\qQ\ge3$. We give a short proof of the two cases $\qQ=2,3$ that we
need, following \cite{Liu26}, together with the equality case for $\qQ=2$.

Let $X$ be a $\QQ$-Fano threefold, $q=\qQ(X)$, $A$ a Weil divisor with $-K_X\sim_\QQ qA$ and $H=-K_X$. Let
\[
0\longrightarrow\OO_X\longrightarrow\mathcal E\longrightarrow T_X\longrightarrow0
\]
be the canonical extension of \cite[Construction~4.1, Definition~4.5]{GKP22}; $\mathcal E$ is reflexive of rank $4$,
$c_1(\mathcal E)=c_1(X)=H$ and $c_2(\mathcal E)\cdot H=c_2(X)\cdot H$. Fix $m$ with $mH$ Cartier and very ample. Since the
singular points of $X$ are isolated, general complete intersection curves $C$ and surfaces $S$ cut out by members of
$|kmH|$, $k\gg0$, lie in the smooth locus, avoid the finite set where $\mathcal E$ is not locally free, and, by
Mehta--Ramanathan \cite{MR82} and Flenner \cite{Fle84}, who treats normal projective varieties, the Harder--Narasimhan filtration of $\mathcal E$ with respect to $H$ restricts to that
of $\mathcal E|_C$ and of $\mathcal E|_S$. In particular $c_2(X)\cdot H$ is the number $24-\sum(r_P-1/r_P)$ of
Theorem~\ref{thm:c2}, also when computed with orbifold Chern classes.

\begin{lemma}\label{lem:canext}
\begin{enumerate}
\item Every Harder--Narasimhan factor $\mathcal G_i$ of $\mathcal E$ has positive slope. Writing $c_1(\mathcal G_i)\equiv
q_iA$ and $r_i=\operatorname{rk}\mathcal G_i$, the $q_i$ are positive integers, $\sum q_i=q$, $\sum r_i=4$ and
$q_1/r_1>q_2/r_2>\dots$.
\item If $r_1=1$, then $2q_1<q$.
\item If $q=2$, then $\mathcal E$ is semistable and $H^3\le\frac83\,c_2(X)\cdot H$.
\item If $q=3$, then $H^3\le\frac{36}{13}\,c_2(X)\cdot H$.
\end{enumerate}
\end{lemma}

\begin{proof}
(i) The restriction $T_X|_C$ is ample \cite[Proposition~3.6]{LL25a}, and the extension restricted to $C$ is not
split, since its class maps to $c_1(H|_C)\ne0$. So $\mathcal E|_C$ is nef. Let $Q$ be its last Harder--Narasimhan
quotient. If $\deg Q=0$, the map $\OO_C\to Q$ is nonzero (otherwise $Q$ is a quotient of $T_X|_C$), its saturated
image is $\OO_C$ by semistability of $Q$, and either $Q=\OO_C$, which splits the extension, or $Q/\OO_C$ is a
degree-zero quotient of $T_X|_C$; both are impossible. Hence all slopes are positive. As $\Cl(X)/\Clt(X)=\ZZ A$, each
$c_1(\mathcal G_i)$ is numerically an integer multiple $q_iA$, so $q_i\ge1$; the other statements are clear.

(ii) The composition $\mathcal G_1\hookrightarrow\mathcal E\to T_X$ is nonzero: otherwise $\mathcal G_1\subset\OO_X$,
contradicting its positive slope. Since $\mathcal G_1$ has rank one and $T_X$ is torsion-free, this composition is
injective. The saturation $L$ of its image in $T_X$ satisfies $(2+1/r_X)\,c_1(L)\cdot H^2\le H^3$ by \cite[Proposition~3.8]{LL25a}, where $r_X$ is the
global index; since $c_1(L)\cdot H^2\ge c_1(\mathcal G_1)\cdot H^2=q_1A\cdot H^2$, we get $(2+1/r_X)q_1\le q$, so $2q_1<q$.

(iii), (iv) On $S$ the Chern character is additive along the filtration, and $\mathcal G_i|_S$ is semistable, so
Bogomolov's inequality \cite{Bog78} gives $\operatorname{ch}_2(\mathcal G_i|_S)\le c_1(\mathcal G_i|_S)^2/(2r_i)$. Summing,
\begin{equation}\label{eq:canext}
2\,c_2(X)\cdot H\ \ge\ \Bigl(1-\frac1{q^2}\sum_i\frac{q_i^2}{r_i}\Bigr)H^3 .
\end{equation}
If $q=2$ and $\mathcal E$ is not semistable, then by (i) there are two factors with $q_1=q_2=1$ and $r_1<r_2$, so
$(r_1,r_2)=(1,3)$, contradicting (ii). Thus $\mathcal E$ is semistable and \eqref{eq:canext} with $\sum q_i^2/r_i=1$
gives $H^3\le\frac83c_2(X)\cdot H$. If $q=3$, the possible types $((q_i,r_i))_i$ allowed by (i) and (ii) are
$((3,4))$, $((1,1),(2,3))$ and $((2,2),(1,2))$, with $\sum q_i^2/r_i=\frac94,\frac73,\frac52$, and \eqref{eq:canext}
gives $H^3\le\frac83,\frac{27}{10},\frac{36}{13}$ times $c_2(X)\cdot H$, respectively.
\end{proof}

\begin{proposition}\label{prop:liu26}
The numerical types No.~$26$, $31$, $54$, $60$, $80$, $100$, $114$, $128$, $151$ and $166$ do not occur.
\end{proposition}

\begin{proof}
Write $c_1^3=H^3=\qQ^3A_X^3$ and $c_2c_1=c_2(X)\cdot H=24-\sum_{P\in\B(X)}(r_P-1/r_P)$ (Theorem~\ref{thm:c2}):
{\renewcommand{\arraystretch}{1.5}
\[
\begin{array}{c|c|c|c|c|c}
\text{No.} & \B(X) & \qQ & c_1^3 & c_2c_1 & c_1^3/c_2c_1\\\hline
26 & [2\times[2,1],[5,2],[14,5]] & 3 & \frac{459}{70} & \frac{159}{70} & \frac{153}{53}>\frac{36}{13}\\
31 & [2\times[2,1],[4,1],[5,2],[8,3]] & 3 & \frac{513}{40} & \frac{183}{40} & \frac{171}{61}>\frac{36}{13}\\
80 & [8\times[2,1],[9,1]] & 2 & \frac{80}{9} & \frac{28}{9} & \frac{20}{7}>\frac83\\
114 & [5\times[2,1],[6,1],[7,2]] & 2 & \frac{232}{21} & \frac{80}{21} & \frac{29}{10}>\frac83\\
151 & [[5,1],2\times[9,4]] & 2 & \frac{176}{45} & \frac{64}{45} & \frac{11}{4}>\frac83\\\hline
54 & [8\times[2,1],[11,1]] & 2 & \frac{32}{11} & \frac{12}{11} & \frac83\\
60 & [8\times[2,1],2\times[5,2]] & 2 & \frac{32}{5} & \frac{12}{5} & \frac83\\
100 & [3\times[2,1],[7,2],[10,3]] & 2 & \frac{256}{35} & \frac{96}{35} & \frac83\\
128 & [8\times[2,1],[5,2],[7,2]] & 2 & \frac{32}{35} & \frac{12}{35} & \frac83\\
166 & [2\times[2,1],4\times[4,1],[5,2]] & 2 & \frac{16}{5} & \frac65 & \frac83
\end{array}
\]
}(supplementary file \texttt{classification/liu26.py}). The first five violate Lemma~\ref{lem:canext}(iii), (iv). The last
five have $\qQ=2$ and attain equality in Lemma~\ref{lem:canext}(iii). Then $\mathcal E$ is semistable with equality in
the Bogomolov inequality, and by \cite[Theorem~1.3]{GKP22} (its other alternative, a quotient of an Abelian variety, is excluded since $-K_X$
is ample)
$X\cong\PP^3/G$ for a finite group $G\subset\operatorname{PGL}_4$ acting freely in codimension one. Then
$K_{\PP^3}$ is the pull-back of $K_X$ and $|G|=64/H^3$, which is $22$, $10$, $\frac{35}4$, $70$, $20$ for Nos.~$54$,
$60$, $100$, $128$, $166$. For No.~$100$ this is not an integer. In the other cases $G$ contains an involution $g$. Up to a
scalar, $g$ has eigenvalue multiplicities $(3,1)$ or $(2,2)$ on $\CC^4$. In the first case $g$ fixes a plane, which is
excluded. In the second it fixes two disjoint lines, and at a general point of them the stabiliser of $G$ contains no
pseudo-reflection (one would fix a plane), so $X$ is singular along the image of these lines. This contradicts the fact
that terminal singularities are isolated.
\end{proof}

No realised numerical type violates the bounds of Lemma~\ref{lem:canext}, and among the open numerical types with
$\qQ\le3$ no other one attains equality.

\section{The table}\label{sec:table}

Each row of Table~\ref{tab:main} is one numerical type. The columns are as follows.
\begin{description}
\item[No.] The number of the numerical type, used throughout the paper.
\item[$\B(X)$] The basket. An entry $[r,a]$ is a point of type $\frac1r(1,a,-a)$, and $m{\times}[r,a]$ stands for
$m$ equal entries. The entries are ordered by $(r,a)$; this order is used in the columns of local types below.
\item[$T_X$] The maximal torsion group of the row, defined before Theorem~\ref{thm:numerical}; $\Clt(X)$ is isomorphic to a subgroup of
$T_X$, and equals it if $|T_X|$ is prime.
\item[$\qW$, $\qQ$] The Weil and $\QQ$-Fano indices.
\item[$A_X^3$] The degree of $A_X$; thus $(-K_X)^3=\qQ^3A_X^3$.
\item[generators of $T_X$] By the injectivity of \eqref{eq:local}, an element $\tau\in T_X$ is determined by the
vector $\bigl(i_P(\tau)\bigr)_P$ of its local types, where $\tau\sim i_P(\tau)K_X$ near $P$, $0\le i_P(\tau)<r_P$,
and $P$ runs over the entries of $\B(X)$ in the order of the second column. We give this vector for a
generator of $T_X$; if $T_X\cong\ZZ/2\times\ZZ/2$ or $\ZZ/3\times\ZZ/3$, we give two generators, separated by a
semicolon. An entry $0$ means that $\tau$ is Cartier at $P$, and the order of $\tau$ is the least common multiple
of the numbers $r_P/\gcd(r_P,i_P(\tau))$.
\item[$A_X$] The vector $\bigl(i_P(A_X)\bigr)_P$ of local types of $A_X$, in the same order. Since $i_P(K_X)=1$,
the class $\tau_0\in T_X$ with $-K_X\sim\qQ A_X+\tau_0$ has local types $i_P(\tau_0)\equiv-1-\qQ\,i_P(A_X)\pmod{r_P}$.
For example, in No.~$1$ we have $\B(X)=[[2,1],[6,1],[10,3]]$; the generator $\sigma$ of $T_X=\ZZ/2$ has local
types $(0,3,5)$, so $\sigma$ is Cartier at the point of index $2$, $\sigma\sim3K_X$ near the point of index $6$ and
$\sigma\sim5K_X$ near the point of index $10$. The local types of $A_X$ are $(1,5,7)$, and
$7\cdot(1,5,7)\equiv(-1,-1,-1)$ modulo $(2,6,10)$, so $\tau_0=0$ and $-K_X\sim7A_X$.

The divisor $A_X$ is normalised as follows. If $\qW=\qQ$, it is the unique divisor with $-K_X\sim\qQ A_X$
(Lemma~\ref{lem:qwqq}). If $\qW=1$, it is the divisor among $A_X+\tau$, $\tau\in T_X$, for which the sequence
$\bigl(h^0(X,kA_X)\bigr)_{k\ge1}$ is lexicographically largest; ties are broken by the smallest vector of local
types. The numerical type does not depend on this choice. This normalises the maximal series over $T_X$: if $\Clt(X)$
is a proper subgroup of $T_X$, the shift used here need not lie in $\Clt(X)$, and an actual fundamental divisor may have
different values of $h^0(X,kA_X)$ (for instance in Nos.~$175$ and $176$, where $\Clt(X)\cong\ZZ/2$ is possible). Some proofs in
\S\ref{sec:exclusions} use another choice, which is then stated.
\item[$h^0(kA_X)$] The numbers $h^0(X,kA_X)$ for $k=1,\dots,5$, for this choice of $A_X$, computed by
Theorem~\ref{thm:RR}; like the column $A_X$, these are data of the maximal series over $T_X$.
\item[$\dim|{-K_X}|$] The number $h^0(X,-K_X)-1$, computed by Theorem~\ref{thm:RR} for the class
$\qQ A_X+\tau_0$; the value $-1$ means that $|-K_X|=\varnothing$. Like the other numerical columns, this is a
value attached to the numerical type: each row of Table~\ref{tab:main} is a numerical candidate, and a
Fano--Enriques threefold with these invariants need not exist (see the last column). If one exists, then
$\dim|-K_X|$ has the stated value.
\item[$\B(Y)$] The basket given by the formula of Lemma~\ref{lem:cover} for the group $T_X$; it is the basket of the cover
$\pi\colon Y\to X$ defined by $T_X$ when $\Clt(X)=T_X$. For instance, for No.~$155$ the column refers to
$T_X=(\ZZ/3)^2$, whereas an example would have $\Clt(X)\cong\ZZ/3$ by Theorem~\ref{thm:nonexistence}(ii), and its cover
defined by $\Clt(X)$ would have basket $6\times[3,1]$.
\item[$\exists$] The status: $+$ if the numerical type is realised and $-$ if it is excluded, followed by the
number of the proposition where this is proved (a link in the electronic version). For the cyclic quotients of
Table~\ref{tab:examples} this is Proposition~\ref{prop:examples}. An empty entry means that the existence is open.
\end{description}

\begin{landscape}
\thispagestyle{empty}\pagestyle{empty}
{\tiny
\setlength{\tabcolsep}{2.5pt}
\renewcommand{\arraystretch}{1.15}
\begin{longtable}{r|l|c|c|c|c|l|l|l|c|l|l}
\caption{Numerical types of Fano--Enriques threefolds with $\qQ\ge2$.}\label{tab:main}\\
No. & $\B(X)$ & $T_X$ & $\qW$ & $\qQ$ & $A_X^3$ & generators of $T_X$ & $A_X$ & $h^0(kA_X)$ & $\dim|{-K_X}|$ & $\B(Y)$ & $\exists$\\\hline
\endfirsthead
No. & $\B(X)$ & $T_X$ & $\qW$ & $\qQ$ & $A_X^3$ & generators of $T_X$ & $A_X$ & $h^0(kA_X)$ & $\dim|{-K_X}|$ & $\B(Y)$ & $\exists$\\\hline
\endhead
\input{table_main}
\end{longtable}
}
\end{landscape}
\pagestyle{headings}

\begin{table}[p]
\centering\scriptsize\setlength{\tabcolsep}{2.5pt}
\renewcommand{\arraystretch}{1.25}
\caption{The possible torsion groups $\Clt(X)$ for the numerical types with $|T_X|$ not prime, up to automorphisms of
$T_X$ and shifts of the fundamental divisor (\S\ref{sec:numerical}). The columns give: the number in
Table~\ref{tab:main} and the maximal group $T_X$; the group $\Clt(X)$ and its generators, by their local types in the order of
the basket in Table~\ref{tab:main}; a compatible fundamental divisor $A_X$, by its local types, normalised as in
Table~\ref{tab:main} among its shifts by $\Clt(X)$; the Weil index; the sums $\sum_{\tau\in\Clt(X)}h^0(X,kA_X+\tau)$ for
$k=1,\dots,5$; the basket of the cover defined by $\Clt(X)$; and the status of the variant, as in Table~\ref{tab:main}
(empty if open). For Nos.~$156$--$167$, $\Clt(X)=T_X$ is forced.}\label{tab:subgroups}
\begin{tabular}{rcccccccc}
No. & $T_X$ & $\Clt(X)$ & generators of $\Clt(X)$ & $A_X$ & $\qW$ & $\sum_\tau h^0(kA_X+\tau)$ & $\B(Y)$ & $\exists$\\\hline
\input{table_subgroups}
\end{tabular}
\end{table}

\end{document}

%% file: table_examples.tex
1 & $Y_{8}\subset\PP(1,2,3,4,5)$ & $\mumu_{2}(0,1,1,1,1)$ & $(0)$ & $[2,1],\,[6,1],\,[10,3]$ & 7 & 7 & $\frac{1}{30}$ \\
2 & $Y_{6}\subset\PP(1,2,3^{2},4)$ & $\mumu_{2}(0,1,0,1,1)$ & $(0)$ & $2{\times}[2,1],\,[3,1],\,[4,1],\,[8,3]$ & 7 & 7 & $\frac{1}{24}$ \\
4 & $Y_{4}\subset\PP(1^{2},2^{2},3)$ & $\mumu_{2}(0,1,1,1,0)$ & $(0)$ & $[2,1],\,2{\times}[4,1],\,[6,1]$ & 5 & 5 & $\frac{1}{6}$ \\
6 & $Y_{6}\subset\PP(1^{2},2,3,4)$ & $\mumu_{2}(0,1,1,0,1)$ & $(0)$ & $2{\times}[2,1],\,[4,1],\,[8,3]$ & 5 & 5 & $\frac{1}{8}$ \\
7 & $Y_{6}\subset\PP(1,2^{2},3^{2})$ & $\mumu_{3}(0,1,2,1,2)$ & $(0)$ & $[2,1],\,2{\times}[9,2]$ & 5 & 5 & $\frac{1}{18}$ \\
9 & $Y_{4}\subset\PP(1^{3},2,3)$ & $\mumu_{2}(0,1,1,1,0)$ & $(0)$ & $5{\times}[2,1],\,[6,1]$ & 1 & 4 & $\frac{1}{3}$ \\
10 & $Y_{6,8}\subset\PP(1,2,3^{2},4,5)$ & $\mumu_{2}(0,1,0,1,1,0)$ & $(0,0)$ & $3{\times}[2,1],\,[3,1],\,[10,3]$ & 1 & 4 & $\frac{1}{15}$ \\
15 & $Y_{6}\subset\PP(1^{2},2,3^{2})$ & $\mumu_{3}(0,1,1,1,2)$ & $(0)$ & $2{\times}[9,4]$ & 4 & 4 & $\frac{1}{9}$ \\
16 & $\PP^3$ & $\mumu_{5}(1,2,3,4)$ & $-$ & $4{\times}[5,2]$ & 4 & 4 & $\frac{1}{5}$ \\
22 & $Y_{4}\subset\PP(1^{3},2^{2})$ & $\mumu_{2}(0,0,1,1,1)$ & $(0)$ & $4{\times}[2,1],\,2{\times}[4,1]$ & 3 & 3 & $\frac{1}{2}$ \\
33 & $Y_{3}\subset\PP(1^{4},2)$ & $\mumu_{3}(0,0,1,2,1)$ & $(0)$ & $4{\times}[3,1],\,[6,1]$ & 1 & 3 & $\frac{1}{2}$ \\
36 & $Y_{4}\subset\PP(1^{3},2^{2})$ & $\mumu_{4}(0,1,3,1,3)$ & $(0)$ & $2{\times}[2,1],\,2{\times}[8,3]$ & 3 & 3 & $\frac{1}{4}$ \\
37 & $Y_{2}\subset\PP^{4}$ & $\mumu_{5}(0,1,2,3,4)$ & $(0)$ & $4{\times}[5,2]$ & 3 & 3 & $\frac{2}{5}$ \\
38 & $Y_{4}\subset\PP(1^{4},2)$ & $\mumu_{2}(0,0,1,1,1)$ & $(0)$ & $8{\times}[2,1]$ & 1 & 2 & $1$ \\
39 & $Y_{2,2}\subset\PP^{5}$ & $\mumu_{2}(0,0,0,1,1,1)$ & $(0,0)$ & $8{\times}[2,1]$ & 1 & 2 & $2$ \\
42 & $Y_{4,6}\subset\PP(1^{2},2^{2},3^{2})$ & $\mumu_{2}(0,1,1,1,0,1)$ & $(0,1)$ & $2{\times}[2,1],\,2{\times}[6,1]$ & 1 & 2 & $\frac{1}{3}$ \\
46 & $Y_{6,8}\subset\PP(1^{2},2,3,4,5)$ & $\mumu_{2}(0,1,1,0,1,0)$ & $(0,0)$ & $3{\times}[2,1],\,[10,3]$ & 1 & 2 & $\frac{1}{5}$ \\
50 & $Y_{4,4}\subset\PP(1^{3},2^{2},3)$ & $\mumu_{2}(0,1,1,1,1,1)$ & $(0,0)$ & $5{\times}[2,1],\,[6,1]$ & 1 & 2 & $\frac{2}{3}$ \\
61 & $Y_{4,6}\subset\PP(1^{2},2^{2},3^{2})$ & $\mumu_{2}(0,1,0,1,0,1)$ & $(0,0)$ & $8{\times}[2,1],\,[3,1]$ & 1 & 2 & $\frac{1}{3}$ \\
109 & $Y_{8,10}\subset\PP(1,2,3,4,5^{2})$ & $\mumu_{2}(0,1,0,1,0,1)$ & $(0,0)$ & $5{\times}[2,1],\,[5,2],\,[6,1]$ & 1 & 2 & $\frac{1}{15}$ \\
134 & $Y_{4,6}\subset\PP(1^{2},2^{2},3^{2})$ & $\mumu_{3}(0,1,1,1,1,2)$ & $(2,0)$ & $2{\times}[9,1]$ & 2 & 2 & $\frac{2}{9}$ \\
141 & $Y_{3}\subset\PP^{4}$ & $\mumu_{3}(0,0,1,1,2)$ & $(0)$ & $6{\times}[3,1]$ & 2 & 2 & $1$ \\
156 & $Y_{4}\subset\PP(1^{4},2)$ & $\mumu_{4}(0,1,2,3,3)$ & $(2)$ & $2{\times}[2,1],\,4{\times}[4,1]$ & 1 & 2 & $\frac{1}{2}$ \\
157 & $Y_{2,2}\subset\PP^{5}$ & $\mumu_{4}(0,0,1,1,2,3)$ & $(0,2)$ & $2{\times}[2,1],\,4{\times}[4,1]$ & 1 & 2 & $1$ \\
163 & $Y_{4,4}\subset\PP(1^{3},2^{2},3)$ & $\mumu_{4}(0,1,2,1,3,0)$ & $(0,2)$ & $2{\times}[2,1],\,[4,1],\,[12,5]$ & 1 & 2 & $\frac{1}{3}$ \\
168 & $Y_{4}\subset\PP(1^{4},2)$ & $\mumu_{5}(0,1,3,4,4)$ & $(3)$ & $4{\times}[5,1]$ & 2 & 2 & $\frac{2}{5}$ \\
170 & $Y_{3}\subset\PP^{4}$ & $\mumu_{5}(0,1,2,3,4)$ & $(0)$ & $4{\times}[5,2]$ & 2 & 2 & $\frac{3}{5}$ \\

%% file: table_main.tex
1 & $[2,1],\,[6,1],\,[10,3]$ & $\ZZ/2$ & 7 & 7 & $1/30$ & $(0,3,5)$ & $(1,5,7)$ & $1,1,1,2,3$ & 6 & $2{\times}[2,1],\,[3,1],\,[5,2]$ & $+$\,\ref{prop:examples} \\
2 & $2{\times}[2,1],\,[3,1],\,[4,1],\,[8,3]$ & $\ZZ/2$ & 7 & 7 & $1/24$ & $(1,1,0,2,4)$ & $(1,1,2,1,1)$ & $0,1,2,3,4$ & 7 & $[2,1],\,2{\times}[3,1],\,[4,1]$ & $+$\,\ref{prop:examples} \\
3 & $2{\times}[4,1],\,[12,5]$ & $\ZZ/2$ & 5 & 5 & $1/12$ & $(2,0,6)$ & $(3,3,7)$ & $1,1,2,4,6$ & 5 & $[2,1],\,2{\times}[4,1],\,[6,1]$ & $-$\,\ref{prop:no3} \\
4 & $[2,1],\,2{\times}[4,1],\,[6,1]$ & $\ZZ/2$ & 5 & 5 & $1/6$ & $(1,2,2,3)$ & $(1,3,3,1)$ & $1,2,4,8,12$ & 11 & $2{\times}[2,1],\,[3,1]$ & $+$\,\ref{prop:examples} \\
5 & $[2,1],\,[4,1],\,[14,3]$ & $\ZZ/2$ & 5 & 5 & $1/28$ & $(1,0,7)$ & $(1,3,11)$ & $1,1,1,2,3$ & 2 & $2{\times}[4,1],\,[7,3]$ & $+$\,\ref{prop:tom} \\
6 & $2{\times}[2,1],\,[4,1],\,[8,3]$ & $\ZZ/2$ & 5 & 5 & $1/8$ & $(1,1,2,4)$ & $(1,1,3,3)$ & $1,2,3,6,9$ & 8 & $[2,1],\,[4,1]$ & $+$\,\ref{prop:examples} \\
7 & $[2,1],\,2{\times}[9,2]$ & $\ZZ/3$ & 5 & 5 & $1/18$ & $(0,6,3)$ & $(1,7,7)$ & $1,1,1,2,4$ & 3 & $3{\times}[2,1],\,2{\times}[3,1]$ & $+$\,\ref{prop:examples} \\
8 & $3{\times}[2,1],\,[10,3]$ & $\ZZ/2$ & 1 & 4 & $2/5$ & $(1,1,1,5)$ & $(0,0,0,6)$ & $2,4,8,14,23$ & 13 & $[5,2]$ & $-$\,\ref{prop:genusbound} \\
9 & $5{\times}[2,1],\,[6,1]$ & $\ZZ/2$ & 1 & 4 & $1/3$ & $(1,1,1,1,1,3)$ & $(1,0,0,0,0,5)$ & $2,4,7,12,19$ & 11 & $[3,1]$ & $+$\,\ref{prop:examples} \\
10 & $3{\times}[2,1],\,[3,1],\,[10,3]$ & $\ZZ/2$ & 1 & 4 & $1/15$ & $(1,1,1,0,5)$ & $(0,0,0,2,6)$ & $1,1,2,3,5$ & 2 & $2{\times}[3,1],\,[5,2]$ & $+$\,\ref{prop:examples} \\
11 & $3{\times}[2,1],\,[5,2],\,[10,3]$ & $\ZZ/2$ & 1 & 4 & $1/5$ & $(1,1,1,0,5)$ & $(0,0,0,1,6)$ & $1,2,4,7,12$ & 6 & $3{\times}[5,2]$ & $-$\,\ref{prop:11} \\
12 & $3{\times}[2,1],\,[7,2],\,[10,3]$ & $\ZZ/2$ & 1 & 4 & $4/35$ & $(1,1,1,0,5)$ & $(0,0,0,5,6)$ & $1,1,2,4,7$ & 3 & $[5,2],\,2{\times}[7,2]$ &  \\
13 & $5{\times}[2,1],\,[5,2],\,[6,1]$ & $\ZZ/2$ & 1 & 4 & $2/15$ & $(1,1,1,1,1,0,3)$ & $(1,0,0,0,0,1,5)$ & $1,2,3,5,8$ & 4 & $[3,1],\,2{\times}[5,2]$ & $-$\,\ref{prop:13} \\
14 & $5{\times}[2,1],\,[6,1],\,[7,2]$ & $\ZZ/2$ & 1 & 4 & $1/21$ & $(1,1,1,1,1,3,0)$ & $(1,0,0,0,0,5,5)$ & $1,1,1,2,3$ & 1 & $[3,1],\,2{\times}[7,2]$ &  \\
15 & $2{\times}[9,4]$ & $\ZZ/3$ & 4 & 4 & $1/9$ & $(6,3)$ & $(2,2)$ & $0,2,2,5,6$ & 4 & $2{\times}[3,1]$ & $+$\,\ref{prop:examples} \\
16 & $4{\times}[5,2]$ & $\ZZ/5$ & 4 & 4 & $1/5$ & $(4,3,2,1)$ & $(1,1,1,1)$ & $0,2,4,7,12$ & 6 & $\varnothing$ & $+$\,\ref{prop:examples} \\
17 & $[2,1],\,[14,3]$ & $\ZZ/2$ & 3 & 3 & $2/7$ & $(1,7)$ & $(1,9)$ & $1,2,5,9,14$ & 4 & $[7,3]$ & $+$\,\ref{prop:17} \\
18 & $2{\times}[2,1],\,[14,5]$ & $\ZZ/2$ & 3 & 3 & $9/14$ & $(1,0,7)$ & $(1,1,9)$ & $2,5,10,18,29$ & 9 & $2{\times}[2,1],\,[7,2]$ &  \\
19 & $3{\times}[2,1],\,[10,3]$ & $\ZZ/2$ & 3 & 3 & $4/5$ & $(1,1,1,5)$ & $(1,1,1,3)$ & $2,6,12,22,36$ & 11 & $[5,2]$ & $-$\,\ref{prop:genusbound} \\
20 & $3{\times}[2,1],\,[14,5]$ & $\ZZ/2$ & 3 & 3 & $1/7$ & $(1,0,0,7)$ & $(1,1,1,9)$ & $1,2,3,5,7$ & 2 & $4{\times}[2,1],\,[7,2]$ &  \\
21 & $4{\times}[2,1],\,[10,3]$ & $\ZZ/2$ & 3 & 3 & $3/10$ & $(1,1,1,0,5)$ & $(1,1,1,1,3)$ & $1,3,5,9,14$ & 4 & $2{\times}[2,1],\,[5,2]$ & $+$\,\ref{prop:tom} \\
22 & $4{\times}[2,1],\,2{\times}[4,1]$ & $\ZZ/2$ & 3 & 3 & $1/2$ & $(1,1,1,1,2,2)$ & $(1,1,1,1,1,1)$ & $1,4,8,15,23$ & 7 & $2{\times}[2,1]$ & $+$\,\ref{prop:examples} \\
23 & $5{\times}[2,1],\,3{\times}[4,1]$ & $\ZZ/2$ & 3 & 3 & $1/4$ & $(1,1,1,1,0,2,2,0)$ & $(1,1,1,1,1,1,1,1)$ & $0,2,4,8,11$ & 3 & $4{\times}[2,1],\,2{\times}[4,1]$ &  \\
24 & $2{\times}[2,1],\,[4,1],\,[14,3]$ & $\ZZ/2$ & 3 & 3 & $1/28$ & $(1,0,0,7)$ & $(1,1,1,9)$ & $0,0,1,2,2$ & 0 & $2{\times}[2,1],\,2{\times}[4,1],\,[7,3]$ & $+$\,\ref{prop:tom} \\
25 & $2{\times}[2,1],\,[4,1],\,[8,3]$ & $\ZZ/2$ & 3 & 3 & $7/8$ & $(1,1,2,4)$ & $(1,1,1,5)$ & $2,6,13,24,39$ & 12 & $[2,1],\,[4,1]$ & $-$\,\ref{prop:genusbound} \\
26 & $2{\times}[2,1],\,[5,2],\,[14,5]$ & $\ZZ/2$ & 3 & 3 & $17/70$ & $(1,0,0,7)$ & $(1,1,3,9)$ & $1,2,4,7,11$ & 3 & $2{\times}[2,1],\,2{\times}[5,2],\,[7,2]$ & $-$\,\ref{prop:liu26} \\
27 & $3{\times}[2,1],\,[4,1],\,[8,3]$ & $\ZZ/2$ & 3 & 3 & $3/8$ & $(1,1,0,2,4)$ & $(1,1,1,1,5)$ & $1,3,6,11,17$ & 5 & $3{\times}[2,1],\,[4,1]$ & $+$\,\ref{prop:tom} \\
28 & $3{\times}[2,1],\,[5,2],\,[10,3]$ & $\ZZ/2$ & 3 & 3 & $2/5$ & $(1,1,1,0,5)$ & $(1,1,1,3,3)$ & $1,3,6,11,18$ & 5 & $3{\times}[5,2]$ &  \\
29 & $4{\times}[2,1],\,2{\times}[4,1],\,[5,2]$ & $\ZZ/2$ & 3 & 3 & $1/10$ & $(1,1,1,1,2,2,0)$ & $(1,1,1,1,1,1,3)$ & $0,1,2,4,5$ & 1 & $2{\times}[2,1],\,2{\times}[5,2]$ & $+$\,\ref{prop:multi} \\
30 & $5{\times}[2,1],\,[4,1],\,[10,3]$ & $\ZZ/2$ & 3 & 3 & $1/20$ & $(1,1,1,0,0,0,5)$ & $(1,1,1,1,1,1,3)$ & $0,1,1,2,2$ & 0 & $4{\times}[2,1],\,2{\times}[4,1],\,[5,2]$ &  \\
31 & $2{\times}[2,1],\,[4,1],\,[5,2],\,[8,3]$ & $\ZZ/2$ & 3 & 3 & $19/40$ & $(1,1,2,0,4)$ & $(1,1,1,3,5)$ & $1,3,7,13,21$ & 6 & $[2,1],\,[4,1],\,2{\times}[5,2]$ & $-$\,\ref{prop:liu26} \\
32 & $4{\times}[2,1],\,2{\times}[4,1],\,[8,3]$ & $\ZZ/2{\times}\ZZ/2$ & 3 & 3 & $1/8$ & $(1,1,1,1,2,2,0);\,(1,1,0,0,2,0,4)$ & $(1,1,1,1,1,1,5)$ & $0,1,2,4,5$ & 1 & $4{\times}[2,1],\,2{\times}[4,1]$ &  \\
33 & $4{\times}[3,1],\,[6,1]$ & $\ZZ/3$ & 1 & 3 & $1/2$ & $(2,2,2,2,2)$ & $(2,0,0,0,5)$ & $2,4,8,14,23$ & 7 & $[2,1]$ & $+$\,\ref{prop:examples} \\
34 & $[2,1],\,2{\times}[3,1],\,[12,5]$ & $\ZZ/3$ & 1 & 3 & $1/4$ & $(0,2,2,8)$ & $(1,0,0,9)$ & $1,2,4,8,12$ & 3 & $3{\times}[2,1],\,[4,1]$ & $+$\,\ref{prop:34} \\
35 & $4{\times}[3,1],\,[5,2],\,[6,1]$ & $\ZZ/3$ & 1 & 3 & $1/10$ & $(2,2,2,2,0,2)$ & $(2,0,0,0,3,5)$ & $1,1,2,3,5$ & 1 & $[2,1],\,3{\times}[5,2]$ &  \\
36 & $2{\times}[2,1],\,2{\times}[8,3]$ & $\ZZ/4$ & 3 & 3 & $1/4$ & $(1,1,6,2)$ & $(1,1,5,5)$ & $1,2,4,7,11$ & 3 & $2{\times}[2,1]$ & $+$\,\ref{prop:examples} \\
37 & $4{\times}[5,2]$ & $\ZZ/5$ & 3 & 3 & $2/5$ & $(4,3,2,1)$ & $(3,3,3,3)$ & $1,2,6,11,19$ & 5 & $\varnothing$ & $+$\,\ref{prop:examples} \\
38 & $8{\times}[2,1]$ & $\ZZ/2$ & 1 & 2 & $1$ & $(1,1,1,1,1,1,1,1)$ & $(0,0,0,0,1,1,1,1)$ & $2,6,12,23,38$ & 4 & $\varnothing$ & $+$\,\ref{prop:examples} \\
39 & $8{\times}[2,1]$ & $\ZZ/2$ & 1 & 2 & $2$ & $(1,1,1,1,1,1,1,1)$ & $(0,0,0,0,1,1,1,1)$ & $3,10,22,43,73$ & 8 & $\varnothing$ & $+$\,\ref{prop:examples} \\
40 & $8{\times}[2,1]$ & $\ZZ/2$ & 1 & 2 & $3$ & $(1,1,1,1,1,1,1,1)$ & $(0,0,0,0,1,1,1,1)$ & $4,14,32,63,108$ & 12 & $\varnothing$ & $-$\,\ref{prop:B24} \\
41 & $8{\times}[2,1]$ & $\ZZ/2$ & 1 & 2 & $4$ & $(1,1,1,1,1,1,1,1)$ & $(0,0,0,0,1,1,1,1)$ & $5,18,42,83,143$ & 16 & $\varnothing$ & $-$\,\ref{prop:B24} \\
42 & $2{\times}[2,1],\,2{\times}[6,1]$ & $\ZZ/2$ & 1 & 2 & $1/3$ & $(1,1,3,3)$ & $(0,1,1,4)$ & $1,2,5,9,14$ & 2 & $2{\times}[3,1]$ & $+$\,\ref{prop:examples} \\
43 & $2{\times}[2,1],\,2{\times}[6,1]$ & $\ZZ/2$ & 1 & 2 & $4/3$ & $(1,1,3,3)$ & $(0,1,1,4)$ & $2,6,15,29,49$ & 6 & $2{\times}[3,1]$ & $+$\,\ref{prop:multi} \\
44 & $2{\times}[2,1],\,2{\times}[6,1]$ & $\ZZ/2$ & 1 & 2 & $7/3$ & $(1,1,3,3)$ & $(0,1,1,4)$ & $3,10,25,49,84$ & 10 & $2{\times}[3,1]$ & $-$\,\ref{prop:largedeg} \\
45 & $2{\times}[2,1],\,2{\times}[6,1]$ & $\ZZ/2$ & 1 & 2 & $10/3$ & $(1,1,3,3)$ & $(0,1,1,4)$ & $4,14,35,69,119$ & 14 & $2{\times}[3,1]$ & $-$\,\ref{prop:largedeg} \\
46 & $3{\times}[2,1],\,[10,3]$ & $\ZZ/2$ & 1 & 2 & $1/5$ & $(1,1,1,5)$ & $(0,0,0,2)$ & $1,2,4,6,10$ & 1 & $[5,2]$ & $+$\,\ref{prop:examples} \\
47 & $3{\times}[2,1],\,[10,3]$ & $\ZZ/2$ & 1 & 2 & $6/5$ & $(1,1,1,5)$ & $(0,0,0,2)$ & $2,6,14,26,45$ & 5 & $[5,2]$ & $-$\,\ref{prop:47} \\
48 & $3{\times}[2,1],\,[10,3]$ & $\ZZ/2$ & 1 & 2 & $11/5$ & $(1,1,1,5)$ & $(0,0,0,2)$ & $3,10,24,46,80$ & 9 & $[5,2]$ & $-$\,\ref{prop:largedeg} \\
49 & $3{\times}[2,1],\,[10,3]$ & $\ZZ/2$ & 1 & 2 & $16/5$ & $(1,1,1,5)$ & $(0,0,0,2)$ & $4,14,34,66,115$ & 13 & $[5,2]$ & $-$\,\ref{prop:largedeg} \\
50 & $5{\times}[2,1],\,[6,1]$ & $\ZZ/2$ & 1 & 2 & $2/3$ & $(1,1,1,1,1,3)$ & $(0,0,0,0,0,4)$ & $2,4,9,16,26$ & 3 & $[3,1]$ & $+$\,\ref{prop:examples} \\
51 & $5{\times}[2,1],\,[6,1]$ & $\ZZ/2$ & 1 & 2 & $5/3$ & $(1,1,1,1,1,3)$ & $(0,0,0,0,0,4)$ & $3,8,19,36,61$ & 7 & $[3,1]$ & $+$\,\ref{prop:51} \\
52 & $5{\times}[2,1],\,[6,1]$ & $\ZZ/2$ & 1 & 2 & $8/3$ & $(1,1,1,1,1,3)$ & $(0,0,0,0,0,4)$ & $4,12,29,56,96$ & 11 & $[3,1]$ & $-$\,\ref{prop:largedeg} \\
53 & $5{\times}[2,1],\,[6,1]$ & $\ZZ/2$ & 1 & 2 & $11/3$ & $(1,1,1,1,1,3)$ & $(0,0,0,0,0,4)$ & $5,16,39,76,131$ & 15 & $[3,1]$ & $-$\,\ref{prop:largedeg} \\
54 & $8{\times}[2,1],\,[11,1]$ & $\ZZ/2$ & 1 & 2 & $4/11$ & $(1,1,1,1,1,1,1,1,0)$ & $(0,0,0,0,1,1,1,1,5)$ & $0,3,3,9,12$ & 1 & $2{\times}[11,1]$ & $-$\,\ref{prop:liu26} \\
55 & $8{\times}[2,1],\,[11,2]$ & $\ZZ/2$ & 1 & 2 & $1/11$ & $(1,1,1,1,1,1,1,1,0)$ & $(0,0,0,0,1,1,1,1,5)$ & $0,1,1,3,3$ & -1 & $2{\times}[11,2]$ &  \\
56 & $8{\times}[2,1],\,[11,4]$ & $\ZZ/2$ & 1 & 2 & $3/11$ & $(1,1,1,1,1,1,1,1,0)$ & $(0,0,0,0,1,1,1,1,5)$ & $0,2,3,6,10$ & 0 & $2{\times}[11,4]$ &  \\
57 & $8{\times}[2,1],\,2{\times}[3,1]$ & $\ZZ/2$ & 1 & 2 & $2/3$ & $(1,1,1,1,1,1,1,1,0,0)$ & $(0,0,0,0,1,1,1,1,1,1)$ & $1,4,8,15,25$ & 2 & $4{\times}[3,1]$ &  \\
58 & $8{\times}[2,1],\,2{\times}[3,1]$ & $\ZZ/2$ & 1 & 2 & $5/3$ & $(1,1,1,1,1,1,1,1,0,0)$ & $(0,0,0,0,1,1,1,1,1,1)$ & $2,8,18,35,60$ & 6 & $4{\times}[3,1]$ & $-$\,\ref{prop:N3} \\
59 & $8{\times}[2,1],\,2{\times}[5,1]$ & $\ZZ/2$ & 1 & 2 & $1/5$ & $(1,1,1,1,1,1,1,1,0,0)$ & $(0,0,0,0,1,1,1,1,2,2)$ & $0,2,2,5,8$ & 0 & $4{\times}[5,1]$ &  \\
60 & $8{\times}[2,1],\,2{\times}[5,2]$ & $\ZZ/2$ & 1 & 2 & $4/5$ & $(1,1,1,1,1,1,1,1,0,0)$ & $(0,0,0,0,1,1,1,1,2,2)$ & $1,4,8,17,29$ & 2 & $4{\times}[5,2]$ & $-$\,\ref{prop:liu26} \\
61 & $8{\times}[2,1],\,[3,1]$ & $\ZZ/2$ & 1 & 2 & $1/3$ & $(1,1,1,1,1,1,1,1,0)$ & $(0,0,0,0,1,1,1,1,1)$ & $1,3,5,9,14$ & 1 & $2{\times}[3,1]$ & $+$\,\ref{prop:examples} \\
62 & $8{\times}[2,1],\,[3,1]$ & $\ZZ/2$ & 1 & 2 & $4/3$ & $(1,1,1,1,1,1,1,1,0)$ & $(0,0,0,0,1,1,1,1,1)$ & $2,7,15,29,49$ & 5 & $2{\times}[3,1]$ & $-$\,\ref{prop:N3} \\
63 & $8{\times}[2,1],\,[3,1]$ & $\ZZ/2$ & 1 & 2 & $7/3$ & $(1,1,1,1,1,1,1,1,0)$ & $(0,0,0,0,1,1,1,1,1)$ & $3,11,25,49,84$ & 9 & $2{\times}[3,1]$ & $-$\,\ref{prop:largedeg} \\
64 & $8{\times}[2,1],\,[3,1]$ & $\ZZ/2$ & 1 & 2 & $10/3$ & $(1,1,1,1,1,1,1,1,0)$ & $(0,0,0,0,1,1,1,1,1)$ & $4,15,35,69,119$ & 13 & $2{\times}[3,1]$ & $-$\,\ref{prop:largedeg} \\
65 & $8{\times}[2,1],\,3{\times}[3,1]$ & $\ZZ/2$ & 1 & 2 & $1$ & $(1,1,1,1,1,1,1,1,0,0,0)$ & $(0,0,0,0,1,1,1,1,1,1,1)$ & $1,5,11,21,36$ & 3 & $6{\times}[3,1]$ &  \\
66 & $8{\times}[2,1],\,4{\times}[3,1]$ & $\ZZ/2$ & 1 & 2 & $1/3$ & $(1,1,1,1,1,1,1,1,0,0,0,0)$ & $(0,0,0,0,1,1,1,1,1,1,1,1)$ & $0,2,4,7,12$ & 0 & $8{\times}[3,1]$ &  \\
67 & $8{\times}[2,1],\,[5,1]$ & $\ZZ/2$ & 1 & 2 & $3/5$ & $(1,1,1,1,1,1,1,1,0)$ & $(0,0,0,0,1,1,1,1,2)$ & $1,4,7,14,23$ & 2 & $2{\times}[5,1]$ &  \\
68 & $8{\times}[2,1],\,[5,1]$ & $\ZZ/2$ & 1 & 2 & $8/5$ & $(1,1,1,1,1,1,1,1,0)$ & $(0,0,0,0,1,1,1,1,2)$ & $2,8,17,34,58$ & 6 & $2{\times}[5,1]$ & $-$\,\ref{prop:N3} \\
69 & $8{\times}[2,1],\,[5,1]$ & $\ZZ/2$ & 1 & 2 & $13/5$ & $(1,1,1,1,1,1,1,1,0)$ & $(0,0,0,0,1,1,1,1,2)$ & $3,12,27,54,93$ & 10 & $2{\times}[5,1]$ & $-$\,\ref{prop:largedeg} \\
70 & $8{\times}[2,1],\,[5,2]$ & $\ZZ/2$ & 1 & 2 & $2/5$ & $(1,1,1,1,1,1,1,1,0)$ & $(0,0,0,0,1,1,1,1,2)$ & $1,3,5,10,16$ & 1 & $2{\times}[5,2]$ &  \\
71 & $8{\times}[2,1],\,[5,2]$ & $\ZZ/2$ & 1 & 2 & $7/5$ & $(1,1,1,1,1,1,1,1,0)$ & $(0,0,0,0,1,1,1,1,2)$ & $2,7,15,30,51$ & 5 & $2{\times}[5,2]$ & $-$\,\ref{prop:N3} \\
72 & $8{\times}[2,1],\,[5,2]$ & $\ZZ/2$ & 1 & 2 & $12/5$ & $(1,1,1,1,1,1,1,1,0)$ & $(0,0,0,0,1,1,1,1,2)$ & $3,11,25,50,86$ & 9 & $2{\times}[5,2]$ & $-$\,\ref{prop:largedeg} \\
73 & $8{\times}[2,1],\,[7,1]$ & $\ZZ/2$ & 1 & 2 & $6/7$ & $(1,1,1,1,1,1,1,1,0)$ & $(0,0,0,0,1,1,1,1,3)$ & $1,5,9,19,31$ & 3 & $2{\times}[7,1]$ &  \\
74 & $8{\times}[2,1],\,[7,1]$ & $\ZZ/2$ & 1 & 2 & $13/7$ & $(1,1,1,1,1,1,1,1,0)$ & $(0,0,0,0,1,1,1,1,3)$ & $2,9,19,39,66$ & 7 & $2{\times}[7,1]$ & $-$\,\ref{prop:N3} \\
75 & $8{\times}[2,1],\,[7,2]$ & $\ZZ/2$ & 1 & 2 & $5/7$ & $(1,1,1,1,1,1,1,1,0)$ & $(0,0,0,0,1,1,1,1,3)$ & $1,4,8,16,26$ & 2 & $2{\times}[7,2]$ &  \\
76 & $8{\times}[2,1],\,[7,2]$ & $\ZZ/2$ & 1 & 2 & $12/7$ & $(1,1,1,1,1,1,1,1,0)$ & $(0,0,0,0,1,1,1,1,3)$ & $2,8,18,36,61$ & 6 & $2{\times}[7,2]$ & $-$\,\ref{prop:N3} \\
77 & $8{\times}[2,1],\,[7,3]$ & $\ZZ/2$ & 1 & 2 & $3/7$ & $(1,1,1,1,1,1,1,1,0)$ & $(0,0,0,0,1,1,1,1,3)$ & $1,3,5,10,16$ & 1 & $2{\times}[7,3]$ &  \\
78 & $8{\times}[2,1],\,[7,3]$ & $\ZZ/2$ & 1 & 2 & $10/7$ & $(1,1,1,1,1,1,1,1,0)$ & $(0,0,0,0,1,1,1,1,3)$ & $2,7,15,30,51$ & 5 & $2{\times}[7,3]$ & $-$\,\ref{prop:N3} \\
79 & $8{\times}[2,1],\,[9,1]$ & $\ZZ/2$ & 1 & 2 & $1/9$ & $(1,1,1,1,1,1,1,1,0)$ & $(0,0,0,0,1,1,1,1,4)$ & $0,2,1,4,4$ & 0 & $2{\times}[9,1]$ &  \\
80 & $8{\times}[2,1],\,[9,1]$ & $\ZZ/2$ & 1 & 2 & $10/9$ & $(1,1,1,1,1,1,1,1,0)$ & $(0,0,0,0,1,1,1,1,4)$ & $1,6,11,24,39$ & 4 & $2{\times}[9,1]$ & $-$\,\ref{prop:liu26} \\
81 & $8{\times}[2,1],\,[9,2]$ & $\ZZ/2$ & 1 & 2 & $7/9$ & $(1,1,1,1,1,1,1,1,0)$ & $(0,0,0,0,1,1,1,1,4)$ & $1,4,8,17,28$ & 2 & $2{\times}[9,2]$ &  \\
82 & $8{\times}[2,1],\,[9,4]$ & $\ZZ/2$ & 1 & 2 & $4/9$ & $(1,1,1,1,1,1,1,1,0)$ & $(0,0,0,0,1,1,1,1,4)$ & $1,3,5,10,16$ & 1 & $2{\times}[9,4]$ &  \\
83 & $2{\times}[2,1],\,2{\times}[6,1],\,[7,1]$ & $\ZZ/2$ & 1 & 2 & $4/21$ & $(1,1,3,3,0)$ & $(0,1,1,4,3)$ & $0,1,2,5,7$ & 1 & $2{\times}[3,1],\,2{\times}[7,1]$ &  \\
84 & $2{\times}[2,1],\,2{\times}[6,1],\,[7,2]$ & $\ZZ/2$ & 1 & 2 & $1/21$ & $(1,1,3,3,0)$ & $(0,1,1,4,3)$ & $0,0,1,2,2$ & 0 & $2{\times}[3,1],\,2{\times}[7,2]$ &  \\
85 & $2{\times}[2,1],\,2{\times}[6,1],\,[7,3]$ & $\ZZ/2$ & 1 & 2 & $16/21$ & $(1,1,3,3,0)$ & $(0,1,1,4,3)$ & $1,3,8,16,27$ & 3 & $2{\times}[3,1],\,2{\times}[7,3]$ & $-$\,\ref{prop:netpencil} \\
86 & $2{\times}[2,1],\,2{\times}[6,1],\,[9,2]$ & $\ZZ/2$ & 1 & 2 & $1/9$ & $(1,1,3,3,0)$ & $(0,1,1,4,4)$ & $0,0,1,3,4$ & 0 & $2{\times}[3,1],\,2{\times}[9,2]$ &  \\
87 & $2{\times}[2,1],\,[3,1],\,2{\times}[6,1]$ & $\ZZ/2$ & 1 & 2 & $2/3$ & $(1,1,0,3,3)$ & $(0,1,1,1,4)$ & $1,3,8,15,25$ & 3 & $4{\times}[3,1]$ & $-$\,\ref{prop:netpencil} \\
88 & $2{\times}[2,1],\,[3,1],\,2{\times}[6,1]$ & $\ZZ/2$ & 1 & 2 & $5/3$ & $(1,1,0,3,3)$ & $(0,1,1,1,4)$ & $2,7,18,35,60$ & 7 & $4{\times}[3,1]$ & $-$\,\ref{prop:8893} \\
89 & $2{\times}[2,1],\,[5,1],\,2{\times}[6,1]$ & $\ZZ/2$ & 1 & 2 & $14/15$ & $(1,1,0,3,3)$ & $(0,1,2,1,4)$ & $1,4,10,20,34$ & 4 & $2{\times}[3,1],\,2{\times}[5,1]$ & $-$\,\ref{prop:netpencil} \\
90 & $2{\times}[2,1],\,[5,2],\,2{\times}[6,1]$ & $\ZZ/2$ & 1 & 2 & $11/15$ & $(1,1,0,3,3)$ & $(0,1,2,1,4)$ & $1,3,8,16,27$ & 3 & $2{\times}[3,1],\,2{\times}[5,2]$ & $-$\,\ref{prop:netpencil} \\
91 & $3{\times}[2,1],\,2{\times}[3,1],\,[10,3]$ & $\ZZ/2$ & 1 & 2 & $13/15$ & $(1,1,1,0,0,5)$ & $(0,0,0,1,1,2)$ & $1,4,10,18,32$ & 3 & $4{\times}[3,1],\,[5,2]$ &  \\
92 & $3{\times}[2,1],\,[3,1],\,[10,3]$ & $\ZZ/2$ & 1 & 2 & $8/15$ & $(1,1,1,0,5)$ & $(0,0,0,1,2)$ & $1,3,7,12,21$ & 2 & $2{\times}[3,1],\,[5,2]$ &  \\
93 & $3{\times}[2,1],\,[3,1],\,[10,3]$ & $\ZZ/2$ & 1 & 2 & $23/15$ & $(1,1,1,0,5)$ & $(0,0,0,1,2)$ & $2,7,17,32,56$ & 6 & $2{\times}[3,1],\,[5,2]$ & $-$\,\ref{prop:8893} \\
94 & $3{\times}[2,1],\,[3,1],\,[10,3]$ & $\ZZ/2$ & 1 & 2 & $38/15$ & $(1,1,1,0,5)$ & $(0,0,0,1,2)$ & $3,11,27,52,91$ & 10 & $2{\times}[3,1],\,[5,2]$ & $-$\,\ref{prop:largedeg} \\
95 & $3{\times}[2,1],\,3{\times}[3,1],\,[10,3]$ & $\ZZ/2$ & 1 & 2 & $1/5$ & $(1,1,1,0,0,0,5)$ & $(0,0,0,1,1,1,2)$ & $0,1,3,4,8$ & 0 & $6{\times}[3,1],\,[5,2]$ &  \\
96 & $3{\times}[2,1],\,[5,1],\,[10,3]$ & $\ZZ/2$ & 1 & 2 & $4/5$ & $(1,1,1,0,5)$ & $(0,0,0,2,2)$ & $1,4,9,17,30$ & 3 & $2{\times}[5,1],\,[5,2]$ &  \\
97 & $3{\times}[2,1],\,[5,2],\,[10,3]$ & $\ZZ/2$ & 1 & 2 & $3/5$ & $(1,1,1,0,5)$ & $(0,0,0,2,2)$ & $1,3,7,13,23$ & 2 & $3{\times}[5,2]$ & $-$\,\ref{prop:97} \\
98 & $3{\times}[2,1],\,[5,2],\,[10,3]$ & $\ZZ/2$ & 1 & 2 & $8/5$ & $(1,1,1,0,5)$ & $(0,0,0,2,2)$ & $2,7,17,33,58$ & 6 & $3{\times}[5,2]$ & $-$\,\ref{prop:98} \\
99 & $3{\times}[2,1],\,[7,1],\,[10,3]$ & $\ZZ/2$ & 1 & 2 & $2/35$ & $(1,1,1,0,5)$ & $(0,0,0,3,2)$ & $0,1,1,2,3$ & 0 & $[5,2],\,2{\times}[7,1]$ &  \\
100 & $3{\times}[2,1],\,[7,2],\,[10,3]$ & $\ZZ/2$ & 1 & 2 & $32/35$ & $(1,1,1,0,5)$ & $(0,0,0,3,2)$ & $1,4,10,19,33$ & 3 & $[5,2],\,2{\times}[7,2]$ & $-$\,\ref{prop:liu26} \\
101 & $3{\times}[2,1],\,[7,3],\,[10,3]$ & $\ZZ/2$ & 1 & 2 & $22/35$ & $(1,1,1,0,5)$ & $(0,0,0,3,2)$ & $1,3,7,13,23$ & 2 & $[5,2],\,2{\times}[7,3]$ &  \\
102 & $5{\times}[2,1],\,2{\times}[3,1],\,[6,1]$ & $\ZZ/2$ & 1 & 2 & $1/3$ & $(1,1,1,1,1,0,0,3)$ & $(0,0,0,0,0,1,1,4)$ & $1,2,5,8,13$ & 1 & $5{\times}[3,1]$ &  \\
103 & $5{\times}[2,1],\,2{\times}[3,1],\,[6,1]$ & $\ZZ/2$ & 1 & 2 & $4/3$ & $(1,1,1,1,1,0,0,3)$ & $(0,0,0,0,0,1,1,4)$ & $2,6,15,28,48$ & 5 & $5{\times}[3,1]$ & $-$\,\ref{prop:104} \\
104 & $5{\times}[2,1],\,[3,1],\,[6,1]$ & $\ZZ/2$ & 1 & 2 & $1$ & $(1,1,1,1,1,0,3)$ & $(0,0,0,0,0,1,4)$ & $2,5,12,22,37$ & 4 & $3{\times}[3,1]$ & $-$\,\ref{prop:104} \\
105 & $5{\times}[2,1],\,[3,1],\,[6,1]$ & $\ZZ/2$ & 1 & 2 & $2$ & $(1,1,1,1,1,0,3)$ & $(0,0,0,0,0,1,4)$ & $3,9,22,42,72$ & 8 & $3{\times}[3,1]$ & $-$\,\ref{prop:N4} \\
106 & $5{\times}[2,1],\,3{\times}[3,1],\,[6,1]$ & $\ZZ/2$ & 1 & 2 & $2/3$ & $(1,1,1,1,1,0,0,0,3)$ & $(0,0,0,0,0,1,1,1,4)$ & $1,3,8,14,24$ & 2 & $7{\times}[3,1]$ &  \\
107 & $5{\times}[2,1],\,[5,1],\,[6,1]$ & $\ZZ/2$ & 1 & 2 & $4/15$ & $(1,1,1,1,1,0,3)$ & $(0,0,0,0,0,2,4)$ & $1,2,4,7,11$ & 1 & $[3,1],\,2{\times}[5,1]$ &  \\
108 & $5{\times}[2,1],\,[5,1],\,[6,1]$ & $\ZZ/2$ & 1 & 2 & $19/15$ & $(1,1,1,1,1,0,3)$ & $(0,0,0,0,0,2,4)$ & $2,6,14,27,46$ & 5 & $[3,1],\,2{\times}[5,1]$ & $-$\,\ref{prop:108} \\
109 & $5{\times}[2,1],\,[5,2],\,[6,1]$ & $\ZZ/2$ & 1 & 2 & $1/15$ & $(1,1,1,1,1,0,3)$ & $(0,0,0,0,0,2,4)$ & $1,1,2,3,4$ & 0 & $[3,1],\,2{\times}[5,2]$ & $+$\,\ref{prop:examples} \\
110 & $5{\times}[2,1],\,[5,2],\,[6,1]$ & $\ZZ/2$ & 1 & 2 & $16/15$ & $(1,1,1,1,1,0,3)$ & $(0,0,0,0,0,2,4)$ & $2,5,12,23,39$ & 4 & $[3,1],\,2{\times}[5,2]$ & $-$\,\ref{prop:108} \\
111 & $5{\times}[2,1],\,[5,2],\,[6,1]$ & $\ZZ/2$ & 1 & 2 & $31/15$ & $(1,1,1,1,1,0,3)$ & $(0,0,0,0,0,2,4)$ & $3,9,22,43,74$ & 8 & $[3,1],\,2{\times}[5,2]$ & $-$\,\ref{prop:N4} \\
112 & $5{\times}[2,1],\,[6,1],\,[7,1]$ & $\ZZ/2$ & 1 & 2 & $11/21$ & $(1,1,1,1,1,3,0)$ & $(0,0,0,0,0,4,3)$ & $1,3,6,12,19$ & 2 & $[3,1],\,2{\times}[7,1]$ &  \\
113 & $5{\times}[2,1],\,[6,1],\,[7,2]$ & $\ZZ/2$ & 1 & 2 & $8/21$ & $(1,1,1,1,1,3,0)$ & $(0,0,0,0,0,4,3)$ & $1,2,5,9,14$ & 1 & $[3,1],\,2{\times}[7,2]$ &  \\
114 & $5{\times}[2,1],\,[6,1],\,[7,2]$ & $\ZZ/2$ & 1 & 2 & $29/21$ & $(1,1,1,1,1,3,0)$ & $(0,0,0,0,0,4,3)$ & $2,6,15,29,49$ & 5 & $[3,1],\,2{\times}[7,2]$ & $-$\,\ref{prop:liu26} \\
115 & $5{\times}[2,1],\,[6,1],\,[7,3]$ & $\ZZ/2$ & 1 & 2 & $2/21$ & $(1,1,1,1,1,3,0)$ & $(0,0,0,0,0,4,3)$ & $1,1,2,3,4$ & 0 & $[3,1],\,2{\times}[7,3]$ &  \\
116 & $5{\times}[2,1],\,[6,1],\,[7,3]$ & $\ZZ/2$ & 1 & 2 & $23/21$ & $(1,1,1,1,1,3,0)$ & $(0,0,0,0,0,4,3)$ & $2,5,12,23,39$ & 4 & $[3,1],\,2{\times}[7,3]$ & $-$\,\ref{prop:108} \\
117 & $5{\times}[2,1],\,[6,1],\,[9,2]$ & $\ZZ/2$ & 1 & 2 & $4/9$ & $(1,1,1,1,1,3,0)$ & $(0,0,0,0,0,4,4)$ & $1,2,5,10,16$ & 1 & $[3,1],\,2{\times}[9,2]$ &  \\
118 & $5{\times}[2,1],\,[6,1],\,[9,4]$ & $\ZZ/2$ & 1 & 2 & $1/9$ & $(1,1,1,1,1,3,0)$ & $(0,0,0,0,0,4,4)$ & $1,1,2,3,4$ & 0 & $[3,1],\,2{\times}[9,4]$ &  \\
119 & $8{\times}[2,1],\,2{\times}[3,1],\,[5,1]$ & $\ZZ/2$ & 1 & 2 & $4/15$ & $(1,1,1,1,1,1,1,1,0,0,0)$ & $(0,0,0,0,1,1,1,1,1,1,2)$ & $0,2,3,6,10$ & 0 & $4{\times}[3,1],\,2{\times}[5,1]$ &  \\
120 & $8{\times}[2,1],\,2{\times}[3,1],\,[5,2]$ & $\ZZ/2$ & 1 & 2 & $1/15$ & $(1,1,1,1,1,1,1,1,0,0,0)$ & $(0,0,0,0,1,1,1,1,1,1,2)$ & $0,1,1,2,3$ & -1 & $4{\times}[3,1],\,2{\times}[5,2]$ &  \\
121 & $8{\times}[2,1],\,[3,1],\,[5,1]$ & $\ZZ/2$ & 1 & 2 & $14/15$ & $(1,1,1,1,1,1,1,1,0,0)$ & $(0,0,0,0,1,1,1,1,1,2)$ & $1,5,10,20,34$ & 3 & $2{\times}[3,1],\,2{\times}[5,1]$ &  \\
122 & $8{\times}[2,1],\,[3,1],\,[5,2]$ & $\ZZ/2$ & 1 & 2 & $11/15$ & $(1,1,1,1,1,1,1,1,0,0)$ & $(0,0,0,0,1,1,1,1,1,2)$ & $1,4,8,16,27$ & 2 & $2{\times}[3,1],\,2{\times}[5,2]$ &  \\
123 & $8{\times}[2,1],\,[3,1],\,[7,1]$ & $\ZZ/2$ & 1 & 2 & $4/21$ & $(1,1,1,1,1,1,1,1,0,0)$ & $(0,0,0,0,1,1,1,1,1,3)$ & $0,2,2,5,7$ & 0 & $2{\times}[3,1],\,2{\times}[7,1]$ &  \\
124 & $8{\times}[2,1],\,[3,1],\,[7,2]$ & $\ZZ/2$ & 1 & 2 & $1/21$ & $(1,1,1,1,1,1,1,1,0,0)$ & $(0,0,0,0,1,1,1,1,1,3)$ & $0,1,1,2,2$ & -1 & $2{\times}[3,1],\,2{\times}[7,2]$ &  \\
125 & $8{\times}[2,1],\,[3,1],\,[7,3]$ & $\ZZ/2$ & 1 & 2 & $16/21$ & $(1,1,1,1,1,1,1,1,0,0)$ & $(0,0,0,0,1,1,1,1,1,3)$ & $1,4,8,16,27$ & 2 & $2{\times}[3,1],\,2{\times}[7,3]$ &  \\
126 & $8{\times}[2,1],\,[3,1],\,[9,2]$ & $\ZZ/2$ & 1 & 2 & $1/9$ & $(1,1,1,1,1,1,1,1,0,0)$ & $(0,0,0,0,1,1,1,1,1,4)$ & $0,1,1,3,4$ & -1 & $2{\times}[3,1],\,2{\times}[9,2]$ &  \\
127 & $8{\times}[2,1],\,[5,1],\,[7,3]$ & $\ZZ/2$ & 1 & 2 & $1/35$ & $(1,1,1,1,1,1,1,1,0,0)$ & $(0,0,0,0,1,1,1,1,2,3)$ & $0,1,0,1,1$ & -1 & $2{\times}[5,1],\,2{\times}[7,3]$ &  \\
128 & $8{\times}[2,1],\,[5,2],\,[7,2]$ & $\ZZ/2$ & 1 & 2 & $4/35$ & $(1,1,1,1,1,1,1,1,0,0)$ & $(0,0,0,0,1,1,1,1,2,3)$ & $0,1,1,3,4$ & -1 & $2{\times}[5,2],\,2{\times}[7,2]$ & $-$\,\ref{prop:liu26} \\
129 & $2{\times}[2,1],\,[3,1],\,[5,1],\,2{\times}[6,1]$ & $\ZZ/2$ & 1 & 2 & $4/15$ & $(1,1,0,0,3,3)$ & $(0,1,1,2,1,4)$ & $0,1,3,6,10$ & 1 & $4{\times}[3,1],\,2{\times}[5,1]$ &  \\
130 & $2{\times}[2,1],\,[3,1],\,[5,2],\,2{\times}[6,1]$ & $\ZZ/2$ & 1 & 2 & $1/15$ & $(1,1,0,0,3,3)$ & $(0,1,1,2,1,4)$ & $0,0,1,2,3$ & 0 & $4{\times}[3,1],\,2{\times}[5,2]$ &  \\
131 & $3{\times}[2,1],\,[3,1],\,[5,1],\,[10,3]$ & $\ZZ/2$ & 1 & 2 & $2/15$ & $(1,1,1,0,0,5)$ & $(0,0,0,1,2,2)$ & $0,1,2,3,6$ & 0 & $2{\times}[3,1],\,2{\times}[5,1],\,[5,2]$ &  \\
132 & $5{\times}[2,1],\,[3,1],\,[5,1],\,[6,1]$ & $\ZZ/2$ & 1 & 2 & $3/5$ & $(1,1,1,1,1,0,0,3)$ & $(0,0,0,0,0,1,2,4)$ & $1,3,7,13,22$ & 2 & $3{\times}[3,1],\,2{\times}[5,1]$ &  \\
133 & $5{\times}[2,1],\,[3,1],\,[5,2],\,[6,1]$ & $\ZZ/2$ & 1 & 2 & $2/5$ & $(1,1,1,1,1,0,0,3)$ & $(0,0,0,0,0,1,2,4)$ & $1,2,5,9,15$ & 1 & $3{\times}[3,1],\,2{\times}[5,2]$ &  \\
134 & $2{\times}[9,1]$ & $\ZZ/3$ & 2 & 2 & $2/9$ & $(6,3)$ & $(4,4)$ & $0,3,2,7,8$ & 2 & $2{\times}[3,1]$ & $+$\,\ref{prop:examples} \\
135 & $2{\times}[9,1]$ & $\ZZ/3$ & 2 & 2 & $11/9$ & $(6,3)$ & $(4,4)$ & $1,7,12,27,43$ & 6 & $2{\times}[3,1]$ & $-$\,\ref{prop:N4} \\
136 & $2{\times}[9,1]$ & $\ZZ/3$ & 2 & 2 & $20/9$ & $(6,3)$ & $(4,4)$ & $2,11,22,47,78$ & 10 & $2{\times}[3,1]$ & $-$\,\ref{prop:largedeg} \\
137 & $2{\times}[9,2]$ & $\ZZ/3$ & 2 & 2 & $5/9$ & $(6,3)$ & $(4,4)$ & $1,3,6,13,21$ & 2 & $2{\times}[3,1]$ & $+$\,\ref{prop:multi} \\
138 & $2{\times}[9,2]$ & $\ZZ/3$ & 2 & 2 & $14/9$ & $(6,3)$ & $(4,4)$ & $2,7,16,33,56$ & 6 & $2{\times}[3,1]$ & $-$\,\ref{prop:largedeg} \\
139 & $2{\times}[9,4]$ & $\ZZ/3$ & 2 & 2 & $8/9$ & $(6,3)$ & $(4,4)$ & $2,5,10,19,32$ & 4 & $2{\times}[3,1]$ & $-$\,\ref{prop:139149} \\
140 & $2{\times}[9,4]$ & $\ZZ/3$ & 2 & 2 & $17/9$ & $(6,3)$ & $(4,4)$ & $3,9,20,39,67$ & 8 & $2{\times}[3,1]$ & $-$\,\ref{prop:largedeg} \\
141 & $6{\times}[3,1]$ & $\ZZ/3$ & 2 & 2 & $1$ & $(2,2,2,1,1,1)$ & $(1,1,1,1,1,1)$ & $1,5,12,21,37$ & 4 & $\varnothing$ & $+$\,\ref{prop:examples} \\
142 & $6{\times}[3,1]$ & $\ZZ/3$ & 2 & 2 & $2$ & $(2,2,2,1,1,1)$ & $(1,1,1,1,1,1)$ & $2,9,22,41,72$ & 8 & $\varnothing$ & $-$\,\ref{prop:B24} \\
143 & $7{\times}[3,1]$ & $\ZZ/3$ & 2 & 2 & $1/3$ & $(2,2,2,1,1,1,0)$ & $(1,1,1,1,1,1,1)$ & $0,2,5,7,13$ & 1 & $3{\times}[3,1]$ & $+$\,\ref{prop:multi} \\
144 & $7{\times}[3,1]$ & $\ZZ/3$ & 2 & 2 & $4/3$ & $(2,2,2,1,1,1,0)$ & $(1,1,1,1,1,1,1)$ & $1,6,15,27,48$ & 5 & $3{\times}[3,1]$ & $-$\,\ref{prop:N4} \\
145 & $2{\times}[3,1],\,2{\times}[9,2]$ & $\ZZ/3$ & 2 & 2 & $2/9$ & $(0,0,6,3)$ & $(1,1,4,4)$ & $0,1,2,5,8$ & 0 & $8{\times}[3,1]$ &  \\
146 & $[3,1],\,2{\times}[9,1]$ & $\ZZ/3$ & 2 & 2 & $5/9$ & $(0,6,3)$ & $(1,4,4)$ & $0,4,5,13,19$ & 3 & $5{\times}[3,1]$ & $-$\,\ref{prop:netpencil} \\
147 & $[3,1],\,2{\times}[9,2]$ & $\ZZ/3$ & 2 & 2 & $8/9$ & $(0,6,3)$ & $(1,4,4)$ & $1,4,9,19,32$ & 3 & $5{\times}[3,1]$ & $-$\,\ref{prop:netpencil} \\
148 & $[3,1],\,2{\times}[9,4]$ & $\ZZ/3$ & 2 & 2 & $2/9$ & $(0,6,3)$ & $(1,4,4)$ & $1,2,3,5,8$ & 1 & $5{\times}[3,1]$ & $+$\,\ref{prop:multi} \\
149 & $[3,1],\,2{\times}[9,4]$ & $\ZZ/3$ & 2 & 2 & $11/9$ & $(0,6,3)$ & $(1,4,4)$ & $2,6,13,25,43$ & 5 & $5{\times}[3,1]$ & $-$\,\ref{prop:139149} \\
150 & $[5,1],\,2{\times}[9,2]$ & $\ZZ/3$ & 2 & 2 & $7/45$ & $(0,6,3)$ & $(2,4,4)$ & $0,1,1,4,6$ & 0 & $2{\times}[3,1],\,3{\times}[5,1]$ &  \\
151 & $[5,1],\,2{\times}[9,4]$ & $\ZZ/3$ & 2 & 2 & $22/45$ & $(0,6,3)$ & $(2,4,4)$ & $1,3,5,10,17$ & 2 & $2{\times}[3,1],\,3{\times}[5,1]$ & $-$\,\ref{prop:liu26} \\
152 & $[5,2],\,2{\times}[9,4]$ & $\ZZ/3$ & 2 & 2 & $13/45$ & $(0,6,3)$ & $(2,4,4)$ & $1,2,3,6,10$ & 1 & $2{\times}[3,1],\,3{\times}[5,2]$ &  \\
153 & $6{\times}[3,1],\,[5,1]$ & $\ZZ/3$ & 2 & 2 & $3/5$ & $(2,2,2,1,1,1,0)$ & $(1,1,1,1,1,1,2)$ & $0,3,7,12,22$ & 2 & $3{\times}[5,1]$ & $-$\,\ref{prop:netpencil} \\
154 & $6{\times}[3,1],\,[5,2]$ & $\ZZ/3$ & 2 & 2 & $2/5$ & $(2,2,2,1,1,1,0)$ & $(1,1,1,1,1,1,2)$ & $0,2,5,8,15$ & 1 & $3{\times}[5,2]$ &  \\
155 & $8{\times}[3,1]$ & $\ZZ/3{\times}\ZZ/3$ & 2 & 2 & $2/3$ & $(2,2,2,1,1,1,0,0);\,(2,1,0,0,1,2,2,1)$ & $(1,1,1,1,1,1,1,1)$ & $0,3,8,13,24$ & 2 & $\varnothing$ &  \\
156 & $2{\times}[2,1],\,4{\times}[4,1]$ & $\ZZ/4$ & 1 & 2 & $1/2$ & $(1,1,3,3,1,1)$ & $(0,1,0,2,1,3)$ & $1,3,6,12,19$ & 2 & $\varnothing$ & $+$\,\ref{prop:examples} \\
157 & $2{\times}[2,1],\,4{\times}[4,1]$ & $\ZZ/4$ & 1 & 2 & $1$ & $(1,1,3,3,1,1)$ & $(0,1,0,0,3,3)$ & $2,5,11,22,37$ & 4 & $\varnothing$ & $+$\,\ref{prop:examples} \\
158 & $2{\times}[2,1],\,4{\times}[4,1]$ & $\ZZ/4$ & 1 & 2 & $3/2$ & $(1,1,3,3,1,1)$ & $(0,1,0,2,1,3)$ & $2,7,16,32,54$ & 6 & $\varnothing$ & $-$\,\ref{prop:B24} \\
159 & $2{\times}[2,1],\,4{\times}[4,1]$ & $\ZZ/4$ & 1 & 2 & $2$ & $(1,1,3,3,1,1)$ & $(0,1,0,0,3,3)$ & $3,9,21,42,72$ & 8 & $\varnothing$ & $-$\,\ref{prop:B24} \\
160 & $2{\times}[2,1],\,2{\times}[3,1],\,4{\times}[4,1]$ & $\ZZ/4$ & 1 & 2 & $1/6$ & $(1,1,0,0,3,3,1,1)$ & $(0,1,1,1,0,2,1,3)$ & $0,1,2,4,6$ & 0 & $8{\times}[3,1]$ &  \\
161 & $2{\times}[2,1],\,[3,1],\,4{\times}[4,1]$ & $\ZZ/4$ & 1 & 2 & $1/3$ & $(1,1,0,3,3,1,1)$ & $(0,1,1,0,0,3,3)$ & $1,2,4,8,13$ & 1 & $4{\times}[3,1]$ &  \\
162 & $2{\times}[2,1],\,[3,1],\,4{\times}[4,1]$ & $\ZZ/4$ & 1 & 2 & $5/6$ & $(1,1,0,3,3,1,1)$ & $(0,1,1,0,2,1,3)$ & $1,4,9,18,30$ & 3 & $4{\times}[3,1]$ & $-$\,\ref{prop:162} \\
163 & $2{\times}[2,1],\,[4,1],\,[12,5]$ & $\ZZ/4$ & 1 & 2 & $1/3$ & $(1,1,3,3)$ & $(0,0,0,4)$ & $1,2,5,8,13$ & 1 & $[3,1]$ & $+$\,\ref{prop:examples} \\
164 & $2{\times}[2,1],\,[4,1],\,[12,5]$ & $\ZZ/4$ & 1 & 2 & $4/3$ & $(1,1,3,3)$ & $(0,0,0,4)$ & $2,6,15,28,48$ & 5 & $[3,1]$ & $-$\,\ref{prop:largedeg} \\
165 & $2{\times}[2,1],\,4{\times}[4,1],\,[5,1]$ & $\ZZ/4$ & 1 & 2 & $1/10$ & $(1,1,3,3,1,1,0)$ & $(0,1,0,2,1,3,2)$ & $0,1,1,3,4$ & 0 & $4{\times}[5,1]$ &  \\
166 & $2{\times}[2,1],\,4{\times}[4,1],\,[5,2]$ & $\ZZ/4$ & 1 & 2 & $2/5$ & $(1,1,3,3,1,1,0)$ & $(0,1,0,0,3,3,2)$ & $1,2,4,9,15$ & 1 & $4{\times}[5,2]$ & $-$\,\ref{prop:liu26} \\
167 & $2{\times}[2,1],\,[3,1],\,[4,1],\,[12,5]$ & $\ZZ/4$ & 1 & 2 & $2/3$ & $(1,1,0,3,3)$ & $(0,0,1,0,4)$ & $1,3,8,14,24$ & 2 & $5{\times}[3,1]$ & $-$\,\ref{prop:netpencil} \\
168 & $4{\times}[5,1]$ & $\ZZ/5$ & 2 & 2 & $2/5$ & $(4,3,2,1)$ & $(2,2,2,2)$ & $0,3,4,9,16$ & 2 & $\varnothing$ & $+$\,\ref{prop:examples} \\
169 & $4{\times}[5,1]$ & $\ZZ/5$ & 2 & 2 & $7/5$ & $(4,3,2,1)$ & $(2,2,2,2)$ & $1,7,14,29,51$ & 6 & $\varnothing$ & $-$\,\ref{prop:B24} \\
170 & $4{\times}[5,2]$ & $\ZZ/5$ & 2 & 2 & $3/5$ & $(4,3,2,1)$ & $(2,2,2,2)$ & $1,3,6,13,23$ & 2 & $\varnothing$ & $+$\,\ref{prop:examples} \\
171 & $4{\times}[5,2]$ & $\ZZ/5$ & 2 & 2 & $8/5$ & $(4,3,2,1)$ & $(2,2,2,2)$ & $2,7,16,33,58$ & 6 & $\varnothing$ & $-$\,\ref{prop:B24} \\
172 & $2{\times}[5,1],\,2{\times}[5,2]$ & $\ZZ/5$ & 2 & 2 & $1$ & $(4,1,4,1)$ & $(2,2,2,2)$ & $1,5,10,21,37$ & 4 & $\varnothing$ & $+$\,\ref{prop:V5} \\
173 & $[3,1],\,4{\times}[5,1]$ & $\ZZ/5$ & 2 & 2 & $11/15$ & $(0,4,3,2,1)$ & $(1,2,2,2,2)$ & $0,4,7,15,27$ & 3 & $5{\times}[3,1]$ & $-$\,\ref{prop:173} \\
174 & $[3,1],\,2{\times}[5,1],\,2{\times}[5,2]$ & $\ZZ/5$ & 2 & 2 & $1/3$ & $(0,4,1,4,1)$ & $(1,2,2,2,2)$ & $0,2,3,7,13$ & 1 & $5{\times}[3,1]$ &  \\
175 & $2{\times}[2,1],\,2{\times}[3,1],\,2{\times}[6,1]$ & $\ZZ/6$ & 1 & 2 & $1$ & $(1,1,2,1,5,1)$ & $(0,1,2,0,5,0)$ & $2,5,11,21,36$ & 4 & $\varnothing$ & $-$\,\ref{prop:175} \\
176 & $2{\times}[2,1],\,3{\times}[3,1],\,2{\times}[6,1]$ & $\ZZ/6$ & 1 & 2 & $1/3$ & $(1,1,2,1,0,5,1)$ & $(0,1,2,0,1,5,0)$ & $1,2,4,7,12$ & 1 & $6{\times}[3,1]$ &  \\

%% file: table_subgroups.tex
32 & $\ZZ/2{\times}\ZZ/2$ & $\ZZ/2$ & $(0,0,1,1,0,2,4)$ & $(1,1,1,1,1,1,5)$ & 3 & $0,2,4,8,11$ & $5{\times}[2,1],\,3{\times}[4,1]$ &  \\
 &  & $\ZZ/2$ & $(1,1,1,1,2,2,0)$ & $(1,1,1,1,1,1,5)$ & 3 & $1,2,4,7,11$ & $2{\times}[2,1],\,2{\times}[8,3]$ &  \\
 &  & $\ZZ/2{\times}\ZZ/2$ & $(0,0,1,1,0,2,4);(1,1,1,1,2,2,0)$ & $(1,1,1,1,1,1,5)$ & 3 & $1,4,8,15,23$ & $4{\times}[2,1],\,2{\times}[4,1]$ &  \\
36 & $\ZZ/4$ & $\ZZ/2$ & $(0,0,4,4)$ & $(1,1,5,5)$ & 3 & $1,4,8,15,23$ & $4{\times}[2,1],\,2{\times}[4,1]$ &  \\
 &  & $\ZZ/4$ & $(1,1,2,6)$ & $(1,1,5,5)$ & 3 & $3,8,16,29,47$ & $2{\times}[2,1]$ & $+$\,\ref{prop:examples} \\
155 & $\ZZ/3{\times}\ZZ/3$ & $\ZZ/3$ & $(1,1,1,2,2,2,0,0)$ & $(1,1,1,1,1,1,1,1)$ & 2 & $2,9,22,41,72$ & $6{\times}[3,1]$ &  \\
 &  & $\ZZ/3{\times}\ZZ/3$ & $(0,1,2,1,0,2,1,2);(1,0,2,1,2,0,2,1)$ & $(1,1,1,1,1,1,1,1)$ & 2 & $8,27,64,125,216$ & -- & $-$\,\ref{thm:nonexistence} \\
156 & $\ZZ/4$ & $\ZZ/4$ & $(1,1,1,1,3,3)$ & $(0,1,0,2,1,3)$ & 1 & $4,11,24,45,76$ & -- & $+$\,\ref{prop:examples} \\
157 & $\ZZ/4$ & $\ZZ/4$ & $(1,1,1,1,3,3)$ & $(0,1,0,0,3,3)$ & 1 & $6,19,44,85,146$ & -- & $+$\,\ref{prop:examples} \\
158 & $\ZZ/4$ & $\ZZ/4$ & $(1,1,1,1,3,3)$ & $(0,1,0,2,1,3)$ & 1 & $8,27,64,125,216$ & -- & $-$\,\ref{prop:B24} \\
159 & $\ZZ/4$ & $\ZZ/4$ & $(1,1,1,1,3,3)$ & $(0,1,0,0,3,3)$ & 1 & $10,35,84,165,286$ & -- & $-$\,\ref{prop:B24} \\
160 & $\ZZ/4$ & $\ZZ/4$ & $(1,1,0,0,1,1,3,3)$ & $(0,1,1,1,0,2,1,3)$ & 1 & $0,3,8,13,24$ & $8{\times}[3,1]$ &  \\
161 & $\ZZ/4$ & $\ZZ/4$ & $(1,1,0,1,1,3,3)$ & $(0,1,1,0,0,3,3)$ & 1 & $2,7,16,29,50$ & $4{\times}[3,1]$ &  \\
162 & $\ZZ/4$ & $\ZZ/4$ & $(1,1,0,1,1,3,3)$ & $(0,1,1,0,2,1,3)$ & 1 & $4,15,36,69,120$ & $4{\times}[3,1]$ & $-$\,\ref{prop:162} \\
163 & $\ZZ/4$ & $\ZZ/4$ & $(1,1,1,9)$ & $(0,0,0,4)$ & 1 & $3,8,17,31,52$ & $[3,1]$ & $+$\,\ref{prop:examples} \\
164 & $\ZZ/4$ & $\ZZ/4$ & $(1,1,1,9)$ & $(0,0,0,4)$ & 1 & $7,24,57,111,192$ & $[3,1]$ & $-$\,\ref{prop:largedeg} \\
165 & $\ZZ/4$ & $\ZZ/4$ & $(1,1,1,1,3,3,0)$ & $(0,1,0,2,1,3,2)$ & 1 & $0,3,4,9,16$ & $4{\times}[5,1]$ &  \\
166 & $\ZZ/4$ & $\ZZ/4$ & $(1,1,1,1,3,3,0)$ & $(0,1,0,0,3,3,2)$ & 1 & $2,7,16,33,58$ & $4{\times}[5,2]$ & $-$\,\ref{prop:liu26} \\
167 & $\ZZ/4$ & $\ZZ/4$ & $(1,1,0,1,9)$ & $(0,0,1,0,4)$ & 1 & $3,12,29,55,96$ & $5{\times}[3,1]$ & $-$\,\ref{prop:netpencil} \\
175 & $\ZZ/6$ & $\ZZ/2$ & $(1,1,0,0,3,3)$ & $(0,1,1,1,1,4)$ & 1 & $2,9,22,41,72$ & $6{\times}[3,1]$ & $-$\,\ref{prop:175} \\
 &  & $\ZZ/6$ & $(1,1,1,2,1,5)$ & $(0,1,2,0,5,0)$ & 1 & $8,27,64,125,216$ & -- & $-$\,\ref{prop:B24} \\
176 & $\ZZ/6$ & $\ZZ/2$ & $(1,1,0,0,0,3,3)$ & $(0,1,1,1,1,1,4)$ & 1 & $0,3,8,13,24$ & $8{\times}[3,1]$ &  \\
 &  & $\ZZ/6$ & $(1,1,1,2,0,1,5)$ & $(0,1,2,0,1,5,0)$ & 1 & $2,9,22,41,72$ & $6{\times}[3,1]$ &  \\